%% file: invitation-CBB.tex
\documentclass[twoside]{book}

\def\thetitle{Lectures on Alexandrov spaces with curvature bounded below}
\def\theauthors{Vitali Kapovitch, Nina Lebedeva, and Anton Petrunin}

\usepackage{lectures,acts}
\usepackage[colorlinks=true,
citecolor=black,
linkcolor=black,
anchorcolor=black,
filecolor=black,
menucolor=black,
urlcolor=black,
pdftitle={\thetitle},
pdfsubject={Geometry},
pdfauthor={\theauthors}
]{hyperref}
\makeindex

\begin{document}

\frontmatter
\title{\parbox{8cm}{\centering\thetitle}}
\author{Vitali Kapovitch, Nina Lebedeva, and Anton Petrunin}
\date{}
\maketitle
\thispagestyle{empty}

\mainmatter
\newpage
\cleardoublepage
\phantomsection
\pdfbookmark[0]{\contentsname}{bm:toc}
\tableofcontents

\include{preface}

\include{prelim}

\include{definitions}
\include{globalization}
\include{calculus}
\include{gradient-flow}
\include{splitting}
\include{dim}
\include{volume}
\include{boundary}
\include{CBB-quotients}
\include{overview}

\backmatter

\include{sols}

{
\cleardoublepage
\phantomsection
\sloppy
\addcontentsline{toc}{chapter}{\indexname}
\input{invitation-CBB.ind}
}

{
\def\emph{\textit}
\sloppy
\printbibliography[heading=bibintoc]
\fussy
}

\end{document}

%% file: preface.tex
\chapter*{Preface}
\addcontentsline{toc}{chapter}{Preface}
\markboth{PREFACE}{PREFACE}

This book is similar to our ``Invitation to Alexandrov geometry'' written jointly with Stephanie Alexander \cite{alexander-kapovitch-petrunin-2019}.
We try to demonstrate the beauty and power of Alexandrov geometry by reaching interesting applications and theorems with minimal preparation.
This time it is about curvature bounded below.

This subject is more technical and requires more preparation, so we had to skip several proofs.
Namely, we do not prove the existence part of the generalized Picard theorem (\ref{thm:glob-exist-grad-curv})
and Perelman's theorem about conic neighborhoods (\ref{thm:spherical-nbhd}).
Apart from the last lecture, the rest is nearly rigorous.
Several important statements are left as exercises (they are marked with \textbf{Exercise!}).
Solutions to all of them are given at the end of the book.

\medskip

Following the preliminaries, Lecture~\ref{chap:defs} introduces the main object of our study --- spaces with curvature bounded below in the sense of Alexandrov.

In Lecture~\ref{chap:globalization}, we formulate and prove the globalization theorem that the local Alexandrov condition implies the global one.
To simplify the presentation, we treat only the compact case; it captures the main ideas of the proof with fewer technicalities.

In Lecture~\ref{chap:derivative}, we develop calculus --- tangent space and space of directions, differential, and gradient.

Lecture~\ref{chap:GF} introduces gradient flow, which will be used as the main technical tool.

Lecture~\ref{chap:splitting} proves the line splitting theorem.
We then introduce and study linear subspaces of tangent spaces.

In Lecture~\ref{chap:dim}, we introduce linear dimension and volume.
Next, we prove the Bishop--Gromov inequality and the right-inverse theorem,
introduce distance charts, and show that all reasonable notions of dimension coincide for Alexandrov spaces.

In Lecture~\ref{chap:lim}, we show that a lower curvature bound survives under Gromov--Hausdorff limits and prove Gromov's selection theorem.
Further, we present Perelman's construction of strictly concave functions and
apply it together with Gromov's selection theorem to prove the homotopy finiteness theorem.
This proof illustrates the main source of applications of Alexandrov geometry.

In Lecture~\ref{chap:bry}, we introduce the boundary of finite-dimensional Alexandrov spaces and prove the doubling theorem.

In Lecture~\ref{chap:L/G}, we show that quotients of Alexandrov spaces by isometric group actions are Alexandrov spaces.
This gives another source of applications of Alexandrov geometry; several examples are given.

Finally, in Lecture~\ref{chap:surfaces} we briefly discuss convex surfaces in Euclidean space; this subject is the main precursor to modern Alexandrov geometry.

Let us list some available texts on Alexandrov spaces with curvature bounded below:
\begin{itemize}
\item The first introduction to Alexandrov geometry of all dimensions was given in the original paper by Yuriy Burago, Michael Gromov, and Grigory Perelman \cite{burago-gromov-perelman} and its extension \cite{perelman1991}, written by Perelman.
\item A brief and reader-friendly introduction was written by Katsuhiro Shiohama \cite[Sections 1--8]{shiohama}.
\item Another reader-friendly introduction was written by Dmitri Burago, Yuriy
Burago, and Sergei Ivanov \cite[Chapter 10]{burago-burago-ivanov}.
\item Survey by Conrad Plaut \cite{plaut:survey}.
\item Survey by the third author \cite{petrunin:survey}.
\item The book by two of us and Stephanie Alexander \cite{alexander-kapovitch-petrunin2024}.
\end{itemize}

\parbf{Acknowledgments.}
This book owes much to Stephanie Alexander: it grew out of our collaboration and draws extensively on our joint book \cite{alexander-kapovitch-petrunin2024}.
Our notes were shaped by a series of lectures given by the authors on various occasions: at Penn State, including the MASS program; at the Summer School ``Algebra and Geometry'' in Yaroslavl; at SPbSU; and at the University of Toronto.
We are grateful to these institutions for their hospitality and support.

The authors were partially supported as follows:
Vitali Kapovitch by an NSERC Discovery Grant;
Nina Lebedeva by the Ministry of Science and Higher Education of the Russian Federation, agreement no. 075-15-2025-343;
Anton Petrunin by NSF grant DMS-2005279.

%% file: prelim.tex
\chapter*{Preliminaries}
\phantomsection
\addcontentsline{toc}{chapter}{Preliminaries}
\markboth{\MakeUppercase{Preliminaries}}{}%

\section{Prerequisites}

The main prerequisite: the reader should know and love elementary geometry, including convexity.

We also assume that the reader is familiar with the following topics in metric geometry:
\begin{itemize}
\item Compactness and proper metric spaces;
recall that a metric space is \index{proper space}\emph{proper} if all its closed balls (of finite radius) are compact.
\item Complete metric spaces and completion.
\item Curves, semicontinuity of length, and rectifiability.
\item Hausdorff and Gromov--Hausdorff convergence.
These are discussed briefly in \ref{sec:Hausdorff convergence}--\ref{sec:Almost isometries}.
The definitions are there, but some prior familiarity with these concepts will be very helpful.
\end{itemize}
All these topics are treated in detail in the book of Dmitri Burago, Yurii Burago, and Sergei Ivanov \cite{burago-burago-ivanov}; see also the book by the third author \cite{petrunin2023pure}.
Occasionally, we use the Baire category theorem and Rademacher's theorem, but these could be used as black boxes.

We use some topology. 
An introductory text in algebraic topology should be sufficient most of the time.
For some examples, we apply more advanced results, but these are used as black boxes.

Since most of the applications come from Riemannian geometry, it is beneficial to be familiar with the Toponogov comparison theorem and related topics.
The classical book by Jeff Cheeger and David Ebin \cite{cheeger-ebin} contains more than one needs.

\section{Notations}

The distance between two points $x$ and $y$ in a metric space $\spc{X}$ will be denoted by \index{9@$\dist{x}{y}{}=\dist{x}{y}{\spc{X}}$ (distance)}$\dist{x}{y}{}=\dist{x}{y}{\spc{X}}$, or \index{dist@$\distfun$ (distance function)}$\distfun_xy=\distfun_x(y)_{\spc{X}}$.\label{page:|x-y|X}
Furthermore, we define the distance function to a subset $A\subset \spc{X}$ as
\[\distfun_Ax=\distfun_A(x)_{\spc{X}}
\df
\inf
\set{\dist{a}{x}{\spc{X}}}{a\in A}.\]
We may use the subscript $\spc{X}$ to emphasize
that the distance is taken in the space~${\spc{X}}$.

Given radius $r\in[0,\infty]$ and center $x\in \spc{X}$, the sets
\begin{align*}
\oBall(x,r)&=\set{y\in \spc{X}}{\dist{x}{y}{}<r},
\\
\cBall[x,r]&=\set{y\in \spc{X}}{\dist{x}{y}{}\le r}
\end{align*}
are called, respectively, the \index{open ball}\emph{open} and the \index{closed ball}\emph{closed balls}.
The notations $\oBall(x,r)_{\spc{X}}$ and $\cBall[x,r]_{\spc{X}}$
might be used if we need to emphasize that these balls are taken in the metric space $\spc{X}$.

The \index{diameter}\emph{diameter} of $\spc{X}$ is defined as
\[\diam \spc{X}=\sup\set{\dist{x}{y}{}}{x,y\in \spc{X}}.\]

The \index{radius}\emph{radius} of $\spc{X}$ is defined as
\[\rad \spc{X}=\inf\set{R>0}{\oBall(x,R)= \spc{X}\ \text{for some}\ x\in\spc{X}}.\]

We will denote by \index{60@$\SSS^m$, $\EE^m$, $\HH^m$, and $\MM^m(\kappa)$}$\SSS^m$, $\EE^m$, and $\HH^m$ the $m$-dimensional sphere (with angle metric),
the Euclidean space, and the Lobachevsky space, respectively.
More generally, $\MM^m(\kappa)$ will denote the \index{model!space}\emph{model $m$-space} of curvature $\kappa$;
that is,
\begin{itemize}
\item if $\kappa>0$, then $\MM^m(\kappa)$ is the $m$-sphere of radius $\tfrac{1}{\sqrt{\kappa}}$, so $\SSS^m\z=\MM^m(1)$,
\item $\MM^m(0)=\EE^m$,
\item if $\kappa<0$, then $\MM^m(\kappa)$ is the Lobachevsky $m$-space $\HH^m$ rescaled by the factor $\tfrac{1}{\sqrt{-\kappa}}$;
in particular, $\MM^m(-1)=\HH^m$.
\end{itemize}

\section{Nets and packing}\label{sec:packing}

Let $S$ be a subset of a metric space $\spc{X}$.
Recall that a set $Z\subset \spc{X}$ is called an \index{net@$\eps$-net}\emph{$\eps$-net} of $S$ if for any point $s\in S$, there is a point $z\in Z$ such that $\dist{s}{z}{}\le\eps$.

The \index{pack@$\eps$-pack}\emph{$\eps$-pack} of $\spc{X}$ (or \index{packing number}\emph{packing number}) is the maximal number (possibly infinite) of points in $\spc{X}$ at distance $>\eps$ from each other;
it is denoted by $\pack_\eps\spc{X}$\index{70@$\pack_\eps\spc{X}$ ($\eps$-pack)}.
If $n=\pack_\eps\spc{X}<\infty$, then a set $\{x_1,x_2,\dots,x_n\}$ in $\spc{X}$
such that $\dist{x_i}{x_j}{}>\eps$ is called a \index{maximal packing}\emph{maximal $\eps$-packing} in $\spc{X}$.

We will use the following characterizations of compact sets.

\begin{thm}[!]{Exercise}\label{ex:net}
Let $X$ be a closed subset of a complete metric space.

\begin{subthm}{ex:net:finite}
Show that $X$ is compact if and only if it admits a finite $\eps$-net for any $\eps>0$.
\end{subthm}

\begin{subthm}{ex:net:compact}
Show that $X$ is compact if and only if it admits a compact $\eps$-net for any $\eps>0$.
\end{subthm}

\end{thm}

\begin{thm}[!]{Exercise}\label{ex:pack-net}
Show that any maximal $\eps$-packing $x_1,\dots,x_n$ is an $\eps$-net.
Conclude that a complete metric space $\spc{X}$ is compact if and only if $\pack_\eps\spc{X}\z<\infty$ for any $\eps>0$.
\end{thm}

\section{Length spaces}\label{sec:length}

Let $\spc{X}$ be a metric space.
If for any $\eps>0$ and any pair of points $x,y\in\spc{X}$, there is a path $\alpha$ connecting $x$ to $y$ such that
\[\length\alpha< \dist{x}{y}{}+\eps,\]
then $\spc{X}$ is called a \index{length space}\emph{length space} and the metric on $\spc{X}$ is called a \index{length metric}\emph{length metric}.\label{page:length metric}

\begin{thm}{Exercise}\label{ex:compact+connceted}
Let $\spc{X}$ be a complete length space.
Show that for any compact subset $K\subset\spc{X}$
there is a compact path-connected subset $K'\subset\spc{X}$ that contains $K$.
\end{thm}

\parbf{Induced length metric.}
From the definition of length, it follows that if $\alpha\:[0,1]\to\spc{X}$ is a path from $x$ to $y$
(that is, $\alpha(0)=x$ and $\alpha(1)=y$), then 
\[\length\alpha\ge \dist{x}{y}{}.\]
Set 
\[\yetdist{x}{y}{}=\inf_\alpha\{\,\length\alpha\,\},\]
where the greatest lower bound is taken for paths $\alpha$ from $x$ to~$y$.
It is straightforward to check that $(x,y)\mapsto \yetdist{x}{y}{}$ is an \index{metric@$\infty$-metric}\emph{$\infty$-metric};
that is, $(x,y)\mapsto \yetdist{x}{y}{}$ is a metric in the extended nonnegative reals $[0,\infty]$.
The metric $\yetdist{*}{*}{}$ is called the \index{induced length metric}\emph{induced length metric}.

\begin{thm}{Exercise}\label{ex:compact=>complete}
Suppose $(\spc{X},\dist{*}{*}{})$ is a complete metric space.
Show that $(\spc{X},\yetdist{*}{*}{})$ is complete;
that is, any Cauchy sequence of points in $(\spc{X},\yetdist{*}{*}{})$ converges in $(\spc{X},\yetdist{*}{*}{})$.
\end{thm}

Let $A$ be a subset of a metric space $\spc{X}$.
Given two points $x,y\in A$,
consider the value
\[\dist{x}{y}{A}=\inf_{\alpha}\{\,\length\alpha\,\},
\]
where the greatest lower bound is taken for all paths $\alpha$ from $x$ to $y$ in~$A$.
In other words, $\dist{*}{*}{A}$ denotes the induced length metric on the subspace $A$.
(The notation $\dist{*}{*}{A}$ conflicts with the previously defined notation for distance $\dist{x}{y}{\spc{X}}$ in a metric space $\spc{X}$.
However, most of the time we will work with ambient length spaces where the meaning will be unambiguous.)

\section{Geodesics}

Let $\spc{X}$ be a metric space,
and let $\II$\index{60@$\II$ (real interval)} be a real interval.
A distance-preserving map $\gamma\:\II\to \spc{X}$ is called a \index{geodesic}\emph{geodesic}%
\footnote{Other texts may refer to geodesics in our sense of the word as  a \textit{shortest path} or a \textit{minimizing geodesic}.
Also, our meaning of the term \textit{geodesic} is closely related to, but different from, the term used in Riemannian geometry.};
in other words, $\gamma\:\II\z\to \spc{X}$ is a geodesic if 
\[\dist{\gamma(s)}{\gamma(t)}{}=|s-t|\]
for any pair $s,t\in \II$.

If $\gamma\:[a,b]\to \spc{X}$ is a geodesic such that $p=\gamma(a)$, $q=\gamma(b)$, then we say that $\gamma$ is a geodesic from $p$ to $q$.
In this case, the image of $\gamma$ is denoted by $[p q]$\index{10@$[pq]=[p q]_{\spc{X}}$ (geodesic)}, and, with abuse of notation, we also call it a \index{geodesic}\emph{geodesic}.
We may write $[p q]_{\spc{X}}$ 
to emphasize that the geodesic $[p q]$ lies in the space ${\spc{X}}$.

In general, a geodesic from $p$ to $q$ need not exist, and if it exists, it need not be unique;
for example, any meridian is a geodesic between poles on the sphere.
However, once we write $[p q]$ we assume that we have chosen such a geodesic.

A \index{geodesic!path}\emph{geodesic path} is a constant-speed parametrization of a geodesic by the unit interval $[0,1]$.

A metric space is called \index{geodesic!space}\emph{geodesic} if any pair of its points can be joined by a geodesic.
A metric space $\spc{X}$ is called \index{geodesic space@$\ell$-geodesic space}\emph{$\ell$-geodesic}
if any two points $x,y\in\spc{X}$ such that $\dist{x}{y}{}<\ell$ can be connected by a geodesic.
For instance, any geodesic space is $\infty$-geodesic.

Evidently, any geodesic space is a length space.

\begin{thm}{Exercise}\label{ex:compact-length}
Show that any proper length space is geodesic.
\end{thm}

\section{Menger's lemma}
Given two points $x,y$ in a metric space $\spc{X}$, a point $z\in \spc{X}$ is called a  \index{midpoint}\emph{midpoint} between them if
\[\dist{x}{z}{}=\tfrac12\cdot\dist{x}{y}{}
\quad\text{and}\quad
\dist{y}{z}{}=\tfrac12\cdot\dist{x}{y}{}.\]

If $\spc{X}$ is geodesic, then for any two points there is a midpoint between them.
The following lemma says that the converse holds for complete metric spaces.

\begin{thm}{Menger's lemma}\label{lem:mid>geod}\index{Menger's lemma}
Let $\spc{X}$ be a complete metric space.
Assume that for any pair of points $x,y\in \spc{X}$, 
there is a midpoint~$z$.
Then $\spc{X}$ is a geodesic space.

\end{thm}

This lemma is due to Karl Menger \cite[Section 6]{menger}.

\parit{Proof.}
Choose $x,y\in \spc{X}$;
set $\gamma(0)=x$ and $\gamma(1)=y$.

\begin{figure}[ht!]
\vskip-0mm
\centering
\includegraphics{mppics/pic-104}
\end{figure}

Let $\gamma(\tfrac12)$ be a midpoint between $\gamma(0)$ and $\gamma(1)$.
Further, let $\gamma(\frac14)$ 
and $\gamma(\frac34)$ be midpoints between the pairs $(\gamma(0),\gamma(\tfrac12))$ 
and $(\gamma(\tfrac12),\gamma(1))$, respectively.
Applying the above procedure recursively,
at the $n$-th step we define $\gamma(\tfrac{k}{2^n})$,
for every odd integer $k$ such that $0<\tfrac k{2^n}<1$, 
as a midpoint of the already defined
$\gamma(\tfrac{k-1}{2^n})$ and $\gamma(\tfrac{k+1}{2^n})$.

This way, we define $\gamma(t)$ for all dyadic rationals $t$ in $[0,1]$.
Moreover, the map $\gamma$ has Lipschitz constant $\dist{x}{y}{}$.
Since $\spc{X}$ is complete, $\gamma$ can be extended continuously to $[0,1]$ with the same Lipschitz constant.
Therefore,
\[
\length\gamma\le \dist{x}{y}{}.
\]
Hence $\gamma$ is a geodesic path from $x$ to $y$.
\qedsf

\begin{thm}{Exercise}\label{ex:menger}
Let $\spc{X}$ be a complete metric space.
Assume that for any pair of points $x,y\in \spc{X}$, 
there is an \index{almost midpoint}\emph{almost midpoint};
that is, given $\eps>0$, there is a point $z$ such that 
\[\dist{x}{z}{}<\tfrac12\cdot\dist{x}{y}{}+\eps 
\quad\text{and}\quad
\dist{y}{z}{}<\tfrac12\cdot\dist{x}{y}{}+\eps.\]
Show that $\spc{X}$ is a length space.
Conversely, if $\spc{X}$ is a length space, then for any two points there exist almost midpoints between them.
\end{thm}

\section{Triangles and model triangles}\label{sec:Triangles and model triangles}

\parbf{Triangles.}
Given a triple of points $p,q,r$ in a metric space $\spc{X}$, a choice of geodesics $([q r], [r p], [p q])$ will be called a \index{triangle}\emph{triangle} and denoted by $\trig p q r=\trig p q r_{\spc{X}}$\index{10@$\trig p q r=\trig p q r_{\spc{X}}$ (triangle)}.

Given a triple of points $p,q,r\in \spc{X}$, there may be no triangle
$\trig p q r$ simply because one of the pairs of these points cannot be joined by a geodesic.
Also, many different triangles with these vertices may exist, any of which can be denoted by $\trig p q r$.
If we write $\trig p q r$, it means that we have chosen such a triangle.

\parbf{Model triangles.}
Given three points $p,q,r$ in a metric space $\spc{X}$,
let us define the \index{model!triangle}\emph{model triangle} $\trig{\tilde p}{\tilde q}{\tilde r}$ 
(briefly, 
$\trig{\tilde p}{\tilde q}{\tilde r}=\modtrig(p q r)_{\EE^2}$%
\index{20@$\modtrig$ (model triangle)}) to be a triangle in the Euclidean plane $\EE^2$ with the same sides;
that is,
\[
\dist{\tilde p}{\tilde q}{\EE^2}=\dist{p}{q}{\spc{X}},
\quad\dist{\tilde q}{\tilde r}{\EE^2}=\dist{q}{r}{\spc{X}},
\quad\dist{\tilde r}{\tilde p}{\EE^2}=\dist{r}{p}{\spc{X}}.
\]

In the same way, we can define the \index{hyperbolic model triangle}\emph{hyperbolic} and the \index{spherical model triangles}\emph{spherical model triangles} $\modtrig(p q r)_{\HH^2}$, $\modtrig(p q r)_{\SSS^2}$
in the Lobachevsky plane $\HH^2$ and the unit sphere~$\SSS^2$.
In the latter case, the model triangle is said to be defined if in addition
\[\dist{p}{q}{}+\dist{q}{r}{}+\dist{r}{p}{}< 2\cdot\pi.\]
Then the model triangle again exists and is unique up to an isometry of~$\SSS^2$.

\parbf{Model angles.}
If 
$\trig{\tilde p}{\tilde q}{\tilde r}=\modtrig(p q r)_{\EE^2}$ 
and $\dist{p}{q}{},\dist{p}{r}{}>0$, 
the angle measure of 
$\trig{\tilde p}{\tilde q}{\tilde r}$ at $\tilde p$ 
will be called the \index{model!angle}\emph{model angle} of the triple $p$, $q$, $r$ and will be denoted by
$\angk p q r_{\EE^2}$%
\index{30@$\angk{p}{q}{r}$ (model angle)}.\label{page:model-angle}

For example, if $\dist{p}{q}{}=\dist{q}{r}{}=\dist{r}{p}{}$, then $\angk p q r_{\EE^2}=\tfrac\pi3$ regardless of the existence and relative position of geodesics $[pq]$ and $[pr]$.

In the same way we define $\angk p q r_{\MM^2(\kappa)}$;
in particular, $\angk p q r_{\HH^2}$ and $\angk p q r_{\SSS^2}$.
We may use the notation $\angk p q r$ if it is evident which of the model spaces is meant.

\begin{thm}[!]{Exercise}\label{ex:k-><mono}
Show that for any triple of points $p$, $q$, and $r$,
the function
\[\kappa\mapsto \angk p q r_{\MM^2(\kappa)}\]
is nondecreasing in the interval of definition.
\end{thm}

\section{Hinges and their angle measure}\label{sec:angles}

\parbf{Hinges.} Let $p,x,y\in \spc{X}$ be a triple of points such that $p$ is distinct from $x$ and~$y$.
A pair of geodesics $([p x],[p y])$ will be called a \index{hinge}\emph{hinge} and will be denoted by
$\hinge p x y=([p x],[p y])$\index{11@$\hinge p x y$ (hinge)}.

\parbf{Angles.}
The angle measure of a hinge $\hinge p x y$ is defined as the following limit \index{31@$\mangle$ (angle measure)}
\[\mangle\hinge p x y=\lim_{\bar x,\bar y\to p} \angk p{\bar x}{\bar y},\]
where $\bar x\in\left]p x\right]$ and $\bar y\in\left]p y\right]$.
The angle is only defined if this limit exists.

If $\mangle\hinge p x y$ is defined, then
\[0\le \mangle\hinge p x y\le \pi.\]

\begin{thm}{Exercise}\label{ex:angkK}
Suppose that in the above definition one uses spherical or hyperbolic model angles instead of Euclidean ones.
Show that this does not change the value $\mangle\hinge p x y$.
\end{thm}

\begin{thm}{Exercise}\label{ex:undefined-angle}
Give an example of a hinge $\hinge p x y$ in a metric space with undefined angle measure $\mangle\hinge p x y$.
\end{thm}

\section{Triangle inequality for angles}

\begin{thm}{Proposition}\label{claim:angle-3angle-inq}
Consider three geodesics $[px_1]$, $[px_2]$, and $[px_3]$ in a metric space.
Suppose all the angle measures $\alpha_{i, j}=\mangle\hinge p {x_i}{x_j}$ are defined.
Then 
\[\alpha_{1,3}\le \alpha_{1,2}+\alpha_{2,3}.\]

\end{thm}

\parit{Proof.}
We can assume that $\alpha_{1,3}>0$.
Furthermore, since $\alpha_{1,3}\le\pi$, we can assume that $\alpha_{1,2}+\alpha_{2,3}< \pi$.

Denote by $\gamma_i$ the unit-speed parametrization of $[px_i]$ from $p$ to $x_i$.
Given a small $\eps>0$, for all sufficiently small $t,\tau,s\in[0,\infty)$ we have
\begin{align*}
\sqrt{t^2+\tau^2-2\cdot t\cdot \tau\cdot \cos(\alpha_{1,3}-\eps)}
&\le
\dist{\gamma_1(t)}{\gamma_3(\tau)}{}\le
\\
\le
\dist{\gamma_1(t)}{\gamma_2(s)}{}+\dist{\gamma_2(s)}{\gamma_3(\tau)}{}&<\\
<
\sqrt{t^2+s^2-2\cdot t\cdot s\cdot \cos(\alpha_{1,2}+\eps)}\;&+
\\
\quad+\;\sqrt{s^2+\tau^2-2\cdot s\cdot \tau\cdot \cos(\alpha_{2,3}+\eps)}&\le
\end{align*}

\begin{wrapfigure}{o}{30 mm}
\vskip-18mm
\centering
\includegraphics{mppics/pic-615}
\vskip0mm
\end{wrapfigure}

Below we define 
$s(t,\tau)$ so that for 
$s=s(t,\tau)$, this chain of inequalities can be continued as
\[\le
\sqrt{t^2+\tau^2-2\cdot t\cdot \tau\cdot \cos(\alpha_{1,2}+\alpha_{2,3}+2\cdot \eps)}.
\]
Thus
\[\alpha_{1,3}\le \alpha_{1,2}+\alpha_{2,3}+3\cdot \eps.\]
Hence the result follows.

To define $s(t,\tau)$, consider three half-lines $\tilde \gamma_1$, $\tilde \gamma_2$, $\tilde \gamma_3$ in a Euclidean plane starting at one point, such that
$\mangle(\tilde \gamma_1,\tilde \gamma_2)\z=\alpha_{1,2}+\eps$,
$\mangle(\tilde \gamma_2,\tilde \gamma_3)\z=\alpha_{2,3}+\eps$,
and $\mangle(\tilde \gamma_1,\tilde \gamma_3)\z=\alpha_{1,2}\z+\alpha_{2,3}\z+2\cdot \eps$.
We parametrize each half-line by the distance from the starting point.
Given two positive numbers $t,\tau$, let $s=s(t,\tau)$ be
the number such that 
$\tilde \gamma_2(s)\z\in[\tilde \gamma_1(t)\ \tilde \gamma_3(\tau)]$. 
Clearly, $s\le\max\{t,\tau\}$, so $t,\tau,s$ may be taken sufficiently small.
\qeds 

\begin{thm}[!]{Exercise}\label{ex:adjacent-angles}
Prove that the sum of adjacent angles is at least $\pi$.

More precisely, suppose two hinges $\hinge pxz$ and $\hinge pyz$ are \index{adjacent hinges}\emph{adjacent};
that is, they share the side $[pz]$, and the union of the two sides $[px]$ and $[py]$ forms a geodesic $[xy]$.
Show that
\[\mangle\hinge pxz+\mangle\hinge pyz\ge \pi\]
whenever each angle on the left-hand side is defined.

Give an example showing that the inequality can be strict.
\end{thm}

\begin{thm}{Exercise}\label{ex:first-var}
Let $\gamma$ be the unit-speed parametrization of $[qx]$ from $q$ to $x$.
Assume that the angle measure $\phi=\mangle\hinge q p x$ is defined.
Show that
\[\dist{p}{\gamma(t)}{}
\le
\dist{q}{p}{}-t\cdot \cos\phi+o(t).\]

\end{thm}

\section{Hausdorff convergence}\label{sec:Hausdorff convergence}

\begin{thm}{Definition}\label{def:gen-Haus-conv}
Let $A_1,A_2,\dots$ be a sequence of closed sets in a metric space $\spc{X}$.
We say that the sequence $A_n$ \index{Hausdorff!limit}\emph{converges} to a closed set $A_\infty$ in the {}\emph{sense of Hausdorff} if, for any $x\in\spc{X}$, we have
$\distfun_{A_n}x\z\to \distfun_{A_\infty}x$ as $n\to\infty$.
\end{thm}

For example, suppose $\spc{X}$ is the Euclidean plane and $A_n$ is the circle with radius $n$ and center at the point $(0,n)$.
Then $A_n$ converges to the $x$-axis as $n\to \infty$.

\begin{figure}[ht!]
\vskip-0mm
\centering
\includegraphics{mppics/pic-415}
\end{figure}

Further, consider the sequence of one-point sets $B_n=\{(n,0)\}$ in the Euclidean plane.
It converges to the empty set;
indeed, for any point $x$ we have $\distfun_{B_n}(x)\to\infty$ as $n\to \infty$ and $\distfun_{\emptyset}(x)= \infty$ for any~$x$.

The following exercise is an extension of the so-called Blaschke selection theorem to our version of Hausdorff convergence.

\begin{thm}{Exercise}\label{ex:generalized-selection}
Show that any sequence of closed sets in a proper metric space has a convergent subsequence in the sense of Hausdorff.
\end{thm}

\begin{thm}{Exercise}\label{ex:two>one}
Construct a metric space with a sequence of points $\{x_0,x_1,\dots\}$ such that $\dist{x_0}{x_n}{}=1$ for any $n>0$ and the two-point sets $\{x_0,x_n\}$ converge to the one-point set $\{x_0\}$ in the sense of Hausdorff.
\end{thm}

\section{Hausdorff metric}\label{sec:Hausdorff metric}

\begin{thm}{Definition}\label{def:hausdorff-convergence}
Let $A$ and $B$ be two nonempty compact subsets of a metric space $\spc{X}$.
Then the \index{Hausdorff!distance}\emph{Hausdorff distance} between $A$ and $B$ is defined as 
$$|A-B|_{\Haus\spc{X}}
\df
\sup\set{|\distfun_Ax-\distfun_Bx|}{x\in \spc{X}}.
$$

\end{thm}

The following observation gives a useful reformulation of the definition:

\begin{thm}{Observation}\label{obs:Haus-nbhds}
Suppose $A$ and $B$ are two nonempty compact subsets of a metric space $\spc{X}$.
Then $|A-B|_{\Haus\spc{X}}< R$ if and only if
$B$ lies in the $R$-neighborhood of $A$, 
and 
$A$ lies in the $R$-neighborhood of~$B$.
\end{thm}

According to the following exercise, Hausdorff convergence of nonempty compact subsets is the convergence in the Hausdorff metric.

\begin{thm}{Exercise}\label{ex:Haus-conv}
Let $A_1,A_2,\dots,$ and $A_\infty$ be compact nonempty sets in a proper metric space $\spc{X}$ such that all $A_i$ are contained in a compact subset of $\spc{X}$.

Show that $\dist{A_n}{A_\infty}{\Haus\spc{X}}\to 0$ as $n\to\infty$
if and only if $A_n\to A_\infty$ in the sense of Hausdorff.

Does the statement remain true without the assumption that all $A_i$ are contained in a compact subset of $\spc{X}$?

\end{thm}

\section{Gromov--Hausdorff convergence}\label{sec:Gromov--Hausdorff}

Let $\spc{X}_1,\spc{X}_2,\dots,$ and $\spc{X}_\infty$ be a sequence of proper metric spaces.
Suppose that there is a metric on the disjoint union 
\[\bm{X}=\bigsqcup_{n\in \NN\cup\{\infty\}} \spc{X}_n\] 
that satisfies the following property:

\begin{thm}{Property}\label{propery:GH}
The restriction of the metric on each $\spc{X}_n$ and $\spc{X}_\infty$ coincides with its original metric,
the space $\bm{X}$ is proper,
and $\spc{X}_n\z\to \spc{X}_\infty$ as subsets of $\bm{X}$ in the sense of Hausdorff.
\end{thm}

In this case we say that the metric on $\bm{X}$ \textit{defines} a \index{Gromov--Hausdorff limit}\emph{convergence} $\spc{X}_n\z\to \spc{X}_\infty$ in the {}\emph{sense of Gromov--Hausdorff}.
The metric on $\bigsqcup \spc{X}_n$ makes it possible to talk about limits of sequences $x_n\in \spc{X}_n$ as $n\to\infty$, as well as weak limits of a sequence of Borel measures $\mu_n$ on $\spc{X}_n$ and so on.

The limit space is not uniquely defined by the sequence.
For example, if each space $\spc{X}_n$ is isometric to the half-line, then the limit may be isometric to the half-line or to the whole line.
The first convergence is evident and the second could be guessed from the picture.

\begin{figure}[ht!]
\vskip-0mm
\centering
\includegraphics{mppics/pic-500}
\end{figure}

Note that any sequence of spaces has an empty space as a limit in some Gromov--Hausdorff convergence.
As we will see later, if the limit is nonempty and compact, then it is unique up to isometry.

\parbf{Pointed convergence.}
The isometry class of the limit can be fixed by marking a point $p_n$ in each space $\spc{X}_n$.
We say that $(\spc{X}_n,p_n)$ converges to $(\spc{X}_\infty,p_\infty)$ if there is a metric on $\bm{X}$ as in \ref{propery:GH} such that $p_n\to p_\infty$.
This is called \index{pointed convergence}\emph{pointed Gromov--Hausdorff convergence}.
For example, the sequence $(\spc{X}_n,p_n)=([0,\infty),0)$ converges to $([0,\infty),0)$, while $(\spc{X}_n,p_n)=([0,\infty),n)$ converges to $(\RR,0)$ as $n\to \infty$.

\section{Gromov--Hausdorff metric}\label{sec:Gromov--Hausdorff-metric}

In this section we define the metric on the set of \textit{isometry classes} of compact metric spaces.
In what follows, the term \textit{metric space} might also stand for its \textit{isometry class}.

Loosely speaking, the distance between (isometry classes of) compact metric spaces is the minimal value $r$ such that  both spaces can be isometrically embedded into a third space in such a way that the Hausdorff distance between the images is not greater than $r$.

Observe that we \textit{know} a real number $x$ if we can give a yes/no answer to the question $x<r$ for any $r$.
Keep this in mind when reading the following definition.

\begin{thm}{Definition}\label{def:GH}
The \index{Gromov--Hausdorff distance}\emph{Gromov--Hausdorff distance} $|\spc{X}-\spc{Y}|_{\GH}$ between compact metric spaces $\spc{X}$ and $\spc{Y}$
is defined by the following
relation.

Given $r > 0$, we say that $|\spc{X}-\spc{Y}|_{\GH} < r$ if and only if there exists a metric
space $\spc{W}$ and subspaces $\spc{X}'$ and $\spc{Y}'$ in $\spc{W}$ that are isometric to $\spc{X}$ and $\spc{Y}$,
respectively, such that $|\spc{X}'-\spc{Y}'|_{\Haus\spc{W}} < r$.
(Here $|\spc{X}'\z-\spc{Y}'|_{\Haus\spc{W}}$ denotes the Hausdorff distance between sets $\spc{X}'$ and $\spc{Y}'$ in $\spc{W}$.)
\end{thm}

These distances actually define a metric on (isometry classes of) compact metric spaces; for a proof we refer to \cite{burago-burago-ivanov,petrunin2023pure}.
The resulting metric is called the Gromov--Hausdorff metric;
the corresponding metric space will be denoted by \index{GH@$\GH$ (space of spaces)}$\GH$.
This means in particular that if $|\spc{X}-\spc{Y}|_{\GH}\z=0$ for compact metric spaces $\spc{X}$ and $\spc{Y}$, then they are isometric.
Moreover, each of the references above provides the following statement.

\begin{thm}{Proposition}\label{prop:complete}
$\GH$ is a complete metric space.
\end{thm}

The proofs of these statements are not quite straightforward.
The following two exercises should help the reader reconstruct them
without the cited textbooks.

\begin{thm}{Exercise}\label{ex:non-contracting-map}
Let $f$ be a distance-noncontracting map from a compact metric space $\spc{K}$ to itself.
Show that $f$ is an isometry; that is, it is a distance-preserving bijection.
\end{thm}

\begin{thm}{Exercise}\label{ex:non-expanding-map}
Show that a surjective \index{short map}\emph{short} map (that is, distance-nonexpanding) from a compact metric space to itself is an isometry.
\end{thm}

\parbf{Two convergences.}
Notice that now we have two notions of convergence of metric spaces.
One is defined in the previous section; it works for proper metric spaces and
its limit might not be uniquely defined.
The other is defined as convergence in $\GH$;
the limit in $\GH$ is unique if it exists.
In particular, these convergences are different, but if we restrict our attention to nonempty uniformly bounded compact spaces, then these two convergences are essentially the same.
Let us sketch why.

Suppose $|\spc{X}_n-\spc{X}_\infty|_{\GH}\to 0$ as $n\to \infty$;
that is, $\spc{X}_n$ converges to $\spc{X}_\infty$ in $\GH$.
Then there is a metric on $\spc{V}_n=\spc{X}_n\sqcup \spc{X}_\infty$ such that the restriction of the metric on each $\spc{X}_n$ and $\spc{X}_\infty$ coincides with its original metric, and $\dist{\spc{X}_n}{\spc{X}_\infty}{\Haus\spc{V}_n}<\eps_n$ for some sequence $\eps_n\to 0$.
Gluing all $\spc{V}_n$ along $\spc{X}_\infty$, we get the required space $\bm{X}$, which defines convergence.

Once the convergence is fixed, we can talk about limits of sequences $x_n\in \spc{X}_n$ in $\spc{X}_\infty$ as $n\to\infty$, as well as weak limits of a sequence of Borel measures $\mu_n$ on $\spc{X}_n$ and so on --- all these limits are defined as the corresponding limits in the ambient space $\bm{X}$.

Now suppose we have a Gromov--Hausdorff convergence $\spc{X}_n\to \spc{X}_\infty$ with a common space $\bm{X}$.
Assume that all spaces $\spc{X}_1,\spc{X}_2,\dots$, as well as $\spc{X}_\infty$, are compact, and $\spc{X}_\infty$ is nonempty.
(The latter condition means that the convergence is nontrivial.)
Observe that in this case $\dist{\spc{X}_n}{\spc{X}_\infty}{\Haus\bm{X}}\to 0$ as $n\to \infty$,
and therefore $\dist{\spc{X}_n}{\spc{X}_\infty}{\GH}\to 0$.
In other words, $\spc{X}_n$ converges to $\spc{X}_\infty$ in $\GH$.

\section{Almost isometries}\label{sec:Almost isometries}

\begin{thm}{Definition}
Let $\spc{X}$ and $\spc{Y}$ be metric spaces.
A map $f\:\spc{X}\to\spc{Y}$
is called an \index{isometry@$\eps$-isometry}\emph{$\eps$-isometry}
if the following two conditions hold:

\begin{subthm}{}
$f(\spc{X})$ is an \index{net@$\eps$-net}\emph{$\eps$-net} in $\spc{Y}$; that is, for any $y\in \spc{Y}$ there is $x\in \spc{X}$ such that
\[\dist{f(x)}{y}{\spc{Y}}\le\eps.\]
\end{subthm}

\begin{subthm}{}
$f(\spc{X})$ is \index{distance-preserving}\emph{$\eps$-distance-preserving}; that is,
\[\dist{f(x)}{f(x')}{\spc{Y}}\lg \dist{x}{x'}{\spc{X}}\pm\eps \]
for any $x,x'\in\spc{X}$.
\end{subthm}

\end{thm}

When dealing with Gromov--Hausdorff convergence, the following lemma allows us to bypass constructing explicit metrics on the disjoint unions of $\spc{X}_1,\spc{X}_2,\dots$, and $\spc{X}_\infty$.

\begin{thm}{Lemma}\label{lem:almost-isom}
Let $\spc{X}_1,\spc{X}_2,\dots$, and $\spc{X}_\infty$ be compact metric spaces,
and let $\eps_n\to\0+$ as $n\to\infty$.
Suppose that either 
\begin{subthm}{lem:almost-isom-a}
for each $n$ there is an $\eps_n$-isometry $f_n\:\spc{X}_n\to\spc{X}_\infty$, or
\end{subthm}
\begin{subthm}{lem:almost-isom-b}
for each $n$ there is an $\eps_n$-isometry $h_n\:\spc{X}_\infty\to\spc{X}_n$.
\end{subthm}
Then there is a Gromov--Hausdorff convergence $\spc{X}_n\to \spc{X}_\infty$.

Furthermore, a partial converse holds.

\begin{subthm}{lem:almost-isom-c}
Suppose we have a Gromov--Hausdorff convergence $\spc{X}_n\to \spc{X}_\infty$, the spaces $\spc{X}_n$ are uniformly bounded, and $\spc{X}_\infty$ is nonempty and compact.
Then for some sequence $\eps_n\to\0+$ as $n\to\infty$ there are $\eps_n$-isometries $f_n\:\spc{X}_n\to\spc{X}_\infty$ (and $h_n\:\spc{X}_\infty\to\spc{X}_n$)
such that $x_n\in \spc{X}_n$ converges to $x_\infty \in  \spc{X}_\infty$ with respect to the convergence $\spc{X}_n\z\to \spc{X}_\infty$ if and only if $f_n(x_n)\to x_\infty$ (respectively, $\dist{h_n(x_\infty)}{x_n}{\spc{X}_n}\to 0$) as $n\to\infty$.
\end{subthm}
\end{thm}

\parit{Proof.}
To prove part \ref{SHORT.lem:almost-isom-a}, let us construct a common space $\bm{X}$ for the spaces $\spc{X}_1,\spc{X}_2,\dots$, and $\spc{X}_\infty$
by taking the metric $\rho$ on the disjoint union $\spc{X}_\infty\sqcup\spc{X}_1\sqcup\spc{X}_2\sqcup\dots$ that is defined as follows:
\begin{align*}
\dist{x_n}{y_n}{\bm{X}}&=\dist{x_n}{y_n}{\spc{X}_n},
\\
\dist{x_\infty}{y_\infty}{\bm{X}}&=\dist{x_\infty}{y_\infty}{\spc{X}_\infty},
\\
\dist{x_n}{x_\infty}{\bm{X}}
&=\inf_{y_n\in \spc{X}_n}\set{\dist{x_n}{y_n}{\spc{X}_n}+\eps_n+\dist{x_\infty}{f_n(y_n)}{\spc{X}_\infty}}{},
\\
\dist{x_n}{x_m}{\bm{X}}
&=\inf_{y_\infty\in\spc{X}_\infty}\set{\dist{x_n}{y_\infty}{\bm{X}}+\dist{x_m}{y_\infty}{\bm{X}}}{},
\end{align*}
where we assume that $x_\infty,y_\infty\in \spc{X}_\infty$, and $x_n,y_n\in \spc{X}_n$ for each $n$. 
It remains to observe that this indeed defines a metric on $\bm{X}$, and $\spc{X}_n\to \spc{X}_\infty$ in the sense of Hausdorff.

The proof of the second part is analogous; one only needs to change one line in the definition of the metric as follows:
\[\dist{x_n}{x_\infty}{\bm{X}}
=
\inf_{y_\infty\in \spc{X}_\infty}\set{\dist{x_n}{h_n(y_\infty)}{\spc{X}_n}+\eps_n+\dist{x_\infty}{y_\infty}{\spc{X}_\infty}}{}.\]

We leave part \ref{SHORT.lem:almost-isom-c} as an exercise.
\qedsf

For two metric spaces $\spc{X}$ and $\spc{Y}$,
we write $\spc{X}\le \spc{Y}+\eps$ if
there is a map $f\:\spc{X}\to \spc{Y}$ such that
\[\dist{x}{x'}{\spc{X}}\le \dist{f(x)}{f(x')}{\spc{Y}}+\eps\]
for any $x,x'\in \spc{X}$.

\begin{thm}{Exercise}\label{ex:GH-po}
Let $\spc{X}_1,\spc{X}_2,\dots,$ and $\spc{X}_\infty$ be nonempty compact metric spaces.
Show that there is a Gromov--Hausdorff convergence $\spc{X}_n\z\to\spc{X}_\infty$ if and only if for some sequence $\eps_n\to 0$,
we have
\[\spc{X}_\infty\le \spc{X}_n+\eps_n\quad\text{and}\quad \spc{X}_n\le \spc{X}_\infty+\eps_n.\]
\end{thm}

Lemma~\ref{lem:almost-isom} has a natural analogue for pointed convergence.
For simplicity, we only state part \ref{SHORT.lem:almost-isom-a} of the lemma.
Parts \ref{SHORT.lem:almost-isom-b} and \ref{SHORT.lem:almost-isom-c} can be rephrased similarly.

The proof of the following lemma is analogous to \ref{lem:almost-isom}.

\begin{thm}{Lemma}\label{lem:almost-isom-pointed}
Let $(\spc{X}_1, p_1),(\spc{X}_2,p_2),\dots$, and $(\spc{X}_\infty, p_\infty)$ be pointed metric spaces, and let $\eps(n,R)\to\0+$ as $n\to\infty$ for any fixed $R>0$.
Suppose that for each $n$ there is a map $f_n\:\spc{X}_n\to\spc{X}_\infty$ such that

\begin{subthm}{lem:almost-isom-pointed-basepoint}
$f_n(p_n)\to p_\infty$;
\end{subthm}

\begin{subthm}{lem:almost-isom-pointed-b}
$\bigl|\dist{f_n(x)}{f_n(x')}{\spc{X}_\infty}-\dist{x}{x'}{\spc{X}_n}\bigr|\le \eps(n,R)$ for any $x,x'\z\in \oBall(p_n,R)$;
\end{subthm}

\begin{subthm}{lem:almost-isom-pointed-c}
For any $x \in \oBall(p_\infty,R)$, there is $x_n\in \oBall(p_n,R)$ such that $\dist{x}{f_n(x_n)}{}\le \eps(n,R)$.
\end{subthm}
Then there is a pointed Gromov--Hausdorff convergence $(\spc{X}_n,p_n)\z\to (\spc{X}_\infty,p_\infty)$.
\end{thm}

\section{Remarks}

In principle, our definition of Gromov--Hausdorff distance can be applied to complete metric spaces that are not necessarily compact.
However, according to the following exercise, it only defines a semimetric; that is, zero Gromov--Hausdorff distance does not imply that the spaces are isometric.
For that reason it is not in use.

\begin{thm}{Exercise}\label{ex:GH-noncompact}
\begin{subthm}{ex:GH-noncompact:proper}
Construct two nonisometric proper (noncompact) metric spaces with vanishing Gromov--Hausdorff distance.
\end{subthm}

\begin{subthm}{ex:GH-noncompact:bounded}

Construct two nonisometric complete geodesic metric spaces of bounded diameter with vanishing Gromov--Hausdorff distance.
\end{subthm}

\end{thm}

%% file: definitions.tex
\chapter{Definitions}\label{chap:defs}

Here we will give several definitions of Alexandrov spaces and prove their equivalence.

\section{Four-point comparison}\label{sec:4-point}

Recall that $\angk pxy$ denotes the model angle; see \ref{sec:Triangles and model triangles}.

Let $p,x,y,z$ be a quadruple of points in a metric space.
If the inequality 
\[\angk pxy_{\EE^2}+\angk pyz_{\EE^2}+\angk pzx_{\EE^2}
\le 
2\cdot\pi
\eqlbl{eq:CBB-comparison}\]
holds, then we say that the quadruple meets \index{comparison}\emph{$\EE^2$-comparison}.
If the left-hand side is undefined, then we also assume that the comparison holds.

\begin{thm}{Exercise}\label{ex:CBB+-}
Suppose $\EE^2$-comparison holds for the quadruple $p$, $x_1$, $x_2$, $x_3$.
Show that it also holds for any quadruple $q$, $y_1$, $y_2$, $y_3$ such that
\[\dist{q}{y_i}{}\ge\dist{p}{x_i}{}\qquad\text{and}\qquad\dist{y_i}{y_j}{}\le\dist{x_i}{x_j}{}\]
for all $i$ and $j$.
\end{thm}

Instead of $\EE^2$, we can use $\SSS^2$ or $\HH^2$.
This way we get the definition of $\SSS^2$- or $\HH^2$-comparisons.
Recall that $\angk pxy_{\EE^2}$ and $\angk pxy_{\HH^2}$ are defined if $p\ne x$, $p\ne y$,
but for $\angk pxy_{\SSS^2}$ we require in addition that
\[\dist{p}{x}{}+\dist{p}{y}{}+\dist{x}{y}{}<2\cdot\pi;\]
if this does not hold for one of the angles, then we assume that $\SSS^2$-comparison holds for this quadruple.

More generally, one may apply this definition to the model $\kappa$-plane $\MM^2(\kappa)$ and define $\MM^2(\kappa)$-comparison for any real $\kappa$.
However, when considering $\MM^2(\kappa)$-comparison, it is safe to assume that $\kappa=-1$, $0$, or $1$;
applying rescaling, the $\MM^2(\kappa)$-comparison can be reduced to these three cases.

\begin{thm}{Definition}\label{def:CBB}
A metric space $\spc{X}$ has {}\emph{curvature $\ge\kappa$} in the sense of Alexandrov
if $\MM^2(\kappa)$-comparison
holds for any quadruple of points in $\spc{X}$.
\end{thm}

While this definition can be applied to any metric space,
we will use it mostly for geodesic spaces that are complete (and often compact or proper). 
If a complete geodesic space has curvature $\ge\kappa$ in the sense of Alexandrov, 
then it will be called an $\Alex\kappa$ space; here $\Alex\kappa$ is an adjective.
If $\spc{X}$ is $\Alex\kappa$ for some $\kappa$, then we say that $\spc{X}$ is an \index{Alexandrov space}\emph{Alexandrov space}.

It is common practice in Alexandrov geometry to write proofs for nonnegative curvature and 
leave the general curvature bound as an exercise. These generalizations are usually straightforward. We will add notes when they are not.
We will also often formulate statements just for $\kappa=0$ even when they admit straightforward generalizations to arbitrary curvature bounds;
see \cite{alexander-kapovitch-petrunin2024} for a more formal treatment.

\begin{thm}{Exercise}\label{ex:Euclid-is-CBB}
Show that $\EE^n$ is $\Alex0$.
\end{thm}

\begin{thm}{Exercise}\label{ex:(3+1)-expanding}
Show that a metric space $\spc{X}$ has nonnegative curvature in the sense of Alexandrov
if and only if, for any quadruple of points $p,x_1,x_2,x_3\in \spc{X}$,
there is a quadruple of points $q,y_1,y_2,y_3\in\EE^2$
such that 
\[\dist{p}{x_i}{\spc{X}}\ge\dist{q}{y_i}{\EE^2} 
\quad \text{and}\quad
\dist{x_i}{x_j}{\spc{X}}\le\dist{y_i}{y_j}{\EE^2}\] 
for all $i$ and $j$.
\end{thm}

\section{Alexandrov's lemma}

Recall that $[xy]$ denotes a geodesic from $x$ to $y$;
set
\index{10@$\left]x y\right]$, $\left[x y\right[$, $\left]x y\right[$ (subintervals)}
\[
\left]x y\right]=[xy]\setminus\{x\},
\quad
\left[x y\right[=[xy]\setminus\{y\},
\quad
\left]x y\right[=[xy]\setminus\{x,y\}.\]

\begin{thm}{Alexandrov's lemma}
\index{Alexandrov's lemma}
\label{lem:alex}
Let $p,x,y,z$ be distinct points in a metric space such that $z\in \left]x y\right[$.
Then 
the following expressions have the same sign:

\begin{subthm}{lem-alex-difference}
$\angk x p y
-\angk x p z$,
\end{subthm} 

\begin{subthm}{lem-alex-angle}
$\angk z p x
+\angk z p y -\pi$.
\end{subthm}

\begin{wrapfigure}{r}{25mm}
\vskip-0mm
\centering
\includegraphics{mppics/pic-730}
\end{wrapfigure}

The same holds for the hyperbolic and spherical model angles, 
but in the latter case, one has to assume in addition that
\[\dist{p}{x}{}+\dist{p}{y}{}+\dist{x}{y}{}< 2\cdot\pi.\]

\end{thm}

In the proof we will apply the following statement from elementary geometry, which is an immediate consequence of the cosine law:

\begin{thm}{Observation}\label{angle-monotonicity}
Increasing one side in a plane triangle increases its opposite angle, and the other way around.
Moreover, the same statement holds for spherical and hyperbolic triangles.
\end{thm}

\parit{Proof of \ref{lem:alex}.}
Consider the model triangle $\trig{\tilde x}{\tilde p}{\tilde z}=\modtrig(x p z)$.
Take 
a point $\tilde y$ on the extension of 
$[\tilde x \tilde z]$ beyond $\tilde z$ so that $\dist{\tilde x}{\tilde y}{}=\dist{x}{y}{}$ (and therefore $\dist{\tilde z}{\tilde y}{}=\dist{z}{y}{}$).

\begin{wrapfigure}{r}{33mm}
\vskip-0mm
\centering
\includegraphics{mppics/pic-740}
\end{wrapfigure}

By \ref{angle-monotonicity},
the following expressions have the same sign:
\begin{enumerate}[(i)]
\item $\mangle\hinge{\tilde x}{\tilde p}{\tilde y}-\angk{x}{p}{y}$,
\item $\dist{\tilde p}{\tilde y}{}-\dist{p}{y}{}$,
\item $\mangle\hinge{\tilde z}{\tilde p}{\tilde y}-\angk{z}{p}{y}$.
\end{enumerate}
Since 
\[\mangle\hinge{\tilde x}{\tilde p}{\tilde y}=\mangle\hinge{\tilde x}{\tilde p}{\tilde z}=\angk{x}{p}{z}\]
and
\[ \mangle\hinge{\tilde z}{\tilde p}{\tilde y}
=\pi-\mangle\hinge{\tilde z}{\tilde x}{\tilde p}
=\pi-\angk{z}{x}{p},\]
the statement follows.

The spherical and hyperbolic cases can be proved along the same lines.
\qeds

\begin{thm}[!]{Exercise}\label{ex:alex-lemma-cat}
Assume that $p,x,y,z$ are as in Alexandrov's lemma (\ref{lem:alex}).
Show that
\[\angk p x y
\ge
\angk p x z + \angk p z y,\]
with equality if and only if the expressions in \ref{SHORT.lem-alex-difference} and \ref{SHORT.lem-alex-angle} in Alexandrov's lemma vanish.
\end{thm}

Note that 
\[p\in\left]x y\right[
\quad\Rightarrow\quad
\angk pxy=\pi.\]
Applying this together with Alexandrov's lemma and $\EE^2$-comparison, we get the following.

\begin{thm}{Claim}\label{clm:angle-mono}
Let $p,x,y,z$ be points in an $\Alex0$ space.
Suppose that $z\in\left]x y\right[$.
Then
\[\angk xyp\le \angk xzp.\]
\end{thm}

\begin{wrapfigure}{r}{25mm}
\vskip-0mm
\centering
\includegraphics{mppics/pic-750}
\end{wrapfigure}

\begin{thm}[!]{Exercise}\label{ex:noncreasing}
Let $\hinge p x y$ be a hinge in an $\Alex0$ space.
Consider the function
\[f\:(\dist{p}{\bar x}{},\dist{p}{\bar y}{})\mapsto \angk p{\bar x}{\bar y},\]
where $\bar x\in\left]p x\right]$ and $\bar y\in\left]p y\right]$.
Show that $f$ is nonincreasing in each argument.
\end{thm}

The statement in the last exercise is called \index{angle-sidelength monotonicity}\emph{angle-sidelength monotonicity}.
It implies the following.

\begin{thm}{Angle existence and comparison}\label{clm:angle-defined}
The angle measure of any hinge in an $\Alex0$ 
space is defined and is at least as large as the corresponding model angle;
that is,
\[\mangle\hinge p x y\ge \angk p x y\]
for any hinge $\hinge p x y$ in an $\Alex0$ space.

\end{thm}

\begin{thm}[!]{Exercise}\label{ex:0-angle}
Let $\hinge p x y$ be a hinge in an $\Alex0$ space.
Suppose $\mangle\hinge p x y=0$; show that $[px]\subset [py]$ or $[py]\subset [px]$.

Conclude that geodesics in an $\Alex0$ space cannot \index{bifurcating geodesics}\emph{bifurcate};
that is, if two geodesics $[px]$ and $[py]$ share a nontrivial arc with an endpoint at $p$, then $[px]\subset [py]$ or $[py]\subset [px]$.
\end{thm}

\begin{thm}{Exercise}\label{ex:pi-angle}
Let $[xy]$ be a geodesic in an $\Alex0$ space.
Suppose $z\in \left]xy\right[$.
Show that there is a unique geodesic $[xz]$ and that $[xz]\subset [xy]$.
\end{thm}

Recall that adjacent hinges are defined in \ref{ex:adjacent-angles}.

\begin{thm}[!]{Exercise}\label{ex:adjacent-CBB}
Let $\hinge pxz$ and $\hinge pyz$ be adjacent hinges in an $\Alex0$ 
space.
Show that
\[\mangle\hinge pxz+\mangle\hinge pyz= \pi.\]
\end{thm}

\begin{thm}{Exercise}\label{ex:pxyvw}
Let $\spc{A}$ be an $\Alex0$ 
space.
Show that
\[
\angk xyp=\angk xvp
\quad\Leftrightarrow\quad
\angk xyp=\angk xwp
\]
for any points
$p,x,y,v,w$ in $\spc{A}$ such that $v,w\in \left]xy\right[$.
\end{thm}

\begin{thm}{Exercise (semicontinuity of angles)}\label{ex:angle-lim}
Let $\spc{A}$ be an $\Alex0$ space.
Suppose hinges $\hinge {x_n}{y_n}{z_n}$ in $\spc{A}$ converge to a hinge $\hinge {x_\infty}{y_\infty}{z_\infty}$;
that is, geodesics $[x_ny_n]$ and $[x_nz_n]$ converge to the geodesics $[x_\infty y_\infty]$ and $[x_\infty z_\infty]$ in the sense of Hausdorff.
Show that 
\[\liminf_{n\to\infty}\mangle \hinge {x_n}{y_n}{z_n}\ge \mangle \hinge {x_\infty}{y_\infty}{z_\infty}.\]
\end{thm}

The last inequality might be strict;
this happens on the surface of a cube in $\EE^3$, which is an $\Alex0$ space by \ref{prop:conv-surf-CBB(0)}.

\section{Hinge comparison}

Let $\hinge pxy$ be a hinge in an $\Alex0$ space $\spc{A}$.
By \ref{ex:noncreasing}, the angle measure $\mangle\hinge pxy$ is defined and
\[\mangle\hinge pxy\ge \angk pxy.\]
Further, according to \ref{ex:adjacent-CBB}, we have 
\[\mangle\hinge pxz+\mangle\hinge pyz=\pi\]
for adjacent hinges $\hinge pxz$ and $\hinge pyz$ in $\spc{A}$.

The following theorem adds a converse.

\begin{thm}{Theorem}\label{thm:angle-cbb}
A complete geodesic space $\spc{A}$ is $\Alex0$ if and only if both of the following conditions hold.

\begin{subthm}{angle-a}
For any hinge $\hinge x p y$ in $\spc{A}$, the angle 
$\mangle\hinge x p y$ is defined and 
\[\mangle\hinge x p y\ge\angk x p y.\]
\end{subthm}

\begin{subthm}{angle-b}
For any two adjacent hinges $\hinge pxz$ and $\hinge pyz$ in $\spc{A}$, we have
\[\mangle\hinge pxz+\mangle\hinge pyz\le\pi.\]
\end{subthm}

\end{thm}

Recall that, by \ref{ex:adjacent-CBB}, we have equality in \ref{SHORT.angle-b}.

\parit{Proof.}
Recall that the only-if part follows from \ref{ex:noncreasing} and \ref{ex:adjacent-CBB}; it remains to show the if part.

Consider a point $w\in \mathopen{]} p z \mathclose{[}$ close to $p$.
From \ref{SHORT.angle-b}, it follows that 
\[\mangle\hinge w x z+ \mangle\hinge w x{p}\le\pi\quad \text{and}\quad \mangle\hinge w y z + \mangle\hinge w y{p}\le\pi.\]

\begin{wrapfigure}{o}{30 mm}
\vskip-0mm
\centering
\includegraphics{mppics/pic-805}
\vskip4mm
\end{wrapfigure}

Since $\mangle\hinge w x y\le \mangle\hinge w x p +\mangle\hinge w y{p}$ (see \ref{claim:angle-3angle-inq}), we get 
\[\mangle\hinge w x z+ \mangle\hinge w y z +\mangle\hinge w x y
\le
2\cdot\pi.\]
Applying \ref{SHORT.angle-a}, 
\[\angk w x z
+ \angk w y z 
+\angk w x y
\le
2\cdot\pi.\]
Passing to the limit as $w\to p$, we have
\[\angk p x z 
+ \angk p y z 
+\angk p x y
\le
2\cdot\pi.\]
\qedsf

\section{Equivalent conditions}

The following theorem summarizes \ref{clm:angle-mono}, \ref{clm:angle-defined}, \ref{ex:adjacent-CBB}, and \ref{thm:angle-cbb}.

\begin{thm}{Theorem}\label{thm:defs_of_alex} 
Let $\spc{A}$ be a complete geodesic space.
Then the following conditions are equivalent.

\begin{subthm}{cbb}
$\spc{A}$ is $\Alex0$.
\end{subthm}

\begin{subthm}{2-sum} 
(adjacent angle comparison\index{comparison!adjacent angle comparison})
\[\angk z p x
+\angk z p y\le \pi\]
for any geodesic $[x y]$ and points $z\in \mathopen{]}x y\mathclose{[}$, $p\ne z$ in $\spc{A}$.
\end{subthm}

\begin{subthm}{point-on-side}
(\index{comparison!point-on-side comparison}point-on-side comparison)
\[\angk x p y\le\angk x p z\]
for any geodesic $[x y]$ and $z\in \mathopen{]}x y\mathclose{[}$ in $\spc{A}$.
\end{subthm}

\begin{subthm}{angle}(hinge comparison\index{comparison!hinge comparison})
\index{hinge!comparison}
the angle $\mangle\hinge x p y$ is defined for any hinge $\hinge x p y$ in $\spc{A}$.
Moreover, 
\[\mangle\hinge x p y\ge\angk x p y\]
for any hinge $\hinge x p y$, and
\[\mangle\hinge z p y + \mangle\hinge z p x\le\pi\]
for any adjacent hinges $\hinge z p y$ and $\hinge z p x$.
\end{subthm}

Moreover, the implications \ref{SHORT.cbb}$\Rightarrow$\ref{SHORT.2-sum}$\Rightarrow$\ref{SHORT.point-on-side}$\Rightarrow$\ref{SHORT.angle} hold in any space, not necessarily a geodesic one.
\end{thm}

\begin{thm}{Advanced exercise}\label{ex:urysohn}
Construct a complete geodesic space $\spc{X}$ that is not $\Alex0$, but satisfies the following weaker version of the adjacent angle comparison (\ref{2-sum}).

For any three points $p,x,y\in \spc{X}$, there is a geodesic $[x y]$ such that for any $z\in \left]x y\right[$
\[\angk{z}{p}{x}+\angk{z}{p}{y}
\le
\pi.\]
\end{thm}

\begin{thm}{Exercise}\label{ex:normCBB}
Let $\spc{W}$ be $\RR^n$ with the metric induced by a norm.
Show that if $\spc{W}$ is $\Alex0$, then $\spc{W}$ is isometric to the Euclidean space~$\EE^n$.
\end{thm}

\section{Function comparison}\label{Function comparison}

\parbf{Real-to-real functions.}
Choose $\lambda\in \RR$.
Let $s\:\II\to\RR$ be a locally Lipschitz function defined on an interval $\II\subset\RR$.
The following statements are equivalent;
if one (and therefore any) of them holds for $s$, then we say that $s$ is \index{concave function}\emph{$\lambda$-concave}.
\begin{itemize}
\item We have the inequality $s''\le \lambda$, where the second derivative $s''$ is understood in the sense of distributions.
\item The function $t\mapsto s(t)-\lambda\cdot\tfrac{t^2}2$ is concave.
\item The \index{Jensen inequality}\emph{Jensen inequality}
\[s(a\cdot t_0+(1-a)\cdot t_1)\ge a\cdot s(t_0)+(1-a)\cdot s(t_1)+\tfrac\lambda2\cdot a\cdot(1-a)\cdot(t_1-t_0)^2 \]
holds for any $t_0,t_1\in \II$ and $a\in[0,1]$.
\item For any $t_0\in \II$ there is a quadratic polynomial $\ell=\tfrac\lambda2\cdot t^2+a\cdot t+b$ (it is called a \index{barrier}\emph{barrier}) that supports $s$ at $t_0$ from above;
that is, $\ell(t_0)\z= s(t_0)$ and $\ell(t)\ge s(t)$ for any $t$.
\item For any $t_0\in \II$ there is a quadratic polynomial $\ell=\tfrac\lambda2\cdot t^2+a\cdot t\z+b$ (it is called a \index{local barrier}\emph{local barrier}) that locally supports $s$ at $t_0$ from above;
that is, $\ell(t_0)\z= s(t_0)$ and $\ell(t)\ge s(t)$ for any $t$ in a neighborhood of $t_0$.
\end{itemize}

To prove the equivalence, approximate $s$ by smooth functions, taking convolutions $s_n=s*k_n$ for a suitable sequence of kernels $k_n$.
Observe that all the conditions are equivalent for $s_n$;
passing to the limit, we get the same for $s$.

\begin{thm}{Exercise}\label{ex:concave'}
Show that $\lambda$-concave functions are {}\emph{one-sided differentiable}.
That is, if the function $f\:\II\to\RR$ is $\lambda$-concave, then
\begin{subthm}{}
the \index{$f^+$, $f^-$ (one-sided derivative)}\index{right derivative}\emph{right derivative}
\[f^+(t_0)=\lim_{t\to t_0^+} \frac{f(t)-f(t_0)}{t-t_0}\]
is defined for any $t_0\in \II$ except at the right endpoint of $\II$.
\end{subthm}
\begin{subthm}{}
the \index{left derivative}\emph{left derivative}
\[f^-(t_0)=\lim_{t\to t_0^-} \frac{f(t)-f(t_0)}{t-t_0}\]
is defined for any $t_0\in \II$ except at the left endpoint of $\II$.
\end{subthm}

Furthermore, show that $f^+(t_0)-f^-(t_0)\le 0$ for any $t_0$ in the interior of $\II$.
\end{thm}

The following exercise implies that if the function is defined on an open interval, then the Lipschitz condition can be dropped from the definition of $\lambda$-concavity.

\begin{thm}{Exercise}\label{ex:concave-open}
Suppose a real-to-real function $f$ is defined on an open interval and, for some $\lambda\in\RR$, it satisfies the Jensen inequality stated above.
Show that $f$ is locally Lipschitz.
\end{thm}

\parbf{Functions on metric spaces.}
A \index{function}\emph{function} on a metric space $\spc{A}$ will usually mean a \textit{locally Lipschitz real-valued function defined on an open subset of $\spc{A}$}.
The domain of a function $f$ will be denoted by $\Dom f$.
Thus, when we write $f\:\spc{A}\to \RR$ we do not assume that $f$ is defined on the whole space $\spc{A}$.

We say that such $f$ is \index{91@$\lambda$-concave function}\emph{$\lambda$-concave} (briefly, $f''\le \lambda$) if
for any geodesic $\gamma\:\II\z\to \Dom f$
the real-to-real function $t\mapsto f\circ\gamma(t)$ is $\lambda$-concave.

The following proposition is simple but conceptual ---
it reduces global comparison to an infinitesimal condition on distance functions.

\begin{thm}{Proposition}\label{comp-kappa}
A complete geodesic space $\spc{A}$ is $\Alex0$ if and only if $f''\le 1$ for any function $f$ of the form
\[f\:x\mapsto \tfrac12\cdot\dist[2]{p}{x}{}.\] 
\end{thm} 

\parit{Proof.}
Choose a geodesic $\gamma$ in $\spc{A}$ and two points $x=\gamma(t_0)$, $y=\gamma(t_1)$ for some $t_0<t_1$.
Consider the model triangle $\trig{\tilde p}{\tilde x}{\tilde y}\z=\modtrig(p x y)$.
Let $\tilde \gamma\:[t_0,t_1]\to\EE^2$ be the unit-speed parametrization of $[\tilde x \tilde y]$ from $\tilde x$ to $\tilde y$.

Set
\begin{align*} 
\tilde r(t)&\df\dist{\tilde p}{\tilde\gamma(t)}{},
& 
r(t)&\df\dist{p}{\gamma(t)}{}.
\end{align*}
Clearly, $\tilde r(t_0)=r(t_0)$ and $\tilde r(t_1)=r(t_1)$.
By the point-on-side comparison (\ref{point-on-side}), the implication
\[t_0\le t\le t_1
\qquad\Rightarrow\qquad
\tilde r(t)\le r(t)
\eqlbl{eq:r=<r}\]
holds for any $\gamma$ and $t_0<t_1$.

Observe that $(\tilde r^2)''\equiv2$.
If $(r^2)''\le 2$, then $t\mapsto r^2(t)-\tilde r^2(t)$ is concave.
Since $\tilde r(t_0)=r(t_0)$ and $\tilde r(t_1)=r(t_1)$, we get \ref{eq:r=<r}, which proves the if part.

On the other hand, \ref{eq:r=<r} implies the Jensen inequality for the function $h\:t\mapsto r(t)^2-t^2$.
Since the subinterval $[t_0,t_1]$ can be chosen arbitrarily, we conclude that $h''\le 0$, or, equivalently, $(r^2)''\le 2$, and the only-if part follows.
\qeds

\parbf{Semiconcave functions.}\label{sec:Semiconcave functions}
Let $f$ be a locally Lipschitz real-valued function defined on an open subset $\Dom f$ of an Alexandrov space $\spc{A}$.
Suppose $\phi$ is a continuous function defined on $\Dom f$.
We will write $f''\le \phi$ if for any point $x\in \Dom f$ and any $\eps>0$ there is a neighborhood $U\ni x$ such that the restriction $f|_U$ is $(\phi(x)+\eps)$-concave.

If $f''\le \phi$ for some continuous function $\phi$, then $f$ is called \index{semiconcave function}\emph{semiconcave}.

\begin{thm}{Exercise}\label{ex:distfun-semiconcave}
Let $f$ be a \index{distance function}\emph{distance function} on an $\Alex0$ space $\spc{A}$;
that is, $f(x)\equiv\dist{p}{x}{}$ for some $p\in \spc{A}$.
Show that $f''\le \tfrac1f$.
In particular, $f$ is semiconcave in $\spc{A}\setminus\{p\}$.
\end{thm}

Proposition~\ref{comp-kappa} admits the following generalization, which can be proved along the same lines, but the formulas get longer.

\begin{thm}{Proposition}
A complete geodesic space $\spc{A}$ is $\Alex{\mp1}$
if and only if $f''\z\le \pm f$ for any function of the type $f=\cosh\circ\distfun_p$ (respectively, $f=-\cos\circ\distfun_p|_{\oBall(p,\pi)}$).
\end{thm}

The geometric meaning of these inequalities remains the same:
\textit{distance functions on an $\Alex{\kappa}$ space are more concave than distance functions in $\MM^2(\kappa)$}.

\section{Remarks}

Alexandrov's lemma has the following useful variation; see \cite[10.2]{alexander-kapovitch-petrunin2024} or \cite[3.3]{alexander-kapovitch-kirszbraun}.

\begin{wrapfigure}{o}{30 mm}
\vskip-4mm
\centering
\includegraphics{mppics/pic-210}
\vskip0mm
\end{wrapfigure}

\begin{thm}{Overlap lemma}\label{lem:extend-overlap}
Let $\tilde x^1$, $\tilde x^2$, $\tilde x^3$, $\tilde p^1$, $\tilde p^2$, and $\tilde p^3$ be points in $\EE^2$, $\SSS^2$, or $\HH^2$.
Assume that, for any permutation $(i,j,k)$ of $(1,2,3)$, we have
\begin{enumerate}[(i)]

\item
\label{no-overlap:px=px}
$\dist{\tilde p^i}{\tilde x^k}{}=\dist{\tilde p^j}{\tilde x^k}{}$,

\item
\label{no-overlap:orient-1}
$\tilde p^i$ and $\tilde x^i$ lie in the same closed half-space determined by $[\tilde x^j\tilde x^k]$,
\item no pair of triangles $\trig{\tilde p^i}{\tilde x^j}{\tilde x^k}$ overlaps.
\end{enumerate}
Then
\[\mangle{\tilde p^1} +\mangle {\tilde p^2}+\mangle{\tilde p^3}> 2\cdot\pi,\]
where $\mangle\tilde p^i\df\mangle\hinge{\tilde p^i}{\tilde x^k}{\tilde x^j}$
for a permutation $(i,j,k)$ of $(1,2,3)$.
\end{thm}

The condition \ref{SHORT.angle-b} in \ref{thm:angle-cbb} might not be necessary.
This is a long-standing open problem possibly dating back to Alexandrov \cite[footnote in 4.1.5]{burago-burago-ivanov}.
Let us state it formally.

\begin{thm}{Open question}\label{open:hinge-}
Let $\spc{A}$ be a complete geodesic space (you can also assume that $\spc{A}$ is homeomorphic to $\mathbb{S}^2$ or $\RR^2$)
such that for any hinge $\hinge x p y$ in $\spc{A}$,
the angle $\mangle\hinge x p y$ is defined and
\[\mangle\hinge x p y\ge\angk x p y.\]
Is it true that $\spc{A}$ is an Alexandrov space?
\end{thm}

One could define Alexandrov spaces as complete \textit{length} spaces with curvature $\ge \kappa$ (in the sense of Alexandrov).
This condition is more natural and general, but many statements can be reduced to the geodesic case.
In particular, any complete length space $\spc{A}$ with curvature $\ge \kappa$ \textit{can be isometrically embedded into an $\Alex\kappa$ space} --- the ultrapower of $\spc{A}$; see \cite[4.11+8.4]{alexander-kapovitch-petrunin2024}.
Also, by Plaut's theorem, any point $p$ in $\spc{A}$ can be connected by geodesics to \textit{most} points in $\spc{A}$; see
\cite[8.11]{alexander-kapovitch-petrunin2024} and compare to \ref{ex:grad-dist:geod}.

Our 4-point comparison in Section~\ref{sec:4-point} is closely related to the so-called $\CAT$ comparison, which defines an \textit{upper} curvature bound in the sense of Alexandrov;
this is the subject of our previous book \cite{alexander-kapovitch-petrunin-2019}.
In both comparisons, we check certain conditions on the 6 distances between pairs of points in a 4-point set.
Such conditions for $n$-point subsets were considered by Michael Gromov \cite[Section 1.19$_+$]{gromov1999}, and they have been actively studied since then \cite{toyoda,lebedeva-petrunin-zolotov,lebedeva2019,petrunin2017,lebedeva-petrunin2024,lebedeva-petrunin2023,lebedeva-petrunin2021,lebedeva-petrunin2025,eskenazis-mendel-naor,gromov2001}.

%% file: globalization.tex
\chapter{Globalization}\label{chap:globalization}

The globalization theorem states that every locally Alexandrov space is globally Alexandrov.
We prove it for compact spaces and indicate the proof of the general case.

\section{Globalization}

A complete geodesic metric space $\spc{A}$ is \index{locally $\Alex0$}\emph{locally $\Alex0$} if any point $p\in\spc{A}$ admits a neighborhood $U\ni p$ such that the $\EE^2$-comparison holds for any quadruple of points in $U$.

It is straightforward to obtain a local version of Theorem~\ref{thm:defs_of_alex} for locally $\Alex0$ spaces;
it gives a list of equivalent properties that hold in sufficiently small neighborhoods of any point.
In particular, the analog of \ref{angle} implies that \textit{the angle measure of any hinge in a locally $\Alex0$ space is well defined.}

\begin{thm}{Globalization theorem}\label{thm:glob}\index{globalization theorem}
Any compact locally $\Alex0$ space is $\Alex0$.
\end{thm}

\parit{Proof of \ref{thm:glob} modulo the key lemma (\ref{key-lem:globalization}).}
Note that condition \ref{angle-b} holds in $\spc{A}$ (the necessity of \ref{angle-b} uses only the local condition).
It remains to check \ref{angle-a};
that is,
\[\mangle\hinge x p q\ge\angk x p q
\eqlbl{eq:mod-angle-CBB-comp-glob}\]
for any hinge $\hinge x p q$ in $\spc{A}$.

Inequality \ref{eq:mod-angle-CBB-comp-glob} holds for hinges in a small neighborhood of any point;
this can be proved the same way as \ref{clm:angle-defined} and \ref{ex:adjacent-CBB}, applying the local version of the $\EE^2$-comparison.
Since $\spc{A}$ is compact, there is $\eps>0$ such that \ref{eq:mod-angle-CBB-comp-glob} holds whenever $\dist{x}{p}{}+\dist{x}{q}{}<\eps$.
Applying the key lemma finitely many times, we get that \ref{eq:mod-angle-CBB-comp-glob} holds for any given hinge.
\qeds

\begin{thm}{Key lemma}\label{key-lem:globalization} 
Let $\spc{A}$ be a locally $\Alex0$ space.
If the comparison
\[\mangle\hinge x p q
\ge\angk x p q\]
holds for any hinge $\hinge x p q$ such that
$\dist{x}{p}{}+\dist{x}{q}{}
<
\frac{2}{3}\cdot\ell$,
then it holds for any hinge $\hinge x p q$ such that
$\dist{x}{p}{}+\dist{x}{q}{}
<
\ell$.
\end{thm}

Let $\hinge x p q$ be a hinge in $\spc{A}$.
Denote by $\side \hinge x p q$ its \index{21@$\side \hinge x p q$ (model side)}\index{model!side}\emph{model side};
this is the opposite side in the plane triangle with the same angle and adjacent side lengths as in $\hinge x p q$.

\begin{wrapfigure}{r}{44mm}
\centering
\includegraphics{mppics/pic-105}
\end{wrapfigure}

More precisely,
\begin{align*}
\side \hinge x p q_{\spc{A}}
&\df
\dist{\tilde p}{\tilde q}{\EE^2}
\intertext{if $\hinge {\tilde x} {\tilde p} {\tilde q}_{\EE^2}$ is the \index{model hinge}\emph{model hinge} of $\hinge {x} {p} {q}_{\spc{A}}$; that is,}
\mangle\hinge {\tilde x} {\tilde p} {\tilde q}_{\EE^2}&=\mangle\hinge x p q_{\spc{A}},
\\
\dist{\tilde x} {\tilde p}{\EE^2}&=\dist{x} {p}{\spc{A}},
\\
\dist{\tilde x} {\tilde q}{\EE^2}&=\dist{x} {q}{\spc{A}}.
\end{align*}

Note that
\[\side \hinge x p q \ge\dist{p}{q}{}
\quad\Leftrightarrow\quad
\mangle\hinge x p q\ge \angk x p q.
\]

\parit{Proof.} 
Consider a hinge $\hinge x p q$ such that 
\[\tfrac{2}{3}\cdot\ell \le\dist{p}{x}{}\z+\dist{x}{q}{}< \ell. \eqlbl{eq:thm:=def-loc-two-sided*}\]
It is sufficient to prove that
\[\side \hinge x p q
\ge\dist{p}{q}{}\eqlbl{eq:thm:=def-loc*}\]

First, let us describe a construction of a new hinge $\hinge{x'}p q$ with several properties described below;
in particular,
\[
\dist{p}{x}{}+\dist{x}{q}{}\ge\dist{p}{x'}{}+\dist{x'}{q}{},
\eqlbl{eq:thm:=def-loc-fourstar}\]
and
\[\side \hinge x p q
\ge\side \hinge{x'}p q.
\eqlbl{eq:thm:=def-loc-fivestar}\]

\parit{Construction.}
Assume $\dist{x}{q}{}\ge\dist{x}{p}{}$; otherwise, switch the roles of $p$ and $q$.
Take $x'\in [x q]$ such that 
\[\dist{p}{x}{}+3\cdot\dist{x}{x'}{}
=\tfrac{2}{3}\cdot\ell.
\eqlbl{3|xx'|}\]
Choose a geodesic $[x' p]$ and consider the hinge $\hinge{x'}p q$ formed by $[x' p]$ and $[x' q]\subset [x q]$.
The triangle inequality implies \ref{eq:thm:=def-loc-fourstar}.
Furthermore, 
\begin{align*}
\dist{p}{x}{}\z+\dist{x}{x'}{}&<\tfrac{2}{3}\cdot\ell,
&
\dist{p}{x'}{}\z+\dist{x'}{x}{}&<\tfrac{2}{3}\cdot\ell.
\end{align*}
In particular, the assumption of the lemma implies that
\[\mangle\hinge x p{x'}
\ge\angk x p{x'}
\quad \text{and}\quad 
\mangle\hinge {x'}p x
\ge\angk {x'}p x.
\eqlbl{eq:thm:=def-loc-threestar}\]

{

\begin{wrapfigure}{r}{30 mm}
\vskip-8mm
\centering
\includegraphics{mppics/pic-820}
\vskip-4mm
\end{wrapfigure}

Let
$\trig{\tilde x}{\tilde x'}{\tilde p}=\modtrig(x x' p)$.
Take $\tilde q$ on the extension of $[\tilde x\tilde x']$ beyond $\tilde x'$ such that $\dist{\tilde x}{\tilde q}{}\z=\dist{x}{q}{}$ (and therefore $\dist{\tilde x'}{\tilde q}{}=\dist{x'}{q}{}$).
By~\ref{eq:thm:=def-loc-threestar},
\[\mangle\hinge x p q
=\mangle\hinge x p{x'}\ge\angk x p{x'}
\ \Rightarrow\
\side \hinge x q p\ge\dist{\tilde p}{\tilde q}{}.\]
Hence
\begin{align*}
\mangle\hinge{\tilde x'}{\tilde p}{\tilde q}&= 
\pi
-\angk{x'}p x
\ge
\\
&\ge
\pi-\mangle\hinge{x'}p x
=
\\
&=
\mangle\hinge{x'}p q,
\end{align*}
and \ref{eq:thm:=def-loc-fivestar} follows.

}

\medskip

Let us continue the proof.
Set $x_0=x$.
Applying the above construction inductively, we obtain a sequence of hinges $\hinge{x_n}p q$ with $x_{n+1}\z=x_n'$.
From \ref{eq:thm:=def-loc-fivestar}, we see that the sequence $s_n\z=\side \hinge{x_n}p q$ is nonincreasing.
\begin{figure}[ht!]
\centering
\includegraphics{mppics/pic-825}
\end{figure}

The process terminates at $x_n$ if
$\dist{p}{x_n}{}+\dist{x_n}{q}{}
\z<\tfrac{2}{3}\cdot\ell$.
In this case, by the assumptions of the lemma, $\side \hinge{x_n}p q\z\ge\dist{p}{q}{}$.
Since the sequence $s_n$ is nonincreasing, inequality \ref{eq:thm:=def-loc*} follows.

From now on, we assume that the process does not terminate.
Let us prove the following claim.
\begin{clm}{}\label{clm:>0}
The distances
$\dist{x_n}{x_{n+1}}{}$, $\dist{x_n}{p}{}$, and $\dist{x_n}{q}{}$ are bounded away from zero for all large $n$.
\end{clm}

Set
\begin{align*}
a_n&=\min\{\dist{p}{x_n}{},\dist{q}{x_n}{}\},
\\
b_n&=\max\{\dist{p}{x_n}{},\dist{q}{x_n}{}\},
\\
r_n&=\dist{p}{x_n}{}+\dist{x_n}{q}{}=a_n+b_n.
\end{align*}
By the triangle inequality, $r_n$ is a nonincreasing sequence.
Since the process does not terminate, we have $\tfrac{2}{3}\cdot\ell\le r_n<\ell$ for all $n$.

By \ref{3|xx'|}, $\dist{x_n}{x_{n+1}}{}=\tfrac13\cdot(\tfrac23\cdot\ell-a_n)$.
Since $a_n+b_n=r_n<\ell$ and $a_n\le b_n$, we have $a_n\le \tfrac12\cdot\ell$.
Hence
\[
\tfrac29\cdot\ell\ge \dist{x_n}{x_{n+1}}{}\ge \tfrac1{18}\cdot\ell,
\eqlbl{eq:|xx'|>=l/18}
\]
which proves the claim for $\dist{x_n}{x_{n+1}}{}$.

Note that
\[a_{n+1}=b_n-\dist{x_n}{x_{n+1}}{}\qquad \text{or}\qquad b_{n+1}=b_n-\dist{x_n}{x_{n+1}}{}.\]
In the latter case,
\[
\begin{split}
a_{n+1}&=r_{n+1}-b_{n+1}=
\\
&=(r_{n+1}-r_n)+(r_n-b_n)+\dist{x_n}{x_{n+1}}{}
\\
&=-(r_n-r_{n+1})+a_n+\dist{x_n}{x_{n+1}}{}.
\end{split}
\]
Since $r_n$ does not increase, it must converge;
so $r_n-r_{n+1}\z\to 0$ as $n\to\infty$.
Therefore, $r_n-r_{n+1}\z\le \tfrac1{100}\cdot \ell$ for all large $n$.
Observe that $b_n\ge\tfrac12 \cdot r_n$, and therefore $b_n\ge \tfrac13\cdot\ell$.
Taking \ref{eq:|xx'|>=l/18} into account, we get that in both cases $a_{n+1}\ge\tfrac1{100}\cdot \ell$ for all large $n$,
which finishes the proof of the claim.
\claimqeds

Since $r_n-r_{n+1}\to 0$, \ref{clm:>0} implies that $\angk{x_n}{p_n}{x_{n+1}}\to \pi$,
where $p_n\z=p$ if $x_{n+1}\in [x_nq]$, and otherwise $p_n=q$ (this describes switching the roles of $p$ and $q$ in the construction of $x'$).
Since $\mangle\hinge{x_n}{p_n}{x_{n+1}}\z\ge\angk{x_n}{p_n}{x_{n+1}}$, we have
$\mangle\hinge{x_n}{p_n}{x_{n+1}}\to \pi$ as $n\to\infty$.

It follows that
\[r_n-s_n=\dist{p}{x_n}{}+\dist{x_n}{q}{}-\side \hinge{x_n}p q\to 0.\] 
By the triangle inequality,
$r_n\ge\dist{p}{q}{}$;
therefore,
\[\lim_{n\to\infty}s_n\ge \dist{p}{q}{}.\]
Finally, the monotonicity of the sequence $s_n=\side \hinge{x_n}p q$ implies \ref{eq:thm:=def-loc*}.
\qeds

\section{On the general case}

\begin{thm}{Theorem}\label{thm:globalization+}\index{globalization theorem}
Any locally $\Alex\kappa$ space is $\Alex\kappa$.
\end{thm}

This is a more general version of our globalization theorem~\ref{thm:glob}.
It holds for an arbitrary curvature bound~$\kappa$ and replaces compactness with completeness.
A proof can be found in~\cite[2F]{alexander-kapovitch-petrunin2024}; it is based on the idea of \ref{thm:glob}, but requires a few additional twists that we are about to discuss.

\parbf{Trading compactness for completeness.}
Assume for a moment that $\kappa=0$.
To replace compactness by completeness, we need the following statement.

\begin{thm}{Exercise}\label{ex:alm-min}
Let $\spc{X}$ be a complete metric space.
Suppose $r\:\spc{X}\to \RR$ is a positive continuous function.
Show that for any $\eps>0$, there exists a point $p\in \spc{X}$ such that
\[
r(x)> (1-\eps)\cdot r(p)
\]
for every $x\in \cBall[p,\tfrac{1}{\eps}\cdot r(p)]$.

\end{thm}

Define $r(x)$ to be the maximal radius of a ball centered at $x$ in which the comparison condition holds (this definition should be made precise, but the details are not important here).
Note that if $r(x)=\infty$ for one (and therefore for every) point $x$, then the theorem follows immediately.
To prove the theorem, we need to choose a point $p$ satisfying the conclusion of the exercise and apply the key lemma.

\parbf{On general curvature bounds.}
The case $\kappa\le 0$ in \ref{thm:globalization+} can be proved in the same way as the case $\kappa=0$, but the case $\kappa>0$ requires extra care.
This is because model triangles may be undefined in $\MM^2(\kappa)$ when $\kappa>0$.%
\footnote{Recall that, according to the definitions in Section~\ref{sec:Triangles and model triangles},
if $\angk p q r_{\MM^2(\kappa)}$, and hence $\modtrig(p q r)_{\MM^2(\kappa)}$, are defined for some $\kappa>0$, then
$\dist{p}{q}{}+\dist{q}{r}{}+\dist{r}{p}{}< 2\cdot\pi/\sqrt{\kappa}$.}

\parbf{Meaning of ``\,curvature $\bm{\ge\kappa}$\,''.}
Recall that an $\Alex\kappa$ space is defined as a complete geodesic space having ``curvature $\ge\kappa$ in the sense of Alexandrov.''
This suggests that the following claim should hold.

\begin{thm}{Claim}\label{clm:K>k}
If $\Kappa>\kappa$, then any $\Alex\Kappa$ space is $\Alex\kappa$.
\end{thm}

The proof of this statement uses \ref{thm:globalization+}.

\parit{Proof.}
By~\ref{ex:k-><mono},
\[
\angk x y z_{\MM^2(\kappa)}\le \angk x y z_{\MM^2(\Kappa)}
\]
whenever both sides are defined.
This completes the proof in the case $\Kappa\le 0$.
In general, however, the comparison is guaranteed only for sufficiently small triangles;
in particular, any $\Alex\Kappa$ space is \index{locally $\Alex\kappa$}\emph{locally} $\Alex\kappa$.
Therefore, \ref{thm:globalization+} completes the proof.
\qeds

While the expression ``curvature bounded below by $\kappa$'' makes sense for geodesic spaces,
it does not make much sense for general metric spaces.
For example, if $\spc{X}$ is the set $\{p,x_1,x_2,x_3\}$ with the metric defined by
$\dist{p}{x_i}{}=\pi$ and $\dist{x_i}{x_j}{}=2\cdot\pi$ for all $i\ne j$, then $\spc{X}$ has no defined spherical model angles.
Therefore, formally speaking, $\spc{X}$ has curvature $\ge 1$.
On the other hand, $\spc{X}$ does not have curvature $\ge 0$, since
\[
\angk p{x_1}{x_2}_{\EE^2}
=\angk p{x_2}{x_3}_{\EE^2}
=\angk p{x_1}{x_3}_{\EE^2}
=\pi.
\]
Use this example as a warning when working with a non-intrinsic metric.

\begin{thm}[!]{Exercise}\label{ex:RisCBB(1)}
Let $p$ and $q$ be points in an $\Alex1$ space $\spc{A}$.
Suppose $\dist{p}{q}{}>\pi$, and let $m$ be the midpoint of $[pq]$.
Show that for any hinge $\hinge mxp$, we have either
$\mangle\hinge mxp=0$ or $\mangle\hinge mxp=\pi$.

Conclude that $\spc{A}$ is isometric to a line segment or a circle.
\end{thm}

\begin{thm}{Exercise}\label{ex:perim-k>0}
Suppose
$\spc{A}$ is $\Alex1$
and $\diam \spc{A}\le \pi$.
Show that 
\[\dist{x}{y}{}+\dist{y}{z}{}+\dist{z}{x}{}\le 2\cdot\pi\]
for any triple of points $x,y,z\in \spc{A}$.
\end{thm}

\section{Remarks}

The globalization theorem is also known as the \textit{generalized Toponogov theorem}.
Its $2$-dimensional case was proved by Paolo Pizzetti \cite{pizzetti} and
rediscovered by Alexandr Alexandrov \cite{alexandrov:devel}.
Victor Toponogov \cite{toponogov-globalization+splitting} proved it for Riemannian manifolds of all dimensions.
For Alexandrov spaces of all dimensions, the theorem first appears in the paper by Michael Gromov, Yuriy Burago, and Grigory Perelman \cite{burago-gromov-perelman}.
They prove globalization for \textit{complete length spaces}.
Another version for \textit{noncomplete but geodesic spaces} is due to the third author \cite{petrunin:globalization}.

We took the proof from our book \cite{alexander-kapovitch-petrunin2024}, but present it in a less general form.
This proof is based on simplifications obtained by Conrad Plaut \cite{plaut:dimension} and Dmitri Burago, Yurii Burago, and Sergei Ivanov \cite{burago-burago-ivanov}.
The same proof was discovered independently by Urs Lang and Viktor Schroeder \cite{lang-schroeder:globalization}.
Another simplified argument was found by Katsuhiro Shiohama \cite{shiohama}.
Alternative proofs were found by Shengqi Hu, Xiaole Su, and Yusheng Wang \cite{hu-su-wang2023,hu-su-wang2024}.

%% file: calculus.tex

\chapter{Calculus}\label{chap:derivative}

This lecture discusses the first-order derivatives in Alexandrov spaces, including spaces of directions, tangent spaces, differentials, and gradients.

\section{Space of directions} 
\label{sec:space+directions}

Let $\spc{A}$ be an Alexandrov space.
By \ref{ex:noncreasing}, the angle measure of any hinge in $\spc{A}$ is defined.
Given $p\in \spc{A}$, consider the set $\mathfrak{S}_p$ of all nontrivial geodesics starting at $p$.
By \ref{claim:angle-3angle-inq}, the triangle inequality holds for $\mangle$ on $\mathfrak{S}_p$.
Thus, $(\mathfrak{S}_p,\mangle)$ forms a \index{semimetric}\emph{semimetric} space:
$\mangle$ behaves like a metric but may vanish for distinct geodesics.

The metric space corresponding to $(\mathfrak{S}_p,\mangle)$ is called the \index{70@$\Sigma_p'$ (geodesic directions)}\index{space of geodesic directions}\emph{space of geodesic directions} at $p$; it will be denoted by $\Sigma'_p$ or $\Sigma'_p\spc{A}$.
The elements of $\Sigma'_p$ are called \index{geodesic!direction}\emph{geodesic directions} at $p$.
Each geodesic direction is an equivalence class of geodesics starting at $p$ under the equivalence relation
\[[px]\sim[py]
\quad\Leftrightarrow\quad
\mangle\hinge pxy=0;\]
the direction of $[px]$ is denoted by $\dir px$.\index{40@$\dir{p}{q}$ (direction)}
By \ref{ex:0-angle}, 
\[[px]\sim[py]
\quad\Leftrightarrow\quad
(\,[px]\subset [py]
\quad\text{or}\quad
[px]\supset[py]\,).
\]
 
The completion of $\Sigma'_p$ is called the \index{space of directions}\emph{space of directions} at $p$ and is denoted by \index{70@$\Sigma_p$ (space of directions)}$\Sigma_p$ or $\Sigma_p\spc{A}$.
The elements of $\Sigma_p$ are called \index{direction}\emph{directions} at $p$.

\begin{thm}{Exercise}\label{ex:dir-compact}
Let $\spc{A}$ be an Alexandrov space.
Assume that a sequence of geodesics $[px_n]$ converges to a geodesic $[px_\infty]$ in the sense of Hausdorff,
and $x_\infty\ne p$.
Suppose $\Sigma_p$ is compact.
Show that $\dir p{x_n}\z\to\dir p{x_\infty}$ as $n\to\infty$.

\end{thm}

\section{Tangent space}\label{sec: tangent space}

The \index{Cone@$\Cone$ (Euclidean cone)}\emph{Euclidean cone} $\spc{V}=\Cone\spc{X}$
over a metric space $\spc{X}$
is defined as the metric space whose underlying set is the quotient of $[0,\infty)\times\spc{X}$ by the equivalence relation generated by $(0,p)\sim (0,q)$ for all points $p,q\in\spc{X}$;
its metric is defined by the cosine rule:
\[
\dist{(s,p)}{(t,q)}{\spc{V}} 
=
\sqrt{s^2+t^2-2\cdot s\cdot t\cdot \cos\theta},
\]
where $\theta= \min\{\pi, \dist{p}{q}{\spc{X}}\}$.

Note that
\[\Cone\SSS^n\iso\EE^{n+1};\]
here ``$\iso$'' stands for ``is isometric to''.

Next we define multiplication by scalar, scalar product, norm, and angles on Euclidean cones, extending the usual definitions from Euclidean space.

The point $v$ in $\spc{V}$ that corresponds to $(t,x)\z\in[0,\infty)\times \spc{X}$ will be denoted by \index{33@$r\cdot u$ (multiplication in cone)}$t\cdot x$;
the space $\spc{X}$ might be identified with the unit sphere in $\spc{V}$; so $x=1\cdot x$.
More generally, for $v,w\in \spc{V}$, we will write $w=s\cdot v$ if $v$ corresponds to $(t,x)\z\in[0,\infty)\times \spc{X}$ and $w$ corresponds to $(s\cdot t,x)\z\in[0,\infty)\times \spc{X}$.
The point in $\spc{V}$ formed by the equivalence class of $\{\0\}\times\spc{X}$ is called the \index{origin}\emph{origin} of the cone and is denoted by $\0$ or $\0_{\spc{V}}$.
For $v\in\spc{V}$ the distance $\dist{\0}{v}{\spc{V}}$ is called the \index{norm}\emph{norm} of $v$ and is denoted by $|v|$ or $|v|_{\spc{V}}$.
The \index{angle measure}\emph{angle measure} \index{31@$\mangle$ (angle measure)} $\mangle(v,w)$ and the \index{32@$\<v,w\>$ (scalar product)} \index{scalar product}\emph{scalar product} $\<v,w\>$
of $v=s\cdot p$ and $w=t\cdot q$
are defined by 
\[
\mangle(v,w)\df\min\{\pi, \dist{p}{q}{\spc{X}}\}
\quad\text{and}\quad
\<v,w\>
\df |v|\cdot|w|\cdot\cos(\mangle(v,w)).
\]
If $v=\0$ or $w=\0$, then $\mangle(v,w)$ is undefined and $\<v,w\>\df0$.

\begin{thm}{Exercise}\label{ex:geodesic-cone}
Show that $\Cone\spc{X}$ is geodesic if and only if $\spc{X}$ is \index{geodesic space@$\ell$-geodesic space}\emph{$\pi$-geodesic};
that is, any two points $x,y\in \spc{X}$ such that $\dist{x}{y}{\spc{X}}<\pi$ can be joined by a geodesic in $\spc{X}$.
\end{thm}

\parbf{Tangent space.}
The Euclidean cone $\Cone\Sigma_p$ over the space of directions $\Sigma_p$ is called the \index{tangent space}\emph{tangent space} at $p$ and is denoted by \index{70@$\T_p$ (tangent space)}$\T_p$ or $\T_p\spc{A}$.
The elements of $\T_p\spc{A}$ will be called \index{tangent vector}\emph{tangent vectors} at $p$
(despite the fact that $\T_p$ is not a vector space).
The space of directions $\Sigma_p$ can be (and will be) identified with the unit sphere in~$\T_p$;
that is, with the set $\set{v\in\T_p}{|v|=1}$.

\begin{thm}{Proposition}\label{prop:Tan-is-CBB(0)}
Any tangent space to an Alexandrov space has nonnegative curvature in the sense of Alexandrov.
\end{thm}

Halbeisen's example \cite[13.6]{alexander-kapovitch-petrunin2024} shows that the tangent space $\T_p$ at some point of an Alexandrov space might fail to be geodesic;
in this case $\T_p$ is \textit{not} $\Alex0$.

\parit{Proof.}
Consider the tangent space $\T_p=\Cone \Sigma_p$ of an Alexandrov space $\spc{A}$ at a point $p$.
We need to show that the $\EE^2$-comparison holds for a given quadruple $v_0$, $v_1$, $v_2$, $v_3\in \T_p$.

Recall that the space of geodesic directions $\Sigma_p'$ is dense in $\Sigma_p$.
It follows that the subcone $\T'_p=\Cone\Sigma_p'$ is dense in $\T_p$.
Therefore, we can assume that $v_0$, $v_1$, $v_2$, $v_3\in \T'_p$.

For each $i$, choose a geodesic $\gamma_i$ from $p$ in the direction of $v_i$;
reparametrize each $\gamma_i$ so that it has speed $|v_i|$.
Since the angles are defined, we have
\[\dist{\gamma_i(\eps)}{\gamma_j(\eps)}{\spc{A}}=\eps\cdot\dist{v_i}{v_j}{\T_p}+o(\eps)
\eqlbl{eq:gamma-v}\]
for $\eps>0$.
The quadruple $\gamma_0(\eps)$, $\gamma_1(\eps)$, $\gamma_2(\eps)$, $\gamma_3(\eps)$ meets the $\MM^2(\kappa)$-comparison.
After rescaling all the distances by $\tfrac1\eps$, it becomes the $\MM^2(\eps^2\cdot\kappa)$-comparison.
Passing to the limit as $\eps\to 0$ and applying \ref{eq:gamma-v}, we get the $\EE^2$-comparison for $v_0$, $v_1$, $v_2$, $v_3$.
\qeds

\begin{thm}[!]{Exercise}\label{ex:GHto-tangent}
Let $p$ be a point in an Alexandrov space $\spc{A}$,
and let $\lambda_n\to\infty$.
Suppose $\Sigma_p$ is compact.
Show that there is a pointed Gromov--Hausdorff convergence $(\lambda_n\cdot \spc{A},p)\z\to (\T_p,0)$.
Moreover, the convergence can be chosen so that for any geodesic $\gamma$ that starts at $p$, we have
\[\iota_n\circ\gamma(t/\lambda_n)\to t\cdot \xi,\]
where $\iota_n$ sends a point in $\spc{A}$ to the corresponding point in $\lambda_n\cdot\spc{A}$ and $\xi\in\Sigma_p$ is the direction of $\gamma$ at $0$.
\end{thm}

\section{Differential}\label{sec:differential}
\index{differential of a function}

Let $f$ be a semiconcave function defined in a neighborhood of a point $p$ in an Alexandrov space $\spc{A}$ (recall that semiconcave functions are locally Lipschitz by definition).
Choose a geodesic $\gamma$ that starts at $p$;
let $\xi\in\Sigma_p$ be its direction.
Define 
\[(\dd_pf)(\xi)\df(f\circ\gamma)^+(0),\]
where $(f\circ\gamma)^+$ denotes the right derivative of $(f\circ\gamma)$;
it is defined since $f$ is semiconcave.

By the following exercise, $\dd_pf$ is a Lipschitz function on $\Sigma'_p$.
It follows that the function $\dd_pf\:\Sigma_p'\to\RR$ can be uniquely extended to a Lipschitz function $\dd_pf\:\Sigma_p\to\RR$.
Furthermore, we can extend it to the tangent space by setting
\[(\dd_pf)(r\cdot \xi)
\df
r\cdot (\dd_pf)(\xi)\]
for any $r\ge 0$ and $\xi\in\Sigma_p$.
The resulting function $\dd_pf\:\T_p\to\RR$ is Lipschitz;
it is called the \index{differential}\emph{differential} of $f$ at $p$.

\begin{thm}{Exercise}\label{ex:df(xi)}
Let $f$ be a semiconcave function on an Alexandrov space.
Suppose $\gamma_1$ and $\gamma_2$ are geodesics that start at $p\z\in \Dom f$;
denote by $\theta$ the angle between $\gamma_1$ and $\gamma_2$ at $p$.
Show that 
\[|(f\circ\gamma_1)^+(0)-(f\circ\gamma_2)^+(0)|\le L\cdot \theta,\]
where $L$ is the Lipschitz constant of $f$ in a neighborhood of $p$.
\end{thm}

\begin{thm}[!]{Exercise (First variation formula)} \label{ex:d(distfun)}
Let $p$ and $q$ be distinct points in an Alexandrov space~$\spc{A}$.

\begin{subthm}{ex:d(distfun):<}
Show that
$(\dd_p\distfun_q)(v)\le -\langle\dir pq,v\rangle$
for any $v\in\T_p$.
\end{subthm}

\begin{subthm}{ex:d(distfun):=}
Suppose $\spc{A}$ is proper.
Let $\Uparrow_p^q$ be the set of all directions of geodesics from $p$ to $q$.
Show that
\[(\dd_p\distfun_q)(v)=-\max_{\xi\in\Uparrow_p^q}\langle\xi,v\rangle\]
for any $v\in\T_p$.
\end{subthm}

\end{thm}

\section{Gradient}\label{sec:grad-def}

The following definition generalizes the notion of gradient of a smooth function to semiconcave functions on Alexandrov spaces.
This generalization is not trivial even for concave functions on Euclidean space.
This point is worth keeping in mind while reading further.

\begin{thm}{Definition}\label{def:grad} 
Let $f$ be a semiconcave function on an Alexandrov space and $p\in\Dom f$.
A tangent vector $g\in \T_p$ is called a 
\index{gradient}\emph{gradient} of $f$ at $p$ 
(briefly, $g\z=\nabla_p f$\index{19@$\nabla$ (gradient)}) if
\begin{subthm}{}
$(\dd_p f)(w)\le \<g,w\>$ for any $w\in \T_p$, and
\end{subthm}

\begin{subthm}{}
$(\dd_p f)(g) = \<g,g\>.$
\end{subthm}
\end{thm}

The following exercise provides a property of gradients that will play a key role in the proof of the first distance estimate (\ref{thm:dist-est}).

\begin{figure}[hbt!]
\centering
\includegraphics{mppics/pic-409}
\end{figure}

\begin{thm}[!]{Exercise}\label{ex:monotonicity}
Let $f$ be a $\lambda$-concave function on an Alexandrov space.
Suppose that gradients $\nabla_xf$ and $\nabla_yf$ are defined.
Show that 
\[\<\dir{x}{y},\nabla_{x}f\>
+
\<\dir{y}{x},\nabla_{y}f\>
+
\lambda\cdot\dist{x}{y}{}\ge 0.\]
\end{thm}

\begin{thm}{Proposition}\label{prop:grad-exist}
Suppose that a semiconcave function $f$ is defined in a neighborhood of a point $p$ in an Alexandrov space.
Then the gradient $\nabla_pf$ is uniquely defined.

Moreover, if $\dd_pf\le 0$, then we have $\nabla_pf=0$;
otherwise, $\nabla_pf\z=s\cdot \overline{\xi}$, where 
$s= \dd_pf(\overline{\xi})$
and
$\overline{\xi}\in \Sigma_p$ is the direction that maximizes the value $\dd_pf(\xi)$ for $\xi\in \Sigma_p$.
\end{thm}

\begin{thm}{Key lemma}\label{lem:ohta} 
Let $f$ be a semiconcave function that is defined in a neighborhood of a point $p$
in an Alexandrov space $\spc{A}$. 
Assume
\[s=\sup\set{(\dd_p f)(\xi)}{\xi\in\Sigma_p}\ge0.\]
Then for any $u,v\in \T_p$, we have
\[s\cdot \sqrt{|u|^2+2\cdot\<u,v\> +|v|^2}
\ge 
(\dd_p f)(u)+(\dd_p f)(v).\]

\end{thm}

Suppose that $\T_p\iso\EE^m$ and the differential $\dd_p f$ is a concave function.
Let $w=\tfrac{u+v}2$.
Then
\[2\cdot |w|=\sqrt{|u|^2+2\cdot\<u,v\> +|v|^2}
\eqlbl{eq:2|w|=|u+v|}\]
and
\[2\cdot(\dd_p f)(w)\ge(\dd_p f)(u)+(\dd_p f)(v),
\eqlbl{eq:2df(w)>=df(u)+df(v)}\]
which implies the lemma.

In general, $\T_p$ is not geodesic (and not even a length space), so concavity of $\dd_p f$ does not make sense.
The key lemma, however, says that in a certain sense $\dd_p f$ behaves as a concave function.
The proof mimics the calculations above inside $\spc{A}$.
Instead of vectors $(u$, $v$, $w)$ we choose a triple of points $(u_t$, $v_t$, $w_t)$ in $\spc{A}$
and prove the approximate versions
of \ref{eq:2|w|=|u+v|} and \ref{eq:2df(w)>=df(u)+df(v)};
see \ref{eq:2|wt|=|ut+vt|} and \ref{eq:2df(wt)>=df(ut)+df(vt)}, respectively.
The main difficulty can be seen in the following exercise; technically it lies in the proof of \ref{eq:2|wt|=|ut+vt|} below.

\begin{thm}[!]{Exercise}\label{ex:d(distfun):==}
Let $p$ and $q$ be distinct points in an Alexandrov space~$\spc{A}$.
Suppose the geodesic $[pq]$ can be extended beyond $q$.

Show that
\[(\dd_p\distfun_q)(v)= -\langle\dir pq,v\rangle\]
for any $v\in\T_p$.
\end{thm}

\parit{Proof of \ref{lem:ohta}.}
We will assume that $\spc{A}$ is $\Alex0$ and $f$ is concave;
the general case requires only minor modifications.
Furthermore, we can assume that $u\ne 0$, $v\ne 0$, and $\alpha=\mangle(u,v)>0$; otherwise, the statement is trivial.

{

\begin{wrapfigure}{r}{34 mm}
\vskip4mm
\centering
\includegraphics{mppics/pic-1205}
\vskip0mm
\end{wrapfigure}

Consider a model configuration of five points: $\tilde p$, $\tilde u$, $\tilde v$, $\tilde q$, $\tilde w\in\EE^2$ such that
\begin{itemize}
\item $\mangle\hinge{\tilde p}{\tilde u}{\tilde v}=\alpha$, 
\item $\dist{\tilde p}{\tilde u}{}=|u|$, 
\item $\dist{\tilde p}{\tilde v}{}=|v|$,
\item $\tilde q$ lies on an extension of $[\tilde p\tilde v]$ so that $\tilde v$ is the midpoint of $[\tilde p\tilde q]$,
\end{itemize}
}
\begin{itemize}
\item $\tilde w$ is the midpoint between $\tilde u$ and ${\tilde v}$.
\end{itemize}
Note that 
\[\dist{\tilde p}{\tilde w}{}
=
\tfrac{1}{2}\cdot\sqrt{|u|^2+2\cdot\<u,v\>+|v|^2}.\eqlbl{eq:|p-w|=}\]

Since the space of geodesic directions $\Sigma'_p$ is dense in $\Sigma_p$,
we can assume that there are geodesics in the directions of $u$ and $v$.
Choose such geodesics $\gamma_u$ and $\gamma_v$
and assume that they are parametrized with speeds $|u|$ and $|v|$, respectively.
For all small $t>0$, consider points $u_t,v_t,q_t,w_t\in \spc{A}$ such that
\begin{itemize}
\item $v_t=\gamma_v(t)$ and $q_t=\gamma_v(2\cdot t)$,
\item $u_t=\gamma_u(t)$,
\item $w_t$ is the midpoint of $[u_t v_t]$.
\end{itemize}
Clearly 
\[\dist{p}{u_t}{}=t\cdot |u|,\qquad \dist{p}{v_t}{}=t\cdot|v|,\qquad \dist{p}{q_t}{}=2\cdot t\cdot|v|.\] 
Since $\mangle(u,v)$ is defined, 
we have 
\[\dist{u_t}{v_t}{}=t\cdot\dist{\tilde u}{\tilde v}{}+o(t),
\qquad
\dist{u_t}{q_t}{}=t\cdot\dist{\tilde u}{\tilde q}{}+o(t).\]

From the point-on-side comparison (\ref{point-on-side}), we have
\[\angk{v_t}p{w_t}
\ge
\angk{v_t}p{u_t}
=
\mangle\hinge{\tilde v}{\tilde p}{\tilde u}+\tfrac{o(t)}t\]
and
\[\angk{v_t}{q_t}{w_t}
\ge
\angk{v_t}{q_t}{u_t}
=
\mangle\hinge{\tilde v}{\tilde q}{\tilde u}+\tfrac{o(t)}t.\]
Clearly, 
$\mangle\hinge{\tilde v}{\tilde p}{\tilde u}+\mangle\hinge{\tilde v}{\tilde q}{\tilde u}=\pi$. 
From the adjacent angle comparison (\ref{2-sum}), 
$\angk{v_t}p{u_t}\z+\angk{v_t}{u_t}{q_t}\le \pi$.
Hence
$\angk{v_t}p{w_t}
\to
\mangle\hinge{\tilde v}{\tilde p}{\tilde w}$ as $t\to0+$
and thus 
\[\dist{p}{w_t}{}=t\cdot\dist{\tilde p}{\tilde w}{}+o(t).
\eqlbl{eq:2|wt|=|ut+vt|}\]

Without loss of generality, we can assume that $f(p)=0$.
Since $f$ is concave, we have 
\[
\begin{aligned}
2\cdot f(w_t)&\ge f(u_t)+f(v_t)=
\\
&=t\cdot [(\dd_p f)(u)+(\dd_p f)(v)]+o(t).
\end{aligned}
\eqlbl{eq:2df(wt)>=df(ut)+df(vt)}
\]
 
Applying concavity of $f$, we have
\begin{align*}
(\dd_p f)(\dir p{w_t})
&\ge 
\frac{f(w_t)}{\dist{p}{w_t}{}}
\ge 
\\
&\ge
\frac{t\cdot[(\dd_p f)(u)+(\dd_p f)(v)]+o(t)}{2\cdot t\cdot\dist{\tilde p}{\tilde w}{}+o(t)}.
\end{align*}
By \ref{eq:|p-w|=}, the key lemma follows.
\qeds

\parit{Proof of \ref{prop:grad-exist}; uniqueness.}\label{uniqueness:prop:grad-exist}
If $g,g'\in \T_p$ are two gradients of $f$,
then 
\begin{align*}
\<g,g\>
&=(\dd_p f)(g)\le \<g,g'\>,
&
\<g',g'\>
&=(\dd_p f)(g')\le \<g,g'\>.
\end{align*}
Therefore,
\[\dist[2]{g}{g'}{}=\<g,g\>-2\cdot\<g,g'\>+\<g',g'\>\le0.\] 
It follows that $g=g'$.

\parit{Existence.} 
If $\dd_p f\le 0$, then one can take $\nabla_p f=\0$.

Suppose $s=\sup\set{(\dd_p f)(\xi)}{\xi\in\Sigma_p}>0$.
It is sufficient to show that there is $\overline{\xi}\in\Sigma_p$ such that
\[
(\dd_p f)\left(\overline{\xi}\right)=s.
\eqlbl{overlinexi}
\]
Indeed, suppose $\overline{\xi}$ exists.
Applying \ref{lem:ohta} for $u=\overline{\xi}$, $v=\eps\cdot w$ with $\eps\to0+$, 
we get
\[(\dd_p f)(w)\le \<w,s\cdot\overline{\xi}\>\] 
for any $w\in\T_p$;
that is, $s\cdot\overline{\xi}$ is the gradient at $p$.

Existence of $\overline{\xi}$ is obvious if $\Sigma_p$ is compact.
In the general case, take a sequence of directions $\xi_n\in \Sigma_p$, such that $(\dd_p f)(\xi_n)\to s$.
Applying \ref{lem:ohta} for $u=\xi_n$ and $v=\xi_m$, we get
\[s
\ge
\frac{(\dd_p f)(\xi_n)+(\dd_p f)(\xi_m)}{\sqrt{2+2\cdot\cos\mangle(\xi_n,\xi_m)}}.\]
Therefore $\mangle(\xi_n,\xi_m)\to0$ as $n,m\to\infty$;
that is, $\xi_1,\xi_2,\dots$ is a Cauchy sequence.
Clearly, $\overline{\xi}=\lim_n\xi_n$ meets \ref{overlinexi}.
\qeds

\begin{thm}[!]{Exercise}\label{ex:convergence-grad}
Let $f$ and $g$ be semiconcave functions defined in a neighborhood of a point $p$ in an Alexandrov space.
Show that 
\[\dist[2]{\nabla_p f}{\nabla_p g}{\T_p}
\le 
s\cdot(|\nabla_p f|+|\nabla_p g|),\]
where
\[s=\sup\set{|(\dd_p f)(\xi)-(\dd_p g)(\xi)|}{\xi\in\Sigma_p}.\]

Conclude that if the sequence of restrictions $\dd_p f_n|_{\Sigma_p}$ converges uniformly, then $\nabla_pf_n$ converges as $n\to\infty$.
Here we assume that all functions $f_1$, $f_2,\dots$ are semiconcave near $p$.
\end{thm}

\begin{thm}[!]{Exercise}\label{ex:semicontinuous-grad}
Let $f$ be a locally Lipschitz $\lambda$-concave function on an Alexandrov space $\spc{A}$.

\begin{subthm}{ex:semicontinuous-grad:>s}
Suppose $s\ge 0$.
Show that $|\nabla_xf|> s$ if and only if for some point $y$ we have
\[f(y)-f(x)>s\cdot \ell+\lambda\cdot \tfrac{\ell^2}2,\]
where $\ell=\dist{x}{y}{}$.
\end{subthm}

\begin{subthm}{ex:semicontinuous-grad:lim} Show that $x\mapsto|\nabla_xf|$ is lower semicontinuous;
that is,
\[|\nabla_{x_\infty}f|\le \liminf_{x_n\to x_\infty} |\nabla_{x_n}f|.\]

\end{subthm}

\end{thm}

%% file: gradient-flow.tex
\chapter{Gradient flow}\label{chap:GF}

Here we define the gradient flows of semiconcave functions on Alexandrov spaces and prove related distance estimates.

\section{Velocity of a curve}\label{Velocity of curve}

Let $\alpha$ be a curve in an Alexandrov space $\spc{A}$, and let $p=\alpha(t_0)$.
If for any choice of 
geodesics $[p\,\alpha(t_0+\eps)]$ the vectors 
\[\tfrac{1}{\eps}\cdot\dist{p}{\alpha(t_0+\eps)}{}\cdot\dir p{\alpha(t_0+\eps)}\]
converge as $\eps\to 0+$, then their limit in $\T_p$ is called the \index{right derivative}\emph{right derivative} of $\alpha$ at $t_0$; it will be denoted by $\alpha^+(t_0)$.
In addition, we set $\alpha^+(t_0)\df0$
if $\tfrac{1}{\eps}\cdot\dist{p}{\alpha(t_0+\eps)}{}\to 0$ as $\eps\to 0+$.

The tangent vector $v=\dist px{}\cdot\dir px$ will be called a \index{logarithm}\emph{logarithm} of $x$ at $p$ (briefly, \index{40@$\log_p x$ (logarithm)}$v=\log_p x$).
We set $\log_p p=0$.
In general the logarithm is a multivalued function from $\spc{A}$ to $\T_p$; so $v=\log_p x$ and $w=\log_p x$ do \textit{not} imply $v=w$.
Note that $\gamma^+(0)=\log_px$ for any geodesic path $\gamma$ from $p$ to $x$.\label{page:log}

\begin{thm}{Claim}\label{clm:fa'=dfa'}
Let $\alpha$ be a curve in an Alexandrov space $\spc{A}$.
Suppose $f$ is a semiconcave function defined in a neighborhood of $p\z=\alpha(0)$, and $\alpha^+(0)$ is defined.
Then $(f\circ\alpha)^+(0)$ exists and 
\[(f\circ\alpha)^+(0)
=
(\dd_pf)(\alpha^+(0)).\]

\end{thm}

\parit{Proof.}
Without loss of generality, we can assume that $f(p)=0$.
Suppose $f$ and therefore $\dd_pf$ are $L$-Lipschitz.

Let $\gamma$ be a constant-speed parametrization of a geodesic that starts at $p$, and
such that the distance
$s=\dist{\alpha^+(0)}{\gamma^+(0)}{\T_p}$
is small.
By the definition of differential,
\[(f\circ\gamma)^+(0)=\dd_pf(\gamma^+(0)).\]

By comparison and the definition of $\alpha^+$,
\[\dist{\alpha(\eps)}{\gamma(\eps)}{\spc{A}}\le s\cdot\eps+o(\eps)\]
for $\eps>0$.
Therefore,
\[|f\circ\alpha(\eps)-f\circ\gamma(\eps)|\le L\cdot s\cdot\eps+o(\eps).\]

Suppose $(f\circ\alpha)^+(0)$ is defined.
Then
\[|(f\circ\alpha)^+(0)-(f\circ\gamma)^+(0)|\le L\cdot s.\]
Since $\dd_pf$ is $L$-Lipschitz, we also get 
\[|\dd_pf(\alpha^+(0))-\dd_pf(\gamma^+(0))|\le L\cdot s.\]
It follows that the needed identity holds up to an error $2\cdot L\cdot s$.
The statement follows since $s>0$ can be chosen arbitrarily.

The same argument is applicable if in place of $(f\circ\alpha)^+(0)$
we use any limit of $\tfrac1{\eps_n}\cdot [f\circ\alpha(\eps_n)-f(p)]$ for a sequence $\eps_n\to 0+$.
It proves that all such limits coincide; in particular, $(f\circ\alpha)^+(0)$ is defined and equals $(\dd_pf)(\alpha^+(0))$.
\qeds

\section{Gradient curves}

\begin{thm}{Definition}\label{def:grad-curve}
Let $f$ be a semiconcave function on an Alexandrov space $\spc{A}$.

A locally Lipschitz curve $\alpha\:[t_{\min},t_{\max})\to\Dom f$ will be called an \index{gradient!curve}\emph{$f$-gradient curve} if
\[\alpha^+(t)=\nabla_{\alpha(t)} f\]
for any $t\in[t_{\min},t_{\max})$.
\end{thm}

A complete proof of the following theorem is given in \cite[16.15]{alexander-kapovitch-petrunin2024};
it mimics the proof of the standard Picard theorem on the existence and uniqueness of solutions of ordinary differential equations.
The uniqueness will follow from the first distance estimate (\ref{thm:dist-est}) proved in the next section.
We omit the proof of existence.

\begin{thm}{Picard theorem}\label{thm:glob-exist-grad-curv}
Let $f\:\spc{A}\subto \RR$ be a $\lambda$-concave function on an Alexandrov space $\spc{A}$.
Then for any $p\z\in \Dom f$, there exist a unique $t_{\max}\in(0,\infty]$ and a unique $f$-gradient curve $\alpha\:[0,t_{\max})\to \spc{A}$ with $\alpha(0)=p$ such that for any sequence $t_n\z\to t_{\max}-$, the sequence $\alpha(t_n)$ does not have a limit point in $\Dom f$.
\end{thm}

According to the theorem, the future of a gradient curve is determined by its present, but its past might \textit{not} be determined by the present.
Indeed, consider the function $f\:x\mapsto -|x|$ on $\RR$.
The tangent space $\T_x\RR$ can be identified with $\RR$.
Note that $\nabla_xf=-\mathrm{sgn}\, x$; that is,
\[\nabla_xf=
\begin{cases}
1&\text{if}\quad x<0,
\\
0&\text{if}\quad x=0,
\\
-1&\text{if}\quad x>0.
\end{cases}
\]
So, the $f$-gradient curves go to the origin with unit speed and then stand there forever.
In particular, if $\alpha$ is an $f$-gradient curve that starts at $x$,
then $\alpha(t)=0$ for any $t\ge |x|$.

Here is a slightly more interesting example;
it shows that gradient curves can merge even in the region where $|\nabla f|\z\ne 0$.

\begin{wrapfigure}[8]{r}{34 mm}
\vskip-0mm
\centering
\includegraphics{mppics/pic-1215}
\vskip0mm
\end{wrapfigure}

\begin{thm}{Example}
Consider the function $f\:(x,y)\mapsto-|x|-|y|$ on the $(x,y)$-plane.
It is concave, and its gradient field is sketched in the figure.

Let $\alpha$ be an $f$-gradient curve that starts at $(x,y)$ for $x>y>0$.
Then 
\[\alpha(t)=
\begin{cases}
(x-t,y-t) &\text{for}\quad 0\le t\le y,
\\
(x-t,0) &\text{for}\quad y\le t\le x,
\\
(0,0) &\text{for}\quad x\le t.
\end{cases}
\]

\end{thm}

\section{Distance estimates}

\begin{thm}{Lemma}\label{eq:fist-var-inq+}
Let $\alpha$ be a gradient curve of a $\lambda$-concave function $f$
defined on the whole Alexandrov space.
Choose a point $p$; let $\ell(t)\z\df\distfun_p\z\circ\alpha(t)$ and $q=\alpha(t_0)$.
Assume $p\ne q$.
Then
\[
\ell^+(t_0)\le -\left({f(p)}-{f(q)}-\tfrac\lambda2\cdot\ell^2(t_0)\right)/\ell(t_0)
\]
\end{thm}

\parit{Proof.}
Let $\gamma$ be the unit-speed parametrization of $[qp]$ from $q$ to $p$, so $q=\gamma(0)$.
Then 
\begin{align*}
\ell^+(t_0)&=(\dd_q\distfun_p)(\nabla_qf)\le\tag{by \ref{clm:fa'=dfa'}}
\\
&\le -\langle\dir qp,\nabla_qf\rangle \le \tag{by \ref{ex:d(distfun):<}}
\\
&\le -\dd_qf(\dir qp)=\tag{by \ref{def:grad}}
\\
&=-(f\circ\gamma)^+(0)\le 
\\
&\le -\left({f(p)}-{f(q)}-\tfrac\lambda2\cdot\ell^2(t_0)\right)/\ell(t_0)
\end{align*}
The last two lines follow from
the definition of differential
and the concavity of $t\z\mapsto f\circ\gamma(t)-\tfrac \lambda2\cdot {t^2}$.
\qeds

The following estimate implies uniqueness in the Picard theorem (\ref{thm:glob-exist-grad-curv}).

\begin{thm}{First distance estimate}\label{thm:dist-est}
Let $f$ be a $\lambda$-concave function defined on the whole Alexandrov space $\spc{A}$.
Then
\[\dist{\alpha(t)}{\beta(t)}{}
\le 
e^{\lambda\cdot t}\cdot\dist[{{}}]{\alpha(0)}{\beta(0)}{}\]
for any $t\ge 0$ and any two $f$-gradient curves $\alpha$ and $\beta$.

Moreover, the statement holds for a $\lambda$-concave function defined in an open domain if there is a geodesic $[\alpha(t)\,\beta(t)]$ in $\Dom f$ for any~$t$.
\end{thm}

\parit{Proof.} 
Fix a choice of a geodesic $[\alpha(t)\,\beta(t)]$ for each $t$.
Let $\ell(t)\z=\dist{\alpha(t)}{\beta(t)}{}$.
Note that
\[\ell^+(t)
\le-
\<\dir{\alpha(t)}{\beta(t)},\nabla_{\alpha(t)}f\>-\<\dir{\beta(t)}{\alpha(t)},\nabla_{\beta(t)}f\>
\le
\lambda\cdot\ell(t).\]
Here one has to apply \ref{eq:fist-var-inq+} for the distance from the midpoint $m$ of $[\alpha(t)\,\beta(t)]$, then apply the triangle inequality and \ref{ex:monotonicity}.

Since $\ell$ is Lipschitz, it is differentiable almost everywhere.
Integrating this inequality, we get the result.
\qeds

The following exercise describes a global geometric property of a gradient curve without direct reference to its function.
It is based on the notion of \index{self-contracting curves}\emph{self-contracting curves} introduced by Aris Daniilidis, Olivier Ley, and Stéphane Sabourau \cite{daniilidis-ley-sabourau}.

\begin{thm}{Exercise}\label{ex:elf-contracting}
Let $\alpha$ be a gradient curve of a concave function on an Alexandrov space.
Show that
\[\dist{\alpha(t_1)}{\alpha(t_3)}{\spc{A}}\ge \dist{\alpha(t_2)}{\alpha(t_3)}{\spc{A}}\]
if $t_1\le t_2\le t_3$.
\end{thm}

\begin{thm}{Exercise}\label{ex:mayer}
Let $f$ be a 
concave function defined on an Alexandrov space $\spc{A}$.
Suppose $\hat\alpha\:[0,\ell]\to\spc{A}$ is an arc-length reparametrization of an $f$-gradient curve.
Show that $f\circ\hat\alpha$ is concave.
\end{thm}

The following exercise implies that gradient curves for a uniformly converging sequence of $\lambda$-concave functions converge to the gradient curves of the limit function.

\begin{thm}[!]{Exercise}\label{lem:fg-dist-est}
Let $f$ and $g$ be $\lambda$-concave 
 functions on an Alexandrov space $\spc{A}$.
Suppose
$\alpha,\beta\:[0,t_{\max})\to \spc{A}$ are respectively $f$- and $g$-gradient curves.
Assume $|f-g|<\eps$; let $\ell\:t\mapsto\dist{\alpha(t)}{\beta(t)}{}$.
Show that
\[\ell^+\le \lambda\cdot\ell+\tfrac{2\cdot\eps}{\ell}.\]

Conclude that if $\alpha(0)=\beta(0)$ and $t_{\max}<\infty$, then
\[\dist{\alpha(t)}{\beta(t)}{}
\le
\Const\cdot\sqrt{\eps\cdot t}\]
for some constant $\Const=\Const(t_{\max},\lambda)$.
\end{thm}

\section{Gradient flow}

Let $f$ be a 
semiconcave function defined on an open subset of an Alexandrov space $\spc{A}$.
If there is an $f$-gradient curve $\alpha$ such that $\alpha(0)\z=x$ and $\alpha(t)=y$,
then we will write 
\[\GF^t_f(x)=y.\]
The partially defined map $\GF^t_f$ from $\spc{A}$ to itself is called the \index{gradient!flow}\emph{$f$-gradient flow} for time $t$.
Note that
\[\GF^{t_1+t_2}_f=\GF_f^{t_1}\circ\GF_f^{t_2}.\]
In other words, the gradient flow is a partial action of the \textit{semigroup} $([0,\infty),+)$ on $\spc{A}$.
 
From the first distance estimate \ref{thm:dist-est}, 
it follows that for any $t\z\ge 0$, the domain of definition of $\GF^t_f$ is an open subset of $\spc{A}$.
For sufficiently nice functions, the gradient flow is globally defined.
For example, if $f$ is a $\lambda$-concave function and it is defined on the whole space $\spc{A}$, then $\GF^t_f(x)$ is defined for all $x\in \spc{A}$ and $t\ge0$;
see \cite[16.19]{alexander-kapovitch-petrunin2024}.

Using this new terminology, we can reformulate several statements about gradient curves.
From the first distance estimate, we have the following.

\begin{thm}{Proposition}\label{prop:GF-is-lip}
Let $f$ be a semiconcave function on an Alexandrov space $\spc{A}$.
Then the map $x\mapsto\GF^t_f(x)$ is locally Lipschitz.

Moreover, if $f$ is $\lambda$-concave, then $\GF^t_f$ is $e^{\lambda\cdot t}$-Lipschitz.
\end{thm}

The next proposition follows from \ref{lem:fg-dist-est}.

\begin{thm}{Proposition}\label{grad-curve-conv}
Let $\spc{A}$ be an Alexandrov space.
Suppose $f_n\:\spc{A}\z\to\RR$ is a sequence of
$\lambda$-concave functions 
that uniformly converges to $f_\infty\:\spc{A}\z\to \RR$. 
Then for any $x\in \spc{A}$ and $t\ge 0$, we have
\[\GF_{f_n}^t(x)\to \GF_{f_\infty}^t(x)\]
as $n\to \infty$.
\end{thm}

This proposition can be generalized to a converging sequence $\spc{A}_n\z\to \spc{A}_\infty$ of spaces and a converging sequence of functions $f_n\:\spc{A}_n\z\to\RR$; see \cite[16.21]{alexander-kapovitch-petrunin2024}.

\section{Gradient exponent}\label{gexp}

One of the technical difficulties in Alexandrov geometry comes from
nonextendability of geodesics. 
In particular, the exponential map, $\exp_p\:\T_p\to \spc{A}$,
when defined in the usual way, may be undefined
in an arbitrarily small neighborhood of the origin.

The \index{gradient!exponential map}\emph{gradient exponential map}
\[\gexp_p\:\T_p\to\spc{A},\]
which we are about to introduce, essentially solves this problem.
It shares many properties with the ordinary exponential map and is better in certain respects,
even in the Riemannian universe.

Let $p$ be a point in an $\Alex0$ space $\spc{A}$.
Consider the function $f\z=\tfrac12\cdot\distfun_p^2$.
Recall that $\GF^t_{f}$ denotes the gradient flow.
Let us define the \textit{gradient exponential map} as the limit
\[\gexp_p(v)=\lim_{n\to\infty}\GF^{t_n}_{f}(x_n),\]
where the sequences $x_n\in \spc{A}$ and $t_n\ge 0$ are chosen so that $t_n\to\infty$ and $e^{t_n}\cdot\log_px_n\to v$ as $n\to\infty$ (more formally, we choose one value of $\log_p x_n$ for each $n$).

More intuitively, suppose $i_{\lambda}\:\lambda\cdot \spc{A}\to \spc{A}$ sends a point in the rescaled copy $\lambda\cdot\spc{A}$ to the corresponding point in $\spc{A}$.
By the first distance estimate (\ref{thm:dist-est}), the map
\[\GF^t_{f}\circ\, i_{e^t}\:e^t\cdot \spc{A}\to \spc{A}
\eqlbl{eq:gexp}\]
is short for any $t\ge 0$.
If we have a pointed Gromov--Hausdorff convergence $(e^{t_n}\cdot \spc{A},p)\to (\T_p,\0)$,
then $\gexp_p\:\T_p\to \spc{A}$ is the limit of $\GF^{t_n}_{f}\circ i_{e^{t_n}}$.
This way, we get that $\gexp_p$ is short as a limit of short maps.
This observation is generalized in the following proposition.

\begin{thm}{Proposition}\label{prop:gexp}
Let $\spc{A}$ be a proper $\Alex0$ space.
Then for any $p\in \spc{A}$ the gradient exponent $\gexp_p\:\T_p\to\spc{A}$ is uniquely defined.
Moreover, $\gexp_p$ is a short map and 
\[\gexp_p(\gamma^+(0))=\gamma(1)\]
for any geodesic path $\gamma$ that starts at $p$.
\end{thm}

The last statement implies that 
\[\gexp_p\circ\log_p=\id,\]
so it is appropriate to use the term \textit{exponent} for $\gexp$.

\parit{Proof.} 
Note that $f''\le 1$.
Since the space is proper we can choose a limit in \ref{eq:gexp}.

Let $\gamma$ be a geodesic that starts at $p$.
Observe that $t\mapsto \gamma\circ\exp(t)$ is an $f$-gradient curve.
By the first distance estimate, we have that $\GF^t_{f}$ is $e^t$-Lipschitz.
This implies that any limit in \ref{eq:gexp} has the same value;
that is, $\gexp_p$ is uniquely defined.

Again, since $\GF^t_{f}$ is $e^t$-Lipschitz, we get that $\gexp_p$ is short.
\qeds

\section{Remarks}

The idea of using gradient flows in Alexandrov geometry was inspired by the success of \index{Sharafutdinov's retraction}\emph{Sharafutdinov's retraction} in comparison geometry \cite{sharafutdinov}.
The gradient flow was introduced by the third author to construct quasigeodesics with given initial data \cite{perelman-petrunin:qg,petrunin:qg, petrunin:survey}.
It turned out that gradient flow and gradient exponent are better tools than quasigeodesics.
These tools quickly found applications in other types of singular spaces \cite{jost,mayer,lytchak:open-map,ohta,sevare,ambrosio-gigli-savare}.

For a general lower curvature bound $\kappa$, the construction of the gradient exponent has to be modified;
it is denoted by $\gexp_p^\kappa$ \cite[16.36]{alexander-kapovitch-petrunin2024}. It is done by taking limits of appropriately reparameterized gradient curves of the modified distance function.

For $\kappa=-1$, we have
that $\gexp_p^{-1}(\gamma^+(0))=\gamma(1)$
for any geodesic path $\gamma$ that starts at $p$
and 
\[\dist{\gexp_p^{-1}v}{\gexp_p^{-1}w}{\spc{A}}\le \side\hinge0vw_{\HH^2}.\]
In other words, $\gexp_p^{-1}$ is short if we equip $\T_p$ with the hyperbolic cone metric.

Similarly, for $\kappa=1$, we have $\gexp_p^1(\gamma^+(0))=\gamma(1)$
for any geodesic path $\gamma$ that starts at $p$ and 
\[\dist{\gexp_p^{1}v}{\gexp_p^{1}w}{\spc{A}}\le \side\hinge0vw_{\SSS^2},\]
but this time all this holds only if $|v|,|w|\le\tfrac\pi2$ and $\length\gamma\le\tfrac\pi2$.

The gradient exponential map in a Riemannian manifold $(M,g)$ coincides with the Riemannian exponential map before the cut locus, but \textit{is different} from the Riemannian exponential beyond the cut locus.
Moreover, beyond the cut locus $\gexp_p^\kappa$ depends on $\kappa$ while the Riemannian exponential map does not.
The following exercise is meant to show that this technique can prove something nontrivial even for Riemannian manifolds.

\begin{thm}{Exercise}\label{ex:short-onto}
Let $(M,g)$ be a complete $m$-dimensional Riemannian manifold with sectional curvature at least $1$.
Assume $M$ is not homeomorphic to $\SSS^m$.
Show that there is a short onto map $\SSS^m\z\to (M,g)$.
\end{thm}

%% file: splitting.tex
\chapter{Line splitting}\label{chap:splitting}

In this lecture, we prove the line splitting theorem and apply it to tangent spaces of Alexandrov spaces.

\section{Busemann function}

A \index{half-line}\emph{half-line} is a geodesic defined on the real half-line $[0,\infty)$;
that is, a distance-preserving map from $[0,\infty)$ to a metric space.

If $\gamma\:[0,\infty)\to \spc{X}$ is a half-line,
then the limit 
\[\bus_\gamma(x)=\lim_{t\to\infty}\dist{\gamma(t)}{x}{}- t\eqlbl{eq:def:busemann*}\]
is called the \index{Busemann function}\emph{Busemann function} of $\gamma$.
It behaves like a shifted distance function from the ideal point of $\gamma$.

\begin{thm}{Proposition}\label{prop:busemann}
For any half-line $\gamma$ in a metric space $\spc{X}$,
its Busemann function $\bus_\gamma\:\spc{X}\to \RR$ 
is defined.
Moreover, $\bus_\gamma$ is $1$-Lipschitz and $\bus_\gamma \circ\gamma(t)=-t$ for any $t$.

\end{thm}

\parit{Proof.}
Since $t=\dist{\gamma(0)}{\gamma(t)}{}$, the triangle inequality implies that
\[t\mapsto\dist{\gamma(t)}{x}{}- t\] 
is a nonincreasing function, and
\[\dist{\gamma(t)}{x}{}- t\ge-\dist{\gamma(0)}{x}{}\]
for any $x\in \spc{X}$.
Therefore, the limit in \ref{eq:def:busemann*} is defined,
and it has to be 1-Lipschitz as a limit of 1-Lipschitz functions.

The last statement follows since 
$\dist{\gamma(t)}{\gamma(t_0)}{}\z=t-t_0$ for $t>t_0$.
\qeds

\begin{thm}[!]{Exercise}\label{ex:busemann-CBB}
Show that any Busemann function on an $\Alex0$ space is concave.
\end{thm}

\section{Splitting theorem}

A \index{line}\emph{line} is a distance-preserving map
from $\RR$ to a metric space.
In other words, a line is a geodesic defined on the whole real line $\RR$.

\begin{thm}[!]{Exercise}\label{ex:bus+bus}
Let $\gamma$ be a line in a metric space $\spc{X}$.
Show that for any point $x$ we have
\[\bus_+(x)+\bus_-(x)\ge 0,\]
where $\bus_+$ and $\bus_-$ are the Busemann functions associated with the half-lines $t\mapsto\gamma(t)$ and $t\mapsto\gamma(-t)$, respectively.
\end{thm}

Let $A$ and $B$ be two subsets of a metric space $\spc{X}$.
We say that $\spc{X}$ is a \index{direct sum}\emph{direct sum} of $A$ and $B$,
or briefly,
\[\spc{X}=A\oplus B\]\index{90@$A\oplus B$ (direct sum)}
if there are retractions $\proj_A\:\spc{X}\to A$
and 
$\proj_B\:\spc{X}\to B$
such that $x\mapsto(\proj_A(x),\proj_B(x))$ defines an onto map $\spc{X}\to A\times B$ and
\[\dist[2]{x}{y}{}=\dist[2]{\proj_A(x)}{\proj_A(y)}{}+\dist[2]{\proj_B(x)}{\proj_B(y)}{}\]
for any two points $x,y\in \spc{X}$.

Note that if $\spc{X}=A\oplus B$, then $\spc{X}$ is isometric to $A\times B$.

\begin{thm}{Exercise}\label{ex:oplus}
Suppose $\spc{X}=A\oplus B$.
Show that
\begin{subthm}{ex:oplus:a}
 $A$ intersects $B$ at a single point, and
\end{subthm}

\begin{subthm}{ex:oplus:b}
both sets $A$ and $B$ are \index{convex set}\emph{convex} in $\spc{X}$, meaning that any geodesic with the endpoints in $A$ (or $B$) lies in $A$ (respectively, $B$).
\end{subthm}

\end{thm}

\begin{thm}{Line splitting theorem}\label{thm:splitting}\index{line splitting theorem}
If $\gamma$ is a line in an $\Alex0$ space~$\spc{A}$, then
\[\spc{A}=\spc{A}'\oplus \gamma(\RR)\]
for some subset $\spc{A}'\subset \spc{A}$.
\end{thm}

\begin{thm}{Corollary}\label{cor:splitting}
Any $\Alex0$ space $\spc{A}$ splits isometrically as
\[
\spc{A}=\spc{A}'\oplus H
\]
where $H\subset \spc{A}$ is a subset isometric to a Hilbert space, and $\spc{A}'\subset \spc{A}$ is a convex subset that contains no lines. 
\end{thm}

The following lemma is closely related to the first distance estimate (\ref{thm:dist-est});
it is also a limit case of \ref{prop:gexp}.
The proof is similar.

\begin{thm}{Lemma}\label{lem:dist-estimate}
Suppose $f\:\spc{A}\to\RR$ is a concave 1-Lipschitz function on an $\Alex0$ space $\spc{A}$.
Consider two $f$-gradient curves $\alpha$ and~$\beta$.
Then for any $t, s\ge 0$ we have
\begin{align*}
&\dist[2]{\alpha(s)}{\beta(t)}{}
\le 
\dist[2]{p}{q}{}+
2\cdot(f(p)-f(q))\cdot(s-t)+ (s-t)^2,
\end{align*}
where $p=\alpha(0)$ and $q=\beta(0)$.
\end{thm}

\parit{Proof.}
Since $f$ is 1-Lipschitz, it follows that $|\nabla f|\le1$.
Therefore 
\[f\circ\beta(t)\le f(q)+t\]
for any $t\ge0$.

Set $\ell(t)=\dist{p}{\beta(t)}{}$.
Applying \ref{eq:fist-var-inq+}, we get
\begin{align*}
(\ell^2)^+(t)
&\le 2\cdot \left(f\circ\beta(t)-f(p)\right)\le 
\\
&\le2\cdot\left(f(q)+t-f(p)\right).
\end{align*}
Therefore 
\[\ell^2(t)-\ell^2(0)\le 2\cdot\left(f(q)-f(p)\right)\cdot t + t^2.\]
This proves the needed inequality in the case $s=0$.
Combining it with the first distance estimate (\ref{thm:dist-est}), we get the result in the case $s\le t$.
The case $s\ge t$ follows by switching the roles of $s$ and $t$.
\qeds

\parit{Proof of \ref{thm:splitting}.} Consider two Busemann functions, $\bus_+$ and $\bus_-$, associated with the half-lines $\gamma\:[0,\infty)\to \spc{A}$ and $\gamma\:(-\infty,0]\to \spc{A}$, respectively; that is,
\[
\bus_\pm(x)
\df
\lim_{t\to\infty}\dist{\gamma(\pm t)}{x}{}- t.
\]
According to \ref{ex:busemann-CBB}, 
both $\bus_+$ and $\bus_-$ are concave.

By \ref{comp-kappa},
$\phi\:t\mapsto\distfun_x^2\circ\gamma(t)$
is $2$-concave.
In particular, $\phi(t)\z\le t^2+at+b$ for some constants $a,b\in\RR$; see \ref{Function comparison}.
Therefore,
\[
\dist{\gamma( t)}{x}{}+\dist{\gamma(- t)}{x}{}- 2\cdot t
\le \sqrt{ t^2+at+b}+\sqrt{ t^2-at+b}-2\cdot t
\]
for all large $t$.
Passing to the limit as $t\to\infty$, we get that $\bus_+(x)\z+\bus_-(x)\le 0$.
By \ref{ex:bus+bus}, $\bus_+(x)+\bus_-(x)\ge0$; hence
\[
\bus_+(x)+\bus_-(x)= 0
\]
for any $x\in \spc{A}$.
In particular, the functions $\bus_+$ and $\bus_-$ are \index{affine function}\emph{affine};
that is, they are convex and concave at the same time.

For any $x$,
\begin{align*}
|\nabla_x \bus_\pm|
&=\sup\set{\dd_x\bus_\pm(\xi)}{\xi\in\Sigma_x}=\\
&=\sup\set{-\dd_x\bus_\mp(\xi)}{\xi\in\Sigma_x}\equiv\\
&\equiv1.
\end{align*}

A curve $\alpha$ is a $\bus_\pm$-gradient curve
if and only if $\alpha$ is a geodesic such that $(\bus_\pm\circ\alpha)^+=1$.
Indeed, if $\alpha$ is a geodesic, then $(\bus_\pm\circ\alpha)^+\le 1$ and the equality holds only if $\nabla_\alpha\bus_\pm=\alpha^+$.
Now suppose $\nabla_\alpha\bus_\pm\z=\alpha^+$.
Then $|\alpha^+|\le 1$ and $(\bus_\pm\circ\alpha)^+=1$; therefore 
\begin{align*}
|t_0-t_1|&\ge \dist{\alpha(t_0)}{\alpha(t_1)}{}\ge
\\
&\ge|\bus_\pm\circ\alpha(t_0)-\bus_\pm\circ\alpha(t_1)|=
\\
&=|t_0-t_1|.
\end{align*}

It follows that for any $t>0$, the $\bus_\pm$-gradient flows are inverse to each other;
that is, 
\[\GF_{\bus_+}^t\circ\GF_{\bus_-}^t=\id_\spc{A}.\]
Setting
\[\GF^t\df
\left[
\begin{matrix}
\GF_{\bus_-}^{t}&\hbox{if}\ t\ge0
\\
\GF_{\bus_+}^{-t}&\hbox{if}\ t\le0
\end{matrix}
\right.\]
defines an $\RR$-action on $\spc{A}$.

Consider the level set $\spc{A}'=\bus_+^{-1}\{0\}=\bus_-^{-1}\{0\}$;
it is a closed convex subset of $\spc{A}$, and therefore forms an Alexandrov space.
Let $h\:\spc{A}'\times \RR\to \spc{A}$ be the map defined by $h\:(x,t)\mapsto \GF^t(x)$.
Note that $\GF^{-t}(y)\in \spc{A}'$ for $t=\bus_-(y)=-\bus_+(y)$; hence $h$ is onto.

Applying \ref{lem:dist-estimate} for $\GF_{\bus_+}^t$ and $\GF_{\bus_-}^t$ shows that $h$ is distance nonexpanding and noncontracting at the same time; that is, $h$ is an isometry.
\qeds

\begin{thm}{Exercise}\label{ex:cone-CBB}
Suppose $\spc{X}$ is a complete geodesic space.
Show that $\Cone\spc{X}$ is $\Alex0$ if and only if $\spc{X}$ is $\Alex1$ and $\diam\spc{X}\le \pi$.
\end{thm}

Recall that according to our definition, any circle or closed real interval is $\Alex1$.
Therefore, the condition $\diam\spc{X}\le \pi$ is necessary.
Nevertheless, according to \ref{ex:RisCBB(1)}, most $\Alex1$ spaces have diameter at most $\pi$.

\section{Anti-sum}

The following lemma is a corollary of \ref{ex:convergence-grad}.
It will be used to prove basic properties of tangent spaces.

\begin{thm}{Anti-sum lemma}\label{lem:minus-sum} 
Let $p$ be a point in an Alexandrov space $\spc{A}$.
Given two vectors $u,v\in \T_p$, there is a unique vector $w\in \T_p$ such that
\[\langle u,x\rangle+\langle v,x\rangle+\langle w,x\rangle\ge0\]
for any $x\in \T_p$, and
\[\langle u,w\rangle+\langle v,w\rangle+\langle w,w\rangle=0.\]

\end{thm}

\begin{thm}[!]{Exercise}\label{ex:|antisum|}
Suppose $u,v, w\in \T_p$ are as in \ref{lem:minus-sum}.
Show that 
\[|w|^2\le |u|^2+|v|^2+2\cdot\langle u,v\rangle.\]

\end{thm}

If the tangent space $\T_p$ was $\Alex0$, then the lemma would follow from the existence of the gradient, applied to the function $\T_p\to \RR$ defined by $x\mapsto -(\langle u,x\rangle +\langle v,x\rangle )$, which is concave by \ref{ex:busemann-CBB}.
This idea cannot work as is since $\T_p$ might fail to be geodesic (see Halbeisen's example \cite[13.6]{alexander-kapovitch-petrunin2024}), but it works after a proper revision.

Applying the above lemma with $u=v$, we have the following statement.

\begin{thm}{Existence of polar vector}\label{cor:polar}
Let $\spc{A}$ be an Alexandrov space 
and $p\in \spc{A}$. 
Given a vector $u\in \T_p$, there is a unique vector $u^*\in\T_p$ such that $\langle u^*,u^*\rangle +\langle u,u^*\rangle = 0$ and
$u^*$ is \index{polar vectors}\emph{polar} to $u$;
that is,
\[\langle u^*,x\rangle +\langle u,x\rangle \ge 0\]
for any $x\in \T_p$.
\end{thm}

\parit{Proof of \ref{lem:minus-sum}.} We may assume that $u$ and $v$ are nonzero.
By \ref{ex:d(distfun):==}, we can choose two sequences of points $a_n,b_n$ such that 
\[
(\dd_p\distfun_{a_n})(x)=-\langle\dir{p}{a_n},x\rangle
\quad\text{and}\quad
(\dd_p\distfun_{b_n})(x)=-\langle\dir{p}{b_n},x\rangle
\]
for any $x\in\T_p$.
Furthermore, we can assume that $\dir{p}{a_n}\to u/|u|$ and $\dir{p}{b_n}\to v/|v|$ as $n\to \infty$.

Consider a sequence of functions 
\[f_n=|u|\cdot\distfun_{a_n}+|v|\cdot\distfun_{b_n}.\]
Note that 
\[(\dd_pf_n)(x)=-|u|\cdot\langle \dir{p}{a_n},x\rangle -|v|\cdot\langle \dir{p}{b_n},x\rangle\]
for any $x\in\T_p$.
Thus, $\dd_pf_n$ uniformly converges to $x\mapsto-\langle u,x\rangle -\langle v,x\rangle$ on $\Sigma_p$.
According to \ref{ex:convergence-grad}, 
the sequence $\nabla_pf_n$ converges as $n\to\infty$.
Let 
\[w=\lim_{n\to\infty}\nabla_pf_n.\]
By the definition of gradient,
\begin{align*}
\langle w,w\rangle &=\lim_{n\to\infty}\langle \nabla_pf_n,\nabla_pf_n\rangle =
&&&
\langle w,x\rangle &=\lim_{n\to\infty}\langle \nabla_pf_n,x\rangle \ge
\\
&=\lim_{n\to\infty}(\dd_p f_n)(\nabla_p f_n)
=
&&&
&\ge
\lim_{n\to\infty}(\dd_pf_n)(x)
=
\\
&=-\langle u,w\rangle -\langle v,w\rangle ,
&&&
&=-\langle u,x\rangle -\langle v,x\rangle .
\end{align*}

The proof of uniqueness is left to the reader;
it is very similar to the proof of uniqueness of gradients (page \pageref{uniqueness:prop:grad-exist}).
\qeds

\section{Linear subspace}

\begin{thm}{Definition}\label{def:opp+Lin}
Let $\spc{A}$ be an Alexandrov space, $p\in \spc{A}$, and $u,v\z\in\T_p$.
We say that vectors $u$ and $v$ are \index{opposite vectors}\emph{opposite}\label{def:opposite:page} to each other (briefly, $u+v=0$) if
\begin{itemize}
 \item $|u|=|v|=0$ or
 \item $\mangle(u,v)=\pi$ and $|u|=|v|$.
\end{itemize}

The subcone
\[\Lin_p=\set{v\in\T_p}{\exists\ w\in\T_p\quad \text{such that}\quad w+v=0}\]
will be called the \index{linear subspace}\emph{linear subspace} of $\T_p$.
\end{thm}

Soon we will introduce a natural linear structure on $\Lin_p$, which will justify its name.

\begin{thm}{Proposition}\label{prop:opposite}
Let $p$ be a point in an Alexandrov space.
Given two vectors $u,v\in\T_p$, the following statements are equivalent:
\begin{subthm}{opposite} $u+v=0$;
\end{subthm}
\begin{subthm}{<x,u>} $\langle u,x\rangle +\langle v,x\rangle =0$ for any $x\in\T_p$;
\end{subthm}
\begin{subthm}{<xi,u>} $\langle u,\xi\rangle +\langle v,\xi\rangle =0$ for any $\xi\in\Sigma_p$.
\end{subthm}
\end{thm}

\parit{Proof.}
The equivalence \ref{SHORT.<x,u>}$\Leftrightarrow$\ref{SHORT.<xi,u>} is trivial.

The condition $u+v=0$ is equivalent to 
$\langle u,u\rangle =-\langle u,v\rangle =\langle v,v\rangle$;
thus,
\ref{SHORT.<x,u>}$\Rightarrow$\ref{SHORT.opposite}.
It remains to prove \ref{SHORT.opposite}$\Rightarrow$\ref{SHORT.<x,u>}.

Suppose \ref{SHORT.opposite} holds.
We can assume that $|u|=|v|\ne 0$, since otherwise \ref{SHORT.<x,u>} is trivial.
Since $\mangle(u,v)=\pi$ and $\T_p$ is a metric cone, the half-lines in the directions of $u$ and $v$ form a line in $\T_p$.
By \ref{prop:Tan-is-CBB(0)}, $\T_p$ has nonnegative curvature, and by \ref{ex:adjacent-CBB} we have $\mangle\hinge 0ux+\mangle\hinge 0vx=\pi$, which implies \ref{SHORT.<x,u>}.
\qeds

\begin{thm}[!]{Exercise}\label{prop:two-opp}
Let $u$, $v$, and $w$ be tangent vectors at a point of an Alexandrov space.
Assume $u+v=0$ and $u+w=0$.
Show that $v=w$.
\end{thm}

Let $u\in \Lin_p$; that is, $u+v=0$ for some $v\in\T_p$.
Given $s<0$, let \index{33@$r\cdot u$ (multiplication in cone)}
\[s\cdot u\df (-s)\cdot v.\]
So we can multiply any vector in $\Lin_p$ by any real number (positive and negative).
By \ref{prop:two-opp}, this multiplication is uniquely defined.
By \ref{prop:opposite}, we have the identity
\[\langle -v,x\rangle=-\langle v,x\rangle.\]

\begin{thm}[!]{Exercise}\label{ex:3<,>=0}
Suppose $u$, $v$, and $w$ are as in \ref{lem:minus-sum}.
Show that
\[\langle u,x\rangle +\langle v,x\rangle +\langle w,x\rangle = 0\]
for any $x\in \Lin_p$.
\end{thm}

\begin{thm}[!]{Exercise}\label{ex:-u}
Let $\spc{A}$ be an Alexandrov space,
$p\in \spc{A}$ and $u\in \T_p$.
Suppose $u^*\in \T_p$ is as in \ref{cor:polar};
that is, 
\[\langle u^*,u^*\rangle +\langle u,u^*\rangle = 0
\quad\text{and}\quad
\langle u^*,x\rangle +\langle u,x\rangle \ge 0
\]
for any $x\in \T_p$.
Show that $u=-u^*$ if and only if $|u|=|u^*|$.
\end{thm}

\begin{thm}{Theorem}\label{thm:lin-subcone}
Let $p$ be a point in an Alexandrov space. 
Then $\Lin_p$ is isometric to a Hilbert space.
\end{thm}

\parit{Proof.}
$\Lin_p$ is a closed subset of $\T_p$;
in particular, it is complete.

If we know that any two vectors in $\Lin_p$ can be connected by a geodesic in $\Lin_p$,
then the statement follows from the splitting theorem (\ref{thm:splitting}).
By Menger's lemma (\ref{lem:mid>geod}), it is sufficient to show that any pair $x,y\in\Lin_p$ has a midpoint $w\in \Lin_p$.

Choose $w\in \T_p$ to be the anti-sum of $u=-\tfrac{1}{2}\cdot x$ and $v=-\tfrac{1}{2}\cdot y$;
see \ref{lem:minus-sum}.
By \ref{ex:|antisum|} and \ref{ex:3<,>=0},
\begin{align*}
|w|^2&\le \tfrac14\cdot |x|^2+\tfrac14\cdot|y|^2+\tfrac12\cdot\langle x,y\rangle,
\\
\langle w,x\rangle&= \tfrac12\cdot|x|^2+\tfrac12\cdot\langle x,y\rangle,
\\
\langle w,y\rangle&= \tfrac12\cdot|y|^2+\tfrac12\cdot\langle x,y\rangle.
\end{align*}
It follows that 
\begin{align*}
|x-w|^2
&= |x|^2+|w|^2-2\cdot\langle w,x\rangle\le
\\
&\le \tfrac14\cdot |x|^2+\tfrac14\cdot|y|^2-\tfrac12\cdot\langle x,y\rangle=
\\
&=\tfrac14\cdot|x-y|^2.
\end{align*}
That is, $|x-w|\le \tfrac12\cdot|x-y|$.
Similarly, we get $|y-w|\le \tfrac12\cdot|x-y|$.
Therefore, $w$ is a midpoint of $x$ and $y$.
In addition, we get the equality 
\[|w|^2= \tfrac14\cdot |x|^2+\tfrac14\cdot|y|^2+\tfrac12\cdot\langle x,y\rangle.\]

It remains to show that $w\in\Lin_p$.
Let $w^*$ be the polar vector provided by \ref{cor:polar}.
Note that 
\[|w^*|\le |w|,
\quad
\langle w^*,x\rangle+\langle w,x\rangle=0,
\quad\text{and}\quad
\langle w^*,y\rangle+\langle w,y\rangle=0.
\]
The same calculation as above shows that $w^*$ is a midpoint of $-x$ and $-y$ and 
\[|w^*|^2= \tfrac14\cdot |x|^2+\tfrac14\cdot|y|^2+\tfrac12\cdot\langle x,y\rangle=|w|^2.\]
By \ref{ex:-u}, $w=-w^*$;
hence $w\in\Lin_p$.
\qeds

\begin{thm}{Lemma}\label{ex:grad-dist:G-delta}
Given a point $p$ in an Alexandrov space $\spc{A}$,
let $f\z=\distfun_p$, and let $S$ be the subset of points $x\in\spc{A}$ such that $|\nabla_xf|=1$.
Then $S$ is a dense G-delta set; that is, a countable intersection of open sets.

\end{thm}

\parit{Proof.}
Given a positive integer $n$, consider the set
\[S_n=\set{x\in\spc{A}}{|\nabla_xf|>1-\tfrac1n}.\]
By \ref{ex:semicontinuous-grad:>s}, $S_n$ is open.

Choose a point $q\ne p$.
Observe that $|\nabla_xf|=1$ for any point $x\in\left]pq\right[$.
It follows that $S_n$ is dense in $\spc{A}$.
Since $S=\bigcap_nS_n$, the lemma follows from the Baire category theorem.
\qeds

\begin{thm}[!]{Exercise}\label{ex:grad-dist}
Let $p$, $f$, and $S$ be as in \ref{ex:grad-dist:G-delta}.

\begin{subthm}{ex:grad-dist:lin}
Show that 
\[\nabla_xf+\dir xp=0\]
for any 
$x\in S$;
in particular, $\dir xp\in \Lin_x$.
\end{subthm}

\begin{subthm}{ex:grad-dist:|grad|=1}
Show that if $|\nabla_xf|=1$, then $(\dd_xf)(w)= \langle\nabla_xf,w\rangle$ for every $w\z\in \T_x$.
\end{subthm}

\begin{subthm}{ex:grad-dist:geod}
Show that for any $x\in S$ there is a unique geodesic $[px]$.
\end{subthm}

\end{thm}

This exercise implies the following.

\begin{thm}{Corollary}\label{cor:euclid-subcone}
Given a countable set of points $X$ in a nontrivial Alexandrov space $\spc{A}$,
there is a G-delta dense set $S\subset\spc{A}$
such that 
$\dir sx\in \Lin_s$
for any $s\in S$ and $x\in X$.
\end{thm}

\section{Remarks}

The history of the splitting theorem starts with Stefan Cohn-Vossen \cite{cohn-vossen_line},
who proved it in the $2$-dimensional case.
For Riemannian manifolds of higher dimensions,
it was proved by Victor Toponogov \cite{toponogov-globalization+splitting}.
Then it was generalized by Anatoliy Milka \cite{milka-line}
to Alexandrov spaces;
historically, it was the first result about Alexandrov spaces of dimension higher than~$2$.
Nearly the same proof is used in \cite[1.5]{burago-burago-ivanov}.

Generalizations to Riemannian manifolds with nonnegative Ricci curvature were obtained by Jeff Cheeger and Detlef Gromoll \cite{cheeger-gromoll-split}.
This was further generalized by Jeff Cheeger and Toby Colding for the limits of Riemannian manifolds with almost nonnegative Ricci curvature \cite{cheeger-colding-alm-rigidity}.
Nicola Gigli generalized it further to the so-called {}\emph{RCD spaces} (spaces with synthetically defined Ricci curvature bound) \cite{gigli2013splitting, gigli-splitting-overview}.
Jost-Hinrich Eschenburg obtained an analogous result for Lorentzian manifolds \cite{eshenburg-split}.

The presented proof is close in spirit to the proof given by Cheeger and Gromoll \cite{cheeger-gromoll-split};
it is taken from our book \cite{alexander-kapovitch-petrunin2024}.

\begin{thm}{Open question}
Let $p$ be a point in an Alexandrov space $\spc{A}$.
Suppose that $0\ne v\in \Lin_p$.
Is it true that the tangent space $\T_p$ splits in the direction of $v$?
\end{thm}

\begin{thm}{Open question}\label{open:Halb-proper}
Let $\spc{A}$ be a proper Alexandrov space.
Is it true that for any $p\in \spc{A}$, the tangent space $\T_p$ is a length space?
\end{thm}

%% file: dim.tex
\chapter{Dimension and volume}\label{chap:dim}

\section{Linear dimension}

Let $\spc{A}$ be an Alexandrov space.
We define its \index{linear dimension}\emph{linear dimension} \index{LinDim@$\LinDim$}$\LinDim \spc{A}$ as the least upper bound of integers $m$ such that
the Euclidean space $\EE^m$ is isometric to a subspace of the tangent space $\T_p\spc{A}$ at \textit{some} point $p\in \spc{A}$.

If not stated otherwise, dimension means linear dimension.
In \ref{sec:packing-numbers} and \ref{sec:all-dim}, we will show that the linear dimension of Alexandrov spaces coincides with all reasonable notions of dimension;
thereafter, we will use the notation \index{dim@$\dim$}$\dim\spc{A}$ for $\LinDim\spc{A}$.

\begin{thm}{(\textit{n}+1)-comparison}\label{thm:n+1}
Let $\spc{A}$ be an $\Alex0$ space.
Then for any finite set of points $p,x_1,\dots,x_n\in \spc{A}$, there exists a model configuration
$\tilde p,\tilde x_1,\dots,\tilde x_n\in \EE^m$ for some $m$ such that
\[|\tilde p-\tilde x_i|_{\EE^m}=| p- x_i|_{\spc{A}}
\quad\text{and}\quad
|\tilde x_i-\tilde x_j|_{\EE^m}\ge |x_i- x_j|_{\spc{A}}\]
for any $i$ and $j$.
Moreover, we can assume that $m\le \LinDim\spc{A}$.

If $\spc{A}$ is an $\Alex{\pm1}$ space, then the same holds if we replace $\EE^m$ by $\SSS^m$ or $\HH^m$, respectively; in the $\Alex{1}$ case one has to assume in addition that $\diam \spc{A}\le \pi$.
\end{thm}

\parit{Proof.}
By \ref{cor:euclid-subcone}, we can choose a point $p'$ arbitrarily close to $p$ so that
$\Lin_{p'}\ni \dir{p'}{x_i}$ for any $i$.
Let us identify $\EE^m$ with a subspace of $\Lin_{p'}$ spanned by $\dir{p'}{x_1},\dots,\dir{p'}{x_n}$.
Note that $m\le \LinDim\spc{A}$.

Set $\tilde p'=0\in \EE^m$ and $\tilde x_i=\dist{p'}{x_i}{}\cdot\dir{p'}{x_i}\in \EE^m$ for every $i$.
Note that
\[|\tilde p'-\tilde x_i|_{\EE^m}=| p'- x_i|_{\spc{A}}\]
for every $i$.
Applying the comparison $\mangle\hinge {p'}{x_i}{x_j}\ge \angk {p'}{x_i}{x_j}$, we get
\[|\tilde x_i-\tilde x_j|_{\EE^m}\ge |x_i- x_j|_{\spc{A}}\]
for all $i$ and $j$.
Passing to a limit configuration as $p'\to p$, we get the result.

If $\spc{A}$ is an $\Alex{-1}$ space, then the same argument gives
$\tilde p,\tilde x_1,\z\dots,\tilde x_n\in \EE^m$ such that
\[|\tilde p-\tilde x_i|_{\EE^m}=| p- x_i|_{\spc{A}}
\quad\text{and}\quad
\side\hinge{\tilde p}{\tilde x_i}{\tilde x_j}_{\HH^2}\ge |x_i- x_j|_{\spc{A}}.\]
It remains to observe that $\EE^m$ with the metric defined by
$\dist{\tilde v}{\tilde w}{}\z\df \side\hinge{\tilde p}{\tilde v}{\tilde w}_{\HH^2}$
is isometric to $\HH^m$.

The $\Alex{1}$ case is analogous, but one has to pull $x_i$ toward $p$
to ensure that $\dist{p}{x_i}{}<\pi$ for each $i$, and then pass to an additional limit.
\qeds

\begin{thm}[!]{Exercise}\label{ex:dim-max}
Let $\spc{A}$ be an $\Alex0$ space.
Assume that $\EE^m$ is isometric to a subspace of the tangent space $\T_p\spc{A}$ at \textit{some} point $p\in \spc{A}$ and $\LinDim \spc{A}\z= m<\infty$.
Show that $\T_p\spc{A}\iso \EE^m$.
\end{thm}

\begin{thm}[!]{Exercise}\label{ex:dim=1}
Show that a 1-dimensional Alexandrov space is homeomorphic to a 1-dimensional manifold, possibly with nonempty boundary.
\end{thm}

\begin{thm}{Exercise}\label{ex:resporka}
Let $\spc{A}$ be an $\Alex0$ space.
Show that $\LinDim \spc{A}\z\ge m$ if and only if there are $m+2$ points $p$, $a_0,\z\dots, a_{m}\in \spc{A}$
such that $\angk p{a_i}{a_j}>\tfrac\pi2$ for any pair $i\ne j$.
\end{thm}

\section{Space of directions}

Recall that a metric space $\spc{X}$ is \index{geodesic space@$\ell$-geodesic space}\emph{$\ell$-geodesic}
if any two points $x,y\in\spc{X}$ such that $\dist{x}{y}{}<\ell$ can be connected by a geodesic.

\begin{thm}{Theorem}\label{thm:finite-space-of-directions}
For any point $p$ in a finite-dimensional Alexandrov space, the space of directions $\Sigma_p$ is compact and $\pi$-geodesic.
\end{thm}

By \ref{ex:GHto-tangent}, this theorem implies the following.

\begin{thm}{Corollary}\label{ex:GHto-tangent-finite-dim}
Let $p$ be a point in a finite-dimensional Alexandrov space $\spc{A}$,
and let $\lambda_n\to\infty$.
Then there is a pointed Gromov--Hausdorff convergence $(\lambda_n\cdot \spc{A},p)\z\to (\T_p,0)$.
\end{thm}

\parit{Proof of \ref{thm:finite-space-of-directions}.}
Choose $\eps>0$; suppose $\spc{A}$ is $m$-dimensional.
Assume $\pack_\eps \Sigma_p\ge n$;
that is, we can choose $n$ directions $\xi_1,\dots, \xi_n\in \Sigma_p$ such that $\mangle(\xi_i,\xi_j)\z>\eps$ for any $i\ne j$.
Without loss of generality, we may assume that each direction is geodesic;
that is, there is a point $x_i\in \spc{A}$ such that $\xi_i=\dir p{x_i}$.

Choose a small $r>0$ and $y_i\in [px_i]$ such that $\dist{p}{y_i}{}=r$ for each $i$.
Since $r$ is small, we can assume that $\angk p{y_i}{y_j}>\eps$ for $i\ne j$.
By \ref{cor:euclid-subcone}, we can choose $p'$ arbitrarily close to $p$ such that $\dir{p'}{y_i}\in \Lin_{p'}$ for any $i$;
in particular, all directions $\dir{p'}{y_i}$ lie in a unit sphere of dimension at most $m-1$.
Since $\dist{p'}{p}{}$ is small, $\angk {p'}{y_i}{y_j}>\eps$ if $i\ne j$.
By comparison,
\[\mangle \hinge{p'}{y_i}{y_j}>\eps.\]
Therefore, $\pack_\eps \Sigma_p\le \pack_\eps\SSS^{m-1}$.

Since $\SSS^{m-1}$ is compact, $\pack_\eps\SSS^{m-1}<\infty$ for any $\eps>0$.
By definition, the space of directions $\Sigma_p$ is complete.
Applying \ref{ex:pack-net}, we get that $\Sigma_p$ is compact.

It remains to prove the following claim.

\begin{clm}{}
If $\Sigma_p$ is compact, then it is $\pi$-geodesic.
\end{clm}

Choose two geodesic directions $\xi=\dir px$ and $\zeta=\dir py$;
let
\[\alpha\z=\tfrac12\cdot \mangle \hinge pxy=\tfrac12\cdot \dist{\xi}{\zeta}{\Sigma_p}.\]
Assume $\alpha<\tfrac\pi2$.

It is sufficient to construct an \index{almost midpoint}\emph{almost midpoint} $\dir pz$ of $\xi$ and $\zeta$ in $\Sigma_p$;
that is, we need to show that for any $\eps>0$ there is a geodesic $[pz]$ such that
\[\mangle\hinge pxz\le \alpha+\eps
\quad\text{and}\quad
\mangle\hinge pyz\le \alpha+\eps.\]
Indeed, once this is done, the compactness of $\Sigma_p$ can be used to get an actual midpoint for any two directions in $\Sigma_p$, and Menger's lemma (\ref{lem:mid>geod}) will finish the proof.

Choose a sequence of small positive numbers $r_n\to0$.
Consider sequences $x_n\z\in [px]$ and $y_n\z\in [py]$ such that $\dist{p}{x_n}{}=\dist{p}{y_n}{}=r_n$.
Let $m_n$ be a midpoint of $[x_n\,y_n]$.
Since $\alpha<\tfrac\pi2$, $m_n\ne p$.

Since $\Sigma_p$ is compact, we can pass to a subsequence of $r_n$ such that
$\dir{p}{m_n}$ converges;
denote the limit by $\mu$.
Choose a geodesic $[pz]$ that runs at a small angle to $\mu$.
It remains to show that $\dir pz$ is the needed almost midpoint.

Evidently, $\angk p{x_n}{m_n}=\angk p{y_n}{m_n}$.
By \ref{ex:alex-lemma-cat}, we have
\[\angk p{x_n}{m_n}+\angk p{y_n}{m_n}\le \angk p{x_n}{y_n}.\]

For each large $n$, choose a point $z_n\in [pz]$ such that $\dist{p}{z_n}{}=\dist{p}{m_n}{}$.
By construction, for all large $n$, we have $\mangle\hinge pz{m_n}\approx0$ with any prescribed error.
By comparison, the value $\frac{\dist{z_n}{m_n}{}}{\dist{p}{z_n}{}}$ can be assumed to be arbitrarily small for all large $n$.
Applying this observation and the definition of angle measure, we also get that the following approximations
\begin{align*}
\angk p{x_n}{y_n}&\approx \mangle\hinge p{x_n}{y_n},
\\
\angk p{x_n}{m_n}&\approx\angk p{x_n}{z_n}\approx\mangle\hinge p{x_n}{z_n},
\\
\angk p{m_n}{y_n}&\approx\angk p{z_n}{y_n}\approx\mangle\hinge p{z_n}{y_n},
\end{align*}
hold with any prescribed error for all sufficiently large $n$.
It follows that $\dir pz$ is an almost midpoint of $\dir px$ and $\dir py$, as required.
\qeds

\parit{Remark.}
In the above proof, the angles $\mangle\hinge pxz$ and $\mangle\hinge pyz$ are bounded below by
comparison, but we needed upper bounds.
These were extracted from the definition of angle measure and the compactness of the space of directions.
Halbeisen's example \cite[13.6]{alexander-kapovitch-petrunin2024} shows that it cannot be done without the compactness condition.

\begin{thm}[!]{Exercise}\label{ex:concave-differential}
Let $f$ be a semiconcave function defined in an open set $\Dom f$ in a finite-dimensional Alexandrov space.
Show that the differential $\dd_pf$ is a concave function $\T_p\to\RR$ for any $p\in \Dom f$.
\end{thm}

\section{Right-inverse theorem}

\begin{thm}{Right-inverse theorem}\label{thm:right-inverse}\index{right-inverse theorem}
Suppose $p,a_0,\dots,a_m$ are points in an Alexandrov space $\spc{A}$ such that
\[\angk p{a_i}{a_j}>\tfrac\pi2\]
for any $i\ne j$.
Then the map $f\:\spc{A}\to\RR^m$ defined by
\[f\:x\mapsto (\dist{a_1}{x}{},\dots,\dist{a_m}{x}{})\]
has a continuous right inverse $\map$ defined in a neighborhood of $f(p)$ such that $p=\map\circ f(p)$.
\end{thm}

In the proof we construct a local right inverse $\map$ of $f$ around $f(p)$.
The construction uses gradient flow for a suitably chosen family of functions.
The proof can be extracted from the following exercise;
more details are given in the hints.

\begin{thm}[!]{Exercise}\label{ex:proof-right-inverse}
Suppose $p,a_0,\dots,a_m\in\spc{A}$ and $f\:\spc{A}\to\RR^m$ are as in \ref{thm:right-inverse}.
Assume $\eps>0$ is sufficiently small.
Given $\bm{y}\z=(y_1,\z\dots,y_m)\in \RR^m$,
consider the function on $\spc{A}$ defined by
\[h_{\bm{y}}(x)\df\min\{\,0, \dist{a_1}{x}{}-y_1,\dots,\dist{a_m}{x}{}-y_m\,\}+\eps\cdot\dist{a_0}{x}{}.\]

\begin{subthm}{ex:proof-right-inverse:grad}
Show that for some fixed $r>0$ and $\lambda$, the function $h_{\bm{y}}$ is $\lambda$-concave in $\oBall(p,r)$, and the following hold at any $x\in \oBall(p,r)$:
\begin{enumerate}[(i)]
\item\label{111} $(\dd_x\distfun_{a_i})(\nabla_x h_{\bm{y}})<-\eps^2$ if $\dist{a_i}{x}{}>y_i$ and
\item\label{222} $(\dd_x\distfun_{a_i})(\nabla_x h_{\bm{y}})>\eps^2$ if
\[\dist{a_i}{x}{}-y_i=\min_{j>0}\{\dist{a_j}{x}{}\z-y_j\}<0.\]
\end{enumerate}

\end{subthm}

\begin{subthm}{ex:proof-right-inverse:alpha}
Let $\alpha_{\bm{y}}$ be an $h_{\bm{y}}$-gradient curve that starts at $p$.
Use \ref{SHORT.ex:proof-right-inverse:grad} to show that
\[
f[\alpha_{\bm{y}}(t_0)]
=
\bm{y}\]
if
$\tfrac1{\eps^2}\cdot|f(p)-\bm{y}|
\le
t_0
\le
\tfrac{r}{2}$.

\end{subthm}

\begin{subthm}{ex:proof-right-inverse:end}
Let $t_0\:\bm{y}\mapsto\tfrac{1}{\eps^2}\cdot\bigl|f(p)-\bm{y}\bigr|$.
Use \ref{lem:fg-dist-est} to show that the map
\[\map\:{\bm{y}}\mapsto \alpha_{\bm{y}}\circ t_0(\bm{y})\]
is continuous (in fact Hölder) on $\Omega=\oBall(f(p),\tfrac{\eps^2}{2}\cdot r)\subset\RR^m$
and $f\circ \Phi(\bm{y})=\bm{y}$ for any $\bm{y}\in \Omega$.
(This finishes the proof of \ref{thm:right-inverse}.)
\end{subthm}

\end{thm}

\begin{thm}{Corollary}\label{cor:tan-g-delta}
Let $\spc{A}$ be an Alexandrov space, and let $m$ be a positive integer.
If $\LinDim\spc{A}\ge m$, then
\[\dim\Lin_s\z\ge m\]
for any point $s$ in a dense G-delta subset of points in $\spc{A}$.

Moreover, if $m=\LinDim\spc{A}$, then the set $S$ of points $s\in \spc{A}$ such that $\T_s\spc{A}\iso \EE^m$ is a dense G-delta set in $\spc{A}$.
\end{thm}

It is worth mentioning that the set $S\subset\spc{A}$ is also convex; see \cite[1.10]{petrunin1998}.

\parit{Proof.}
Assume $\spc{A}$ is $\Alex0$.
By the right-inverse theorem and \ref{ex:resporka} there is a ball $B=\cBall[p,r]_{\spc{A}}$ and $\Const>0$ such that
\[\pack_\eps B\ge \frac{\Const}{\eps^m}\]
for all small $\eps>0$.

In particular, given a radius $R$ we can choose $\eps>0$ and an integer $N$ such that
\[\pack_\eps \cBall(0,R+r)_{\EE^{m-1}}< N\le \pack_\eps B.\]

Applying \ref{cor:euclid-subcone} to an $\eps$-packing $x_1,\dots,x_N$ in $B$, we get that \[\dim\Lin_s\z\ge m\]
for any point $s$ in a dense G-delta subset of $\oBall(p,R)$.
Since $R>0$ is arbitrary, the first statement follows.

To prove the last statement, apply \ref{ex:dim-max}.

The general case can be reduced to $\Alex{-1}$ using rescaling,
and it can be done essentially in the same way; one only has to change $\EE^{m-1}$ to $\HH^{m-1}$.
\qeds

\section{Distance chart}

\begin{thm}{Theorem}\label{thm:dist-chart}
Let $p,a_0,\dots,a_m$ be points in an $m$-dimensional Alexandrov space $\spc{A}$ such that
\[\angk p{a_i}{a_j}>\tfrac\pi2\]
for any $i\ne j$.
Then the map $f\:\spc{A}\to\RR^m$ defined by
\[f\:x\mapsto (\dist{a_1}{x}{},\dots,\dist{a_m}{x}{})\]
gives a bi-Lipschitz embedding of a neighborhood $\Omega$ of $p$;
the restriction $f|_\Omega$ is called a \index{distance chart}\emph{distance chart} at $p$.
\end{thm}

The following exercise provides a guide to the proof of the theorem.

\begin{thm}[!]{Exercise}\label{ex:proof-dist-chart}
Suppose $p,a_0,\dots,a_m\in\spc{A}$ and $f\:\spc{A}\to\RR^m$ are as in \ref{thm:right-inverse}.
Show that there is $\eps>0$ such that for any two points $x,y$ in a sufficiently small neighborhood of $p$,
one of the inequalities holds:
\begin{align*}
\mangle\hinge xy{a_1}&<\tfrac\pi2-\eps,\ \dots,\ \mangle\hinge xy{a_m}<\tfrac\pi2-\eps,
\\
\mangle\hinge yx{a_1}&<\tfrac\pi2-\eps,\ \dots,\ \mangle\hinge yx{a_m}<\tfrac\pi2-\eps.
\end{align*}
Use this together with the right-inverse theorem (\ref{thm:right-inverse}) to prove \ref{thm:dist-chart}.
\end{thm}

\section{Volume}

Fix a positive integer $m$.
The $m$-dimensional Hausdorff measure of a Borel set $B$ in a metric space will be called its \index{volume}\emph{$m$-volume}; it will be denoted by $\vol_m B$.
We assume that the Hausdorff measure is calibrated so that the unit cube in $\EE^m$ has unit volume.

This definition will be applied mostly to subsets of $m$-dimensional Alexandrov spaces.
In this case, we may write $\vol B$ instead of $\vol_m B$.

\begin{thm}{Bishop--Gromov inequality}\label{inq:BG}\index{Bishop--Gromov inequality}
Let $\spc{A}$ be an $\Alex0$ space.

\begin{subthm}{inq:BG:a}
For any integer $k\ge 0$, any point $p\in \spc{A}$ and any radius $r>0$, we have
\[\vol_k \oBall(p,r)_{\spc{A}}\le \vol_k \oBall(0,r)_{\T_p\spc{A}}.\]
\end{subthm}

\smallskip

\noindent Further, suppose $\LinDim \spc{A}=m<\infty$.

\smallskip

\begin{subthm}{inq:BG:b}
For any point $p\in \spc{A}$ and any radius $r>0$, we have
\[\vol_m \oBall(p,r)\le \omega_m\cdot r^m,\]
where $\omega_m$ denotes the volume of the unit ball in $\EE^m$.
\end{subthm}

\begin{subthm}{inq:BG:c}
Moreover, the function
\[r\mapsto \frac{\vol_m \oBall(p,r)}{r^m}\]
is nonincreasing.
\end{subthm}

\end{thm}

\parit{Proof.}
Given $x\in\spc{A}$, choose a geodesic path $\gamma_x$ from $p$ to $x$.
Recall that $\gamma_x^+(0)=\log_p x$ (here $\log_p\:\spc{A}\to \T_p$ is a multivalued function; so, $\gamma_x^+(0)$ might not be the only value of $\log_p x$).
By comparison, $\log_p$ is distance-noncontracting.
Note that $\log_p(\oBall(p,r)_{\spc{A}})\subset \oBall(0,r)_{\T_p}$;
this proves part \ref{SHORT.inq:BG:a}.

\begin{wrapfigure}{r}{44 mm}
\vskip-4mm
\centering
\includegraphics{mppics/pic-803}
\vskip1mm
\end{wrapfigure}

If $\T_p\iso \EE^m$, then \ref{SHORT.inq:BG:a} implies \ref{SHORT.inq:BG:b}.
In the general case, by \ref{cor:tan-g-delta}, we can find a point $p'$ arbitrarily close to $p$ such that $\T_{p'}\iso \EE^m$.
If $\eps>\dist{p}{p'}{}$, then $\oBall(p,r)\z\subset \oBall(p',r+\eps)$.
Therefore,
\[\vol \oBall(p,r)\le \omega_m\cdot (r+\eps)^m\]
for any $\eps>0$.
Hence \ref{SHORT.inq:BG:b} follows.

Now, suppose $0<r_1<r_2$.
Consider the map $w\: \spc{A}\selfmap$ defined by $w\:x\mapsto \gamma_x(\tfrac {r_1}{r_2})$; it mimics the dilation with center at $p$ and coefficient $\tfrac {r_1}{r_2}$.
By comparison,
\[\dist{w(x)}{w(y)}{}\ge \tfrac {r_1}{r_2}\cdot \dist{x}{y}{}.\]
Observe that $\oBall(p,r_1) \supset w[\oBall(p,r_2)]$.
Therefore,
\[\vol \oBall(p,r_1)\ge (\tfrac {r_1}{r_2})^m\cdot\vol \oBall(p,r_2).\]
\qedsf

The following exercise generalizes the Bishop--Gromov inequality to the $\Alex{-1}$ case;
it is sufficient for most applications.
A sharper statement for $\Alex\kappa$ spaces is given in \ref{inq:BG+}.

\begin{thm}[!]{Exercise}\label{ex:BG}
Let $\spc{A}$ be a finite-dimensional $\Alex{-1}$ space;
so $\LinDim\spc{A}=m<\infty$.
Show that
\[\vol \oBall(p,r)\le \omega_m\cdot(\sinh r)^m,\]
where $\omega_m$ denotes the volume of the unit ball in $\EE^m$.
Moreover, the function
\[r\mapsto \frac{\vol \oBall(p,r)}{(\sinh r)^m}\]
is nonincreasing.
\end{thm}

\begin{thm}[!]{Exercise}\label{ex:diam-compact:proper}
Show that any finite-dimensional Alexandrov space is proper.

\end{thm}

\begin{thm}{Exercise}\label{ex:vol>0}
Show that any nonempty open set in an Alexandrov space of finite dimension $m$ has positive $m$-volume.
\end{thm}

\section{Packings}\label{sec:packing-numbers}

\begin{thm}{Exercise}\label{ex:tangent=Em}
Let $\spc{A}$ be an $m$-dimensional $\Alex\kappa$ space and $m\z<\infty$.

\begin{subthm}{ex:tangent=Em:2balls}
Show that for any two balls $B_1$ and $B_2$ in $\spc{A}$ with positive radii,
there is $\delta>0$ such that the inequality
$\pack_{\delta\cdot\eps} B_1\ge \pack_\eps B_2$
holds for all sufficiently small $\eps>0$.
\end{subthm}

\begin{subthm}{ex:tangent=Em:ball}
Show that there are constants $c_1=c_1(m)$ and $c_2=c_2(m,r,\kappa)$ such that
\[c_1\cdot \vol B\le\eps^m\cdot\pack_\eps B\le c_2\cdot \vol B\]
for any ball $B=\cBall[p,r]_{\spc{A}}$ with $r>0$ and all small $\eps>0$.
\end{subthm}

\end{thm}

The exercise implies that for Alexandrov spaces, \index{Minkowski dimension}\emph{Minkowski dimension} (also known as \index{box-counting dimension}\emph{box-counting dimension}) coincides with linear dimension.
Namely, let $\spc{A}$ be an Alexandrov space.
Given a point $p\in \spc{A}$ and $r>0$, let
\[m(p,r)\df \inf\set{\alpha\in\RR}{\eps^\alpha\cdot\pack_\eps \cBall[p,r]\to 0\quad\text{as}\quad \eps\to 0}.\]
Then $m(p,r)=\LinDim \spc{A}$; in particular, $m(p,r)$ is an integer, and it does not depend on the choice of $p$ and $r$.

\section{Other dimensions}\label{sec:all-dim}

Now we want to show that \textit{all reasonable definitions of dimension give the same result for Alexandrov spaces}.
More precisely, we have the following theorem; compare with \cite[15.16]{alexander-kapovitch-petrunin2024}.
We refer to \cite{hurewicz-wallman} for definitions of \index{Lebesgue covering dimension}\emph{Lebesgue covering dimension} and \index{Hausdorff!dimension}\emph{Hausdorff dimension}, which will be denoted by \index{TopDim@$\TopDim$}$\TopDim$ and
\index{HausDim@$\HausDim$}$\HausDim$, respectively.

\begin{thm}{Theorem}\label{thm:dim=dim}
For any Alexandrov space $\spc{A}$, we have
\[\LinDim \spc{A}=\TopDim \spc{A}=\HausDim \spc{A}.\]
\end{thm}

\parit{Proof.}
Suppose $\LinDim \spc{A}=\infty$.
By the right-inverse theorem (\ref{thm:right-inverse}), $\spc{A}$ contains a compact subset $K$ with arbitrarily large $\TopDim K$.
In particular,
\[\TopDim\spc{A}=\infty.\]
By Szpilrajn's theorem,
$\HausDim K\ge \TopDim K$.
Thus we also have
\[\HausDim\spc{A}=\infty.\]

Now suppose $\LinDim \spc{A}=m<\infty$.
By the Bishop--Gromov inequality (\ref{inq:BG} and \ref{ex:BG}),
\[\HausDim\spc{A}\le m.\]

Since $\spc{A}$ is proper (\ref{ex:diam-compact:proper}),
Szpilrajn's theorem implies that
\[\TopDim\spc{A}\le \HausDim\spc{A}\le m.\]
Finally, by the right-inverse theorem (\ref{thm:right-inverse}), $m\le\TopDim\spc{A}$.
\qeds

\begin{thm}{Exercise}\label{ex:dim=dim}
Let $\Omega$ be an open subset of an Alexandrov space $\spc{A}$.
Show that
\[\LinDim \spc{A}=\LinDim \Omega=\TopDim \Omega=\HausDim \Omega.\]
\end{thm}

\begin{thm}[!]{Exercise}\label{ex:finite-tan}
Let $p$ be a point in a finite-dimensional Alexandrov space $\spc{A}$.
Prove the following:
\begin{subthm}{ex:finite-tan:tan}
The tangent space $\T_p$ is a proper $\Alex0$ space.
\end{subthm}

\begin{subthm}{ex:finite-space-of-directions-dim}
$\LinDim\Sigma_p+1=\LinDim\T_p=\LinDim\spc{A}$.
\end{subthm}

\begin{subthm}{ex:finite-tan:sigma}
If $\LinDim \spc{A}>1$, then $\Sigma_p$ is geodesic.
\end{subthm}

\end{thm}

Part \ref{SHORT.ex:finite-space-of-directions-dim} of the exercise opens a way to use induction on dimension, which can be used for finite-dimensional Alexandrov spaces;
in Lecture~\ref{chap:bry}, we will use it to define boundary.
The condition $\LinDim \spc{A}>1$ in \ref{SHORT.ex:finite-tan:sigma} is necessary;
indeed $\Sigma_0\RR$ is a two-point space which is not geodesic.

\section{Remarks}

Halbeisen's example \cite{alexander-kapovitch-petrunin2024,halbeisen} shows that compactness of the space of directions is essential in the proof that the space of directions is $\pi$-geodesic (\ref{thm:finite-space-of-directions}).

The proof of the following version of the Bishop--Gromov inequality additionally requires the so-called \textit{coarea formula} for Alexandrov spaces.
The weaker inequality from \ref{ex:BG} is sufficient for the sequel.

\begin{thm}{Optimal Bishop--Gromov inequality}\label{inq:BG+}\index{Bishop--Gromov inequality}
Given a point $p$ in an $m$-dimensional $\Alex\kappa$ space,
consider the function
\[v\:r\mapsto\vol_m\oBall(p,r);\]
denote by $\tilde v(r)$ the volume of an $r$-ball in $\MM^m(\kappa)$.
Then
\[v(r)\le \tilde v(r)\]
for any $r>0$.
Moreover, the function
\[r\mapsto \frac{v(r)}{\tilde v(r)}\] is nonincreasing.
If $\kappa>0$, then one has to assume that $r<\tfrac\pi{\sqrt\kappa}$.
\end{thm}

This inequality was originally proved for Riemannian manifolds with Ricci curvature bounded below.
The first part is also known as \index{Bishop's inequality}\emph{Bishop's inequality}; it is due to Richard Bishop \cite{bishop1964}, \cite[Corollary 4, p. 256]{bishop-crittenden}.
The second part is due to Michael Gromov \cite{gromov1981}.

Theorem~\ref{thm:dim=dim} was essentially proved by Conrad Plaut \cite{plaut:dimension}.
At that time, it was not known whether
\[\LinDim\spc{A}=\infty\quad \Rightarrow\quad \TopDim\spc{A}=\infty\]
for any Alexandrov space $\spc{A}$.
The latter implication was proved by Grigory Perelman and the third author \cite{perelman-petrunin:qg}; our proof is slightly different.

%% file: volume.tex
\chapter{Limit spaces}\label{chap:lim}\label{chap:stability}

This lecture discusses the main source of applications of Alexandrov geometry.
We start with an observation that the curvature bound of Alexandrov spaces survives under a Gromov--Hausdorff limit.
Further, we present Perelman's construction of strictly concave functions and prove Gromov's selection theorem.
We then use all this to prove the homotopy stability theorem (\ref{thm:h-stability}) and the homotopy finiteness theorem (\ref{thm:h-finiteness}).

\section{Survival of curvature bounds}

\begin{thm}{Theorem}\label{thm:CBB-closed}
Let $\spc{X}_n\z\to \spc{X}_\infty$ be a convergence in the sense of Gromov--Hausdorff.
Suppose that for each $n$, the space $\spc{X}_n$ has curvature $\ge\kappa$ in the sense of Alexandrov.
Then the same holds for~$\spc{X}_\infty$.
\end{thm}

\parit{Proof.}
Let us verify the 4-point comparison for a quadruple of points $p_\infty$, $x_\infty$, $y_\infty$, $z_\infty\in \spc{X}_\infty$.
We can assume that all the comparison angles $\angk {p_\infty}{x_\infty} {y_\infty}_{\MM^2(\kappa)}, \angk {p_\infty}{y_\infty}{z_\infty}_{\MM^2(\kappa)}, \angk {p_\infty} {z_\infty}{x_\infty}_{\MM^2(\kappa)}$ are defined, since otherwise there is nothing to check.

By the definition of Gromov--Hausdorff convergence, we can choose points $p_n$, $x_n$, $y_n$, $z_n\in \spc{X}_n$ that converge to $p_\infty$, $x_\infty$, $y_\infty$, $z_\infty\in \spc{X}_\infty$, respectively.
In particular, the 6 distances between $p_n$, $x_n$, $y_n$, $z_n$ converge to the distances between the corresponding pairs of $p_\infty$, $x_\infty$, $y_\infty$, $z_\infty$.

Since $\MM^2(\kappa)$-comparison holds for $p_n$, $x_n$, $y_n$, $z_n\z\in \spc{X}_n$,
passing to the limit, we get the $\MM^2(\kappa)$-comparison for $p_\infty$, $x_\infty$, $y_\infty$, $z_\infty$.
\qeds

\begin{thm}[!]{Exercise}\label{ex:dim-lim}
A sequence $\spc{A}_1,\spc{A}_2,\dots$ of $\Alex\kappa$ spaces converges to $\spc{A}_\infty$ in the sense of Gromov--Hausdorff.
Show that $\spc{A}_\infty$ is $\Alex\kappa$ and
\[\dim \spc{A}_\infty\le \liminf_{n\to\infty} \dim \spc{A}_n.\]
\end{thm}

\section{Gromov's selection theorem}

\begin{thm}{Theorem}\label{thm:gromov-compactness}\index{Gromov's selection theorem}
Let $m$ be a positive integer, and let $D,\kappa\in\RR$.
Then any sequence of $m$-dimensional $\Alex\kappa$ spaces with diameters at most $D$
has a converging subsequence in the sense of Gromov--Hausdorff.
\end{thm}

Recall that $\GH$ denotes the space of compact metric spaces with the Gromov--Hausdorff metric; see Section~\ref{sec:Gromov--Hausdorff-metric}.

\parit{Proof of \ref{thm:gromov-compactness}.}
Suppose $\kappa=0$.
Denote by $\bm{K}$ the set of all isometry classes of $\Alex0$ spaces with dimension $\le m$ and diameter $\le D$.
By \ref{ex:dim-lim}, $\bm{K}$ is a closed subset of $\GH$.

Choose a space $\spc{A}\in \bm{K}$;
suppose $x_1,\dots,x_n\in \spc{A}$ are points such that $\dist{x_i}{x_j}{}> \eps$ for all $i\ne j$.
Note that the balls $B_i=\oBall(x_i,\tfrac\eps2)$ do not overlap.

By \ref{thm:right-inverse} (or \ref{ex:vol>0}), $\vol \spc{A}>0$.
By the Bishop--Gromov inequality, $\vol \spc{A}\z<\infty$,
and if $\eps<D$, then 
\[\vol B_i\ge (\tfrac\eps{2\cdot D})^m\cdot\vol \spc{A}\]
for any $i$.
It follows that $n\le (\tfrac{2\cdot D}\eps)^m$;
that is, 
\[\pack_\eps\spc{A}\le N(\eps)\df(\tfrac{2\cdot D}\eps)^m\]
for all small $\eps>0$.

Choose a maximal $\eps$-packing $x_1,\z\dots,x_n\in \spc{A}$.
By \ref{ex:pack-net}, $\spc{F}_\eps\z\df\{x_1,\z\dots,x_n\}$ is an $\eps$-net of $\spc{A}$.
Observe that $\dist{\spc{F}_\eps}{\spc{A}}{\GH}\le \eps$.

Let $\bm{F}_\eps$ be the set of all metric spaces with
at most $N(\eps)$ points
and
diameter $\le D$.
Note that the set $\bm{F}_\eps$ is a compact subset of $\GH$.

Summarizing, for any $\eps>0$ we can find a compact $\eps$-net $\bm{F}_\eps\subset \GH$ of $\bm{K}$.
Since $\GH$ is complete (\ref{prop:complete}), it remains to apply \ref{ex:net:compact}.

We finished the proof of the case $\kappa=0$.
In the general case, applying rescaling, we can assume that $\kappa=-1$ and then argue as before, using \ref{ex:BG} instead of \ref{inq:BG}.
\qeds

\begin{thm}{Exercise}\label{ex:diam-compact:GH}
Show that any sequence of $m$-dimensional $\Alex\kappa$ spaces with marked points contains a convergent subsequence in the sense of pointed Gromov--Hausdorff convergence (see Section~\ref{sec:Gromov--Hausdorff}).

\end{thm}

\begin{thm}[!]{Exercise}\label{ex:pack-vol:dim}
Show that if $\spc{A}_n$ is a sequence of $m$-dimensional $\Alex0$ spaces with diameter $\le D$, and volume $\ge v_0>0$, then its Gromov--Hausdorff limit $\spc{A}_\infty$ (if it exists) has dimension~$m$.
\end{thm}

\section{Controlled concavity}

Alexandrov spaces have plenty of semiconcave functions, including squared distance functions and their convex combinations.
The following theorem provides a source of strictly concave functions defined on small open sets of finite-dimensional Alexandrov spaces.

\begin{thm}{Theorem}
\label{thm:strictly-concave}
Let $\spc{A}$ be a complete finite-dimensional Alexandrov space.
Then for any point $p\in \spc{A}$, there is a strictly concave function $f$ defined in an
open neighborhood of $p$.

Moreover, given a nonzero vector $v\in \T_p$, the differential, $\dd_p f$, can be chosen
arbitrarily close to $x\mapsto -\<v,x\>$.
\end{thm}

\parit{Idea of the proof.}
Suppose $f$ is a semiconcave function.
One can improve concavity of $f$ by passing to a composition $\phi\circ f$ with a very concave real-to-real function $\phi$ such that $\phi'\approx 1$.
Namely, if $\gamma$ is a geodesic such that $|(f\z\circ\gamma)^\pm|$ is bounded away from zero, then $\phi\circ f\circ\gamma$ can be made as concave as you wish by choosing an appropriate $\phi$.
For smooth real-to-real functions $f$ and $\phi$ this follows from the chain rule $(\phi\circ f)''(x)=\phi''(f(x))\cdot (f'(x))^2+\phi'(f(x)) \cdot f''(x)$.
This trick does not help if $|(f\circ\gamma)^\pm|\approx 0$, but it does not make concavity worse;
so $\phi\circ f$ is very concave in most directions and at least as concave as $f$ in all directions.

In the proof we consider a sum $\phi\circ f_1+\ldots+\phi\circ f_N$ for a suitably chosen sequence of functions $f_1,\ldots,f_N$, so that for any geodesic $\gamma$ in a small neighborhood of $p$ the value $|(f_i\circ\gamma)^\pm|$ is bounded away from zero for some $i$.

\medskip

Before getting into the proof, we need to refresh material from \ref{Function comparison} and \ref{sec:packing-numbers}.

\begin{thm}{Lemma}\label{lem:barrier}
Let $f$ be a locally Lipschitz function defined on a real interval $\II$.
Then $f$ is $\lambda$-concave if and only if for any $t_0\in\II$ there is a smooth (local) upper barrier
$\bar f$ of $f$ at $t_0$ such that $\bar f''(t_0)\le \lambda$.

Moreover, in the only-if part, we can assume, in addition, that $\bar f'(t_0)$ takes a given value in the interval $[f^+(t_0),f^-(t_0)]$ (by \ref{ex:concave'}, this interval is not empty).
\end{thm}

\parit{Proof.}
Passing to the function $t\mapsto f(t)-\lambda\cdot\tfrac{t^2}2$, we can assume that $\lambda=0$.

The only-if part follows since a concave function is supported by the linear function $t\mapsto f(t_0)+ \Const\cdot (t-t_0)$ for any $\Const\in [f^+(t_0),f^-(t_0)]$.

On the other hand, suppose $f$ is not concave.
Then adding a linear function to $f$ we can assume that there is a subinterval $[a,b]\subset\II$ such that
\[
f(a)>0,
\quad
f(b)>0,
\quad\text{but}\quad
f(t)<0
\quad\text{for some}\quad
t\in [a,b].
\]
Let $s$ be the minimal value such that $f(t)\ge \tfrac s2\cdot (t-a)\cdot(t-b)$ for any $t\in [a,b]$.
Note that $s>0$ and $f(t_0)=\tfrac s2\cdot (t_0-a)\cdot(t_0-b)$ for some $t_0\in (a,b)$;
that is, $t\mapsto \tfrac s2\cdot (t-a)\cdot(t-b)$ is a local lower barrier of $f$ at $t_0$.
Therefore, $\bar f''(t_0)\ge s>0$ for any upper barrier $\bar f$ of $f$ at $t_0$,
and the if part follows.
\qeds

\begin{thm}[!]{Exercise}\label{ex:barrier}
Let $f$ be a Lipschitz real-to-real function, and let $\phi$ be a smooth increasing function defined in a neighborhood of $x_0\z=f(t_0)$.
Suppose that $f$ is $\lambda_1$-concave for some $\lambda_1\ge 0$ and there are some constants $c_1,c_2,c_3\in\RR$ such that (1) $|f^+(t_0)|>c_1>0$ (or $|f^-(t_0)|>c_1>0$), (2) $c_2>\phi'(x_0)\ge 0$, and (3) $0>c_3\ge \phi''(x_0)$.
Let $\lambda_2=c_2\cdot\lambda_1+c_1^2\cdot c_3$.
Show that $\phi\circ f$ has a smooth upper barrier $u$ at $t_0$ with
$u''(t_0)\le \lambda_2$.

\end{thm}

\begin{thm}[!]{Exercise}\label{ex:pack-sphere}
Let $\Omega$ be a nonempty open set in an $m$-dimensional Alexandrov space $\spc{A}$; assume $1\le m<\infty$.
Given a point $p\in\spc{A}$, there are $r>0$ and $\Const>0$ such that for any sufficiently small $\eps>0$ there are at least $\Const/\eps^{m-1}$ points $q_1,\dots,q_N\in \Omega$ such that $\dist{p}{q_i}{}=r$ for any $i$ and $\mangle\hinge p{q_i}{q_j}>\eps$ if $i\ne j$.
\end{thm}

\parit{Proof of \ref{thm:strictly-concave}.}
Suppose $r>0$ is small and $c$ is large (these values will be specified a bit later).
Consider the real-to-real function
$$\phi_{r,c}(x)=(x-r)- c\cdot(x-r)^2/r;$$
so,
\begin{align*}
\phi_{r,c}(r)&=0,
&
\phi_{r,c}'(r)&=1,
&
\phi_{r,c}''(r)&=-\frac{2\cdot c}{r}.
\end{align*}

{

\begin{wrapfigure}{o}{44 mm}
\vskip-4mm
\centering
\includegraphics{mppics/pic-901}
\vskip1mm
\end{wrapfigure}

Let $\gamma$ be a geodesic.
Given points $q_1,\dots, q_N$, set
\begin{align*}
\alpha_i(t)&=\mangle(\gamma^+(t),\dir{\gamma(t)}{q_i}),
\\
\lambda_i(t)&=\frac{3-c\cdot \cos^2[\alpha_i(t)]}{r},
\\
h_i&=\phi_{r,c}\circ\distfun_{q_i}\circ\gamma.
\end{align*}

}

Suppose that $\dist{q_i}{\gamma(t_0)}{}$ is sufficiently close to
$r$.
By \ref{ex:barrier}, $h_i$ has a smooth upper barrier $u_i$ at $t_0$ such that $u_i''(t_0)\le \lambda_i(t_0)$.

Assume that $\dist{p}{q_i}{}=r$ for each $i$, and $\eps_1>0$ is small.
If the geodesic $\gamma$ lies in $\oBall(p,\eps_1)$, then from above we have that $\sum h_i$ has a smooth upper barrier $u=\sum u_i$ at $t_0$ such that $u''(t_0)\le \sum \lambda_i(t_0)$.

Suppose that there is $\eps_2>0$ such that
\[\max_i\{|\alpha_i(t_0)-\tfrac\pi2|\}\ge\eps_2>0\eqlbl{eq:angle-ne-pi/2}\]
for any geodesic $\gamma$ in $\oBall(p,\eps_1)$ and every $t_0$
in the domain of $\gamma$.
If $c>3\cdot N/\sin^2\eps_2$, then $\sum\lambda_i(t_0)<0$ for any $t_0$.
By \ref{lem:barrier}, $\sum h_i$ is strictly concave.
It follows that
\[f=\sum_i \phi_{r,c}\circ\distfun_{q_i}\]
is strictly concave in $\oBall(p,\eps_1)$.

It remains to choose $r,\eps_1,\eps_2>0$ and points $q_1,\dots, q_N$ such that $\dist{p}{q_i}{}=r$ for each $i$, and \ref{eq:angle-ne-pi/2} holds.

Let us assume that $\spc{A}$ is $\Alex0$.
Applying \ref{ex:pack-sphere} to a small open ball around $p$ and small $\eps_3>0$,
we can choose points $q_1,\dots,q_N$ such that $N\ge \Const_1/\eps_3^{m-1}$ and $\dist{p}{q_i}{}=r$ for some small $r>0$
and $\angk p{q_j}{q_i}>\eps_3$ (here $\Const_1>0$).
On the other hand, suppose $\gamma$ runs from $x$ to~$y$.
Let us apply the ($n$+1)-point comparison (\ref{thm:n+1}) to $\gamma(t)$, $x$, $y$ and $q_1,\dots,q_N$;
denote the obtained points by $\tilde x, \tilde y,\tilde q_1,\dots,\tilde q_N$.
If $\eps_1\ll \eps_3$, $\eps_2\ll \eps_3$, and
$|\alpha_i(t)-\tfrac\pi2|<\eps_2$ for each $i$, then the points $\tilde q_1,\dots,\tilde q_N$ lie very close to an $r$-sphere in a hyperplane in $\EE^m$ that is perpendicular to $[\tilde x\tilde y]$;
in particular, $N\le \Const_2(m)/\eps_3^{m-2}$.
Therefore, for small $\eps_3>0$ and yet smaller $\eps_1$ and $\eps_2$, points $q_1,\dots,q_N$ meet \ref{eq:angle-ne-pi/2}.

To prove the last statement, choose a geodesic $[pz]$ in the direction sufficiently close to $v$ and apply \ref{ex:pack-sphere} to a small ball centered at a point $\bar z\in [pz]$ close to $p$.

The $\Alex\kappa$-case can be reduced to the $\Alex{-1}$-case by rescaling, and it is done the same way.
The only difference is that the points $\tilde q_1,\dots,\tilde q_N$ lie very close to an $r$-sphere in a hyperplane in $\HH^m$.
\qeds

Given a point $p$, a strictly concave function $f$ can be chosen so that
$f(p)=0$ and $(\dd_p f)(x)\approx-|x|$.
In particular, $p$ is the maximum point of $f$.
It can be constructed by taking the minimum of the functions in the theorem.
Then the set $K=\set{x\in\spc{A}}{f(x)\ge -\eps}$ forms a closed convex neighborhood of $p$ for any small $\eps>0$.
Replacing $f$ by $f+\eps$, we get the following.

\begin{thm}{Corollary}\label{cor:convex-nbhd}
Any point $p$ of a finite-dimensional Alexandrov space admits an arbitrarily small closed convex neighborhood $K$ and a strictly concave function $f$ defined in a neighborhood of $K$ such that $p$ is the maximum point of $f$
and $f|_{\partial K}=0$.
\end{thm}

\begin{thm}{Exercise}\label{ex:no-conc}
Construct an Alexandrov space $\spc{A}$ such that there is no nonempty open set $\Omega\subset \spc{A}$ with a strictly concave function $f$ defined on $\Omega$.
\end{thm}

\begin{thm}{Exercise}\label{ex:dist-chart+strictly-concave}
Let $p,a_0,\dots,a_m$ be points in an $m$-dimensional Alexandrov space $\spc{A}$ such that $\dist{p}{a_i}{}=r$ for any $i$ and $\angk p{a_i}{a_j}>\tfrac\pi2$
for any $i\ne j$, and let $\phi_{r,c}$ be as in the proof of \ref{thm:strictly-concave}, with $c$ large.
Show that the function $f=\phi_{r,c}\circ\distfun_{a_1}+\dots+\phi_{r,c}\circ\distfun_{a_m}$
is strictly concave in a neighborhood of $p$.
\end{thm}

\section{Liftings}\label{sec:Liftings}

Suppose that the Gromov--Hausdorff distance $\dist{\spc{A}}{\spc{A}'}{\GH}$ is sufficiently small, so we may think that both spaces $\spc{A}$ and $\spc{A}'$ lie at a small Hausdorff distance in an ambient metric space $\spc{W}$.
In particular, we may choose a small $\eps>0$, so that for any point $p\in \spc{A}$, there is a point $p'\in \spc{A}'$ such that $\dist{p}{p'}{\spc{W}}<\eps$;
the point $p'$ will be called a \index{lifting}\emph{lifting} (or {}\emph{$\eps$-lifting}) of $p$ in $\spc{A}'$.
We may choose a lifting $p'\in\spc{A}'$ for every point $p\in\spc{A}$;
in this case, the map $p\mapsto p'$ is called a {}\emph{($\eps$-)lifting map}.

Let us emphasize that liftings are not uniquely defined, and the lifting map is not assumed to be continuous.
Also, to talk about liftings, we have to choose $\eps>0$, as well as the inclusions $\spc{A},\spc{A}'\hookrightarrow\spc{W}$.

Let $\spc{A}$ be a compact $m$-dimensional $\Alex\kappa$ space.
Suppose $\spc{A}'$ is another compact $m$-dimensional $\Alex\kappa$ space such that $\dist{\spc{A}}{\spc{A}'}{\GH}$ is sufficiently small --- smaller than some $\eps=\eps(\spc{A})>0$.
If the construction of the previous section is performed in $\spc{A}$,
then it can be repeated in $\spc{A}'$ for the liftings of all points and the same function $\phi$.
It produces a strictly concave function defined in a controlled neighborhood of the lifting $p'$ of $p$.

The results of this and related constructions will also be called \index{lifting}\emph{liftings};
for example, we can talk about a lifting from $\spc{A}$ to $\spc{A}'$ of a function provided by \ref{thm:strictly-concave} (if the Gromov--Hausdorff distance $\dist{\spc{A}}{\spc{A}'}{\GH}$ is small, then these liftings are strictly concave)
and a lifting of a convex neighborhood from \ref{cor:convex-nbhd}.
Here one cannot use \ref{thm:strictly-concave} and \ref{cor:convex-nbhd} as black boxes --- one has to understand the construction, but it is straightforward.

\begin{thm}{Exercise}\label{ex:no-cont-lifting}
Give an example of Gromov--Hausdorff convergence of spaces $\spc{A}_n\z\to \spc{A}_\infty$ such that all $\spc{A}_n$ and $\spc{A}_\infty$ are compact finite-dimensional $\Alex0$ spaces and for any small $\eps>0$ there is no continuous $\eps$-lifting map $\spc{A}_\infty\to \spc{A}_n$ for any large $n$.
\end{thm}

\section{Nerves}

Let $\{\Omega_1,\dots,\Omega_k\}$ be a finite open cover of a compact metric space $\spc{X}$.
Consider an abstract simplicial complex $\spc{N}$ with one vertex $v_i$ for each set $\Omega_i$;
a simplex with vertices $v_{i_1},\dots,v_{i_m}$ is included in $\spc{N}$ if the intersection $\Omega_{i_1}\cap\dots\cap \Omega_{i_m}$ is nonempty.
\begin{figure}[ht!]
\vskip-0mm
\centering
\includegraphics{mppics/pic-1402}
\end{figure}
The obtained simplicial complex $\spc{N}$ is called the \index{nerve}\emph{nerve} of the covering $\{\Omega_i\}$.
Evidently, $\spc{N}$ is a finite simplicial complex ---
it is a subcomplex of a simplex with the vertices $\{v_1,\dots,v_k\}$.
Recall that $\Star_{v_i}$ denotes the union of all simplices in $\spc{N}$ that share the vertex $v_i$.

The next statement follows from \cite[4G.3]{hatcher}.

\begin{thm}{Nerve theorem}\label{thm:nerve}
Let $\{\Omega_1,\dots,\Omega_k\}$ be an open cover of a compact metric space $\spc{X}$
and let $\spc{N}$ be the corresponding nerve with vertices $\{v_1,\dots,v_k\}$.
Suppose that every nonempty finite intersection $\Omega_{\alpha_1}\cap\z\dots\cap\Omega_{\alpha_m}$ is contractible.
Then $\spc{X}$ is homotopy equivalent to the nerve $\spc{N}$ of the cover.

Moreover, homotopy equivalences $a\:\spc{X}\to \spc{N}$ and $b\:\spc{N}\to\spc{X}$ can be chosen so that
if $x\in \Omega_i$, then $a(x)\in \Star_{v_i}$,
and if $y\in\spc{N}$ lies in the simplex with vertices $v_{i_1},\dots, v_{i_m}$, then $b(y)\in \Omega_{i_1}\cup\dots\cup \Omega_{i_m}$.
\end{thm}

\section{Homotopy stability}

\begin{thm}{Theorem}\label{thm:h-stability}
Let $\spc{A}_1,\spc{A}_2,\dots$, and $\spc{A}_\infty$ be compact $m$-dimensional $\Alex\kappa$ spaces, and $m<\infty$.
Suppose $\spc{A}_n\z\to \spc{A}_\infty$ as $n\to \infty$ in the sense of Gromov--Hausdorff.
Then $\spc{A}_\infty$ is homotopy equivalent to $\spc{A}_n$ for all large $n$.

Moreover, given $\eps>0$ there are maps $h_n\:\spc{A}_\infty\to \spc{A}_n$ that are homotopy equivalences and $\eps$-liftings for all large $n$.
\end{thm}

Combining this theorem with Gromov's selection theorem (\ref{thm:gromov-compactness}) and Exercise~\ref{ex:pack-vol:dim}, we get the following.

\begin{thm}{Theorem}\label{thm:h-finiteness}
Given $\kappa\in \RR, D>0, v_0>0$ and a positive integer $m$,
there are only finitely many homotopy types of $m$-dimensional $\Alex\kappa$ spaces with diameter $\le D$, and volume $\ge v_0$.
\end{thm}

\parit{Proof of \ref{thm:h-finiteness} modulo \ref{thm:h-stability}.}
Assume the contrary.
Then we can choose a sequence of $m$-dimensional $\Alex\kappa$ spaces $\spc{A}_1,\spc{A}_2,\dots$ that have different homotopy types and satisfy the assumptions of the theorem.
By Gromov's selection theorem, we can assume that $\spc{A}_n$ converges to some space $\spc{A}_\infty$ in the sense of Gromov--Hausdorff.

By \ref{ex:pack-vol:dim}, $\dim \spc{A}_\infty=m$.
It remains to apply \ref{thm:h-stability}.
\qeds

\parit{Proof of \ref{thm:h-stability}.}
Since $\spc{A}_\infty$ is compact, applying \ref{cor:convex-nbhd}, we can find a finite open cover of $\spc{A}_\infty$ by convex open sets $\Omega_1,\dots, \Omega_k$ such that
for each $\Omega_i$ there is a strictly concave function $f_i$ that is defined in a neighborhood of the closure $\bar \Omega_i$ and such that $f_i|_{\partial \Omega_i}=0$.

By shifting $f_i$ down by small value $\eps>0$,
we can ensure that $\bigcap_{i\in S}\Omega_{i}\ne \emptyset$ if and only if $\bigcap_{i\in S}\bar\Omega_{i}\ne \emptyset$.

Suppose that $W=\bigcap_{i\in S}\Omega_{i}\ne \emptyset$.
Then $W$ is contractible.
In other words, the cover $\{\Omega_1,\dots, \Omega_k\}$ meets the assumptions of the nerve theorem (\ref{thm:nerve}).

Indeed, the function
\[f_S\df\min_{i\in S} f_i\]
is strictly concave and it vanishes on the boundary of $W$.
The $f_S$-gradient flow $(t,x)\mapsto \GF_{f_S}^t(x)$ defines a homotopy
$[0,\infty)\times W\z\to W$.
By the first distance estimate (\ref{thm:dist-est}), $\GF_{f_S}^t(x)$ converges to the (necessarily unique) maximum point of $f_S$.
Therefore, in the obtained homotopy we can reparametrize $[0,\infty)$ by $[0,1)$ and extend it continuously to $[0,1]$;
thus we get that $W$ is contractible.

The functions $f_i$ and sets $\Omega_i$ can be lifted to $\spc{A}_n$ while keeping their properties for all large $n$.
More precisely, there are liftings $f_{i,n}$ of all $f_i$ to $\spc{A}_n$ which are strictly concave for all large $n$ and such that $\bar\Omega_{i,n}=\set{x\in \spc{A}_n}{f_{i,n}(x)\ge 0}$ is a compact convex set and $\Omega_{i,n}\z=\set{x\in \spc{A}_n}{f_{i,n}(x)> 0}$ is an open convex set for each $i$.

Notice that $\{\Omega_{1,n},\dots,\Omega_{k,n}\}$ is an open cover of $\spc{A}_n$ for all large~$n$.
Indeed, suppose we have $p_n\in \spc{A}_n\setminus(\Omega_{1,n}\cup\dots\cup\Omega_{k,n})$ for arbitrarily large $n$.
Since $\spc{A}_\infty$ is compact, a subsequence of $p_n$ converges to some $p_\infty\in\spc{A}_\infty$.
But $p_\infty\in\Omega_i$ for some $i$ and therefore $p_n\in \Omega_{i,n}$ for arbitrarily large $n$ --- a contradiction.

In a similar fashion we can show that if $n$ is large then any collection $\{\Omega_{i,n}\}_{i\in S}$ has a common point in $\spc{A}_n$
if and only if $\{\Omega_{i}\}_{i\in S}$ has a common point in $\spc{A}_\infty$.
Here we have to use that $\bigcap_{i\in S}\Omega_{i}\ne \emptyset$ if and only if $\bigcap_{i\in S}\bar\Omega_{i}\ne \emptyset$.

It follows that for any large $n$ the covers
\begin{itemize}
\item $\{\Omega_{1},\dots,\Omega_{k}\}$ of $\spc{A}_\infty$ and
\item $\{\Omega_{1,n},\dots,\Omega_{k,n}\}$ of $\spc{A}_n$
\end{itemize}
have the same nerve.
By the nerve theorem (\ref{thm:nerve}), $\spc{A}_n$ and $\spc{A}_\infty$ are homotopy equivalent for all large $n$.
\qeds

Note that the proof above implies the following.

\begin{thm}{Theorem}\label{thm:finite-dim-hom-simplicial}
Any compact finite-dimensional Alexandrov space is homotopy equivalent to a finite simplicial complex.
\end{thm}

\section{Remarks}

Originally, Gromov's selection theorem was proved for Riemannian manifolds with a lower bound on Ricci curvature \cite{gromov1981}.
It motivates the study of limits of manifolds with lower Ricci curvature bounds and their synthetic generalizations, the so-called $\mathrm{CD}(K,m)$ spaces; $\mathrm{CD}$ stands for curvature-dimension condition.
This theory has significant applications in Alexandrov geometry;
in particular, it provides a version of the Liouville theorem about phase-space volume of geodesic flow in Alexandrov spaces \cite{brue-mondino-semola}.

The construction of a strictly concave function (\ref{thm:strictly-concave}) is due to Grigory Perelman \cite{perelman1993,perelman-petrunin}.
A sharper version of this construction is given in the survey of the third author \cite[7.2]{petrunin:survey}.
Let us mention the following closely related result of Artem Nepechiy \cite{nepechiy2019}:
\textit{given a point $p$ in an Alexandrov space, there is a $2$-convex function that approximates $\distfun_p^2$ up to second order at $p$}.

The fact that the controlled concavity survives under lifting (see Section \ref{sec:Liftings}) was observed and used by the first author \cite[4.2]{kapovitch2002}.
Functions with controlled concavity can be lifted to a collapsing sequence; in this case the resulting function is not concave but it has some interesting properties.
In particular, if the limit space has dimension $m$, then the Morse index of the lifting is at least $m$;
this property implies a partial answer to an old question of \textit{which finite-dimensional Alexandrov spaces can be approximated by Riemannian manifolds with a fixed lower curvature bound}; see \cite{kapovitch2005}.

Furthermore, such a lifted function has the following property (let us call it \index{gap-convexity}\emph{gap-convexity}) which, as far as we know, has not yet found an application:
\textit{If $f$ is a controllably concave function and $f_n$ are liftings of $f$ to a collapsing sequence of Alexandrov spaces with a uniform lower curvature bound, then for any $\eps>0$ and all large $n$, the restriction of $f_n$ to any geodesic of length at least $\eps$ is concave.}
(Recall that our geodesics are always length-minimizing.)

Let us list some results that can be proved by applying Gromov's selection theorem
in the same fashion as in the proof of the homotopy-type finiteness theorem (\ref{thm:h-finiteness}).

\begin{thm}{Betti-number theorem}\label{thm:betti}
For any $\kappa\in \RR, D>0$ and a positive integer $m$,
there is a constant $\Const=\Const(m,D,\kappa)$ such that
\[\beta_0(M)+\dots+\beta_m(M)\le \Const\]
for any closed $m$-dimensional Riemannian manifold $M$ with sectional curvature $\ge \kappa$ and diameter $\le D$.
Here $\beta_i(M)$ denotes the $i$-th Betti number of $M$.
\end{thm}

Gromov's original proof \cite{gromov-1981} of the Betti-number theorem did not use Alexandrov geometry directly.
However, it is easy to produce a proof following the steps above (no one bothered to write it so far).
The following result was proved by the third author \cite{petrunin2008}, and it uses the same technique.

\begin{thm}{Scalar curvature bound}
Given $\kappa\in \RR, D>0$ and a positive integer $m$, there is a constant $\Const=\Const(m,D,\kappa)$ such that
\[\int_M\Sc\le \Const\]
for any closed $m$-dimensional Riemannian manifold $M$ with sectional curvature $\ge \kappa$ and diameter $\le D$.
Here $\Sc$ denotes the scalar curvature.
\end{thm}

The following theorem is a stronger version of \ref{thm:h-stability}.
Its close relative (\ref{thm:spherical-nbhd}) will play an important role in the following lecture.

\begin{thm}{Stability theorem}\label{thm:stability}
Let $\spc{A}_1,\spc{A}_2,\dots$, and $\spc{A}_\infty$ be compact $m$-dimensional $\Alex\kappa$ spaces, and $m<\infty$.
Suppose $\spc{A}_n\to \spc{A}_\infty$ as $n\to \infty$ in the sense of Gromov--Hausdorff.
Then $\spc{A}_\infty$ is homeomorphic to $\spc{A}_n$ for all large $n$.

Moreover, given $\eps>0$ there are maps $h_n\:\spc{A}_\infty\to \spc{A}_n$ that are homeomorphisms and $\eps$-liftings for all large $n$.
\end{thm}

This theorem was proved by Grigory Perelman \cite{perelman1991}.
The proof was rewritten with more details by the first author \cite{kapovitch};
for more exact statements, see \cite{kapovitch2002,kapovitch2005}.
In private conversations, Perelman claimed that the homeomorphisms in the theorem can be assumed to be bi-Lipschitz with constants that depend on $\spc{A}_\infty$.
Since no proof has been written, this statement should be considered as a conjecture;
partial results in this direction were obtained by Mohammad Alattar \cite{alattar}.

Theorem \ref{thm:h-finiteness} was originally proved by Karsten Grove and Peter Petersen \cite{grove-petersen1988}.
Perelman's stability theorem (\ref{thm:stability}) implies the following stronger statement.

\begin{thm}{Homeomorphism-type finiteness}
For any real numbers $\kappa$, $D\z>0$, $v>0$, and a positive integer $m$, there are only finitely many homeomorphism types of closed $m$-dimensional manifolds that admit a Riemannian metric with sectional curvature $\ge \kappa$, volume $\ge v$, and diameter $\le D$.
\end{thm}

This statement can be improved to diffeomorphism-type finiteness in all dimensions $m\ne 4$.
Indeed, for $m\ne4$, a closed topological $m$-manifold admits only finitely many smooth structures; see \cite{kirby-siebenmann} and \cite{moise,thurston} for cases $m\ge 5$ and $m\le 3$, respectively.

In all these theorems we prove a statement about $m$-dimensional manifolds with lower sectional curvature bound following the following steps:
\begin{itemize}
\item Arguing by contradiction, we assume the existence of a sequence of manifolds that violates the statement.
\item Use Gromov's selection theorem to choose a converging subsequence.
\item Use curvature survival to conclude that the limit space is Alexandrov.
\item Try to arrive at a contradiction using Alexandrov geometry.
\end{itemize}

In principle, the same strategy might work for a sequence of Riemannian manifolds with dimension approaching infinity, but no applications of this kind have been found yet.
The following question gives an example where such a strategy might be successful in principle.

\begin{thm}{Question}
Is it true that no simply connected closed manifold of arbitrarily large dimension admits a Riemannian metric with two-sided sectional-curvature bound and diameter bounded by a fixed sufficiently small positive constant?
\end{thm}

If the dimension is bounded, then Gromov's theorem \cite{gromov1978} implies that such a manifold can be covered by a compact nil-manifold; in particular, the manifold cannot be simply connected.
However, if the dimension is unbounded, then Riemannian manifolds satisfying these conditions may not be covered by compact nil-manifolds; such examples were found by Galina Guzhvina \cite{guzhvina2008}.

%% file: boundary.tex
\chapter{Boundary}\label{chap:bry}

Here we define the boundary of a finite-dimensional Alexandrov space, discuss its properties, and prove the doubling theorem (\ref{thm:doubling:doubling}).

\section{Definition}

The definition of \index{boundary}\emph{boundary} will use induction on dimension, and it will work \textit{only} for finite-dimensional Alexandrov spaces.

If $\spc{A}$ is a $0$-dimensional Alexandrov space, we set
$\partial\spc{A}=\emptyset$.
Suppose $\spc{A}$ is a 1-dimensional Alexandrov space.
By~\ref{ex:dim=1},
$\spc{A}$ is homeomorphic to a 1-dimensional manifold (possibly with nonempty boundary).
This allows us to define the boundary $\partial\spc{A}\subset \spc{A}$ as the boundary of the manifold.

Now assume that the notion of boundary is defined in dimensions $1,\dots,m-1$, where $m\ge 2$.
Let $\spc{A}$ be an $m$-dimensional Alexandrov space.
We say that $p\in \spc{A}$ belongs to the boundary (briefly $p\in \partial \spc{A}$) if 
$\partial\Sigma_p\ne\emptyset$.
By \ref{thm:finite-space-of-directions} and \ref{ex:finite-space-of-directions-dim}, $\Sigma_p$ is an $(m-1)$-dimensional Alexandrov space;
therefore its boundary is already defined and hence this inductive definition makes sense.

\begin{thm}{Exercise}\label{ex:bry}
Prove the following statements:

\begin{subthm}{ex:bry:Em}
For a closed convex set $K\subset \EE^m$ with nonempty interior, the topological boundary of $K$ as a subset of $\EE^m$ coincides with the boundary of $K$ described above.
\end{subthm}

\begin{subthm}{ex:bry:sum}
Suppose a finite-dimensional Alexandrov space $\spc{A}$ is a direct sum of $\spc{A}_1$ and $\spc{A}_2$.
Then
\[\partial \spc{A}=(\partial\spc{A}_1\oplus\spc{A}_2)\,\cup\,(\spc{A}_1\oplus\partial\spc{A}_2)\]
\end{subthm}

\begin{subthm}{ex:bry:cone}
If $\Cone\Sigma$ is an $\Alex0$ space of dimension $\ge 2$ (by \ref{ex:finite-space-of-directions-dim} this implies that $\Sigma$ is $\Alex1$), then
\[\partial \Cone\Sigma=\Cone\partial\Sigma,\]
where $\Cone\partial\Sigma\df\set{s\cdot \xi\in\Cone\Sigma }{\xi\in \partial\Sigma}$.
\end{subthm}

\end{thm}

\section{Conic neighborhoods}

We are going to use the following result without proof;
it is a close relative of Perelman's stability theorem \ref{thm:stability}, and it is proved in \cite{perelman1993}.

Recall that the logarithm $\log_p\:\spc{A}\to \T_p$ is a multivalued function defined in \ref{Velocity of curve}.

\begin{thm}{Theorem}\label{thm:spherical-nbhd}
For any point $p$ in a finite-dimensional Alexandrov space $\spc{A}$ and all sufficiently small $\eps>0$,
there is a homeomorphism $h_\eps\:\oBall(p,\eps)_{\spc{A}}\to \oBall(0,\eps)_{\T_p}$ such that $0=h_\eps(p)$.

Moreover, we may assume that
\[
\sup_{x\in \oBall(p,\eps)}\{\,\tfrac1\eps\cdot\dist{\log_px}{h_\eps(x)}{\T_p}\,\}\to 0
\quad\text{as}\quad
\eps\to 0.\]
\end{thm}

Note that the last condition implies that $h_\eps$ is an $o(\eps)$-isometry.

The above theorem is often used together with the \textit{uniqueness of conic neighborhoods} stated below.

Suppose that an open neighborhood $U$ of a point $x$ in a metric space $\spc{X}$
admits a homeomorphism to $\Cone\Sigma$ such that $x$ is mapped to the origin of the cone.
In this case, we say that $U$ is a \index{conic neighborhood}\emph{conic neighborhood} of~$x$.

\begin{thm}{Uniqueness of conic neighborhoods}\label{lem:kwun}
Any two conic neighborhoods of a given point in a metric space are \index{pointed homeomorphic}\emph{pointed homeomorphic}; that is, there is a homeomorphism between neighborhoods that maps the origin of one cone to the origin of the other.
\end{thm}

\begin{thm}{Advanced exercise}\label{ex:conic}
Prove \ref{lem:kwun} or read the proof in \cite{kwun1964}.
\end{thm}

\begin{thm}[!]{Exercise}\label{ex:conic-tangent}
Suppose $x\mapsto x'$ is a homeomorphism between finite-dimensional Alexandrov spaces $\spc{A}$ and $\spc{A}'$. Show that 

\begin{subthm}{ex:conic-tangen:tangent}
$\T_x\cong \T_{x'}$ (here and below $\cong$ means homeomorphic),
\end{subthm}

\begin{subthm}{ex:conic-tangen:dir}
$\Susp\Sigma_x\cong \Susp\Sigma_{x'}$,
\end{subthm}

\begin{subthm}{ex:conic-tangen:example}
but in general $\Sigma_x\ncong\Sigma_{x'}$.
\end{subthm}

\end{thm}

\section{Topology}

The following theorem asserts that the boundary is a topological invariant, despite the heavy use of geometry in its definition.

\begin{thm}{Theorem}\label{thm:top-bry}
Let $\spc{A}$ and $\spc{A}'$ be homeomorphic finite-dimensional Alexandrov spaces.
Then $\dim \spc{A}=\dim\spc{A}'$ and
\[\partial\spc{A}\ne \emptyset
\quad\Leftrightarrow\quad
\partial\spc{A}'\ne \emptyset.
\]
\end{thm}

While working on the proof, keep in mind that there are pairs of spaces $\spc{K}_1$ and $\spc{K}_2$ such that $\spc{K}_1\ncong \spc{K}_2$, but $\RR\times \spc{K}_1\cong \RR\times \spc{K}_2$.
For example, one can let $\spc{K}_1=\SSS^4$ and let $\spc{K}_2$ be the suspension over the Poincaré homology sphere; compare with \ref{ex:conic-tangen:example}.

Let $\spc{A}$ be an $m$-dimensional Alexandrov space and $m<\infty$.
Define the \index{corank}\emph{corank} of $\spc{A}$ (briefly, $\corank\spc{A}$) as the minimal value $k$ such that $\spc{A}$ splits isometrically as $\RR^{m-k}\times \spc{K}$;
here $\spc{K}$ is a $k$-dimensional Alexandrov space.

In the following proof we will apply induction on the corank of $\spc{A}$.

\parit{Proof of \ref{thm:top-bry}.}
The first statement follows from \ref{thm:dim=dim}.

Suppose we have a counterexample, say $\partial \spc{A}\ne \emptyset$, but $\partial \spc{A}'=\emptyset$.
Let $k\df\corank \spc{A}$ and $k'\df\corank \spc{A}'$.
We can assume that the pair $(k,k')$ is minimal in lexicographic order among all counterexamples;
in particular, $k$ is minimal.
Let $x\mapsto x'$ be a homeomorphism from $\spc{A}$ to $\spc{A}'$.

Choose $x\in \partial \spc{A}$.
Since $\partial \spc{A}'=\emptyset$, we have $x'\notin \partial \spc{A}'$.
Note that 
\[\corank \T_x\le k
\quad\text{and}\quad
\corank \T_{x'}\le k'.
\]
By \ref{ex:conic-tangen:tangent}, $\T_x\cong\T_{x'}$.
Note that $\partial \T_x\ne\emptyset$ and $\partial \T_{x'}=\emptyset$.
Therefore, we may assume that $\spc{A}$ and $\spc{A}'$ are Euclidean cones
and the homeomorphism sends the origin to the origin.
The remaining part of the proof is divided into three cases.

\parit{Case 1.}
Suppose $k>1$.
Let $\spc{A}\iso \RR^{m-k}\times \spc{C}$, where $\spc{C}$ is a $k$-dimensional $\Alex0$ cone.
Observe that $\corank\T_y\le\corank\spc{A}$ for any $y\in\spc{A}$ and the equality holds only if $y$ projects to the origin of $\spc{C}$.

Since $k>1$, we can choose $z\ne 0$ that lies in $\partial\spc{C}$ (compare to \ref{ex:bry:cone}).
Choose $y\in \spc{A}$ that projects to $z$;
in particular, $\corank\T_y<\corank\spc{A}$.
By \ref{ex:conic-tangen:tangent}, $\T_y\cong\T_{y'}$,
$\partial \T_y\ne\emptyset$ and $\partial \T_{y'}=\emptyset$.
The latter contradicts the minimality of $k$.

\parit{Case 2.} Suppose $k\le1$ and $k'>1$.
Since $\partial \spc{A}\ne \emptyset$, we get that $k=1$;
therefore, $\spc{A}=\RR^{m-1}\times[0,\infty)$.

Let $\spc{A}'\iso \RR^{m-k'}\times \spc{C}'$, where $\spc{C}'$ is a $k'$-dimensional $\Alex0$ cone.
Since $\partial\spc{A}\cong\RR^{m-1}$,
the image of $\partial\spc{A}$ in $\spc{A}'$ does not lie in $\RR^{m-k'}\z\times\{0\}$.
In other words, we can choose $y\in \partial \spc{A}$ such that its image $y'\in \spc{A}'$ has a nonzero projection in $\spc{C}'$.
Observe that $\T_y\cong\T_{y'}$,
\[
\corank\T_y\le k=1,
\quad
\corank\T_{y'}< k',
\quad
\partial \T_y\ne\emptyset,
\quad\text{and}\quad
\partial \T_{y'}= \emptyset,\]
which contradicts minimality of $k'$.

\parit{Case 3.}
Suppose $k\le 1$ and $k'\le 1$.
Since $\partial \spc{A}\ne \emptyset$, we have $k=1$.
By \ref{ex:dim=1}, $\spc{A}\z\cong \RR^{m-1}\times[0,\infty)$.
Note that $\spc{A}'\cong\RR^m$, and $\spc{A}\ncong\spc{A}'$ --- a contradiction.
\qeds

\begin{thm}[!]{Exercise}\label{ex:bry2bry}
Let $x\mapsto x'$ be a homeomorphism $\Omega\to\Omega'$
between open subsets of finite-dimensional Alexandrov spaces $\spc{A}$ and $\spc{A}'$.
Show that $x\in \partial \spc{A}$ if and only if $x'\in \partial \spc{A}'$.
\end{thm}

\begin{thm}[!]{Exercise}\label{ex:bry-closed}
Show that the boundary of a finite-dimensional Alexandrov space is a closed subset.
\end{thm}

\section{Tangent space}

Spaces of directions and tangent spaces of an Alexandrov space have already been defined in \ref{sec:space+directions} and \ref{sec: tangent space}, respectively.
Let us extend these definitions to subsets of an Alexandrov space.

Let $X$ be a subset of a finite-dimensional Alexandrov space $\spc{A}$.
Choose $p\in \spc{A}$ and $\xi\in \Sigma_p$.
Suppose $\xi$ is a limit of directions $\dir{p}{x_n}$ for a sequence $x_1,x_2,\dots{}\in X$ that converges to $p$.
Then we say that $\xi$ is in the \index{space of directions}\emph{space of directions} from $p$ to $X$;
briefly \index{70@$\Sigma_p$ (space of directions)}$\xi\in\Sigma_pX$.

Further, $\Cone(\Sigma_pX)$ will be called the \index{tangent space}\emph{tangent space} to $X$ at $p$;
it will be denoted by \index{70@$\T_p$ (tangent space)}$\T_pX$.

Note that $\Sigma_pX\subset \Sigma_p\spc{A}$ and $\T_pX\subset\T_p\spc{A}$.

\begin{thm}{Theorem}\label{thm:partial-Sigma}
For any finite-dimensional Alexandrov space $\spc{A}$, we have
\[\partial (\Sigma_p\spc{A})=\Sigma_p(\partial\spc{A})
\quad\text{and}\quad
\partial(\T_p\spc{A})=\T_p(\partial\spc{A}).\]
\end{thm}

\parit{Proof.}
Choose $\xi\in\Sigma_p(\partial\spc{A})$ and a sequence $x_n\in \partial \spc{A}$ such that $x_n\to p$ and $\dir p{x_n}\to\xi$.

Let $\eps_n=2\cdot \dist{p}{x_n}{}$,
and let $h_{\eps_n}\:\oBall(p,\eps_n)_{\spc{A}}\to \oBall(0,\eps_n)_{\T_p}$ be the homeomorphisms provided by \ref{thm:spherical-nbhd};
in particular, $\tfrac2{\eps_n}\cdot h_{\eps_n}(x_n)\to \xi$ as $n\to\infty$.
By \ref{ex:bry2bry}, $h_{\eps_n}(x_n)\in \partial \T_p$.
By \ref{ex:bry-closed}, $\xi\in \partial \T_p$.
Therefore,
\[\partial (\Sigma_p\spc{A})\supset\Sigma_p(\partial\spc{A})
\quad\text{and}\quad
\partial(\T_p\spc{A})\supset\T_p(\partial\spc{A}).\]

Similarly, choose $\xi\in\partial\Sigma_p$.
Let $h_{\eps_n}\:\oBall(p,\eps_n)_{\spc{A}}\to \oBall(0,\eps_n)_{\T_p}$ be the homeomorphisms provided by \ref{thm:spherical-nbhd} for a sequence $\eps_n\to 0$ as $n\to\infty$.
By \ref{ex:bry2bry}, $x_n=h_{\eps_n}^{-1}(\tfrac{\eps_n}2\cdot\xi)\in \partial \spc{A}$.
By \ref{thm:spherical-nbhd}, $\dir p{x_n}\to \xi$.
Hence
\[\partial (\Sigma_p\spc{A})\subset\Sigma_p(\partial\spc{A})
\quad\text{and}\quad\partial(\T_p\spc{A})\subset\T_p(\partial\spc{A}).\]
\qedsf

\section{Doubling}

Let $A$ be a closed subset of a metric space $\spc{X}$.
The \index{doubling}\emph{doubling} $\spc{W}$ of $\spc{X}$ across $A$ is two copies of $\spc{X}$ glued along $A$.
More precisely, the underlying set of $\spc{W}$ is the quotient $\spc{X}\times\{0,1\}/\sim$, where $(a,0)\sim (a,1)$ for any $a\in A$ and $\spc{W}$ is equipped with the maximal metric such that both maps $\spc{X}\to \spc{W}$ defined by $x\mapsto (x,0)$ and $x\mapsto (x,1)$ are distance-preserving.

Alternatively, one may say that $\spc{W}$ is equipped with the maximal metric such that the projection $\proj\:\spc{W}\to\spc{X}$ defined by $(x,i)\mapsto x$ is a short map.
The metric on $\spc{W}$ can also be defined directly:
\[\dist{(x,i)}{(y,j)}{\spc{W}}=
\begin{cases}
\dist{x}{y}{\spc{X}}&\text{if}\quad i= j.
\\
\inf\set{\dist{x}{a}{\spc{X}}+\dist{y}{a}{\spc{X}}}{a\in A}&\text{if}\quad i\ne j.
\end{cases}
\]

\begin{thm}{Theorem}\label{thm:doubling}\index{doubling theorem}
Let $\spc{A}$ be a finite-dimensional Alexandrov space with nonempty boundary,
and let $f\z=\tfrac12\cdot\distfun_p^2$ for some $p\in \spc{A}$.
Then

\begin{subthm}{thm:doubling:concave}
If $\dim \spc{A}\ge 2$, then
$\distfun_{\partial \Sigma_x}\xi\le \tfrac\pi2$ for any $x\in\partial \spc{A}$ and $\xi\in \Sigma_x$.
Moreover, if $\distfun_{\partial \Sigma_x}\xi= \tfrac\pi2$, then $\mangle(\xi,\zeta)=\tfrac\pi2$ for any $\zeta\in \partial\Sigma_x$.
\end{subthm}

\begin{subthm}{thm:partial-grad:grad}
$\nabla_xf\in \partial\T_x$ for any $x\in\partial \spc{A}$.
\end{subthm}

\begin{subthm}{thm:partial-grad:flow}
If $\alpha$ is an $f$-gradient curve that starts at $x\in \partial \spc{A}$, then $\alpha(t)\in \partial \spc{A}$ for any $t$.
Moreover, if $p\in \partial \spc{A}$, then $\gexp_p(v)\in \partial \spc{A}$ for any $v\in\partial\T_p$.
\end{subthm}

\begin{subthm}{thm:doubling:doubling}
The doubling $\spc{W}$ of $\spc{A}$ across $\partial \spc{A}$ is an Alexandrov space with the same curvature bound as $\spc{A}$.
\end{subthm}

\end{thm}

Part \ref{SHORT.thm:doubling:doubling} is called the {}\emph{doubling theorem}.

\parit{Proof.}
We will denote by 
\ref{SHORT.thm:doubling:concave}$_m,\dots,$\ref{SHORT.thm:doubling:doubling}$_m$ the corresponding statement assuming $m=\dim\spc{A}$.

The proof goes by induction on $m$.
Statement \ref{SHORT.thm:doubling:doubling}$_1$ follows from \ref{ex:dim=1} --- this is the base.
The induction step is a combination of the implications below.

\parit{\ref{SHORT.thm:doubling:doubling}$_{m-1}\Rightarrow$\ref{SHORT.thm:doubling:concave}$_m$.}
If $m=2$, then $\dim\Sigma_x=1$; see \ref{ex:finite-space-of-directions-dim}.
By \ref{ex:dim=1}, $\Sigma_x$ is isometric to a line segment $[0,\ell]$;
we need to show that $\ell\le\pi$.

Recall that the tangent space $\T_x=\Cone\Sigma_x$ is an Alexandrov space; see \ref{ex:finite-tan}.
If $\ell>\pi$, then $\T_x$ has bifurcating geodesics, which is impossible by \ref{ex:0-angle}.

Now suppose $m>2$, so $\dim \Sigma_x>1$.
Assume $\distfun_{\partial \Sigma_x}\xi>\tfrac\pi2$ for some $\xi$.
By \ref{SHORT.thm:doubling:doubling}$_{m-1}$, the doubling $\Xi$ of $\Sigma_x$ is $\Alex1$.
Since $\dim\Xi>1$, we have $\diam\Xi\le\pi$ by \ref{ex:RisCBB(1)}.
Therefore, $\Cone\Xi$ is $\Alex0$ by \ref{ex:cone-CBB}.
Denote by $\xi_0$ and $\xi_1$ the points in $\Xi$ that correspond to $\xi$.
Observe that $\dist{\xi_0}{\xi_1}{\Xi}>\pi$.
The latter contradicts \ref{ex:RisCBB(1)}.

Finally, if $\distfun_{\partial \Sigma_x}\xi= \tfrac\pi2$, then $\dist{\xi_0}{\xi_1}{\Xi}=\pi$.
Therefore, $\Cone \Xi$ contains a line in the directions of $\xi_0$ and $\xi_1$.
By the splitting theorem (\ref{thm:splitting}), $\Cone \Xi$ is a direct product of the line with some subcone;
in other words, $\Xi$ is a spherical suspension with poles $\xi_0$ and $\xi_1$.
In particular, every point of $\Xi$ lies at distance at most $\tfrac\pi2$ from $\xi_0$ or $\xi_1$.
The natural projection $\Xi\to \Sigma_x$ does not increase distances and sends both $\xi_0$ and $\xi_1$ to $\xi$.
Therefore, the second statement of \ref{SHORT.thm:doubling:concave}$_m$ follows.

\parit{\ref{SHORT.thm:doubling:concave}$_{m-1}+$\ref{SHORT.thm:doubling:concave}$_m\Rightarrow$\ref{SHORT.thm:partial-grad:grad}$_m$.}
We can assume that $s=|\nabla_xf|\ne 0$.
By \ref{prop:grad-exist}, $\nabla_xf\z=s\cdot \overline{\xi}$, where $s=\dd_xf(\overline{\xi})>0$ and $\overline{\xi}\in\Sigma_x$ is the direction that maximizes $\dd_xf(\overline{\xi})$.

Let $\zeta\in \partial\Sigma_x$ be a direction that minimizes the angle $\mangle(\overline{\xi},\zeta)$.
It is sufficient to show that $\zeta=\overline{\xi}$.

Assume $\zeta\ne \overline{\xi}$;
let $\eta=\dir[\Sigma_x]\zeta{\overline{\xi}}$.
By \ref{SHORT.thm:doubling:concave}$_m$, $\mangle(\overline{\xi},\zeta)\le \tfrac\pi2$; further,
\ref{SHORT.thm:doubling:concave}$_{m-1}$ implies that 
\[\mangle(\eta,\nu)\le \tfrac\pi2\eqlbl{eq:<pi/2}\]
for any $\nu\in \Sigma_\zeta\Sigma_x$ (if $m=2$, then the last statement is evident). 

Let $\phi\:\Sigma_x\to\RR$ be the restriction of $\dd_xf$ to $\Sigma_x$.
Applying \ref{ex:d(distfun):=} and \ref{eq:<pi/2}, we get that $(\dd_{\zeta}\phi)(\eta)\le 0$.
Recall that $\dd_xf$ is concave (\ref{ex:concave-differential}), so $(\dd_xf)''\le 0$;
by straightforward computations, we have that $\phi''+\phi\le 0$.
If $\phi(\zeta)\le 0$, this implies that $\phi(\overline{\xi})\le 0$ --- a contradiction to the fact that $s>0$.
If $\phi(\zeta)> 0$, then $\phi(\overline{\xi})<\phi(\zeta)$ --- again a contradiction.

\parit{\ref{SHORT.thm:partial-grad:grad}$_m\Rightarrow$\ref{SHORT.thm:partial-grad:flow}$_m$.}
Let $\alpha$ be an $f$-gradient curve and $\ell(t)=\distfun_{\partial \spc{A}}\circ\alpha(t)$.

Choose $t$;
let $x=\alpha(t)$ and $y\in \partial\spc{A}$ be a closest point to $x$.
By \ref{SHORT.thm:partial-grad:grad}$_m$, we have that $\nabla_y f\in\partial \T_y$.
Since the distance $\dist{x}{y}{}$ is minimal, 
we get $\langle \dir yx,v\rangle\le 0$ for any $v\in \partial \T_y$.
In particular,
\[\langle \dir yx,\nabla_y f\rangle\le 0.\]
Applying \ref{ex:monotonicity} to $x$ and $y$,
we get
\[\ell'(t)\le \ell(t)\]
if the left-hand side is defined.
Since $\ell$ is Lipschitz, $\ell'$ is defined almost everywhere.
Integrating the inequality, we get 
\[\ell(t)\le e^t\cdot\ell(0)\]
for any $t\ge 0$.
In particular, if $\ell(0)=0$, then $\ell(t)=0$ for any $t\ge 0$.
Since $\partial\spc{A}$ is closed (\ref{ex:bry-closed}), the first statement follows.

It remains to show that $\gexp_p(v)\in \partial \spc{A}$ for any $v\in\partial\T_p$.
We can assume $v\ne0$; let $v=t\cdot u$ and $|u|=1$.
Note that $u\in\partial\T_p$.
Therefore, we can choose a sequence of points $x_n\in \partial \spc{A}$ such that $x_n\to p$ and $\dir{p}{x_n}\to u$; let $r_n=\dist{p}{x_n}{}$.
From above, the $f$-gradient curve $\alpha_n$ that starts at $x_n$ remains in $\partial\spc{A}$.
By the definition of gradient exponent, $\gexp_p(v)=\gexp_p(t\cdot u)$ is a limit of $\alpha_n(\ln\tfrac{t}{r_n})$.
Hence the last statement follows.

\parit{\ref{SHORT.thm:partial-grad:flow}$_{m}+$\ref{SHORT.thm:doubling:doubling}$_{m-1}\Rightarrow$\ref{SHORT.thm:doubling:doubling}$_m$.}
We will consider the case $\kappa=0$;
other cases can be done in the same way, but formulas get more complicated.

Denote by $\spc{A}_0$ and $\spc{A}_1$ the two copies of $\spc{A}$ in $\spc{W}$;
let us keep the notation $\partial \spc{A}$ for the common boundary of $\spc{A}_0$ and $\spc{A}_1$.

\begin{clm}{}\label{clm:crossing}
A geodesic in $\spc{W}$ either has at most one interior point in $\partial\spc{A}$ or lies entirely in $\partial\spc{A}$.
\end{clm}

\begin{wrapfigure}{r}{45mm}
\vskip-2mm
\centering
\includegraphics{mppics/pic-1315}
\end{wrapfigure}

Indeed, assume $\gamma$ shares at least two points with $\partial \spc{A}$, say $x=\gamma(t_1)$ and $y=\gamma(t_2)$ and these are not endpoints of $\gamma$.
Remove from $\gamma$ the set $\gamma\cap \spc{A}_1$ and replace it by its reflection across $\partial\spc{A}$;
denote the obtained curve by $\hat\gamma$.

Any arc of $\hat\gamma$ with one endpoint in $\partial \spc{A}$
is a geodesic in $\spc{A}_0$.
Since $x,y\in \partial \spc{A}$, the arc of $\hat\gamma$ beyond $y$ lies in the image of the map $t\z\mapsto \GF^t_{f_x}(y)$, where $f_x=\tfrac12\cdot\distfun^2_x$.
By \ref{SHORT.thm:partial-grad:flow}, this arc lies in $\partial\spc{A}$.

Now choose a point $z$ on this arc, so $z\in \partial\spc{A}$.
Applying the same argument, we get that the arc of $\hat\gamma$ up to $y$ lies in $\partial\spc{A}$.
Hence the claim follows.\claimqeds

Choose a point $p$ in $\spc{W}$;
let $f\df\tfrac12\cdot\distfun_p^2$.
It is sufficient to show that $(f\circ\gamma)''\le 1$.
If $p\in \partial\spc{A}$, then the statement follows from function comparison in $\spc{A}_0$ and $\spc{A}_1$.
So, we can assume that $p\in \spc{A}_0\setminus \partial\spc{A}$.
Also, we can assume that $\gamma$ does not lie in $\partial \spc{A}$;
otherwise, the inequality again follows from the comparison in~$\spc{A}_0$.

\begin{wrapfigure}{r}{55mm}
\vskip-2mm
\centering
\includegraphics{mppics/pic-1325}
\end{wrapfigure}

Choose $y=\gamma(t_0)$; without loss of generality we can assume that $t_0=0$.

Assume $y\z\in \spc{A}_0\setminus\partial\spc{A}$.
Then $(f\z\circ\gamma)''(0)\le 1$ in the barrier sense;
it follows from the comparison in $\spc{A}_0$.

Assume $y\in \spc{A}_1\setminus\partial\spc{A}$.
Suppose $[py]$ crosses $\partial\spc{A}$ at $x$;
by \ref{clm:crossing}, this point is unique.
Let $\Sigma_x$ be the space of directions of $\spc{A}$ at $x$,
and let $\Xi$ be its doubling.
As before, we denote by $\Sigma_0$ and $\Sigma_1$ the two copies of $\Sigma_x$ in $\Xi$
and keep the notation $\partial\Sigma_x$ for their common boundary.
By \ref{SHORT.thm:doubling:doubling}$_{m-1}$, $\Xi$ is $\Alex1$; note that $\diam\Xi\le\pi$.

The directions $\dir x{y}$ and $\dir xp$ lie on opposite sides of $\partial\Sigma_x$ and
\[\dist{\dir xy}{\dir xp}\Xi\ge \pi.\]
(Thus equality holds.)
Otherwise, we could choose a direction $\xi\in\partial\Sigma_x$ such that
\[\dist{\dir x{y}}{\xi}\Xi+\dist{\xi}{\dir xp}\Xi<\pi.\]
Furthermore, we could consider the curve $\alpha(t)=\gexp_x(t\cdot \xi)$.
By \ref{SHORT.thm:partial-grad:flow}$_m$, $\alpha$ lies in $\partial \spc{A}$.
By \ref{prop:gexp}
\[\dist{p}{\alpha(s)}{\spc{A}_0}
+\dist{y}{\alpha(s)}{\spc{A}_1}
<\dist{p}{y}{\spc{W}}\]
for small values $s>0$
--- a contradiction.

$\Cone \Xi$ contains a line with directions $\dir x{y}$ and $\dir xp$.
By the splitting theorem, $\Cone \Xi$ splits in these directions;
in particular, 
\[\dist{\dir x{y}}{\xi}{}+\dist{\xi}{\dir xp}{}=\pi\]
for any $\xi\in\Xi$.
It follows that for any $\xi\in\Xi$ there is $\xi'\in\partial\Sigma_x$ such that 
$\xi$ and $\xi'$ lie on some geodesic $[\dir x{y} \dir xp]_\Xi$.

Fix $t\approx 0$ such that $t\ne 0$.
Let $z=\gamma(t)$ and $\xi=\dir xz$;
choose a direction $\xi'\in\partial\Sigma_x$ as above for $\xi$.
Consider the curve $\alpha(s)\z\df\gexp_x(s\cdot\xi')$.
Let us show that 
\[
\begin{aligned}
\dist{p}{z}{\spc{W}}
&\le \dist{p}{\alpha(s)}{\spc{A}_0}+ \dist{\alpha(s)}{z}{\spc{A}_1}\le
\\
&\le\side\hinge yp{z}.
\end{aligned}
\eqlbl{eq:gamma''}
\]
for a suitable value $s$.

The first inequality in \ref{eq:gamma''} is evident.
Set $\phi=\mangle\hinge{x}{y}{z}$ and $\psi\z=\mangle(\dir xp,\xi')$.
The choice of $s$ comes from the model configuration $\tilde p$, $\tilde x$, $\tilde y$, $\tilde w$, $\tilde z\in \EE^2$ such that
\begin{align*}
\tilde x&\in [\tilde p\tilde y],
&
\dist{\tilde p}{\tilde x}{}&=\dist{ p}{x}{},
&
\dist{\tilde p}{\tilde y}{}&=\dist{p}{y}{},
&
\dist{\tilde x}{\tilde z}{}&=\dist{x}{z}{},
\\
\tilde w&\in [\tilde p\tilde z],
&
\mangle\hinge{\tilde x}{\tilde y}{\tilde z}&=\phi,
&
\mangle\hinge{\tilde x}{\tilde p}{\tilde w}&=\psi, 
&
s&=\dist{\tilde x}{\tilde w}{}.
\end{align*}
\begin{figure}[ht!]
\vskip-0mm
\centering
\includegraphics{mppics/pic-1014}
\end{figure}

\noindent
By \ref{prop:gexp}, we get 
\begin{align*}
\dist{p}{\alpha(s)}{\spc{A}_0}&\le \dist{\tilde p}{\tilde w}{},
\\
\dist{\alpha(s)}{z}{\spc{A}_1}&\le\dist{\tilde w}{\tilde z}{};
\end{align*}
by the comparison, 
\[\dist{\tilde p}{\tilde z}{}\le \side\hinge ypz.\]

\begin{thm}{Exercise}\label{ex:pz<ypz}
Prove the last inequality.
\end{thm}

Hence we get $(f\circ\gamma)''(0)\le 1$ in the barrier sense.

Finally, if $\gamma(0)\in\partial\spc{A}$, then a splitting argument shows that
\[(f\circ\gamma)^+(0)-(f\circ\gamma)^-(0)\le 0.\]

Summarizing, we get that $(f\circ\gamma)''\le 1$ on every arc of $\gamma$ that lies entirely in $\spc{A}_0$ or $\spc{A}_1$.
If $\gamma$ crosses $\partial \spc{A}$, then we know that it happens only once and at the crossing moment $t_0$ 
we have $(f\circ\gamma)^+(t_0)\z-(f\circ\gamma)^-(t_0)\z\le 0$.
All this implies that $(f\circ\gamma)''\le 1$.
\qeds

\begin{thm}[!]{Exercise}\label{ex:bry-connected}
Let $\spc{A}$ be an $m$-dimensional $\Alex1$ space.
Assume that $\partial\spc{A}\ne \emptyset$ and $2\le m<\infty$.
Show that $\partial\spc{A}$ is connected.
\end{thm}

\begin{thm}{Exercise}\label{ex:dist-to-bry}
Let $\spc{A}$ be a finite-dimensional $\Alex0$ space with nonempty boundary $\partial\spc{A}$.

\begin{subthm}{ex:dist-to-bry:geod}
Suppose a geodesic $\gamma$ in $\spc{A}$ has an interior point in $\partial\spc{A}$.
Show that $\gamma$ lies in $\partial\spc{A}$.
\end{subthm}

\begin{subthm}{ex:dist-to-bry:dist}
Show that the distance function to the boundary is concave.
\end{subthm}

\end{thm}

\begin{thm}[!]{Exercise}\label{ex:liberman}
Let $\spc{A}$ be a finite-dimensional $\Alex0$ space with nonempty boundary $\partial\spc{A}$.
Suppose $\gamma$ is a geodesic in $\partial\spc{A}$ with the induced length metric.
Show that the function $t\mapsto \tfrac12\cdot\distfun_p^2\circ\gamma(t)$ is 1-concave for any point $p$. 
\end{thm}

\begin{thm}{Exercise}\label{ex:native}
Let $\proj\:\spc{W}\to\spc{A}$ be the natural projection to a finite-dimensional Alexandrov space $\spc{A}$ from its doubling $\spc{W}$ across the boundary.
Suppose $f\:\spc{A}\to\RR$ is a $\lambda$-concave function.
Show that $f\circ\proj\:\spc{W}\to\RR$ is $\lambda$-concave if and only if $\nabla_xf\in \partial \T_x$ 
for any $x\in\partial \spc{A}$.
\end{thm}

\section{Remarks}\label{sec:bry-remarks}

The doubling theorem was generalized by the third author \cite{petrunin1997} to allow the gluing of nonidentical spaces, and was generalized further by Jian Ge and Nan Li \cite{ge-li} to allow the gluing of more than two spaces.

It easily follows by induction on dimension that the doubling of a finite-dimensional Alexandrov space across its boundary results in an Alexandrov space without boundary.
This observation often helps to reduce a statement about general finite-dimensional Alexandrov spaces to Alexandrov spaces without boundary.

Let us discuss several additional tools that are available for spaces with empty boundary.

\begin{thm}{Fundamental-class lemma}\label{lem:fund-class}
Any compact finite-dimensional Alexandrov space $\spc{A}$ without boundary has a fundamental class with $\ZZ/2$ coefficients;
that is, if $\spc{A}$ is $m$-dimensional, then
\[H^m(\spc{A},\ZZ/2)=\ZZ/2.\]

\end{thm}

This lemma was proved by Karsten Grove and Peter Petersen \cite{grove-petersen1993}.
Originally it was stated for Alexander--Spanier cohomology.
We do not make this distinction because for compact Alexandrov spaces it is the same as singular cohomology.
Indeed, both cohomology theories are homotopy invariant \cite[Chapter 6]{Spanier}, compact Alexandrov spaces are homotopy equivalent to finite simplicial complexes (see \ref{thm:finite-dim-hom-simplicial}), and for paracompact CW complexes, Alexander--Spanier cohomology is isomorphic to Čech and singular cohomology \cite[Chapter 6]{Spanier}.

The lemma implies, for example, that on finite-dimensional Alexandrov spaces without boundary
the gradient flow for a $\lambda$-concave function is an onto map;
in other words, gradient curves can be extended into the past.
It is also used in the proof of the following version of the invariance-of-domain theorem \cite[Theorem 3.2]{kapovitch-zhu}.

\begin{thm}{Invariance of domain}\label{thm-inv-domain}
Let $\spc{A}_1$ and $\spc{A}_2$ be two $m$-dimensional Alexandrov spaces with empty boundary; $m$ is finite.
Suppose $\Omega_1$ is an open subset of $\spc{A}_1$ and $f\:\Omega_1\to \spc{A}_2$ is an injective continuous map.
Then $f(\Omega_1)$ is open in $\spc{A}_2$.
\end{thm}

Theorem~\ref{thm:spherical-nbhd} can be used to prove the following. 

\begin{thm}{Topological stratification}\label{thm:top-stratification}
Any $m$-dimensional Alexandrov space with $2\le m<\infty$ can be subdivided into topological manifolds $S_0,\z\dots,S_m$ such that for every $i$ we have $\dim S_i=i$ or $S_i=\emptyset$.
Moreover,
\begin{subthm}{}
the closure of $S_{m-1}$ is the boundary of the space, and
\end{subthm}

\begin{subthm}{}
$S_{m-2}=\emptyset$.
\end{subthm}

\end{thm}

This statement implies that every compact finite-dimensional Alexandrov space has the homotopy type of a finite CW complex.
It remains unknown whether such a space must be homeomorphic to a CW complex.
However, Tadashi Fujioka proved that every finite-dimensional Alexandrov space admits a CS stratification \cite{fujioka2025a}.

The stratification theorem~\ref{thm:top-stratification} can be sharpened as follows.

\begin{thm}{Boundary characterization}
Let $\spc{A}$ be an $m$-dimensional Alexandrov space with $m<\infty$.
Given $p\in\spc{A}$, the following statements are equivalent.

\begin{subthm}{item-boundary}
$p\in \partial \spc{A}$;
\end{subthm}

\begin{subthm}{item-contractible}
$\Sigma_p$ is contractible;
\end{subthm}

\begin{subthm}{item-space-dir-homology}
$\tilde H_{m-1}(\Sigma_p,\ZZ/2)= 0$;
\end{subthm}

\begin{subthm}{item-local-homology}
$H_m(\spc{A},\spc{A}\setminus \{p\},\ZZ/2)= 0$.
\end{subthm}

\end{thm}

Let $f$ be a semiconcave function.
A point $p\in \Dom f$ is called a \index{critical point}\emph{critical} point of $f$ if $\dd_pf\le 0$;
otherwise it is called \index{regular point}\emph{regular}.

The following statement plays a role in the proof of the stability theorem,
but it is also a useful technical tool on its own.

\begin{thm}{Morse lemma}\index{Morse lemma}
Let $f$ be a semiconcave function on a finite-dimensional Alexandrov space without boundary.
Suppose $K$ is a compact set of regular points of $f$ in its level set $f=a$.
Then an open neighborhood $\Omega$ of $K$ admits a homeomorphism $x\mapsto (h(x),f(x))$ to a product space $\Lambda\times (a-\eps,a+\eps)$.
\end{thm}

The Morse lemma was originally proved by Grigory Perelman \cite{perelman1993} for a special class of semiconcave functions built from distance functions.
For such functions it holds for all Alexandrov spaces (possibly with nonempty boundary).
For general semiconcave functions on spaces without boundary, it follows from \cite{perelman1994} where the technical tools from \cite{perelman1993} (the notion of inner product between differentials of semiconcave functions and the notion of admissible maps) were extended to general semiconcave functions.
With those definitions all the main Morse theory results of \cite{perelman1993} follow by the same proof, using \cite[2.2]{perelman1994} instead of \cite[Lemma 1]{perelman1993}.
On Alexandrov spaces with boundary, all the Morse theory results hold if one requires in addition that the semiconcave functions remain semiconcave when canonically extended to the doubling.%
\footnote{Semiconcave functions are often defined to meet this property, but we do not follow this convention.}
This is automatic for the special class of semiconcave functions considered in \cite{perelman1993}; therefore this distinction is not made there.

Subsets of Alexandrov spaces that satisfy the condition in \ref{thm:partial-grad:flow} are called extremal.
More precisely, a subset $E$ is \index{extremal set}\emph{extremal} if any $f$-gradient curve that starts in $E$ remains in $E$;
here $f$ is an arbitrary function of the form $\tfrac12\cdot \distfun_p^2$.

Extremal subsets were introduced by Grigory Perelman and the third author \cite{perelman-petrunin}.
They will pop up in the next lecture.

The following conjecture is one of the oldest questions in Alexandrov geometry that remains open.

\begin{thm}{Conjecture}\label{conj:bry}
Let $S$ be a component of the boundary of a finite-dimensional Alexandrov space.
Then $S$, equipped with the induced length metric, is an Alexandrov space with the same curvature bound.
\end{thm}

%% file: CBB-quotients.tex
\chapter{Quotients}\label{chap:L/G}

This lecture discusses applications of Alexandrov geometry to isometric group actions.

\section{Quotient space}

Suppose that a group $G$ acts isometrically on a metric space $\spc{X}$.
Note that
\[\dist{G\cdot x}{G\cdot y}{\spc{X}/G}
\df
\inf
\set{\dist{x}{g\cdot y}{\spc{X}}}{g\in G}\]
defines a semimetric on the orbit space $\spc{X}/G$.
Moreover, if the orbits of the action are closed,
then it is a genuine metric.

\begin{thm}{Theorem}\label{thm:CBB/G}
Suppose that a group $G$ acts isometrically on a proper $\Alex0$ space $\spc{A}$, and $G$ has closed orbits.
Then the quotient space $\spc{A}/G$ is $\Alex0$.

\end{thm}

A more general statement will be given in \ref{thm:submetry-CBB-1}.

\parit{Proof.}
Denote by $\sigma\:\spc{A}\to \spc{A}/G$ the quotient map.

Fix a quadruple of points $p,x_1,x_2,x_3\in \spc{A}/G$.
Choose $\hat p\in \spc{A}$ such that $\sigma(\hat{p})=p$.
Since $\spc{A}$ is proper and the $G$-orbits are closed, we can choose points $\hat{x}_i\in \spc{A}$ such that $\sigma(\hat x_i)=x_i$ and
\[\dist{p}{x_i}{\spc{A}/G}
=
\dist{\hat{p}}{\hat{x}_i}{\spc{A}}\]
for all $i$.

Note that 
\[\dist{x_i}{x_j}{\spc{A}/G}
\le 
\dist{\hat{x}_i}{\hat{x}_j}{\spc{A}}
\]
for all $i$ and $j$.
By \ref{angle-monotonicity},
\[\angk p{x_i}{x_j}
\le
\angk {\hat{p}}{\hat{x}_i}{\hat{x}_j}
\eqlbl{eq:angles-M-L}\]
for all $i$ and $j$.

By the $\EE^2$-comparison in $\spc{A}$,
we have
\[\angk {\hat{p}}{\hat{x}_1}{\hat{x}_2}
+\angk {\hat{p}}{\hat{x}_2}{\hat{x}_3}
+\angk {\hat{p}}{\hat{x}_3}{\hat{x}_1}
\le 
2\cdot\pi.\]
Applying \ref{eq:angles-M-L},
we get 
\[\angk p{x_1}{x_2}
+\angk p{x_2}{x_3}
+\angk p{x_3}{x_1}\le 2\cdot\pi;\]
that is,
the $\EE^2$-comparison holds for any quadruple in $\spc{A}/G$.

Since $\spc{A}$ is proper and geodesic and the orbits are closed, we get that $\spc{A}/G$ is geodesic, and the statement follows.
\qeds

\begin{thm}{Very advanced exercise}\label{ex:Hilbert/G}
Let $G$ be a compact Lie group with a bi-invariant Riemannian metric.
Show that $G$ is isometric to a quotient of a Hilbert space by an isometric group action.

Conclude that $G$ is $\Alex0$.
\end{thm}

\section{Submetries}

A map $\sigma\:\spc{X}\to\spc{Y}$ between metric spaces
is called a \index{submetry}\emph{submetry} if 
\[\sigma(\oBall(p,r)_\spc{X})=\oBall(\sigma(p),r)_{\spc{Y}}\]
for any $p\in \spc{X}$ and $r\ge 0$.

Suppose $G$ and $\spc{A}$ are as in \ref{thm:CBB/G}.
Observe that the quotient map $\sigma\:\spc{A}\to \spc{A}/G$ is a submetry.
The following two exercises show that this is not the only source of submetries. 

\begin{thm}{Exercise}\label{ex:sumbetries(S^2)}
Construct submetries
\begin{subthm}{ex:sumbetries(S^2):1}
$\sigma_1\:\mathbb{S}^2\to[0,\pi]$,
\end{subthm}
\begin{subthm}{ex:sumbetries(S^2):2}
$\sigma_2\:\mathbb{S}^2\to[0,\tfrac\pi2]$,
\end{subthm}
\begin{subthm}{ex:sumbetries(S^2):n}
$\sigma_n\:\mathbb{S}^2\to[0,\tfrac\pi n]$ (for an integer $n\ge 1$)
\end{subthm}
such that the fibers $\sigma_n^{-1}\{x\}$ are connected for any $x$.
\end{thm}

\begin{thm}{Exercise}\label{ex:sumbetries(E^2)}
Let $\sigma\:\EE^2\to [0,\infty)$ be a submetry.
Show that $K\z=\sigma^{-1}\{0\}$ is a closed convex set with empty interior and $\sigma(x)\z=\distfun_Kx$.
\end{thm}

The proof of \ref{thm:CBB/G} works for submetries;
that is, \textit{if $\sigma\:\spc{A}\to\spc{B}$ is a submetry and $\spc{A}$ is a proper $\Alex0$ space, then so is $\spc{B}$}.
Theorem \ref{thm:CBB/G} admits a straightforward generalization to the $\Alex{-1}$ case.

In the $\Alex1$ case, the proof produces a weaker statement: \textit{if $\sigma\:\spc{A}\to\spc{B}$ is a submetry and $\spc{A}$ is $\Alex1$ then $\SSS^2$-comparison holds for any quadruple $p$, $x_1$, $x_2$, $x_3$ in $\spc{B}$ if $\dist{p}{x_i}{}<\tfrac\pi 2$ for each $i$}.
In particular, $\spc{B}$ is \textit{locally} $\Alex1$.
As above, assuming orbits are closed, the quotient of a proper geodesic space is geodesic.
Therefore, the globalization theorem implies that it is globally $\Alex1$.
The same holds for the targets of submetries from $\Alex1$ spaces.
With a bit of extra work, one can extend the statement to nonproper spaces \cite[8.34]{alexander-kapovitch-petrunin2024}.
Thus, we have the following.

\begin{thm}{Theorem}\label{thm:submetry-CBB-1}
Let $\sigma\:\spc{A}\to\spc{B}$ be a submetry.
If $\spc{A}$ is $\Alex\kappa$, then so is $\spc{B}$.

In particular, if $G$ acts isometrically on an $\Alex\kappa$ space $\spc{A}$, and $G$ has closed orbits, then the quotient space $\spc{A}/G$ is $\Alex\kappa$.

\end{thm}

\section{Hopf's conjecture}

\textit{Does $\mathbb{S}^2\times\mathbb{S}^2$ admit a Riemannian metric with positive sectional curvature?} \index{Hopf's conjecture}\emph{Hopf's conjecture} says that the answer should be negative.
Let us take a close look at the following partial result obtained by Wu-Yi Hsiang and Bruce Kleiner \cite{hsiang-kleiner}.

\begin{thm}{Theorem}\label{thm:hsiang-kleiner}
There is no Riemannian metric on $\SSS^2\times\SSS^2$ with sectional curvature $\ge 1$ and a nontrivial isometric $\SSS^1$-action.
\end{thm}

Recall that a group action $G\acts\spc{X}$ is called \index{effective action}\emph{effective} if for any nonidentity element $g\in G$ there is $x\in\spc{X}$ such that $g\cdot x\ne x$.

\begin{thm}{Key lemma}\label{lem:S^3/S^1}
Suppose $\SSS^1\acts\SSS^3$ is an effective isometric action without fixed points
and $\Sigma=\SSS^3/\SSS^1$ is its quotient space.
Then there is a distance-noncontracting map $\Sigma\to \tfrac12\cdot \SSS^2$, where $\tfrac12\cdot \SSS^2$ denotes the standard 2-sphere rescaled by a factor of $\tfrac12$.
\end{thm}

The proof of the lemma is guided by the following exercise.

\begin{thm}[!]{Exercise}\label{ex:S^3/S^1}
Suppose $\SSS^1\acts\SSS^3$ is an effective isometric action without fixed points.
Let us think of $\SSS^3$ as the unit sphere in $\RR^4$.

\begin{subthm}{ex:S^3/S^1:pq}
Show that one can identify $\RR^4$ with $\CC^2$ so that the action
is given by matrix multiplication
\[\left(\begin{matrix}
u^p&0\\
0& u^q
\end{matrix}
\right),\]
where $(p,q)$ is a pair of relatively prime positive integers and $u\in \SSS^1=\set{z\in\CC}{|z|=1}$.
In particular, our $\SSS^1$ is a subgroup of the torus that acts by
matrix multiplication
\[\left(\begin{matrix}
v&0\\
0& w
\end{matrix}
\right),\]
where $v,w\in \SSS^1$.
\end{subthm}

\smallskip

\noindent Fix $p$ and $q$ as above.
Let $\Sigma_{p,q}=\SSS^3/\SSS^1$ be the quotient space.

\smallskip

\begin{subthm}{ex:S^3/S^1:sphere}
Show that $\Sigma_{p,q}=\SSS^3/\SSS^1$ is a topological sphere with $\SSS^1$-symmetry.
This symmetry has two fixed points, say the north pole and the south pole, that correspond to the orbits of $(1,0)$ and $(0,1)$ in $\SSS^3$.
\end{subthm}

\smallskip

\noindent Denote by $S(r)$ the circle of radius $r$ centered at the north pole of $\Sigma_{p,q}$.

\begin{subthm}{ex:S^3/S^1:a}
Denote by $T(r)$ the inverse image of $S(r)$ in $\SSS^3$, and let $a(r)$ be its area.
Show that $T(r)$ is an orbit of the torus action and
\[a(r)=4\cdot \pi^2\cdot\sin r\cdot \cos r.\]

\end{subthm}

\smallskip

\begin{subthm}{ex:S^3/S^1:b}
Let $b_{p,q}(r)$ be the length of the $\SSS^1$-orbit in $\SSS^3$ that corresponds to a point on $S(r)$. 
Show that
\[b_{p,q}(r)=2\cdot \pi\cdot\sqrt{(p\cdot \cos r)^2+(q\cdot \sin r)^2}.\]
\end{subthm}

\smallskip

\begin{subthm}{ex:S^3/S^1:c}
Let $c_{p,q}(r)$ be the length of $S(r)$.
Show that
\[a(r)=c_{p,q}(r)\cdot b_{p,q}(r).\]
\end{subthm}

\smallskip

\begin{subthm}{ex:S^3/S^1:cc}
Show that $c_{p,q}(r)\le c_{1,1}(r)$ for any pair $(p,q)$ of relatively prime positive integers.
Use it to construct a distance-noncontracting map $\Sigma_{p,q}\to \tfrac12\cdot \SSS^2\iso\Sigma_{1,1}$.
\end{subthm}

\end{thm}

\parit{Proof of \ref{thm:hsiang-kleiner}.}
Assume $\spc{B}=(\SSS^2\times\SSS^2,g)$ is a counterexample.
By the Toponogov theorem, $\spc{B}$ is $\Alex1$.
By \ref{thm:submetry-CBB-1}, the quotient space $\spc{A}\z=\spc{B}/\SSS^1$ is $\Alex1$;
evidently, $\spc{A}$ is 3-dimensional.

Denote by $F\subset \spc{B}$ the fixed point set of the $\SSS^1$-action.
Then $\chi(\spc{B})\z=\chi(F)$.
Each connected component of $F$ is either an isolated point or a 2-dimensional geodesic submanifold in $\spc{B}$;
the latter has to have positive curvature, and therefore it is homeomorphic to $\SSS^2$ or $\RP^2$.
Notice that 
\begin{itemize}
 \item each isolated point contributes 1 to the Euler characteristic of~$\spc{B}$,
 \item each sphere contributes 2 to the Euler characteristic of $\spc{B}$, and
 \item each projective plane contributes 1 to the Euler characteristic of~$\spc{B}$.
\end{itemize}
Since $\chi(\spc{B})=4$, we are in one of the following three cases:
\begin{enumerate}
 \item\label{case1} $F$ has exactly 4 isolated points,
 \item\label{case2} $F$ has one 2-dimensional submanifold and at least 2 isolated points,
 \item\label{case3} $F$ has at least two 2-dimensional submanifolds.
\end{enumerate}
In each case we will arrive at a contradiction.

\parit{Case \ref{case1}.}
Suppose $F$ has exactly 4 isolated points $x_1$, $x_2$, $x_3$, and $x_4$.
Denote by $y_1$, $y_2$, $y_3$, and $y_4$ the corresponding points in $\spc{A}$.
Note that $\Sigma_{y_i}\spc{A}$ is isometric to a quotient of $\SSS^3$ by an isometric $\SSS^1$-action without fixed points.

By \ref{lem:S^3/S^1}, the angle $\mangle\hinge{y_i}{y_j}{y_k}$ is at most $\tfrac\pi2$
for any three distinct points $y_i$, $y_j$, and $y_k$.%
\footnote{It also follows since $y_i$ must be an extremal point, and hence $\diam\Sigma_{y_i}\spc{A}\le\tfrac\pi2$.}
In particular, all four triangles $[y_1y_2y_3]$, $[y_1y_2y_4]$, $[y_1y_3y_4]$, and $[y_2y_3y_4]$ are nondegenerate.
By the comparison, the sum of the angles in each triangle is strictly greater than $\pi$.

Denote by $\omega$ the sum of all 12 angles in the 4 triangles $[y_1y_2y_3]$, $[y_1y_2y_4]$, $[y_1y_3y_4]$, and $[y_2y_3y_4]$.
From above,
\[\omega>4\cdot\pi.\]

On the other hand, by \ref{lem:S^3/S^1} any triangle in $\Sigma_{y_1}\spc{A}$ has perimeter at most $\pi$.
In particular, 
\[\mangle\hinge{y_1}{y_2}{y_3}+\mangle\hinge{y_1}{y_3}{y_4}+\mangle\hinge{y_1}{y_4}{y_2}\le \pi.\]
Apply the same argument in $\Sigma_{y_2}\spc{A}$, $\Sigma_{y_3}\spc{A}$, and $\Sigma_{y_4}\spc{A}$;
adding the results, we get 
\[\omega\le 4\cdot\pi\]
--- a contradiction.

\parit{Case \ref{case2}.}
Suppose $F$ contains exactly one surface $S$.
Then the projection of $S$ to $\spc{A}$ is its boundary $\partial \spc{A}$.
The doubling $\spc{W}$ of $\spc{A}$ across its boundary has at least 4 singular points --- each singular point of $\spc{A}$ corresponds to two singular points of $\spc{W}$.

By the doubling theorem, $\spc{W}$ is an $\Alex1$ space.
Therefore we arrive at a contradiction in the same way as in the first case.

\parit{Case \ref{case3}.} This case is impossible by \ref{ex:bry-connected}.
\qeds

\section{Erdős' problem revisited}

A point $p$ in an Alexandrov space is called \index{extremal point}\emph{extremal} if $\mangle\hinge pxy\le \tfrac\pi2$ for any hinge $\hinge pxy$ with the vertex at $p$; equivalently, $\diam \Sigma_p\le \pi/2$.

\begin{thm}{Theorem}\label{thm:extr-point}
Let $\spc{A}$ be a compact $m$-dimensional $\Alex0$ space.
Then it has at most $2^m$ extremal points.
\end{thm}

\parit{Proof of \ref{thm:extr-point}.}
Let $\{p_1,\dots,p_N\}$ be extremal points in $\spc{A}$.
For each $p_i$ consider its open \index{Voronoi domain}\emph{Voronoi domain} $V_i$; that is, 
\[V_i=\set{x\in \spc{A}}{\dist{p_i}{x}{}<\dist{p_j}{x}{}\ \text{for any}\ j\ne i}.\]
Clearly $V_i\cap V_j=\emptyset$ if $i\not=j$.

Choose $0<\alpha\le 1$.
Given a point $x\in\spc{A}$, choose a geodesic $[p_ix]$ and denote by $x_i$ the point on $[p_ix]$ such that $\dist{p_i}{x_i}{}=\alpha\cdot\dist{p_i}{x}{}$;
let $\map_i\:x\mapsto x_i$ be the corresponding map.
By the comparison, 
\[\dist{x_i}{y_i}{}\ge\alpha\cdot \dist{x}{y}{}\]
for any $x$, $y$, and $i$.
Therefore 
\[\vol(\map_i \spc{A})\ge\alpha^m\cdot\vol \spc{A}.\]

Now, suppose $\alpha<\tfrac12$.
Then $x_i\in V_i$ for any $x\in \spc{A}$.
Indeed, assume $x_i\notin V_i$;
then there is $p_j$ such that $\dist{p_i}{x_i}{}\ge\dist{p_j}{x_i}{}$.
By comparison, we have $\angk{p_j}{p_i}{x}_{\EE^2}>\tfrac\pi2$;
that is, $p_j$ is not an extremal point.

It follows that $\vol V_i\ge\alpha^m\cdot\vol \spc{A}$
for any $0<\alpha<\tfrac12$; hence 
\[\vol V_i\ge\tfrac1{2^m}\cdot\vol \spc{A}.\]
Since $V_1,\dots,V_N$ are disjoint subsets of $\spc{A}$, we have $N\le 2^m$.
\qeds

\begin{thm}{Exercise}\label{ex:extr-point(-1)}
Find a bound on the number of extremal points in an $\Alex{-1}$ space in terms of its dimension $m$ and diameter $D$.
\end{thm}

\section{Crystallographic actions}

An isometric action $\Gamma\acts \EE^m$ is called \index{crystallographic action}\emph{crystallographic} if it is 
\index{properly discontinuous}\emph{properly discontinuous} (that is, for any compact set $K\subset \EE^m$ and any $x\z\in \EE^m$ there are only finitely many elements $g\in \Gamma$ such that $g\cdot x\in K$) and \index{cocompact action}\emph{cocompact} (that is, the quotient space $\spc{A}=\EE^m/\Gamma$ is compact).

Let $F$ be a maximal finite subgroup of $\Gamma$;
that is, if $F<H<\Gamma$ for a finite group $H$, then $F=H$.
Denote by $\mathfrak{M}(\Gamma)$ the number of maximal finite subgroups of $\Gamma$ up to conjugation.

\begin{thm}{Open question}
Let $\Gamma\acts \EE^m$ be a crystallographic action.
Is it true that $\mathfrak{M}(\Gamma)\le 2^m$?
\end{thm}

Note that any finite subgroup $F$ of $\Gamma$ fixes an affine subspace $A_F$ in $\EE^m$.
If $F$ is maximal, then $A_F$ completely describes $F$.
Indeed, since the action is properly discontinuous, the subgroup of $\Gamma$ that fixes $A_F$ pointwise has to be finite.
This subgroup must contain $F$, but since $F$ is maximal, it must coincide with $F$.

Denote by $\mathfrak{M}_k(\Gamma)$ the number of maximal finite subgroups $F<\Gamma$ (up to conjugation) such that $\dim A_F=k$.

Choose a finite subgroup $F<\Gamma$; consider a conjugate subgroup $F'=g \cdot F \cdot g^{-1}$.
Note that $A_{F'}=g\cdot A_F$.
In particular, the subspaces $A_F$ and $A_{F'}$ have the same image in the quotient space $\spc{A}=\EE^m/\Gamma$.
Therefore, to count subgroups up to conjugation, we need to count the images of their fixed sets.
By the lemma below (\ref{lem:extr/G}), $\mathfrak{M}_0(\Gamma)$ cannot exceed the number of extremal points in $\spc{A}=\EE^m/\Gamma$.
Combining this observation with \ref{thm:extr-point}, we get the following.

\begin{thm}{Proposition}\label{prop:2m}
Let $\Gamma\acts \EE^m$ be a crystallographic action.
Then $\mathfrak{M}_0(\Gamma)\le 2^m$.
\end{thm}

\begin{thm}{Lemma}\label{lem:extr/G}
Let $\Gamma\acts \EE^m$ be a crystallographic action, and let $F$ be a maximal finite subgroup of $\Gamma$ that fixes an isolated point $p$.
Then the image of $p$ in the quotient space $\spc{A}=\EE^m/\Gamma$ is an extremal point.
\end{thm}

\parit{Proof.}
Let $q$ be the image of $p$.
Suppose $q$ is not extremal;
that is, $\mangle \hinge q{y_1}{y_2}>\tfrac\pi2$ for some hinge $\hinge q{y_1}{y_2}$ in $\spc{A}$.

Choose the inverse images $x_1,x_2\in \EE^m$ of $y_1,y_2\in \spc{A}$ such that $\dist{p}{x_i}{\EE^m}=\dist{q}{y_i}{\spc{A}}$.
Note that $\mangle \hinge p{x_1}{x_2}\ge \mangle \hinge q{y_1}{y_2}>\tfrac\pi2$.
Moreover, since $p$ is fixed by $F$, we have
\[\mangle \hinge p{x_1}{g\cdot x_2}>\tfrac\pi2
\eqlbl{eq:>pi/2}\]
for any $g\in F$.

Denote by $z\in\EE^m$ the barycenter of the orbit $F\cdot x_2$.
Note that $z$ is a fixed point of $F$.
By \ref{eq:>pi/2}, $z\ne p$;
so $F$ must fix the line $pz$.
But $p$ is an isolated fixed point of $F$ --- a contradiction.
\qeds

\begin{thm}{Exercise}\label{ex:number(m-1)}
Let $\Gamma\acts \EE^m$ be a crystallographic action.
Show that
\begin{subthm}{ex:number(m-1):2}
$\mathfrak{M}_{m-1}(\Gamma)\le 2$, and
\end{subthm}

\begin{subthm}{ex:number(m-1):1}
if $\mathfrak{M}_{m-1}(\Gamma)=1$, then $\mathfrak{M}_0(\Gamma)\le 2^{m-1}$.
\end{subthm}

Construct crystallographic actions with equalities in \ref{SHORT.ex:number(m-1):2} and \ref{SHORT.ex:number(m-1):1}.
\end{thm}

\section{Remarks}

Submetries were introduced by Valerii Berestovskii \cite{berestovskii1987} and have attracted attention in various contexts of differential and metric geometry.

A more general form of Theorem \ref{thm:hsiang-kleiner} was proved by Karsten Grove and Burkhard Wilking \cite{grove-wilking};
it classifies isometric $\SSS^1$ actions on 4-dimensional manifolds with nonnegative sectional curvature.
Their proof is as beautiful as the original work of Wu-Yi Hsiang and Bruce Kleiner.

It is expected that \textit{no $\Alex1$ space with a nontrivial isometric $\SSS^1$-action can be homeomorphic to $\SSS^2\times\SSS^2$};
that is, \ref{thm:hsiang-kleiner} holds for a general $\Alex1$ space.
The proof of \ref{thm:hsiang-kleiner} would work if we had the following generalization of \ref{lem:S^3/S^1};
see \cite{harvey-searle}.

\begin{thm}{Open question}
Let $\Sigma$ be an $\Alex1$ space homeomorphic to $\SSS^3$.
Suppose $\SSS^1$ acts on $\Sigma$ isometrically and without fixed points.
Is it true that any triangle in $\Sigma/\SSS^1$ has perimeter at most $\pi$?
(And if the answer is yes, is there a distance-noncontracting map $\Sigma/\SSS^1\z\to \tfrac12\cdot\SSS^2$?)
\end{thm}

{\sloppy

\begin{thm}{Advanced exercise}\label{ex:S1actsS3}
Suppose $\SSS^1$ acts isometrically on an $\Alex1$ space $\spc{A}$ that is homeomorphic to $\SSS^3$.
Assume its fixed-point set is a closed local geodesic $\gamma$.
Show that
\[\length\gamma\le2\cdot\pi.\]
\end{thm}

}

An analogous question for a $\ZZ_2$-action is open \cite{petrunin-involution}.

Theorem \ref{thm:extr-point} is a translation of the following classical problem in discrete geometry to Alexandrov's language.

\begin{thm}{Problem}\label{erdos-problem}
Let $F$ be a set of points in $\EE^m$ such that any triangle formed by three distinct points in $F$ has no obtuse angles.
Then $|F|\z\le2^m$.
Moreover, if $|F|=2^m$, then $F$ consists of the vertices of an $m$-dimensional rectangle.
\end{thm}

This problem was posed by Paul Erdős \cite{erdos} and solved by Ludwig Danzer and Branko Grünbaum \cite{danzer-gruenbaum}.
Grigory Perelman noticed that, after proper definitions, the same proof works in Alexandrov spaces \cite{perelman-Erdos}; thus, it proves \ref{thm:extr-point}.
Applying our argument to the convex hull of $F$ in \ref{erdos-problem} proves that $|F|\le 2^m$;
classifying the case of equality requires more work.

Compact $m$-dimensional $\Alex0$ spaces with the maximal number of extremal points include $m$-dimensional rectangles and the quotients of flat tori by reflections across a point
(this action has $2^m$ isolated fixed points; each corresponds to an extremal point in the quotient space $\spc{A}=\TT^m/\ZZ_2$).
The second author has proved \cite{lebedeva2015} that \textit{every $m$-dimensional $\Alex0$ space with $2^m$ extremal points is a quotient of Euclidean space by a crystallographic action}.

Extremal subsets of Alexandrov spaces were briefly discussed in \ref{sec:bry-remarks}.
The following definition is more relevant to isometric group actions.

A closed subset $E$ of a finite-dimensional Alexandrov space is called
\index{extremal set}\emph{extremal} if $\mangle\hinge pxy\z\le \tfrac\pi2$ for any hinge $\hinge pxy$ such that $x\notin E$ and $p\in E$ minimizes $\dist{x}{p}{}$.
A nonempty extremal set is called \index{minimal extremal set}\emph{minimal} if it contains no proper extremal subsets.

For example, the whole space and the empty set are extremal.
The vertices, edges, and faces of a cube, as well as their unions, are extremal subsets.
Vertices of the cube are its only minimal extremal subsets.

Counting maximal finite subgroups in a crystallographic group $\Gamma$ (up to conjugation) is equivalent to counting the minimal extremal subsets in the quotient space $\spc{A}=\EE^m/\Gamma$.
So, \ref{prop:2m} would follow from the next conjecture.

\begin{thm}{Conjecture}
Any $m$-dimensional compact $\Alex0$ space has at most $2^m$ minimal extremal subsets.
\end{thm}

Let us mention another related conjecture.
A nonempty extremal set is called \index{primitive extremal set}\emph{primitive} if it contains no proper extremal subsets with nonempty relative interior.
(Given an Alexandrov space and a point $p$ in it, the set of points reachable from $p$ by a recursive application of the gradient flow for squares of distance functions forms a primitive extremal set.)
For example, each face of an $m$-dimensional cube is a primitive extremal subset;
therefore the cube has exactly $3^m$ primitive extremal subsets, including the whole cube.

\begin{thm}{Conjecture}
Any $m$-dimensional compact $\Alex0$ space has at most $3^m$ primitive extremal subsets.
\end{thm}

Some nonsharp explicit estimates on the number of extremal subsets
were obtained by Tadashi Fujioka \cite{fujioka2025};
these results are closely related to Gromov's Betti number theorem \ref{thm:betti}.

%% file: overview.tex
\chapter{Surfaces}\label{chap:surfaces}

This lecture is less rigorous; it offers a glimpse into the geometry of convex surfaces, a major precursor of modern Alexandrov geometry.
For a deeper dive into this theory, we recommend the classic and brilliantly written books by Alexandr Alexandrov \cite{alexandrov, alexandrov-1948}.
Additionally, the book by Alexey Pogorelov \cite{pogorelov1969} is highly recommended, despite being a challenging read.

\section{Polyhedral surfaces}

A \index{polyhedral surface}\emph{polyhedral surface} is defined as a 2-dimensional manifold (possibly with boundary) equipped with a length metric such that there is a finite triangulation in which each triangle is isometric to a Euclidean triangle.
A \index{triangulation}\emph{triangulation} of a polyhedral surface is always assumed to satisfy this condition.

Note that, according to our definition, every polyhedral surface is compact.

Consider a point $p$ on a polyhedral surface $\spc{P}$.
We can assume that $p$ is a vertex of a triangulation of $\spc{P}$,
which can be achieved by subdividing the triangulation.
Let $\theta_p$ denote the {}\emph{total angle} around $p$, which is the sum of all angles at $p$ in the triangles that have $p$ as a vertex.

Note that $\theta_p$ is independent of the choice of triangulation.
If $p$ is an interior point, then the value $2\cdot\pi-\theta_p$ is called the \index{curvature}\emph{curvature} at $p$.
If $p$ lies on the boundary of $\spc{P}$, then the value $\pi-\theta_p$ is called the \index{inner turn}\emph{inner turn} at $p$.

A point with nonzero curvature or inner turn will be called an \index{essential vertex}\emph{essential vertex} of the surface.
Note that an essential vertex is a vertex in any triangulation.

\begin{thm}[!]{Exercise}\label{ex:geodesic-vertex}
Show that geodesics on a polyhedral surface with nonnegative curvature and nonnegative inner turns can have essential vertices only at their endpoints.
\end{thm}

The following statement is an analog (and, properly interpreted, a special case of) of the Gauss--Bonnet formula.

\begin{thm}{Exercise}\label{ex:gauss-bonnet}
Let $\spc{P}$ be a polyhedral surface.
Denote by $K(\spc{P})$ and $T(\spc{P})$ the sum of the curvatures at all interior points
and the sum of the inner turns at all boundary points of $\spc{P}$, respectively.
Show that
\[
K(\spc{P}) + T(\spc{P}) = 2 \cdot \pi \cdot \chi(\spc{P}),
\]
where $\chi(\spc{P})$ denotes the Euler characteristic of $\spc{P}$.
\end{thm}

The following proposition states that the new definition of curvature is consistent with the $\Alex0$ comparison.

\begin{thm}{Proposition}\label{prop:poly-CBB}
A polyhedral surface is $\Alex0$ if and only if it has nonnegative curvature at every interior point and nonnegative inner turn at each boundary point.
\end{thm}

\parit{Proof.}
To prove the if direction, observe that near an interior point the surface is locally isometric to a cone over a circle of length $\le 2\cdot\pi$ and near a boundary point it is locally isometric to the cone over an interval of length $\le \pi$.
Then apply \ref{ex:cone-CBB} and globalization.




\begin{thm}{Exercise}\label{ex:poly-CBB}
The converse is left to the reader.\qeds
\end{thm}

\section{Approximation}

The following theorem is a key additional tool in the Alexandrov geometry of surfaces.
We will use this result in the proof of \ref{cor:Alex0-convex} to reduce questions about $\Alex0$ surfaces to polyhedral surfaces with nonnegative curvature.

\begin{thm}{Theorem}\label{thm:approximation}
Any closed $\Alex0$ surface is a Gromov--Hausdorff limit of homeomorphic polyhedral surfaces with nonnegative curvature.
\end{thm}

We give only a sketch of the proof, but the sketch is not short.
First we apply the following exercise to construct an approximation $\tilde{\spc{P}}_\delta$ of a given surface $\spc{P}$.
It reduces the theorem to saying that this is indeed an approximation (\ref{clm:approximation}).
Further, we use a couple of deep results (\ref{thm:cont-vol} and \ref{thm:vol-short}) to prove the claim.

\begin{thm}{Exercise}\label{ex:construction}
Let $\spc{P}$ be a closed $\Alex0$ surface.

\begin{subthm}{ex:approximation:nbhd}
Show that any point $p$ admits an arbitrarily small closed convex polygonal neighborhood $N$;
that is, $N$ is convex and bounded by a broken geodesic.

\end{subthm}

\begin{subthm}{ex:approximation:triangulation}
Given $\delta>0$, show that $\spc{P}$ admits a triangulation $\tau$ by convex triangles
with positive inner turn at each vertex and diameter smaller than $\delta$.
\end{subthm}

\begin{subthm}{ex:approximation:poly}
Suppose that $v$ is a vertex of a triangulation $\tau$ of $\spc{P}$ by convex triangles.
Let $\theta_v$ be the sum of angles at $v$ in all the triangles of~$\tau$.
Show that $\theta_v\le 2\cdot\pi$.
\end{subthm}

\end{thm}

\parit{Construction.}
Let $\spc{P}$ be a closed $\Alex0$ surface.
By part \ref{SHORT.ex:approximation:triangulation}, we can triangulate $\spc{P}$ by small convex triangles, say, the diameter of each triangle is less than a given $\delta>0$.
Now, replace each triangle in the triangulation with its corresponding model solid triangle, and denote the resulting polyhedral surface by $\tilde{\spc{P}}_\delta$.

Note that $\tilde{\spc{P}}_\delta$ is homeomorphic to $\spc{P}$;
moreover, there exists a homeomorphism $\spc{P} \to \tilde{\spc{P}}_\delta$ (it will be denoted by $x\mapsto \tilde x$) that maps each edge in the triangulation of $\spc{P}$ isometrically to the corresponding edge of $\tilde{\spc{P}}_\delta$, and each triangle in $\spc{P}$ is mapped to the corresponding model triangle in $\tilde{\spc{P}}_\delta$.

By the hinge comparison (\ref{angle}) and part \ref{SHORT.ex:approximation:poly} of the exercise, the total angle around each vertex in $\tilde{\spc{P}}_\delta$ does not exceed $2\cdot\pi$.
Thus, the resulting polyhedral surface $\tilde{\spc{P}}_\delta$ has nonnegative curvature.

\medskip

Observe that Theorem \ref{thm:approximation} follows from the following statement.

\begin{thm}{Claim}\label{clm:approximation}
If $\tilde{\spc{P}}_\delta$ is provided by the construction, then $\tilde{\spc{P}}_\delta\to\spc{P}$ as $\delta\to 0$ in the sense of Gromov--Hausdorff.
\end{thm}

This claim looks self-evident, but it is not.
A remarkably elegant proof was given in \cite[VII §~6]{alexandrov-1948}.
We will outline an alternative proof based on the following exercise and two theorems, which will be stated without proof.

The first theorem is due to Yuri Burago, Mikhael Gromov, and Grigori Perelman \cite[10.8]{burago-gromov-perelman};
it generalizes Alexandrov's theorem for surfaces \cite[X §~2]{alexandrov-1948}.
The second theorem is due to Nan Li \cite{li}.

\begin{thm}{Theorem}\label{thm:cont-vol}
Let $\spc{X}_1, \spc{X}_2,\dots$ be $m$-dimensional $\Alex\kappa$ spaces that converge to $\spc{X}_\infty$ in the sense of Gromov--Hausdorff.
Then $m$-volumes on $\spc{X}_i$ weakly converge to the $m$-volume on $\spc{X}_\infty$.
\end{thm}

\begin{thm}{Theorem}\label{thm:vol-short}
Let $\spc{X}$ be an Alexandrov space without boundary, and let $\spc{Y}$ be an arbitrary Alexandrov space of the same dimension.
Then any short volume-preserving map $\spc{X}\to\spc{Y}$ is an isometry.
\end{thm}

Suppose a convex solid triangle $\Delta$ in an $\Alex0$ surface has angles $\alpha$, $\beta$, and $\gamma$.
Let us define its excess by
\[\excess\Delta=\alpha+\beta+\gamma-\pi.\]
Since the angles of a model triangle sum to $\pi$,
the hinge comparison (\ref{angle}) implies that the excess is nonnegative.

\begin{thm}{Exercise}\label{ex:approximation}
Let $\tau$ be a triangulation of a closed $\Alex0$ surface $\spc{P}$ by convex triangles $\Delta_1,\dots,\Delta_n$.


\begin{subthm}{ex:approximation:excess}
Show that
\[\excess\Delta_1+\dots+\excess\Delta_n \le 2\cdot\pi\cdot \chi(\spc{P}),\]
where $\chi(\spc{P})$ denotes the Euler characteristic of $\spc{P}$.
\end{subthm}

\begin{subthm}{ex:approximation:length}
Let $x$ and $y$ be points on the sides of a triangle $\Delta_i$, and let $\tilde x$ and $\tilde y$ be the corresponding points in the corresponding triangle $\tilde \Delta$ in~$\tilde{\spc{P}}$.
Show that
\[\dist{\tilde x}{\tilde y}{\tilde \Delta}\le\dist{x}{y}{\Delta}\le \dist{\tilde x}{\tilde y}{\tilde \Delta}+\excess\Delta\cdot \diam \Delta.\]

\end{subthm}

\begin{subthm}{ex:approximation:area}
Let $\Delta$ be a solid triangle in the triangulation $\tau$ of $\spc{P}$, and let $\tilde \Delta$ be the corresponding triangle in $\tilde{\spc{P}}$.
Show that
\[\area \tilde \Delta\le \area \Delta\le \area \tilde \Delta+\tfrac12\cdot\excess\Delta\cdot (\diam \Delta)^2.\]

\end{subthm}

\end{thm}

Note that part \ref{SHORT.ex:approximation:excess} implies that $\chi(\spc{P})\ge 0$.
Therefore, $\spc{P}$ is homeomorphic to a sphere, projective plane, torus, or Klein bottle.
In the latter two cases, the construction produces a flat surface $\tilde{\spc{P}}_\delta$, which has to be isometric to $\spc{P}$.
Therefore the cases of the sphere and the projective plane are more interesting.

\parit{Proof of \ref{clm:approximation}.}
Choose a sequence of positive numbers $\delta_n\to0$;
let $\tilde{\spc{P}}_{\delta_n}$ be polyhedral surfaces provided by the construction and let $\tau_n$ be the corresponding triangulation.

According to part \ref{SHORT.ex:approximation:length} of the exercise, the spaces $\tilde{\spc{P}}_{\delta_n}$ have bounded diameter.
Therefore by Gromov's selection theorem, we can pass to a converging sequence of $\tilde{\spc{P}}_{\delta_n}$;
denote its Gromov--Hausdorff limit by $\tilde{\spc{P}}$.
If $\tilde{\spc{P}}_{\delta_n}$ does not converge to $\spc{P}$, then we can assume that $\tilde{\spc{P}}$ is not isometric to $\spc{P}$.

Choose two points $x,y\in \spc{P}$, and connect them by a geodesic.
Denote by $s_1,\dots,s_m$ the points of the geodesic on the sides of the triangulation $\tau_n$;
we assume that these points appear in the same order on the geodesic.
Denote by $\tilde x$, $\tilde y$, and $\tilde s_1,\dots,\tilde s_m$ the corresponding points in $\tilde{\spc{P}}_{\delta_n}$.
By part \ref{SHORT.ex:approximation:length} of the exercise,
\[\dist{\tilde s_{i-1}}{\tilde s_i}{\tilde{\spc{P}}_{\delta_n}}\le\dist{s_{i-1}}{s_i}{\spc{P}},\]
Note also that
\[\dist{\tilde x}{\tilde s_1}{\tilde{\spc{P}}_{\delta_n}}\le \delta_n
\quad\text{and}\quad
\dist{\tilde s_m}{\tilde y}{\tilde{\spc{P}}_{\delta_n}}\le \delta_n
\]
Therefore
\[\dist{\tilde x}{\tilde y}{\tilde{\spc{P}}_{\delta_n}}\le\dist{x}{y}{\spc{P}}+2\cdot\delta_n.\]

Passing to the limit, we get a short onto map $\spc{P}\to \tilde{\spc{P}}$.
On the other hand, applying parts \ref{SHORT.ex:approximation:excess} and \ref{SHORT.ex:approximation:area}, we get that
\[
\area \spc{P} -\pi\cdot\chi(\spc{P})\cdot \delta_n^2
\le
\area \tilde{\spc{P}}_{\delta_n}
\le
\area \spc{P}.
\]
By Theorem \ref{thm:cont-vol}, $\area \spc{P}= \area \tilde{\spc{P}}$.
Applying Theorem \ref{thm:vol-short}, we get that the short map $\spc{P}\to \tilde{\spc{P}}$ is an isometry --- a contradiction.
\qeds

\parit{Remark.}
The main difficulty in the proof comes from the nonconvexity of triangles in the triangulation of $\tilde{\spc{P}}_\delta$.
If these triangles were convex, then the estimates in \ref{SHORT.ex:approximation:length} and \ref{SHORT.ex:approximation:excess} would imply that $\tilde{\spc{P}}_\delta$ is close to $\spc{P}$ in the sense of Gromov--Hausdorff.

\section{Surface of polyhedrons and bodies}

Let us define a \index{convex body}\emph{convex body} as a compact convex subset of $\EE^3$ with nonempty interior.
The \index{surface of a convex body}\emph{surface} of a convex body is defined as its boundary equipped with the induced length metric.

\begin{thm}[!]{Exercise}\label{ex:surf-S2}
Show that the surface of a convex body is homeomorphic to $\SSS^2$.
\end{thm}

A \index{convex polyhedron}\emph{convex polyhedron} is a convex body
that is the convex hull of a finite set of points; the essential points in this set are called its \index{vertex}\emph{vertices}.

Note that the surface of a convex polyhedron $K$ is a closed polyhedral surface.

\begin{thm}{Exercise}\label{pr:tetrahedron}
Assume that the surface of a nondegenerate tetrahedron $T$ has curvature $\pi$ at each vertex.
Show that

\begin{subthm}{pr:tetrahedron:=}
all faces of $T$ are congruent;
\end{subthm}

\begin{subthm}{pr:tetrahedron:perp} the line containing the midpoints of opposite edges of $T$ intersects these edges at right angles.
\end{subthm}

\end{thm}

\begin{thm}{Claim}\label{clm:total-angle}
The surface $\spc{P}$ of any convex polyhedron $K$ has nonnegative curvature.
Moreover, a point $v$ is a vertex of $K$ if and only if
$v$ is an essential vertex of $\spc{P}$.
\end{thm}

A proof of this claim can be found in a school textbook \cite[§~48]{kiselev-stereo-en};
one can also deduce it from \ref{claim:angle-3angle-inq}.

\begin{thm}[!]{Exercise}\label{ex:surface-covergence}
Let $K_1,K_2,\dots,$ and $K_\infty$ be convex bodies in $\EE^m$.
Denote by $\spc{P}_n$ the surface of $K_n$.
Suppose $K_n\z\to K_\infty$ in the sense of Hausdorff.
Show that $\spc{P}_n\to \spc{P}_\infty$ in the sense of Gromov--Hausdorff.
\end{thm}

Since any convex body is a Hausdorff limit of a sequence of convex polyhedrons, the next proposition follows from \ref{prop:poly-CBB}, \ref{ex:surface-covergence}, and \ref{thm:CBB-closed}.

\begin{thm}{Proposition}\label{prop:conv-surf-CBB(0)}
The surface of a convex body in $\EE^3$ is $\Alex0$.

\end{thm}

This can also be shown by approximating the defining convex function of a convex body by a smooth convex function using a convolution.
By the Gauss formula, the level sets of the smooth function have nonnegative Gauss curvature and hence are $\Alex0$.
Now the result again follows by \ref{ex:surface-covergence}.

\section{Uniqueness theorem}

\begin{thm}{Theorem}\label{thm:alexandrov-uni'}
Any two convex polyhedrons in $\EE^3$ with isometric surfaces are congruent.

Moreover, any isometry between the surfaces of convex polyhedrons can be extended to an isometry of the whole $\EE^3$.
\end{thm}

Assuming the isometry between the surfaces to be face-to-face gives a reformulation of Cauchy's theorem; his argument, with a small addition, proves \ref{thm:alexandrov-uni'}.

First, let us recall Cauchy's proof, assuming the reader knows it.
If not, then read it in one of the classical texts; see, for example, \cite{aigner-zigler,dolbilin,tabacnikov-fuks}.

\parit{Sketch of Cauchy's proof.}
Suppose $K$ and $K'$ are convex polyhedrons in $\EE^3$;
denote their surfaces
by $\spc{P}$ and $\spc{P}'$.
Suppose there is an isometry $\iota\:\spc{P}\to \spc{P}'$ that sends each face of $K$ to a face of $K'$.

Let us mark an edge of $K$ with ``$+$'' (or ``$-$'') if the dihedral angle at this edge in $K$ is smaller (respectively, larger) than the corresponding angle in $K'$.
Further, we consider the graph $\Gamma$ that is formed by all marked edges.
If $\Gamma$ is empty, then Cauchy's theorem follows; assume the contrary.

The graph $\Gamma$ is embedded into $\spc{P}$, which is homeomorphic to the sphere.
In particular, the edges coming from one vertex have a natural cyclic order.
Given a vertex $v$ of $\Gamma$, we can count the \textit{number of sign changes} around $v$;
that is, the number of consecutive pairs of edges with different signs.

We need two statements:

\begin{thm}{Local lemma}
The number of sign changes at each vertex of $\Gamma$ is at least $4$.
\end{thm}

\begin{thm}{Global lemma}
No (nonempty) planar graph meets the condition of the local lemma.
\end{thm}

Once the lemmas are proved, Cauchy's theorem follows.
\qeds

Once more, the argument above is written only to make sure we are on the same page;
it will not work without reading the actual proof.

\parit{Alexandrov's addition.}
We need to remove the assumption that the isometry $\iota\:\spc{P}\z\to \spc{P}'$ is face-to-face.
Mark in $\spc{P}$ all the edges of $K$ as we did above.
In addition, if an edge in $K'$ does not correspond to an edge of $K$, then mark its inverse image in $K$ with ``$-$''; these lines on $K$ will be referred to as \index{fake edges and vertices}\emph{fake edges}.

The marked lines divide $\spc{P}$ into regions, and the restriction of $\iota$ to each region is a rigid motion.
A vertex of the obtained graph can be a vertex of $K$, or it can be a {}\emph{fake vertex};
that is, it might be an intersection of an edge and a fake edge (as $v$ in the picture).

\begin{figure}[ht!]
\vskip-0mm
\centering
\includegraphics{mppics/pic-50}
\vskip-0mm
\end{figure}

For an actual vertex, the local lemma can be proved the same way as in Cauchy's theorem.
For a fake vertex $v$, it is easy to see that both parts of the edge coming thru $v$ are marked with minus
while both of the fake edges at $v$ are marked with plus.
Therefore, we still have at least four sign changes at $v$.
The remaining argument works as before.
\qeds

\section{Existence theorem}

By \ref{prop:poly-CBB}, \ref{clm:total-angle}, and \ref{ex:surf-S2}, the surface of a convex polyhedron is $\Alex0$ and is homeomorphic to the sphere.
Alexandrov's theorem states that the converse holds if one also allows \textit{twice covered polygons}.
In other words, we have to consider a plane polygon as a degenerate polyhedron;
in this case, its surface is defined as the doubling of the polygon across its boundary.

From now on, we assume that \textit{a polyhedron can degenerate to a plane polygon}.

\begin{thm}{Theorem}\label{thm:alexandrov-first}
A polyhedral metric on the two-sphere is isometric to the surface of a convex polyhedron (possibly degenerate) if and only if it has nonnegative curvature.

\end{thm}

Applying this theorem with the approximation theorem (\ref{thm:approximation}) and \ref{ex:surface-covergence}, we get the following statement.
Here we again assume that a convex body can degenerate to a convex plane figure,
and, in this case, its surface is defined as the doubling of the figure across its boundary.

\begin{thm}{Corollary}\label{cor:Alex0-convex}
A metric on the two-sphere is $\Alex0$ if and only if it is isometric to the surface of a convex body (possibly degenerate).

\end{thm}

The proof of the existence theorem will be discussed in the following two sections.
It is instructive to solve the following exercise before going further.

\begin{thm}[!]{Exercise}\label{ex:alexandrov=<4}
Let $\spc{P}$ be the 2-sphere equipped with a polyhedral metric with nonnegative curvature.

\begin{subthm}{ex:alexandrov=<4:>=3}
Prove that $\spc{P}$ has at least 3 essential vertices.
\end{subthm}

\begin{subthm}{ex:alexandrov=<4:=3}
If $\spc{P}$ has exactly 3 essential vertices $u$, $v$, and $w$, then it is isometric to the doubling of the solid model triangle $\modtrig(uvw)$.
\end{subthm}

\begin{subthm}{ex:alexandrov=<4:4}
If $\spc{P}$ has exactly 4 essential vertices, then it is isometric to the surface of a tetrahedron (possibly degenerate to a quadrangle).
\end{subthm}

\end{thm}

\section{Reformulation}

In this section, we introduce several notions and use them to reformulate the existence theorem (\ref{thm:reformulation}).

\paragraph{Space of polyhedrons.}
Let us denote by $\bm{K}$ the space of all convex polyhedrons in $\EE^3$,
including polyhedrons that degenerate to a plane polygon.
Polyhedrons in $\bm{K}$ will be considered up to a motion of the space; we will not distinguish between a convex polyhedron and its congruence class.

The space $\bm{K}$ will be considered with the topology induced by the {}\emph{Hausdorff metric up to a motion};
that is, the distance between (congruence classes of) polyhedrons $K$ and $L$ is defined by
\[\dist{K}{L}{}\df \inf_\mu \{\,\dist{K}{\mu(L)}{\Haus}\,\},\]
where $\mu$ runs among all motions of $\EE^3$ and $\dist{K}{L}{\Haus}$ denotes Hausdorff distance; see Section~\ref{sec:Hausdorff metric}.

We say that a polyhedron $K$ in $\bm{K}$ has \index{no symmetries}\emph{no symmetries} if $K\z\ne \mu(K)$ for any nontrivial motion $\mu$ of $\EE^3$.
The set of all polyhedrons without symmetries in $\bm{K}$ will be denoted by $\bm{K}^\circ$.
Observe that $\bm{K}^\circ$ is an open set in $\bm{K}$.

Further, denote by $\bm{K}_n$ the polyhedrons in $\bm{K}$ with exactly $n$ vertices, and let $\bm{K}_n^\circ=\bm{K}_n\cap \bm{K}^\circ$.
Since any polyhedron has at least 3 vertices, the space $\bm{K}$ admits a subdivision into a countable number of subsets $\bm{K}_3,\bm{K}_4,\dots$

\paragraph{Space of surfaces.}
The space of polyhedral surfaces with nonnegative curvature that are homeomorphic to the 2-sphere will be denoted by $\bm{P}$.
The surfaces in $\bm{P}$ will be considered up to an isometry, and the whole space $\bm{P}$ will be equipped with the natural topology induced by the Gromov--Hausdorff metric.

We say that a surface $\spc{P}$ in $\bm{P}$ has \index{no symmetries}\emph{no symmetries} if there is no nontrivial isometry
$\mu\:\spc{P}\selfmap$.
The set of all surfaces without symmetries in $\bm{P}$ will be denoted by $\bm{P}^\circ$.
Observe that $\bm{P}^\circ$ is an open set in $\bm{P}$.

The subset of $\bm{P}$ of all surfaces with exactly $n$ essential vertices will be denoted by $\bm{P}_n$; let $\bm{P}_n^\circ=\bm{P}_n\cap \bm{P}^\circ$.
By \ref{ex:alexandrov=<4:>=3}, any surface in $\bm{P}$ has at least 3 essential vertices.
Therefore $\bm{P}$ is subdivided into countably many subsets $\bm{P}_3,\bm{P}_4,\dots$

\paragraph{From a polyhedron to its surface.}
Recall that the surface of a convex polyhedron is a sphere with nonnegative curvature.
Therefore, passing from a polyhedron to its surface, we define a map
\[\iota\:\bm{K}\to \bm{P}.\]

The existence theorem (\ref{thm:alexandrov-first}) follows from the next statement.

\begin{thm}{Theorem}\label{thm:reformulation}
For any integer $n\ge 3$,
the map $\iota$ defines a bijection from $\bm{K}_n$ to~$\bm{P}_n$.
\end{thm}

\section{About the proof of existence}

By \ref{ex:surface-covergence}, the map $\iota\:\bm{K}\to\bm{P}$ is continuous.
Combining \ref{clm:total-angle} with the uniqueness theorem (\ref{thm:alexandrov-uni'}), we get that $\iota(\bm{K}_n)\subset \bm{P}_n$ and the map $\iota\:\bm{K}_n\to\bm{P}_n$ is injective.
It remains to prove the following.

\begin{thm}{Claim}\label{clm:surjective}
For any $n\ge 3$, the map $\iota\:\bm{K}_n\to\bm{P}_n$ is surjective.
\end{thm}

The proof is based on the construction of a one-parameter family of polyhedrons that starts at an arbitrary polyhedron
and ends at a polyhedron with its surface isometric to the given surface $\spc{P}$.
This is called the \index{continuity method}\emph{continuity method}; it is often used in the theory of differential equations.

\medskip

Now let us get into details.
First, observe that the second part of the uniqueness theorem (\ref{thm:alexandrov-uni'}) implies that $\iota(\bm{K}_n^\circ)\subset \bm{P}_n^\circ$.

\begin{thm}{Lemma}\label{lem:connected}
For any integer $n\ge 4$, the space $\bm{P}_n^\circ$ is connected and dense in $\bm{P}_n$.
\end{thm}

Note that $\bm{P}_3^\circ=\emptyset$;
indeed the surface of a triangle admits a reflection symmetry.
The case $n=4$ can be deduced from \ref{ex:alexandrov=<4:4}; thus, we can assume that $n\ge 5$.

The density of $\bm{P}_n^\circ$ in $\bm{P}_n$ follows from a general-position-type argument.

The proof of the first statement is not complicated, but it requires ingenuity;
it can be done by the direct construction of a one-parameter family of surfaces in $\bm{P}_n^\circ$ that connects two given surfaces.
Such a family can be obtained by concatenating the following elementary deformations (direct or reversed).

Start with a surface $\spc{P}$ from $\bm{P}_n^\circ$.
Suppose $v$ and $w$ are essential vertices in $\spc{P}$.
Let us cut $\spc{P}$ along a geodesic from $v$ to~$w$.
This way we obtain a sphere with a hole.
The hole can be patched by a family of patches so that the obtained surface remains in $\bm{P}_n$;
in other words, the surface still has exactly $n$ essential vertices.
(After the patching, the vertices $v$ and $w$ may become inessential.
There is a three-parameter family of such patches, so we have something to choose from; see the picture.)
Choosing a one-parameter family of such patches, we can get an elementary deformation of~$\spc{P}$.

\begin{figure}[t!]
\vskip-0mm
\centering
\includegraphics{mppics/pic-1172}
\vskip-0mm
{\captionsetup{width=.85\linewidth}
\caption*{
The surface of a round cone is shown from above.
The gray annular region is squeezed between a pair of geodesic bigons.
A patch (which must be homeomorphic to a disk) can be obtained from the region by gluing the two sides of one bigon.
The vertices of the glued bigon become essential vertices after patching.
The vertices of the unglued bigon must be attached to the vertices $v$ and $w$.}}
\end{figure}

Again, applying a general-position-type argument to the above construction, we get a path in $\bm{P}_n^\circ$, assuming that the starting and ending surfaces are in $\bm{P}_n^\circ$.

\begin{thm}{Lemma}\label{lem:open}
The map $\iota\:\bm{K}_n^\circ\to\bm{P}_n^\circ$ is open;
that is, it maps any open set in $\bm{K}_n^\circ$ to an open set in $\bm{P}_n^\circ$.

In particular, for any $n\ge 3$, the image $\iota(\bm{K}_n^\circ)$ is open in~$\bm{P}_n^\circ$.
\end{thm}

This statement follows from the so-called \index{invariance of domain}\emph{invariance of domain theorem},
which states that a \textit{continuous injective map between manifolds of the same dimension is open}.

Recall that $\iota$ defines a continuous and injective map $\bm{K}_n^\circ\to\bm{P}_n^\circ$.
It remains to check that both spaces $\bm{K}_n^\circ$ and $\bm{P}_n^\circ$ are $(3\cdot n-6)$-dimensional manifolds.

Choose a polyhedron $K$ from $\bm{K}_n^\circ$.
It is uniquely determined by the $3\cdot n$ coordinates of its $n$ vertices.
We can assume that the first vertex is at the origin,
the second has a positive $x$-coordinate
and the remaining two coordinates vanish,
and the third has a vanishing $z$-coordinate and a positive $y$-coordinate.
Therefore, all polyhedrons in $\bm{K}_n$ that lie sufficiently close to $K$ can be described by $3\cdot n-6$ parameters.
If $K$ has no symmetries, then this description is one-to-one;
in this case, a neighborhood of $K$ in $\bm{K}_n$ admits a parametrization by an open set in $\RR^{3\cdot n-6}$.

The case of surfaces is analogous.
We need to subdivide a given $\spc{P}\in\bm{P}_n^\circ$ into plane triangles using only essential vertices.
By Euler's formula, there are exactly $3\cdot n-6$ edges in this subdivision.
The lengths of the edges completely describe the surface $\spc{P}$ and any surface nearby.
Since $\spc{P}$ has no symmetries, this description is one-to-one, and a neighborhood of $\spc{P}$ in $\bm{P}_n$ admits a parametrization by an open set in $\RR^{3\cdot n-6}$.

\begin{thm}{Lemma}\label{lem:closed}
The map $\iota\:\bm{K}_n\to\bm{P}_n$ is closed;
that is, the image of a closed set in $\bm{K}_n$ is closed in $\bm{P}_n$.

In particular, for any $n\ge 3$, the set $\iota(\bm{K}_n)$ is closed in~$\bm{P}_n$.
\end{thm}

Choose a sequence of polyhedrons $K_1,K_2,\ldots\in\bm{K}_n$.
Assume that the sequence $\spc{P}_i=\iota(K_i)$ converges in $\bm{P}_n$ as $i\to \infty$;
denote its limit by $\spc{P}_\infty$.
We need to construct a polyhedron $K_\infty\in \bm{K}_n$ such that $\iota(K_\infty)=\spc{P}_\infty$;
let us do it.

Passing to a subsequence, we can assume that $K_i$ converges in $\bm{K}$ as $i\to\infty$;
denote the limit polyhedron by $K_\infty$.
Since $\iota$ is continuous, $\iota(K_i)$ converges to $\iota(K_\infty)$ in~$\bm{P}$; so, $\iota(K_\infty)=\spc{P}_\infty$.
Recall that $\iota(\bm{K}_m)\z\subset\bm{P}_m$ for each $m$; therefore, $K_\infty\in \bm{K}_n$.

\parit{Proof of \ref{clm:surjective}.}
The case $n\le 4$ is already solved in \ref{ex:alexandrov=<4}; so we assume that $n\ge 5$.
By \ref{lem:closed} and \ref{lem:open},
$\iota(\bm{K}_n^\circ)$ is a nonempty closed and open set in $\bm{P}_n^\circ$, and $\bm{P}_n^\circ$ is connected.
Therefore, $\iota(\bm{K}_n^\circ)=\bm{P}_n^\circ$.

By \ref{lem:closed}, $\iota(\bm{K}_n)$ is closed in $\bm{P}_n$.
By \ref{lem:connected}, $\bm{P}_n^\circ$ is dense in $\bm{P}_n$.
Since $\iota(\bm{K}_n^\circ)=\bm{P}_n^\circ$, we have $\bm{P}_n^\circ\subset \iota(\bm{K}_n)$;
therefore, $\iota(\bm{K}_n)=\bm{P}_n$;
that is, $\iota\:\bm{K}_n\z\to\bm{P}_n$ is surjective.
\qeds

\section{Ambient space}

On one hand, Alexandrov's surface theory is simpler than that of general Alexandrov spaces since it has extra tools.
On the other hand, these tools come with extra structure, which makes the theory more complicated.
The following result of Joseph Liberman \cite{liberman} gives an example.

\begin{thm}{Theorem}
Any geodesic in the surface of a convex body is one-sided differentiable as a curve in $\EE^3$.
\end{thm}

\parit{Proof.}
Let $\gamma$ be a geodesic on the surface of a convex body $K$.
By \ref{ex:liberman} (which is Liberman's lemma),
the function $f_p\:t\z\mapsto \distfun_p\circ\gamma(t)$ is semiconcave for any $p\in K$.
In particular, one-sided derivatives $f_p^\pm(t)$ are defined for every $t$.

Given $x=\gamma(t)$, choose three points $p_1,p_2,p_3\in K$ in general position;
that is, the points $x$, $p_1$, $p_2$, and $p_3$ do not lie in one plane.
Observe that the distance functions $\distfun_{p_i}$ give smooth coordinates in a neighborhood of $x$.
By the above, the functions $f_{p_i}$ have one-sided derivatives at $t$.
Since the coordinates are smooth, we get that $\gamma^+(t)$ is defined as well. The same argument shows that $\gamma^-(t)$ is defined.
\qeds

{

\begin{wrapfigure}{r}{27 mm}
\vskip-0mm
\centering
\includegraphics{mppics/pic-1129}
\end{wrapfigure}

\begin{thm}{Exercise}\label{ex:convex}
Suppose a plane $\Pi$ cuts off a disk $\Delta$ from the surface of a convex body $K$, and the reflection of $\Delta$ across $\Pi$ lies in $K$.
Show that $\Delta$ is a convex subset of the surface;
that is, if a geodesic has endpoints in $\Delta$, then it lies completely in $\Delta$.
\end{thm}

}

The following exercise gives a sharper version of comparison for convex surfaces;
it is due to Anatolii Milka \cite[Theorem 2]{milka1982}.

\begin{thm}{Very advanced exercise}\label{ex:milka}
Let $\spc{P}$ be the surface of a nondegenerate convex body $K\subset\EE^3$,
and let $\gamma_1,\gamma_2\: [0,1]\to\spc{P}$ be geodesic paths that start at one point $p\z=\gamma_1(0)\z=\gamma_2(0)$.
Suppose $x_i=\gamma_i(1)$, and $y_i\z=p+\gamma_i^+(0)$.
Show that
\[\dist{x_1}{x_2}{\spc{P}}\le \dist{y_1}{y_2}{W},\]
where $W$ denotes the complement of the interior of $K$.

\end{thm}

\section{Remarks}

The excess of solid triangles can be extended to a sigma-additive measure on Borel subsets of an $\Alex\kappa$ surface.
This is the so-called \index{curvature measure}\emph{curvature measure} \cite[V]{alexandrov-1948}.
There have been several attempts to generalize this to higher dimensions, but so far, we only have partial results; see \cite{gigli2018, lebedeva-petrunin2024b} and the references therein.

Much of the theory described in this lecture has been generalized to surfaces with bounded integral curvature \cite{alexandrov-zalgaller,reshetnyak1993}.

The statement of Cauchy's theorem was conjectured by Adrien-Marie Legendre in the first edition of his geometry textbook \cite{legendre}.
It was motivated by a vague definition in Euclid's Elements, which could be interpreted as
\textit{polyhedrons are equal if their faces are pairwise equal}.
The local lemma was already known to Legendre.
Joseph-Louis Lagrange, a colleague of Legendre, suggested this problem to Augustin-Louis Cauchy, who soon solved it \cite{cauchy}.

The observation that the face-to-face condition can be removed was made by
Alexandr Alexandrov in 1941; in the same paper he proved the uniqueness theorem \cite{alexandrov-1941}.
A quite different proof was found by Yurii Volkov in his thesis \cite{volkov}; it uses a deformation of three-dimensional polyhedral space.
(Be aware that the proof of this theorem given in the book by Igor Pak contains an essential mistake \cite{petrunin-2023}.)

In Cauchy's proof \cite{cauchy}, the theorem was deduced from an analog of the following lemma.
The proof of the lemma contained a small mistake, which was corrected a century later \cite{sabitov}.
Several proofs of the arm lemma can be found in the letters between Isaac Schoenberg and Stanisław Zaremba \cite{schoenberg-zaremba}.

\begin{thm}{Arm lemma}\label{lem:arm}
Assume that $A=[a_0 a_1\dots a_n]$ is a convex polygon in $\EE^2$
and $A'=[a'_0 a'_1\dots a'_n]$ is a polygonal line in $\EE^3$
such that
$|a_i-a_{i+1}|=|a'_i-a'_{i+1}|$ for any $i\in\{0,\dots,n-1\}$
and
$\measuredangle a_i\le \measuredangle a'_i$
for each $i\in\{1,\dots,n-1\}$.
Then
\[|a_0-a_n|\le |a'_0-a'_n|\]
and equality holds if and only if $A$ is congruent to $A'$.
\end{thm}

The following variation of the arm lemma is due to Viktor Zalgaller \cite{zalgaller};
it makes sense for nonconvex spherical polygons and can be used instead of the standard arm lemma.

\begin{thm}{Another arm lemma}
Let $A=[a_1\dots a_n]$ and $A'\z=[a'_1\z\dots a'_n]$ be two spherical $n$-gons (not necessarily convex).
Assume that $A$ lies in a half-sphere,
the corresponding sides of $A$ and $A'$ are equal,
and each angle of $A$ is at least the corresponding angle in $A'$.
Then $A$ is congruent to~$A'$.
\end{thm}

Another close relative of the arm lemma is Reshetnyak's majorization theorem \cite{reshetnyak}.

Alexandrov gave two proofs of the global lemma \cite[2.1.2 and 2.1.3]{alexandrov}.
The first is combinatorial, and the second is more visual.
The argument in the second proof was reused by Anton Klyachko in his \index{car-crash lemma}\emph{car-crash lemma} \cite{klyachko}.

Proposition \ref{prop:conv-surf-CBB(0)} generalizes to the boundaries of convex bodies in $\EE^m$ for any $m\ge 2$.
It could be considered as a special case of the conjecture about the boundary of an Alexandrov space; see \ref{conj:bry}.
Another partial case, for convex sets in Riemannian manifolds, is proved by Stephanie Alexander and the two of us \cite{alexander-kapovitch-petrunin-2008}.

\begin{wrapfigure}{o}{30mm}
\vskip-0mm
\centering
\includegraphics{mppics/pic-15}
\vskip-0mm
\end{wrapfigure}

According to the uniqueness theorem, a convex polyhedron is completely determined by the intrinsic metric of its surface.
In particular, knowing the metric, we could find the position of the edges.
However, in practice, it is not easy to do.
For example, the surface glued from a rectangle, as shown in the picture, is the surface of a tetrahedron.
Some of the glued lines appear inside the faces of the tetrahedron, and some edges (dashed lines) do not follow the sides of the rectangle.

Let us also state the following result of Alexey Pogorelov \cite[chapter III]{pogorelov};
an alternative proof was found by Yurii Volkov \cite{volkov1968}.

\begin{thm}{Theorem}
Any two convex bodies in $\EE^3$ with isometric surfaces are congruent.

Moreover, any isometry between the surfaces of convex bodies can be extended to an isometry of the whole $\EE^3$.
\end{thm}

Notice that the statement is equivalent to the following:
\textit{Two polyhedra $K$ and $K'$ with close surfaces in the sense of Gromov--Hausdorff are almost congruent};
that is, there is a motion $\mu$ of $\EE^3$ such that the Hausdorff distance between $K$ and $\mu(K')$ is small.
At first glance, this theorem might look like a small improvement of Alexandrov's uniqueness,
but this improvement is huge, and the proof is quite hard.

%% file: sols.tex

\chapter{Semisolutions}

\parbf{\ref{ex:net}};
\ref{SHORT.ex:net:finite}.
Suppose $X$ is compact.
Then for any $\eps>0$, any cover of $X$ by open $\eps$-balls has a finite subcover.
Note that the centers of these balls form an $\eps$-net of $X$.

Now suppose $X$ has a finite $\eps$-net.
Show that any sequence $x_n$ of points in $X$ has a subsequence such that all of its points lie in one $\eps$-ball.
Apply this statement for $\eps=\tfrac1n$ together with a diagonal procedure.

\parit{\ref{SHORT.ex:net:compact}.}
Let $Z$ be a compact $\eps$-net of $X$.
By \ref{SHORT.ex:net:finite}, $Z$ has a finite $\eps$-net, say~$F$.
Note that $F$ is a $2\cdot\eps$-net of $X$.
Since $\eps>0$ is arbitrary, we get the result.

\parbf{\ref{ex:pack-net}.} If $x_1,\dots,x_n$ is not an $\eps$-net, then there is a point $y$ such that $\dist{x_i}{y}{}\ge\eps$ for any $i$.
Therefore $x_1,\dots,x_n$ is not a maximal packing --- a contradiction.

\parbf{\ref{ex:compact+connceted}.}
Choose a sequence of positive numbers $\varepsilon_n$ that converges to $0$ very fast.
Choose a finite $\varepsilon_n$-net $N_n$ of $K$ for each $n$.
We can assume that $\eps_0>\diam K$, and $N_0$ is a one-point set.
If $\dist{x}{y}{}<\eps_k$ for some $x\in N_{k+1}$ and $y\in N_{k}$, then connect them by a curve of length at most $\eps_k$.

Let $K'$ be the union of all these curves and $K$.
Show that $K'$ is compact and path-connected.

\parit{Source:} This problem is due to Eugene Bilokopytov \cite{bilokopytov}.

\parbf{\ref{ex:compact=>complete}.}
Choose a Cauchy sequence $x_n$ in $(\spc{X},\|*\z-*\|)$; it is sufficient to show that a subsequence of $x_n$ converges.

Observe that the sequence $x_n$ is Cauchy in $(\spc{X},|*-*|)$;
denote its limit by $x_\infty$.

Passing to a subsequence, we can assume that $\|x_{n-1}-x_{n}\|\z<\tfrac1{2^{n}}$ for each $n$.
It follows that there is a 1-Lipschitz path $\gamma$ in $(\spc{X},\|*-*\|)$ such that $x_n=\gamma(\tfrac1{2^n})$ for each $n$ and $x_\infty=\gamma(0)$.
Therefore,
\begin{align*}
\|x_\infty-x_n\|&\le \length\gamma|_{[0,\frac1{2^n}]}\le \tfrac1{2^n}.
\end{align*}
In particular, $x_n$ converges to $x_\infty$ in $(\spc{X},\|*\z-*\|)$.

\parit{Source:} \cite[Corollary]{hu-kirk}; see also \cite[Lemma 2.3]{petrunin-stadler}.

\parbf{\ref{ex:compact-length}.}
Given a pair of points $p$ and $q$, choose a sequence of paths $\gamma_n$ from $p$ to $q$ such that
\[\length\gamma_n\to \dist pq{}
\quad\text{as}\quad
n\to\infty;\]
these paths exist since we are in a length space.
We can assume that each $\gamma_n$ is parametrized proportionally to the arc length;
in particular, the paths $\gamma_n$ are equicontinuous.
Show that paths $\gamma_n$ lie in a closed ball, say $\cBall[p,r]$ of some radius $r<\infty$.
Since the space is proper, $\cBall[p,r]$ is compact.
By the Arzelà--Ascoli theorem, we can pass to a converging subsequence of $\gamma_n$.
Show that its limit is a geodesic path from $p$ to $q$.

\parbf{\ref{ex:menger}.}
Choose a sequence $\eps_n>0$ that converges to zero very fast, say such that $\sum_n10^n\cdot \eps_n$ is small.
Follow the argument in the proof of Menger's lemma, taking $\eps_n$-midpoints at the $n$-th stage.

\parbf{\ref{ex:k-><mono}.}
Let us write the Riemannian metric on $\MM^2(\kappa)$ in polar coordinates $(\theta,r)$;
it has the form 
$(\begin{smallmatrix}
h^2&0
\\
0&1
\end{smallmatrix})$, where $h=h(\kappa,r)\ge 0$.
Calculate $h(\kappa,r)$.
Show that for fixed $r$, the function $\kappa\mapsto h(\kappa,r)$ is nonincreasing in the domain of definition.
Suppose $\kappa<\Kappa$, and consider the partially defined map $\MM^2(\kappa)\to\MM^2(\Kappa)$ that sends a point to the point with the same polar coordinates.
Show that this map is short in its domain of definition.
Use it to prove the statement in the exercise.

\parbf{\ref{ex:angkK}.} Show and use that 
$\angk p{x}{y}_{\SSS^2}-\angk p{x}{y}_{\EE^2}=O(\dist[2]{p}{x}{}+\dist[2]{p}{y}{})$
and 
$\angk p{x}{y}_{\EE^2}-\angk p{x}{y}_{\HH^2}=O(\dist[2]{p}{x}{}+\dist[2]{p}{y}{})$.

\parbf{\ref{ex:undefined-angle}.}
Consider a hinge in the plane $\RR^2$ with a metric defined by a norm, say by the $\ell^\infty$-norm.

\parbf{\ref{ex:adjacent-angles}.}
Show that $\mangle\hinge pxy=\pi$ and apply \ref{claim:angle-3angle-inq}.

\parbf{\ref{ex:first-var}.}
Denote by $\alpha$ the arc-length parametrization of $[qp]$ from $q$ to $p$.
Choose $\eps>0$.
Observe that 
\[\dist[2]{\gamma(t)}{\alpha(\tfrac1\eps\cdot t)}{}\le t^2\cdot(1-\tfrac2\eps\cdot\cos\phi+\tfrac1{\eps^2})+o(t^2).\]
By the triangle inequality
\[\dist{p}{\gamma(t)}{}\le \dist{\gamma(t)}{\alpha(\tfrac1\eps\cdot t)}{}+\dist{q}{p}{}-\tfrac1\eps\cdot t.\]
Conclude that
\[\dist{p}{\gamma(t)}{}
\le\dist{q}{p}{}-t\cdot \cos \phi+\delta(\eps)\cdot t+o(t),\]
where $\delta(\eps)\to 0$ as $\eps\to0$.
The statement follows since $\eps>0$ is arbitrary.

\parbf{\ref{ex:generalized-selection}.}
Since the space is proper, it is separable; 
that is, we can choose a countable everywhere dense set $\{x_1,x_2,\dots\}$.

Let $A_1,A_2,\dots$ be a sequence of closed sets.
Applying a diagonal procedure, we can pass to a subsequence such that for each $i$ the sequence $\distfun_{A_n}x_i$ converges as $n\to\infty$; denote its limit by $f(x_i)$.

Suppose $f(x_i)<\infty$ for some $i$.
Since $\distfun_{A_n}$ is $1$-Lipschitz for any $n$, we have 
\[|f(x_i)-f(x_j)|\le \dist{x_i}{x_j}{}\]
for all $i$ and $j$.
Hence $f(x_i)<\infty$ for any $i$ and the function $f$ can be extended to a continuous function defined on the ambient space.
Show that $A_\infty=f^{-1}\{0\}$ is the limit of $A_n$ in the sense of Hausdorff.

If $f(x_i)=\infty$ for some $i$, then the same holds for any $i$.
Show that in this case $A_n\to\emptyset$ in the sense of Hausdorff.

\parbf{\ref{ex:two>one}.}
Set $\dist{x_i}{x_j}{}=2$ for $0\ne i\ne j\ne 0$.

\parbf{\ref{ex:Haus-conv}.}
Apply the definition of Hausdorff distance (\ref{def:hausdorff-convergence}).
Without the assumption that all $A_i$ are contained in a compact subset of $\spc{X}$, the statement is false, as shown by the example $\spc{X}=\RR$, $A_n=\{0,n\}$, and $A_\infty=\{0\}$.

\parbf{\ref{ex:non-contracting-map}.}
Given a pair of points $x_0,y_0\in \spc{K}$, 
consider two sequences $x_0,x_1,\dots$ and $y_0,y_1,\dots$
such that $x_{n+1}=f(x_n)$ and $y_{n+1}\z=f(y_n)$ for each $n$.

Since $\spc{K}$ is compact, 
we can choose an increasing sequence of integers $n_i$
such that both sequences $(x_{n_i})_{i=1}^\infty$ and $(y_{n_i})_{i=1}^\infty$
converge.
In particular, both are Cauchy;
that is,
\[
|x_{n_i}-x_{n_j}|_{\spc{K}}\to 0 
\quad\text{and}\quad
|y_{n_i}-y_{n_j}|_{\spc{K}}\to 0
\]
as $\min\{i,j\}\to\infty$.

Since $f$ is distance-noncontracting, 
\[
|x_0-x_{|n_i-n_j|}|
\le 
|x_{n_i}-x_{n_j}|
\]
for any $i$ and $j$.
Therefore, there is a sequence $m_i\to\infty$ such that
\[
x_{m_i}\to x_0\quad\text{and}\quad y_{m_i}\to y_0
\leqno({*})\]
as $i\to\infty$.

Since $f$ is distance-noncontracting, the sequence $\ell_n=|x_n-y_n|_{\spc{K}}$ is nondecreasing.
By $({*})$, $\ell_{m_i}\to\ell_0$ as $m_i\to\infty$.
It follows that 
\[\ell_0=\ell_1=\dots\]
In particular, 
\[|x_0-y_0|_{\spc{K}}=\ell_0=\ell_1=|f(x_0)-f(y_0)|_{\spc{K}}\]
for any pair of points $(x_0,y_0)$ in $\spc{K}$.
That is, the map $f$ is distance-preserving and hence injective.
From $({*})$, we also get that $f(\spc{K})$ is everywhere dense.
Since $\spc{K}$ is compact, $f$ is surjective --- hence the result.

\parit{Remarks.}
This is a basic lemma in the introduction to Gromov--Hausdorff distance \cite[see 7.3.30 in][]{burago-burago-ivanov}.
The presented proof was suggested by Travis Morrison.

\parbf{\ref{ex:non-expanding-map}.}
Apply \ref{ex:non-contracting-map} to a right inverse of the map.

\parbf{\ref{ex:GH-po}.}
The only-if part is trivial.
Let us prove the if part.

Let $f_n\:\spc{X}_n\to \spc{X}_\infty$ and $h_n\:\spc{X}_\infty\to \spc{X}_n$ be the maps in the definition of the inequalities $\spc{X}_n\le \spc{X}_\infty+\eps_n$ and $\spc{X}_\infty\le \spc{X}_n+\eps_n$, respectively.
Apply \ref{ex:non-contracting-map} to show that any partial limit of $f_n\circ h_n$ is an isometry of $\spc{X}_\infty$.
Conclude that $f_n$ is an $\eps_n'$-isometry for some converging-to-zero sequence $\eps_n'$ and apply \ref{lem:almost-isom}.

\parbf{\ref{ex:GH-noncompact}}; \ref{SHORT.ex:GH-noncompact:proper}
 Consider the graphs of the following functions with the induced metric from $\RR^2$.
\[
x\mapsto \cos x+\cos \tfrac x\pi
\quad\text{and}\quad
x\mapsto \cos x+\sin \tfrac x\pi.
\]

\parit{\ref{SHORT.ex:GH-noncompact:bounded}}
For every rational number $q\in[1,2]$ consider an interval of length~$q$.
Let $\spc{X}$ be obtained by identifying all endpoints of the intervals.

Let $\spc{Y}$ be constructed in the same way but skipping the interval of length $1.5$.

\parbf{\ref{ex:CBB+-}.} Note that we can assume that we have equality in all inequalities except one.
In other words, we can assume that equality holds in all inequalities except $\dist{y_1}{y_2}{}\le\dist{x_1}{x_2}{}$ or $\dist{q}{y_2}{}\ge\dist{p}{x_2}{}$.

In the first case,
\[
\angk q{y_1}{y_2}\le \angk p{x_1}{x_2},
\quad
\angk q{y_2}{y_3}= \angk p{x_2}{x_3},
\quad
\text{and}
\quad
\angk q{y_3}{y_1}= \angk p{x_3}{x_1}.
\]
Therefore, the comparison for $p,x_1,x_2,x_3$ implies the comparison for $q,y_1,y_2,y_3$.

In the second case, let us argue by contradiction.
Assume that the comparison does not hold for $q,y_1,y_2,y_3$.
Show that
\[\angk q{y_1}{y_2}+\angk q{y_2}{y_3}> \pi.\]
Then show and use that
\[\angk p{x_1}{x_2}+\angk p{x_2}{x_3}\ge\angk q{y_1}{y_2}+\angk q{y_2}{y_3}.\]

\parbf{\ref{ex:Euclid-is-CBB}.}
The 4-point comparison (\ref{def:CBB}) reduces our question to the following.
\textit{Any spherical triangle has perimeter at most $2\cdot\pi$.}
Choose a spherical triangle $[xyz]$.
Let $x'$ be the antipode of $x$; that is, $x'=-x$.
The spherical triangle inequality (\ref{claim:angle-3angle-inq}) implies that
\[\dist{y}{z}{\mathbb{S}^2}\le \dist{y}{x'}{\mathbb{S}^2}+\dist{x'}{z}{\mathbb{S}^2}.\]
Observe that 
\[
\dist{x}{y}{\mathbb{S}^2}+\dist{y}{x'}{\mathbb{S}^2}=\pi,
\quad\text{and}\quad
\dist{x}{z}{\mathbb{S}^2}+\dist{z}{x'}{\mathbb{S}^2}=\pi.
\]
Hence
\[\dist{x}{y}{\mathbb{S}^2}+\dist{x}{z}{\mathbb{S}^2}+\dist{y}{z}{\mathbb{S}^2}\le2\cdot \pi.\]

\parbf{\ref{ex:(3+1)-expanding}.}
For the only-if part consider the following two cases.

If $\angk p{x_1}{x_2}+\angk p{x_2}{x_3}\ge \pi$, then choose two model triangles $[qy_1y_2]\z=\modtrig(px_1x_2)$ and $[qy_2y_3]=\modtrig(px_2x_3)$ that lie on the opposite sides of $[qy_2]$.
By the comparison, $\dist{y_1}{y_3}{}\ge \dist{x_1}{x_3}{}$.
Therefore the obtained configuration meets all the conditions.

If $\angk p{x_1}{x_2}+\angk p{x_2}{x_3}\le \pi$, then choose a model triangle $[qy_1y_2]\z=\modtrig(px_1x_2)$
and take $y_3$ on the extension of $[y_1q]$ behind $q$ such that $\dist{q}{y_3}{}=\dist{p}{x_3}{}$.
Then $\mangle \hinge q{y_2}{y_3}\ge \angk p{x_2}{x_3}$, therefore $\dist{y_2}{y_3}{}\ge \dist{x_2}{x_3}{}$.
Further, $\dist{y_1}{y_3}{}=\dist{x_1}{p}{}+\dist{p}{x_3}{} \ge \dist{x_1}{x_3}{}$,
and again, the obtained configuration meets all the conditions.

To prove the if part, choose a configuration $q,y_1,y_2,y_3$ that meets all the conditions and maximize the sum
\[\dist{y_1}{y_2}{}+\dist{y_2}{y_3}{}+\dist{y_3}{y_1}{}.\]
Show that $\dist{q}{y_i}{}=\dist{p}{x_i}{}$ for each $i$ and $q$ lies in the solid triangle $y_1y_2y_3$;
in particular 
\[\mangle \hinge q{y_1}{y_2}+\mangle \hinge q{y_2}{y_3}+ \mangle \hinge q{y_3}{y_1}=2\cdot\pi.\]
Applying \ref{angle-monotonicity}, we get
\[\angk p{x_1}{x_2}+\angk p{x_2}{x_3}+\angk p{x_3}{x_1}
\le 
2\cdot\pi.
\]

\parbf{\ref{ex:alex-lemma-cat}.}
Consider model triangles $[\tilde p\tilde x\tilde z]=\modtrig(pxz)$ and $[\tilde p\tilde y\tilde z]=\modtrig(pyz)$
that share side $[\tilde p\tilde z]$ and lie on its opposite sides.
Note that 
\begin{align*}
\dist{\tilde x}{\tilde y}{\EE^2}
&\le \dist{\tilde x}{\tilde z}{\EE^2}+\dist{\tilde z}{\tilde y}{\EE^2}=
\\
&=\dist{x}{z}{\spc{X}}+\dist{z}{y}{\spc{X}}=
\\
&=\dist{x}{y}{\spc{X}},
\end{align*}
where $\spc{X}$ is our metric space.
It remains to apply \ref{angle-monotonicity}.

\parbf{\ref{ex:noncreasing}.}
Apply \ref{clm:angle-mono}.

\parbf{\ref{ex:0-angle}.}
Without loss of generality, we can assume that $\dist{p}{x}{}\le \dist{p}{y}{}$.
Choose $\bar x\in [px]$;
let $\bar y\in [py]$ be such that $\dist{p}{\bar x}{}=\dist{p}{\bar y}{}$.
Apply \ref{clm:angle-mono} to show that $\bar x=\bar y$.
Conclude that $[px]\subset [py]$.

\parbf{\ref{ex:pi-angle}.}
Assume that there are two distinct geodesics from $z$ to $x$.
Then we can choose distinct points $p$ and $q$ on these geodesics such that $\dist{z}{p}{}=\dist{z}{q}{}$.
Observe that $\angk zpq>0$.
By comparison,
\[\angk zpq+\angk zpy+\angk zqy\le 2\cdot\pi.\]
Therefore, one of the angles $\angk zpy$ or $\angk zqy$ is strictly less than $\pi$.
The latter contradicts the triangle inequality.

Alternatively, if there are two distinct geodesics from $z$ to $x$ then their concatenations with $[yz]$ give bifurcating geodesics from $y$ to $x$. This contradicts \ref{ex:0-angle}.

\parbf{\ref{ex:adjacent-CBB}.}
By \ref{ex:adjacent-angles}, we have
\[\mangle\hinge pxz+\mangle\hinge pyz\ge \pi.\]
Since $p\in \left]xy\right[$ we have
\[\angk p{\bar x}{\bar y}=\pi\]
for any $\bar x\in \left]px\right[$ and $\bar y\in \left]py\right]$.
By comparison, we have that 
\[\angk p{\bar x}{\bar z}+\angk p{\bar z}{\bar y}\le\pi\]
for any $\bar z\in \left]pz\right]$.
Passing to the limit as
$\dist{p}{\bar x}{}\to 0$,
$\dist{p}{\bar y}{}\to 0$, and
$\dist{p}{\bar z}{}\to 0$,
we get \[\mangle\hinge pxz+\mangle\hinge pyz\le \pi.\]

\parbf{\ref{ex:pxyvw}.} 
Without loss of generality, we can assume that $x$, $v$, $w$, and $y$ appear on 
$[xy]$ in this order.
By \ref{clm:angle-mono},
\[
\angk xyp\ge \angk xwp \ge\angk xvp.
\]
Hence, $\Rightarrow$ follows.

By Alexandrov's lemma,
\begin{align*}
\angk xyp=\angk xvp
\quad&\Leftrightarrow\quad
\angk yxp=\angk yvp,
\\
\angk xyp=\angk xwp
\quad&\Leftrightarrow\quad
\angk yxp=\angk ywp.
\end{align*}
Whence, $\Leftarrow$ follows.

\parbf{\ref{ex:angle-lim}.} Suppose $\mangle \hinge {x_\infty}{y_\infty}{z_\infty}>\alpha$.
Then we can choose $\bar y_\infty\in\left]x_\infty y_\infty\right]$
and $\bar z_\infty\in\left]x_\infty z_\infty\right]$ such that 
$\angk{x_\infty}{\bar y_\infty}{\bar z_\infty}>\alpha$.
Now choose $\bar y_n\in\left]x_n y_n\right]$ and $\bar z_n\in\left]x_n z_n\right]$ such that $\bar y_n\to \bar y_\infty$ and $\bar z_n\to \bar z_\infty$.
Observe that 
\[\liminf_{n\to\infty}\mangle \hinge {x_n}{y_n}{z_n}\ge\liminf_{n\to\infty}\angk{x_n}{\bar y_n}{\bar z_n} = \angk{x_\infty}{\bar y_\infty}{\bar z_\infty},\]
hence the result.

\parbf{\ref{ex:urysohn}.}
The Urysohn space provides an example;
see, for example, \cite[Lecture 2]{petrunin2023pure}.

\parbf{\ref{ex:normCBB}.}
Choose a triangle $[0vw]$.
Note that $m=\tfrac12(v+w)$ is the midpoint of $[vw]$.

Use comparison to show that
\[2\cdot |\tfrac12(v+w)|^2+2\cdot |\tfrac12(v-w)|^2\ge |v|^2+|w|^2.\]

This inequality implies the opposite one;
it follows if we rewrite it via $x=\tfrac12(v+w)$ and $y=\tfrac12(v-w)$.
Hence we have 
\[2\cdot |\tfrac12(v+w)|^2+2\cdot |\tfrac12(v-w)|^2= |v|^2+|w|^2\]
for any $v,w$.
That is, the norm meets the parallelogram identity.
It is well known that any such norm is quadratic, and the statement follows.

\parbf{\ref{ex:concave'}.}
Choose a $\lambda$-concave function $f\:\II\to\RR$.
Note that
\[h\:t\mapsto f(t)-\lambda\cdot\tfrac{t^2}2\]
is concave,
and it is sufficient to show that $h$ is one-sided differentiable.

Choose $t_0\in \II$.
Show and use that the function
$t\mapsto \tfrac {h(t)-h(t_0)}{t-t_0}$ is nonincreasing.

\parbf{\ref{ex:concave-open}.}
We can assume that $\lambda=0$; otherwise pass to the function $t\mapsto f(t)-\lambda\cdot\tfrac{t^2}2$.

Suppose that $f$ is defined at $t_0$.
Since $f$ is defined on an open interval, we can choose $a$ and $b$ such that $a<t_0<b$ and $f$ is defined on $[a,b]$.
Show and use that
\begin{align*}
\frac{f(b)-f(t_0)}{b-t_0}&\le \frac{f(t)-f(t_0)}{t-t_0}\le \frac{f(a)-f(t_0)}{a-t_0}
\end{align*}
for any $t\in [a,b]\setminus\{t_0\}$.

\parbf{\ref{ex:distfun-semiconcave}.}
From \ref{comp-kappa}, this inequality follows in the sense of distributions, and hence in any other sense.

\parbf{\ref{ex:alm-min}.}
Suppose such a point does not exist;
that is, for any $p\in \spc{X}$ there is a point $p'$ such that $r(p')\le (1-\eps)\cdot r(p)$ and $\dist p{p'}{}\le\tfrac{1}{\eps}\cdot r(p)$.
Construct a sequence of points $p_0,p_1,\dots$ such that $p_n=p_{n-1}'$ for any~$n$.
Show that this sequence is Cauchy; denote its limit by $p_\infty$.
Arrive at a contradiction by showing that $r(p_\infty)\le0$.

\parit{Comment.}
A slightly more general version of this statement was proved by Ivar Ekeland \cite{ekeland}.

\parbf{\ref{ex:RisCBB(1)}.}
Suppose $\mangle\hinge mxp\ne 0$ and $\mangle\hinge mxp\ne\pi$;
equivalently, $\mangle\hinge mxp\ne 0$ and $\mangle\hinge mxq\ne0$.

We can assume that $\dist pq{}$ only slightly exceeds $\pi$,
so $\dist pm{}<\pi$ and $\dist qm{}<\pi$.
We can also assume that $\dist xm{}<\pi$.
Use the comparison to show that 
\[\dist px{}+\dist qx{} < \dist pq{},\]
and arrive at a contradiction with the triangle inequality.

Extend $[pq]$ to a maximal local geodesic $\gamma$.
Argue as above to show that any point lies on $\gamma$.
If $\gamma$ is closed, then the space is isometric to a circle;
otherwise, it is isometric to a line segment.

\parbf{\ref{ex:perim-k>0}.}
Arguing by contradiction, suppose 
\[\dist{p}{q}{}+\dist{q}{r}{}+\dist{r}{p}{}> 2\cdot\pi\eqlbl{eq:perimeter-of-triange<2pi}\] 
for $p,q,r\in \spc{A}$. 
Rescaling the space slightly by a constant $<1$, we can assume that $\diam\spc{A}<\pi$,
but the inequality \ref{eq:perimeter-of-triange<2pi} still holds.
By \ref{clm:K>k},
after rescaling $\spc{A}$ is still $\Alex1$.

\begin{wrapfigure}{o}{26mm}
\vskip-3mm
\centering
\includegraphics{mppics/pic-380}
\vskip-0mm
\end{wrapfigure}

Take $z_0\in [q r]$ at maximal distance from $p$.
Consider the following model configuration:
two geodesics $[\tilde p\tilde z_0]$, $[\tilde q\tilde r]$ in $\mathbb{S}^2$ such that 
\begin{align*}
\dist{\tilde p}{\tilde z_0}{}&=\dist{p}{z_0}{},
&
\dist{\tilde q}{\tilde r}{}&=\dist{q}{r}{},
\\ 
\dist{\tilde z_0}{\tilde q}{}&=\dist{z_0}{q}{},
&
\dist{\tilde z_0}{\tilde r}{}&=\dist{z_0}{r}{},
\end{align*}
and 
\[\mangle\hinge{\tilde z_0}{\tilde q}{\tilde p}
=\mangle\hinge{\tilde z_0}{\tilde r}{\tilde p}
=\tfrac\pi2.\]

Choose $\tilde z\in [\tilde q\tilde r]$,
and let $z\in [q r]$ be the corresponding point.
By comparison, $\dist pz{}\le\dist {\tilde p}{\tilde z}{}$ if $z$ lies near $z_0$.
Moreover, this inequality holds as long as 
\[\dist{\tilde p}{\tilde z_0}{}+\dist{\tilde z_0}{\tilde z}{}+\dist{\tilde p}{\tilde z}{}<2\cdot\pi.\]
But this inequality holds for all $\tilde z\in [\tilde q\tilde r]$ since $\dist{\tilde p}{\tilde z_0}{}<\pi$, $\dist{\tilde z_0}{\tilde q}{}<\pi$, and $\dist{\tilde z_0}{\tilde r}{}<\pi$.
Hence $\dist pq{}\le\dist {\tilde p}{\tilde q}{}$ and $\dist pr{}\le\dist {\tilde p}{\tilde r}{}$.
Since $\dist{q}{r}{}=\dist{\tilde q}{\tilde r}{}$, the latter contradicts \ref{eq:perimeter-of-triange<2pi}.

\parbf{\ref{ex:dir-compact}.}
Suppose $\dir p{x_n}\not\to\dir p{x_\infty}$.
Since $\Sigma_p$ is compact, we may pass to a converging subsequence of $\dir p{x_n}$;
denote its limit by $\xi$.
We may assume that $\mangle (\dir p{x_\infty},\xi)>0$.

Denote by $\gamma_n$ and $\gamma_\infty$ the arc-length parametrizations of $[px_n]$ and $[px_\infty]$ from $p$.
For a geodesic $\alpha$ that starts from $p$ and goes in a direction sufficiently close to $\xi$, we have
\[\dist{\alpha(t)}{\gamma_n(t)}{}<\eps\cdot t
\quad\text{and}\quad
\dist{\alpha(t)}{\gamma_\infty(t)}{}>a\cdot t
\]
for some $\eps \ll a$, all large $n$, and all sufficiently small $t$.
These two inequalities imply 
that 
\[\dist{\gamma_n(t)}{\gamma_\infty(t)}{}>\tfrac a2\cdot t\]
for all small $t$ and all large $n$.
On the other hand, by assumption, $\dist{\gamma_n(t)}{\gamma_\infty(t)}{}\to0$ as $n\to\infty$ --- a contradiction.

\parit{Comments.}
The compactness of $\Sigma_p$ is necessary.
A counterexample with noncompact $\Sigma_p$ can be built using an iterated warped product of line segments and applying \cite[Theorem 1.2]{alexander-bishop2004}.
The space $\spc{A}$ can be assumed to be compact.

\parbf{\ref{ex:geodesic-cone}.}
Any point of $\Cone \spc{X}$ can be connected to the origin by a geodesic.
Given a nonzero element $v\in\Cone \spc{X}$, denote by $v'$ its projection in~$\spc{X}$; so, $v=|v|\cdot v'$.

Suppose $\spc{X}$ is $\pi$-geodesic.
Choose two nonzero elements $v,w\z\in\Cone \spc{X}$; let $\alpha=\dist{v'}{w'}{\spc{X}}$.
If $\alpha\ge \pi$, then $\mangle(v,w)=\pi$, and the product of geodesics $[v0]$ and $[0w]$ is a geodesic $[vw]$.
If $\alpha<\pi$, then $\alpha=\mangle(v,w)$, and there is a geodesic $\gamma\:[0,\alpha]\to \spc{X}$ from $v'$ to $w'$.
Consider the hinge $\hinge {\tilde o}{\tilde v}{\tilde w}$ in the plane
such that $\mangle\hinge {\tilde o}{\tilde v}{\tilde w}=\alpha$, $\dist{\tilde o}{\tilde v}{}=|v|$, and $\dist{\tilde o}{\tilde w}{}=|w|$.
Let $t\mapsto (\phi(t),r(t))$ be the geodesic $[\tilde v\tilde w]$ written in polar coordinates with origin at $\tilde o$, so that $\phi(0)=0$.
Show that $t\mapsto r(t)\cdot\gamma\circ\phi(t)$ is a geodesic from $v$ to $w$;
here we identify $\spc{X}$ with the unit sphere in $\Cone \spc{X}$.

To prove the converse, reverse the steps in the argument.

\parbf{\ref{ex:GHto-tangent}.}
Let $\spc{A}_n=\lambda_n\cdot\spc{A}$.
Note that for any $n$ the space $\Sigma_p\spc{A}$ is identical to $\Sigma_{\iota_n(p)} \spc{A}_n$.
In particular, we can canonically and isometrically identify $\T_p\spc{A}$ with $\T_{\iota_n(p)}\spc{A}_n$.
So for any geodesic $\gamma$ that starts at $p$, the vector $\gamma^+(0)$ corresponds to $\frac{1}{\lambda_n}\cdot(\iota_n\circ\gamma)^+(0)$.

Consider the logarithm maps $f_n=\log_{\iota_n(p)}\:\spc{A}_n\to \T_p\spc{A}$.
We claim that this sequence of maps satisfies the assumptions in \ref{lem:almost-isom-pointed}.

The condition in \ref{SHORT.lem:almost-isom-pointed-basepoint} is evident, and it is sufficient to check the conditions in \ref{SHORT.lem:almost-isom-pointed-b} and \ref{SHORT.lem:almost-isom-pointed-c} only for $R=1$.

Choose $\eps>0$.
By compactness of $\Sigma_p$ we can find a finite $\eps$-net $\xi_1,\dots,\xi_N$ in $\Sigma_p$. Moreover, without loss of generality we can assume that these directions are geodesic;
that is, there exist geodesics $\gamma_1,\ldots, \gamma_N$ starting at $p$ such that $\xi_i=\gamma_i^+(0)$ for each $i$.

Choose $T>0$ such that all $\gamma_i$ are defined on $[0,T]$.
Show that for any $\lambda_n>\frac{1}{T}$ the image under $f_n$ of the union $\bigcup_i\gamma_i([0,T])$ is an $\eps$-net in $\oBall(0,1)_{T_p}$.
This proves \ref{SHORT.lem:almost-isom-pointed-c}.

By comparison, we have that
\[\dist{\xi_i}{\xi_j}{\Sigma_p}\ge \angk p{\gamma_i(t_i)}{\gamma_j(t_j)}\]
for all $i\ne j$ and any $t_i,t_j\in (0,T]$.
By the definition of angle, we can assume that $T$ has been chosen so that in addition
\[\dist{\xi_i}{\xi_j}{\Sigma_p}\le \angk p{\gamma_i(t)}{\gamma_j(t)}+\eps\]
for all $i\ne j$ and any $t\in (0,T]$.

By construction of the map $f_n$ this implies that 
\[\left|\dist{x}{x'}{\spc{A}_n}-\dist{f_n(x)}{f_n(x')}{\T_p}\right|<\eps\]
for all $\lambda_n>\frac{1}{T}$ and all points $x,x'$ in $\bigcup_i\gamma_i([0,\frac{1}{\lambda_n}])\subset \oBall(p,1)_{\spc{A}_n}$.

Now hinge comparison and the triangle inequality imply that the same holds for arbitrary points $x,x'$ in $\oBall(p,1)_{\spc{A}_n}$ with $\eps$ replaced by $3\cdot\eps$.
This verifies \ref{SHORT.lem:almost-isom-pointed-b}.

\parbf{\ref{ex:df(xi)}.}
Since angles are defined, it follows that 
\[\dist{\gamma_1(t)}{\gamma_2(t)}{}\le \theta\cdot t\]
for all small $t>0$.
Since $f$ is $L$-Lipschitz, we get 
\[|f(\gamma_1(t))-f(\gamma_2(t))|\le L\cdot \theta\cdot t,\]
and the statement follows.

\parbf{\ref{ex:d(distfun)}}; \ref{SHORT.ex:d(distfun):<}
Show that we can assume there is a geodesic in the direction of $v$, and apply \ref{ex:first-var}.

\parit{\ref{SHORT.ex:d(distfun):=}.}
By \ref{SHORT.ex:d(distfun):<}, $(\dd_p\distfun_q)(v)\le-\max_{\xi\in\Uparrow_p^q}\langle\xi,v\rangle$.
Suppose this inequality is strict for some $v$.
We can assume that $|v|=1$ and there is a geodesic, say $\gamma$ in the direction of $v$.
Let $(\dd_p\distfun_q)(v)=-\cos\alpha_0$ for some $\alpha_0\in [0,\pi]$.
Note that any geodesic from $p$ to $q$ makes angle bigger than $\alpha_0$ with $\gamma$.

The function $f=\distfun_q\circ\gamma$ is Lipschitz.
By Rademacher's theorem it is differentiable almost everywhere;
moreover, 
\[f(t)-f(0)=\int_0^t f'(s)\cdot ds.\]
Suppose $f'(t)$ is defined.
Use \ref{SHORT.ex:d(distfun):<} to show that 
$f'(t)=-\cos\alpha(t)$, where $\alpha(t)$ is the angle between $\gamma$ and any geodesic from $\gamma(t)$ to $q$.
We can choose a sequence $t_n\to 0$ such that 
\[\lim_{n\to\infty}\alpha(t_n) \le \alpha_0.\]
Consider a sequence of geodesics $[\gamma(t_n)\, q]$.
Since the space is proper, we can pass to a convergent subsequence.
Its limit is a geodesic from $p$ to $q$; denote it by $[pq]$.

Use \ref{ex:angle-lim} to show that $[pq]$ makes an angle at most $\alpha_0$ with $\gamma$ --- a contradiction.
 
\parbf{\ref{ex:monotonicity}.}
Let $\gamma\:[0,\ell]\to \spc{A}$ be the geodesic $[xy]$ parametrized from $x$ to $y$,
and let $\phi=f\circ\gamma$.
Observe that 
\[\phi'(0)=\dd_xf(\dir xy)\le \<\dir{x}{y},\nabla_{x}f\>.\]
The same way we get $-\phi'(\ell)\le \<\dir{y}{x},\nabla_{y}f\>$.
Since $f$ is $\lambda$-concave, we have
\begin{align*}
f(y)&\le f(x)+\phi'(0)\cdot \ell+\tfrac\lambda2\cdot\ell^2,
\\
f(x)&\le f(y)-\phi'(\ell)\cdot \ell+\tfrac\lambda2\cdot\ell^2.
\end{align*}
It remains to add these inequalities.

\parbf{\ref{ex:d(distfun):==}.}
If the space is proper, then the statement follows from \ref{ex:d(distfun):=} and \ref{ex:pi-angle}.

To do the general case, let us argue by contradiction.
By the assumption, we can choose a point $z$ on the extension of $[pq]$ beyond~$q$.
We can assume that $|v|=1$ and it is a direction of a geodesic, say $[px]$.

Show that there is a sequence $x_n\in \left]px\right]$ such that $\dist{p}{x_n}{}\to0$ and
$\mangle \hinge q{x_n}p>\eps$ for each $n$ and some fixed $\eps>0$.
Observe that $\mangle\hinge q{x_n}z\z<\pi-\eps$; therefore
\[\dist{z}{x_n}{}<\dist{x_n}{q}{}+\dist{q}{z}{}-\delta\]
for each $n$ and some fixed $\delta>0$.
Pass to the limit as $x_n\to p$ and arrive at a contradiction to the triangle inequality.

\parbf{\ref{ex:convergence-grad}.}
Note that
$|(\dd_p f)(v)-(\dd_p g)(v)|\le s\cdot|v|$
for any $v\in \T_p$.
From the definition of gradient (\ref{def:grad}) we have:
\begin{align*}
&(\dd_p f)(\nabla_p g)\le\<\nabla_p f,\nabla_p g\>,
&&(\dd_p g)(\nabla_p f)\le\<\nabla_p f,\nabla_p g\>,
\\
&(\dd_p f)(\nabla_p f)=\<\nabla_p f,\nabla_p f\>,
&&(\dd_p g)(\nabla_p g)=\<\nabla_p g,\nabla_p g\>.
\end{align*}
Therefore,
\begin{align*}
&\dist[2]{\nabla_pf}{\nabla_pg}{}
=\<\nabla_p f,\nabla_p f\>+\<\nabla_p g,\nabla_p g\>-2\cdot\<\nabla_p f,\nabla_p g\>
\le
\\
&\qquad\le (\dd_p f)(\nabla_p f)+(\dd_p g)(\nabla_p g)-
(\dd_p f)(\nabla_p g)-(\dd_p g)(\nabla_p f)
\le
\\
&\qquad\le s\cdot(|\nabla_p f|+|\nabla_p g|).
\end{align*}

\parbf{\ref{ex:semicontinuous-grad}.}
Suppose $|\nabla_xf|> s$.
Then we can choose a geodesic $\gamma$ that starts at $x$ such that 
$(f\circ\gamma)^+(0)>s$.
In particular, there is $\eps>0$ such that
\[f\circ\gamma(t)-f(x)>(s+\eps)\cdot t+o(t),\]
and the only-if part follows.

Now suppose $f(y)-f(x)>s\cdot \ell+\lambda\cdot \tfrac{\ell^2}2$,
where $\ell=\dist{x}{y}{}$.
Let $\gamma\:[0,\ell]\to \spc{A}$ be a geodesic from $x$ to $y$.
Since $f\circ\gamma$ is $\lambda$-concave, we have
\[f\circ\gamma(\ell)\le f\circ\gamma(0)+(f\circ\gamma)^+(0)\cdot\ell+\lambda\cdot \tfrac{\ell^2}2.\]
It follows that 
\[(\dd_xf)(\dir xy)=(f\circ\gamma)^+(0)>s,\]
and by \ref{prop:grad-exist}, $|\nabla_x f|>s$.

This proves \ref{SHORT.ex:semicontinuous-grad:>s}.
To prove \ref{SHORT.ex:semicontinuous-grad:lim}, argue by contradiction and apply \ref{SHORT.ex:semicontinuous-grad:>s}.

\parbf{\ref{ex:elf-contracting}.}
Note that $f\circ\alpha$ is a nondecreasing function.
Apply \ref{ex:d(distfun):<} and the definition of gradient to show that
\begin{align*}
-(\dd_{\alpha(t)}\distfun_{\alpha(t_3)})(\nabla_{\alpha(t)}f)
&\ge
\langle \nabla_{\alpha(t)}f,\dir{\alpha(t)}{\alpha(t_3)}\rangle
\ge
\\
&\ge
\dd_{\alpha(t)}f(\dir{\alpha(t)}{\alpha(t_3)})\ge
\\&\ge0
\end{align*}
for any $t<t_3$.
Conclude that the function $t\mapsto \distfun_{\alpha(t_3)}\circ\alpha(t)$ is nonincreasing for $t\le t_3$.

\parbf{\ref{ex:mayer}.}
Suppose $s>s_0$.
Then
\begin{align*}
(f\circ\hat\alpha)^+(s_0)&=|\nabla_{\hat\alpha(s_0)}f|
\ge
\\
&\ge
(\dd_{\hat\alpha(s_0)}f)(\dir{\hat\alpha(s_0)}{\hat\alpha(s)})
\ge
\\
&\ge
\frac{f\circ\hat\alpha(s)-f\circ\hat\alpha(s_0)}{\dist{\hat\alpha(s)}{\hat\alpha(s_0)}{}}.
\end{align*} 
Since $s-s_0\ge\dist{\hat\alpha(s)}{\hat\alpha(s_0)}{}$,
\[(f\circ\hat\alpha)^+(s_0)\ge
\frac{f\circ\hat\alpha(s)-f\circ\hat\alpha(s_0)}{s-s_0},\]
which implies the statement.

\parbf{\ref{lem:fg-dist-est}.}
Fix $t$, and let $p=\alpha(t)$ and $q=\beta(t)$.
Apply \ref{eq:fist-var-inq+} to get
\begin{align*}
\ell^+&\le -\<\dir{p}{q},\nabla_{p}f\>
-\<\dir{q}{p},\nabla_{q}g\>
\le
\\&\le -{\left({f(q)}-{f(p)}-\lambda\cdot\tfrac{\ell^2}2\right)}/{\ell}
-{\left({g(p)}-{g(q)}-\lambda\cdot\tfrac{\ell^2}2\right)}/{\ell}\le
\\&\le \lambda\cdot\ell+\tfrac{2\cdot\eps}{\ell}.
\end{align*}
Integrating this inequality, we get the second statement.

\parbf{\ref{ex:short-onto}.}
Choose a point $p\in M$.
Observe that $\gexp_p^1$ provides a short map from the unit hemisphere $\SSS_+^m$ onto $\cBall[p,\tfrac\pi2]_{(M,g)}$.
Composing this map with the quotient map $\SSS^m\to \SSS_+^m$ gives us a short map $s\:\SSS^m\z\to(M,g)$.

Since $M$ is not homeomorphic to the sphere,
the diameter-sphere theorem implies that $\diam (M,g)\le \tfrac\pi2$.
In particular, $\cBall[p,\tfrac\pi2]_{(M,g)}\z=(M,g)$.
Hence, the map $s$ is onto.

\parbf{\ref{ex:busemann-CBB}.}
Apply \ref{ex:distfun-semiconcave}.

\parbf{\ref{ex:bus+bus}.}
By the triangle inequality,
\[\dist{\gamma(-t)}{x}{}+\dist{\gamma(t)}{x}{}-2\cdot t\ge 0\]
for any $t\ge 0$.
Passing to the limit as $t\to\infty$, we get the result.

\parbf{\ref{ex:oplus}}; \ref{SHORT.ex:oplus:a}.
Show and use that $\proj_AB=\proj_BA$ is a one-point set.

\parit{\ref{SHORT.ex:oplus:b}.}
Suppose $\alpha$ is a rectifiable curve with endpoints in $A$.
Show and use that the curve $\beta=\proj_A\circ\alpha$ has the same endpoints and
\[\length\beta\le \length\alpha;\]
moreover, in case of equality, $\beta=\alpha$.

\parbf{\ref{ex:cone-CBB}.}
Suppose $\Cone\spc{X}$ is $\Alex0$.
Observe that two half-lines in $\Cone\spc{X}$ that start from the origin and go in directions $x,y\in\spc{X}$ form a line if and only if $\dist{x}{y}{\spc{X}}\ge \pi$.
Apply the splitting theorem or non-branching of geodesics in $\Cone\spc{X}$ to show that for any $x\in \spc{X}$ there is at most one point $y$ such that $\dist{x}{y}{\spc{X}}\ge \pi$; if such a point exists, then
$\dist{x}{y}{\spc{X}}=\pi$.
Conclude that $\diam \spc{X}\z\le \pi$.

Now choose a quadruple of points $p,x_1,x_2,x_3\in \spc{X}$;
we will identify $\spc{X}$ with the unit sphere in $\Cone\spc{X}$.
Suppose $\dist{p}{x_i}{\spc{X}}<\tfrac\pi2$ for any $i$.
Consider the following points in the cone: $y_i=\tfrac1{\cos \dist{p}{x_i}{\spc{X}}}\cdot x_i$, and $q=p$.
Show that $\EE^2$-comparison for $q,y_1,y_2,y_3$ in $\Cone\spc{X}$ implies $\SSS^2$-comparison for $p,x_1,x_2,x_3$ in $\spc{X}$.
Conclude that $\spc{X}$ is locally $\Alex1$. 
Apply the globalization theorem (\ref{thm:globalization+}).

Now assume $\spc{X}$ is $\Alex1$ and $\diam\spc{X}\le \pi$.
By \ref{ex:perim-k>0}, the perimeter of any triangle in $\spc{X}$ is at most $2\cdot\pi$.
We need to check $\EE^2$-comparison for a given quadruple of points $q,y_1,y_2,y_3$ in $\Cone\spc{X}$.
We can assume that none of these points is the origin; otherwise perturb them a bit.

Set $x_i=y_i/|y_i|$ for each $i$ and $p=q/|q|$; we can assume that $p,x_1,x_2,x_3$ are distinct in $\spc{X}$, which is the unit sphere in $\Cone\spc{X}$.

Assume the model triangles $\modtrig(px_1x_2)$, $\modtrig(px_2x_3)$, and $\modtrig(px_3x_1)$ are defined;
that is, the perimeters of the triangles $[px_1x_2]$, $[px_2x_3]$, and $[px_3x_1]$ are strictly less than $2\cdot\pi$.
Note that $\EE^3\iso\Cone\SSS^2$.
Use this together with the $\SSS^2$-comparison for $p,x_1,x_2,x_3$ in $\spc{X}$ to show that $\EE^2$-comparison holds for $q,y_1,y_2,y_3$ in $\Cone\spc{X}$.

Finally, if one of the model triangles is undefined, consider a rescaling of $\spc{X}$ by a factor $\lambda$ slightly smaller than $1$.
Apply the argument above to show that the comparison holds for the corresponding points in $\Cone(\lambda\cdot\spc{X})$ and pass to the limit as $\lambda\to 1$.

\parit{Comment.}
The last part of the proof is close to the argument in \ref{thm:CBB-closed}.

\parbf{\ref{ex:|antisum|}.}
Observe that
\begin{align*}
\langle u,u\rangle+\langle v,u\rangle+\langle w,u\rangle &\ge 0,
\\
\langle u,v\rangle+\langle v,v\rangle+\langle w,v\rangle &\ge 0,
\\
\langle u,w\rangle+\langle v,w\rangle+\langle w,w\rangle &= 0.
\end{align*}
Add the first two inequalities and subtract the last identity.

\parbf{\ref{prop:two-opp}.}
Apply \ref{prop:opposite} to show that 
$\langle v,v\rangle =\langle v,w\rangle=\langle w,w\rangle$,
and use it.

\parbf{\ref{ex:3<,>=0}.} Show and use that
\[\langle u,x\rangle +\langle v,x\rangle +\langle w,x\rangle \ge 0\]
and
\[\langle u,-x\rangle +\langle v,-x\rangle +\langle w,-x\rangle \ge 0.\]

\parbf{\ref{ex:-u}.} Part $\Rightarrow$ is evident.
To prove part $\Leftarrow$, observe that 
\[\langle u^*,u^*\rangle =-\langle u,u^*\rangle\le \langle u,u\rangle\]
and since $|u|=|u^*|$, we have equality.

\parbf{\ref{ex:grad-dist}.}
Apply \ref{ex:-u}.

\parbf{\ref{ex:dim-max}.}
Let $L\subset \T_p$ be a subspace isometric to $\EE^m$.
Assume $L\ne \T_p$.

Let $N=\pack_\eps\SSS^{m-1}$.
Observe that the unit sphere in $L$ is isometric to $\SSS^{m-1}$ and
use it to show that $\pack_\eps\Sigma_p>N$ for some $\eps>0$.

It follows that we can choose $N+1$ directions in $\Sigma_p$ at angles larger than $\eps$ from each other.
Use it to construct points $x_1,\dots,x_{N+1}$ such that
$\angk p{x_i}{x_j}>\eps$ for $i\ne j$, and apply \ref{thm:n+1}.

\parbf{\ref{ex:dim=1}.}
Argue as in \ref{ex:RisCBB(1)}.

\parbf{\ref{ex:resporka}.}
The only-if part is trivial.
Suppose the configuration $p$, $a_0,\z\dots, a_{m}\in \spc{A}$ meets the condition.
By \ref{ex:grad-dist} the directions $\dir q{a_0},\z\dots,\dir q{a_m}\in \Lin_q$ for a G-delta dense set of points $q\in \spc{A}$.
If $q$ is sufficiently close to $p$, then $\angk q{a_i}{a_j}>\tfrac\pi2$,
and therefore, $\mangle\hinge q{a_i}{a_j}>\tfrac\pi2$ for $i\ne j$.
Conclude that $\dim\Lin_q\ge m$ in this case.
To see this, show and use that given $m+1$ vectors in $\RR^k$ making pairwise obtuse angles with each other, any $m$ of them must be linearly independent.

\parbf{\ref{ex:concave-differential}.}
We can assume that $f(p)=0$ and $f$ is $1$-concave.

Choose a sequence of positive numbers $\lambda_n\to \infty$.
Let $\iota_n\:\spc{A}\z\to \lambda_n\cdot \spc{A}$ be the rescaling maps.
Consider functions $f_n\:\lambda_n\cdot \spc{A}\to \RR$ defined by
\[f_n(\iota_n(x))=\lambda_n\cdot f(x).\]
Observe that $f_n$ is $\tfrac1{\lambda_n}$-concave and $f_n\to \dd_pf$ as $n\to\infty$ for the convergence $(\lambda_n\cdot \spc{A},p)\z\to (\T_p,0)$ from \ref{ex:GHto-tangent-finite-dim}.

Choose a geodesic $\gamma$ in $\T_p$; assume it can be extended beyond its endpoints.
Construct a sequence of geodesics $\gamma_n$ in $\lambda_n\cdot \spc{A}$ that converge to $\gamma$ (here you should use that $\gamma$ can be extended).
Since $f_n\circ\gamma_n\z\to \dd_pf\circ\gamma$, we get that $\dd_pf$ is concave on $\gamma$, which implies the statement.

\parit{Remark.}
The same argument shows that given a Gromov--Hausdorff convergence $\spc{X}_n\to\spc{X}_\infty$ of $\Alex{\kappa}$ spaces and $\lambda$-concave functions $f_n$ on $\spc{X}_n$ that pointwise converge to a function $f_\infty$ on $\spc{X}_\infty$, the limit function $f_\infty$ is also $\lambda$-concave.

\parbf{\ref{ex:proof-right-inverse}}; \ref{SHORT.ex:proof-right-inverse:grad}.
By \ref{ex:distfun-semiconcave}, each function $\distfun_{a_i}$ is semiconcave in a small neighborhood of $p$.
Therefore we can choose $\lambda$ and $r>0$ so that $f_{\bm{y}}$ is $\lambda$-concave in $\oBall(p,r)$; further we will assume that $r$ is sufficiently small.
Choose $\alpha>0$ such that $\angk{x}{a_i}{a_j}>\tfrac\pi2+\alpha$ for all $i\ne j$;
we may assume that $\alpha<\tfrac{1}{10}$;
in particular,
\[(\dd_x\distfun_{a_j}{}{})(\dir{x}{a_i})
\ge
-\cos\angk{x}{a_i}{a_j}
>
\tfrac\alpha2\eqlbl{inq-a_j}\]
for $j\ne i$.

By the definition of gradient and \ref{ex:d(distfun):<}, we have
\begin{align*}
-(\dd_x\distfun_{a_i}{}{})(\nabla_x f_{\bm{y}})
&\ge
\<\dir x{a_i},\nabla_x f_{\bm{y}}\>
\ge
\\
&\ge
(\dd_xf_{\bm{y}})(\dir x{a_i}).
\end{align*}
If $\dist{a_i}{x}{}>y_i$, then 
\[\dd_xf_{\bm{y}}=\sigma+\eps\cdot \dd_x\distfun_{a_0},\]
where $\sigma$ is a minimum of a subset of the following functions:
$0$, and $\dd_x\distfun_{a_j}$ for $0\ne j\ne i$.
By \ref{inq-a_j}, 
\[(\dd_x\distfun_{a_i}{}{})(\nabla_x f_{\bm{y}})< -\tfrac\alpha2\cdot\eps.\]
Hence (\ref{111}) holds for all sufficiently small $\eps>0$.

Now assume that $\dist{a_i}{x}{}-y_i=\min_j\{\dist{a_j}{x}{}\z-y_j\}<0$.
Then
\begin{align*}
\dd_x f_{\bm{y}}
&=
\min_{j\in S} \{\,\dd_x\distfun_{a_j}\,\}+\eps\cdot \dd_x\distfun_{a_0}
\le
\\
&\le
\dd_x \distfun_{a_i}{}{}+\eps\cdot\dd_x\distfun_{a_0}{}{},
\end{align*}
where $j\in S$ if and only if $\dist{a_i}{x}{}-y_i=\dist{a_j}{x}{}-y_j$.
Applying \ref{inq-a_j}, we get
\begin{align*}
(\dd_x \distfun_{a_i}{}{})(\nabla_x f_{\bm{y}})
&\ge 
\dd_xf_{\bm{y}}(\nabla_x f_{\bm{y}}) -\eps\cdot(\dd_x \distfun_{a_0}{}{})(\nabla_x f_{\bm{y}}) 
\ge 
\\
&\ge
\bigl[(\dd_xf_{\bm{y}})(\dir x{a_0})\bigr]^2-2\cdot \eps
\ge
\\
&\ge
\bigl[\tfrac\alpha2-\eps\bigr]^2-2\cdot \eps.
\end{align*}
Thus, (\ref{222}) holds for all sufficiently small $\eps>0$. 

\parit{\ref{SHORT.ex:proof-right-inverse:alpha}}
Consider the following real-to-real functions:
\[\begin{aligned}
\phi(t)
&\df
\max_{i}\{\dist{a_i}{\alpha_{\bm{y}}(t)}{}-y_i\},
\\
\psi(t)
&\df
\min_{i}\{\dist{a_i}{\alpha_{\bm{y}}(t)}{}-y_i\}.
\end{aligned}\eqlbl{eq:xy-def}\]
Use \ref{SHORT.ex:proof-right-inverse:grad} to show that for $t\in[0,t_0]$, we have $\phi^+(t)<-\tfrac{1}{10}\cdot\eps^2$ if $\phi(t)>0$
and $\psi^+(t)>\tfrac{1}{10}\cdot\eps^2$ if $\psi(t)<0$.
Conclude that $\phi(t_0)=\psi(t_0)=0$; hence the result.

\parit{\ref{SHORT.ex:proof-right-inverse:end}}
A straightforward application of \ref{lem:fg-dist-est} and a reformulation of \ref{SHORT.ex:proof-right-inverse:alpha}.

\parit{Remark.}
By \ref{lem:fg-dist-est}, the constructed map $\map$ is bi-Hölder with the exponent $\tfrac12$.
In particular, \textit{every infinite-dimensional Alexandrov space $\spc{A}$ contains a bi-Hölder copy of a Euclidean ball of arbitrary dimension}.
It seems plausible that $\spc{A}$ should contain a bi-Lipschitz copy of a Euclidean ball of arbitrary dimension,
but this question is open.

\parbf{\ref{ex:proof-dist-chart}.}
Apply the (\textit{n}+1)-comparison (\ref{thm:n+1}) to show that at least one of the inequalities
\[
\mangle\hinge xy{a_0}<\tfrac\pi2-\eps,\ \dots,\ \mangle\hinge xy{a_m}<\tfrac\pi2-\eps,
\]
holds.
Similarly, we get that at least one of the inequalities
\[
\mangle\hinge yx{a_0}<\tfrac\pi2-\eps,\ \dots,\ \mangle\hinge yx{a_m}<\tfrac\pi2-\eps,
\]
holds.

Suppose our statement does not hold for $x$ and $y$ in a sufficiently small neighborhood of $p$.
It follows that 
\[\mangle\hinge xy{a_0}<\tfrac\pi2-\eps
\quad\text{and}\quad
\mangle\hinge yx{a_0}<\tfrac\pi2-\eps.
\eqlbl{eq:a0}
\]
Note that $\dist{x}{y}{}$ is small compared to $\dist{a_0}{x}{}$, $\dist{a_0}{y}{}$ and $\eps$.
Therefore, the comparison contradicts \ref{eq:a0}. 

By the construction, $f$ is Lipschitz.
By the above, we can choose an index $i>0$ so that $\mangle\hinge xy{a_i}<\tfrac\pi2-\eps$ (if $\mangle\hinge yx{a_i}<\tfrac\pi2-\eps$, then swap $x$ and $y$).
By comparison, there is $c>0$ such that $\dist{a_i}{y}{}\le \dist{a_i}{x}{}-c\cdot \dist{x}{y}{}$.
Hence $f$ is bi-Lipschitz, and now \ref{thm:right-inverse} implies \ref{thm:dist-chart}.

\parbf{\ref{ex:BG}.} 
Follow the proof of the Bishop--Gromov inequality, plus prove the following two inequalities
\begin{align*}
\frac{\sinh r_2}{r_2}\cdot \dist{\log_p x}{\log_p y}{\T_p} &\ge\dist{x}{y}{\spc{A}},
\\
\sinh r_2\cdot\dist{w(x)}{w(y)}{\spc{A}} &\ge \sinh r_1\cdot\dist{x}{y}{\spc{A}}
\end{align*}
for any $x,y\in\oBall(p,r_2)$.

\parbf{\ref{ex:diam-compact:proper}.}
Reuse the argument from the first part of the proof of the Bishop--Gromov inequality.

\parbf{\ref{ex:vol>0}.} Applying the right-inverse theorem, show that some ball has positive volume.
Further, use the Bishop--Gromov inequality to show that any ball with positive radius has positive volume.

Alternatively, use \ref{thm:dist-chart} applied to a neighborhood of any point $p$ with $\T_p\iso \EE^m$.

\parbf{\ref{ex:tangent=Em}}; \ref{SHORT.ex:tangent=Em:2balls}.
Apply the dilation map that was used in the proof of the Bishop–Gromov inequality (\ref{inq:BG} and \ref{ex:BG}).

\parit{\ref{SHORT.ex:tangent=Em:ball}.}
Apply the Bishop--Gromov inequality (\ref{inq:BG} and \ref{ex:BG}).

\parbf{\ref{ex:dim=dim}.}
Suppose $K$ is a compact set in $\spc{A}$ such that $\HausDim K\ge m$.
Use the map $w$ from the proof of the Bishop--Gromov inequality (\ref{inq:BG} and \ref{ex:BG}) to show that any open ball in $\spc{A}$ contains a compact set $K'$ such that $\HausDim K'\ge m$.

Use this in addition to the arguments in \ref{thm:dim=dim}. 

\parbf{\ref{ex:finite-tan}};
\ref{SHORT.ex:finite-tan:tan}. Apply \ref{ex:geodesic-cone}, \ref{prop:Tan-is-CBB(0)}, and \ref{thm:finite-space-of-directions}.

\parit{\ref{SHORT.ex:finite-space-of-directions-dim}.}
Apply \ref{ex:resporka} to show that $\LinDim\T_p\le\LinDim\spc{A}$ (argue as in \ref{prop:Tan-is-CBB(0)}).
Argue as in the finite-dimensional case in the proof of \ref{thm:dim=dim} to obtain that $\LinDim\T_p\ge\LinDim\spc{A}$.

\parit{\ref{SHORT.ex:finite-tan:sigma}.}
By \ref{thm:finite-space-of-directions} for any two points $\xi,\zeta\in\Sigma_p$ such that $\dist{\xi}{\zeta}{\Sigma_p}<\pi$ there is a geodesic $[\xi\zeta]_{\Sigma_p}$.
Suppose $\dist{\xi}{\zeta}{\Sigma_p}\ge\pi$.
Then $\T_p$ contains a line thru the origin in the directions $\xi$ and $\zeta$.
By \ref{SHORT.ex:finite-tan:tan} we can apply the splitting theorem (\ref{thm:splitting}) to $\T_p$.
We get that $\Sigma_p$ is a spherical suspension (over a nonempty space!) with poles $\xi$ and $\zeta$.
Hence, $\dist\xi\zeta{}=\pi$ and there is a geodesic $[\xi\zeta]$.

\parbf{\ref{ex:dim-lim}.}
Apply the 4-point condition to prove that $\spc{A}_\infty$ is $\Alex\kappa$.
To prove the dimension bound apply \ref{ex:resporka}.

\parbf{\ref{ex:diam-compact:GH}.}
Argue as in \ref{thm:gromov-compactness} to construct a Gromov--Hausdorff convergence of $\cBall(p_n,R)_{\spc{A}_n}$ for a given $R>0$, then apply a diagonal procedure to construct the needed convergence.

\parbf{\ref{ex:pack-vol:dim}.}
By \ref{ex:dim-lim}, $\dim\spc{A}_\infty\le m$.
To show that $\dim\spc{A}_\infty\ge m$,
apply \ref{cor:euclid-subcone} and \ref{ex:tangent=Em:ball}.

\parit{Comment.}
The following stronger statement holds
\[\vol_m\spc{A}_\infty=\lim_{n\to\infty} \vol_m\spc{A}_n.\]
In other words, if $\bm{K}_m\subset \GH$ denotes the set of isometry classes of all compact $\Alex\kappa$ spaces with dimension $\le m$, then the function $\vol_m\:\bm{K}_m\to \RR$ is continuous \cite[10.8]{burago-gromov-perelman}.

\parbf{\ref{ex:barrier}.}
We can assume that $|f^+(t_0)|>c_1$.
Choose a smooth upper barrier $\bar f$ of $f$ at $t_0$ with $\bar f''\le \lambda_1$ and $\bar f'(t_0)=f^+(t_0)$, provided by Lemma \ref{lem:barrier}.
Observe that the composition $\phi\circ\bar f$ is an upper barrier of $\phi\circ f$ at $t_0$.
It remains to apply the chain rule to $\phi\circ\bar f$.

\parbf{\ref{ex:pack-sphere}.}
Choose a ball $B=\cBall[x,R]$ such that $p\notin\cBall[x,2\cdot R]\subset \Omega$.
By \ref{ex:tangent=Em:ball}, we can choose a constant $\Const_1>0$ such that for every sufficiently small $\eps>0$ the ball $B$ contains at least $\Const_1/{\eps^m}$ points $q_1,\dots,q_N$ such that $\dist{q_i}{q_j}{}>10\cdot\eps$.
Slice $B$ into sets
\[B_k=\set{z\in B}{k\cdot \eps\le \dist{p}{z}{}< {k+1}\cdot \eps}.\]
There are at most $1+2\cdot R/\eps$ nonempty sets $B_k$.
Therefore, for a fixed constant $\Const_2>0$ we can choose $B_k$ that contains at least $\Const_2/\eps^{m-1}$ points from $\{q_1,\dots,q_N\}$.
Let $r=k\cdot \eps$.
For each $q_i\in B_k$, choose $q_i'\in [pq_i]$ such that $\dist{p}{q_i'}{}=r$.
It remains to observe that (1) $q_i'\in\Omega$ for each $i$ and (2) $\dist{q_i'}{q_j'}{}>8\cdot\eps$ and $\mangle\hinge p{q_i'}{q_j'}>8\cdot \eps/r$ if $i\ne j$.

\parbf{\ref{ex:no-conc}.}
Consider the infinite product $\SSS^1\times ({\tfrac 12}\cdot \SSS^1)\times ({\tfrac 14}\cdot \SSS^1)\times\dots$

\parbf{\ref{ex:dist-chart+strictly-concave}.} Combine \ref{ex:proof-dist-chart} and the proof of \ref{thm:strictly-concave}.

\parbf{\ref{ex:no-cont-lifting}.}
Consider the canonical metric $g$ on the round unit sphere $\SSS^3$ and the Hopf bundle $\SSS^3\to \SSS^2$.
Let $g_n$ be obtained from $g$ by rescaling the Hopf fibers by $\frac{1}{n}$.
Then all $(\SSS^3,g_n)$ have nonnegative curvature and the sequence converges to the round metric on $\SSS^2$ of radius $\frac{1}{2}$.
Verify that this convergence provides the desired counterexample.

\parit{Comment.}
The spaces $(\SSS^3,g_n)$ are called \index{Berger sphere}\emph{Berger spheres}.
Checking the curvature bounds is easy once you realize that $(\SSS^3,g_n)$ is isometric to a quotient space $(\SSS^3\times (\eps_n\cdot \SSS^1))/\SSS^1$, where $\SSS^1$ acts diagonally by the Hopf action on $\SSS^3$ and the standard action on $\SSS^1$.

\parbf{\ref{ex:conic}.}
Let $V$ and $W$ be two conic neighborhoods of a point~$p$.
Without loss of generality, we may assume that $V\Subset W$;
that is, the closure of $V$ lies in $W$.

Construct a sequence of embeddings $f_n\:V\to W$
such that 
\begin{itemize}
\item 
For any compact set $K\subset V$ 
there is a positive integer $n=n_K$ such that 
$f_n(k)=f_m(k)$ for any $k\in K$ and $m, n \ge n_K$.
\item For any point $w\in W$ there is a point $v\in V$ such that $f_n(v)=w$ for all large $n$.
\end{itemize}

Once such a sequence is constructed, $f\:V\to W$ can be defined by $f(v)=f_n(v)$ for all large values of $n$; this gives the needed homeomorphism.

The sequence $f_n$ can be constructed recursively by the formula
\[f_{n+1}=\Psi_n\circ f_n\circ \Phi_n,\]
where $\Phi_n\:V\selfmap$
and $\Psi_n\:W\selfmap$
are homeomorphisms
of the form 
\[\Phi_n(x)=\phi_n(x)\ast x\quad \text{and}\quad \Psi_n(x)=\psi_n(x)\star x,\]
where $\phi_n\:V\to \RR_{\ge 0}$, $\psi_n\:W\to \RR_{\ge 0}$ are suitable continuous functions;
``$\ast$'' and ``$\star$'' denote the multiplications in the cone structures of $V$ and $W$ respectively.

\parit{Comment.}
One may also read the original proof by Kyung Whan Kwun \cite{kwun1964}.

\parbf{\ref{ex:conic-tangent}}; \ref{SHORT.ex:conic-tangen:tangent}.
Apply \ref{thm:spherical-nbhd} and \ref{lem:kwun}.

\parit{\ref{SHORT.ex:conic-tangen:dir}.}
Apply \ref{SHORT.ex:conic-tangen:tangent}.

\parit{\ref{SHORT.ex:conic-tangen:example}.}
Recall that the Poincaré homology sphere can be obtained as a quotient space $\Sigma=\SSS^3/\Gamma$ by an isometric action of a finite group $\Gamma$ --- the so-called binary icosahedral group.
By the double suspension theorem
$\Susp^2\Sigma\cong\SSS^5$; see, for example, \cite{daverman}.
Note that $\Susp^2\Sigma$ is an Alexandrov space and it has a point with space of directions isometric to $\Susp\Sigma$.
Observe that $\Susp\Sigma$ is not a manifold; in particular $\Susp\Sigma\ncong\SSS^4$.
Therefore the pair $\Susp^2\Sigma$ and $\SSS^5$ provides the needed example.

\parbf{\ref{ex:bry2bry}.}
Apply \ref{thm:spherical-nbhd}, \ref{lem:kwun}, and \ref{thm:top-bry}.

\parbf{\ref{ex:bry-closed}.}
Let $\spc{A}$ be a finite-dimensional Alexandrov space.
Choose $x\in\spc{A}$.
By \ref{thm:spherical-nbhd}, a neighborhood $U\ni x$ is homeomorphic to $\T_x$.
Therefore, \ref{ex:bry2bry} implies that $U\cap\partial\spc{A}=\emptyset\ \Leftrightarrow\ x\notin \partial\spc{A}$;
that is, the complement $\spc{A}\setminus\partial\spc{A}$ is open, and therefore, $\partial\spc{A}$ is closed.

\parbf{\ref{ex:pz<ypz}.}
Consider the model triangle $[\tilde x\tilde y\tilde z']=\modtrig(xyz)$.
\begin{figure}[ht!]
\vskip-0mm
\centering
\includegraphics{mppics/pic-1015}
\end{figure}

Show that 
\[\dist{\tilde p}{\tilde z}{}\le \dist{\tilde p}{\tilde z'}{}\le\side\hinge yp{z}.\]

\parbf{\ref{ex:bry-connected}.}
Assume that $\partial\spc{A}$ has at least two connected components, say $A$ and $B$.
Let $\gamma$ be a geodesic that minimizes the distance from $A$ to $B$.

Consider a two-sided infinite sequence of copies of $\spc{A}$
\[\dots,\spc{A}_{-1},\spc{A}_{0},\spc{A}_{1},\dots\]
Let us glue $\spc{A}_{i}$ to $\spc{A}_{i+1}$ along $A$ if $i$ is even and along $B$ if $i$ is odd.

By the doubling theorem, every point in the obtained space $\spc{N}$ has a neighborhood that is isometric to a neighborhood of the corresponding point in $\spc{A}$ or its doubling.
By the globalization theorem, $\spc{N}$ is $\Alex1$.

The copies of $\gamma$ in $\spc{A}_{i}$ form a line in $\spc{N}$.
By the splitting theorem, $\spc{N}$ is isometric to a product $\spc{N}'\oplus \RR$.
Since $\dim\spc{N}>1$, Exercise~\ref{ex:RisCBB(1)} implies that $\diam\spc{N}\le \pi$ --- a contradiction.

\parbf{\ref{ex:dist-to-bry}}; \ref{SHORT.ex:dist-to-bry:geod}
Apply the gradient flow as in the proof of \ref{thm:doubling} (part \ref{SHORT.thm:partial-grad:flow}$_{m}\z+$\ref{SHORT.thm:doubling:doubling}$_{m-1}\z\Rightarrow$\ref{SHORT.thm:doubling:doubling}$_m$).

\parit{\ref{SHORT.ex:dist-to-bry:dist}.}
We can assume that $\dim \spc{A}\ge 2$; otherwise the statement is trivial.

Choose an interior point $x$ on $\gamma$;
we can assume that $x=\gamma(0)$.
Let $y\in \partial \spc{A}$ be a closest point to $x$,
and let $\alpha=\mangle(\dir xy,\gamma^+(0))$.

By \ref{SHORT.ex:dist-to-bry:geod}, we can assume that $x\notin \partial \spc{A}$.
Show that $\T_y=[0,\infty)\times\T_y\partial \spc{A}$
and $\dir yx\perp \T_y\partial \spc{A}$.

Given a vector $v\in \T_y$, denote by $\bar v$ its projection to $\T_y\partial \spc{A}$.
Apply the comparison and \ref{prop:gexp} to show that 
\[\dist{\gamma(t)}{\gexp_y(\overline{\log_y\gamma(t)})}{}\le \dist{x}{y}{}+t\cdot\cos\alpha.\]
Conclude that $(\distfun_{\partial \spc{A}}\circ\gamma)''(0)\le 0$ in the barrier sense.

\parbf{\ref{ex:liberman}.}
Suppose $\gamma$ is defined on the interval $[0,\ell]$.
Assume that the function $\rho\:t\mapsto \tfrac12\cdot\distfun_p^2\circ\gamma(t)$ is not $1$-concave.
Let $\bar\rho\:[0,\ell]\to\RR$ be the minimal $1$-concave function such that $\bar\rho\ge \rho$.
Note that $\bar\rho=\rho$ at the ends of $[0,\ell]$ but $\bar\rho>\rho$ at some point on $(0,\ell)$.

Consider the curve $\bar\gamma(t)\df \GF_f^{s(t)}\gamma(t)$, where $f=\tfrac12\cdot\distfun_p^2$ and $2\cdot s(t)\z=\ln\circ\bar\rho(t)-\ln\circ\rho(t)$.
Use the first distance estimate to show that $\length\bar\gamma<\length\gamma$ and arrive at a contradiction.

\parit{Comment.}
This is a version of \index{Liberman's lemma}\emph{Liberman's lemma} \cite{liberman}; it was proved by Grigory Perelman and the third author \cite{perelman-petrunin}.

\parbf{\ref{ex:native}.}
The only-if direction follows since $f\z\circ\proj$ is invariant under the canonical isometric involution of $\spc{W}$.

Let us prove the if direction.
Choose a geodesic $\gamma$ in $\spc{W}$.
Arguing as in the proof of \ref{thm:doubling:doubling}, we get 
that $\gamma$ can cross the common boundary of two halves $\spc{A}_0$ and $\spc{A}_1$ of $\spc{W}$ at most once, or it lies in the common boundary.

In the latter case, $\lambda$-concavity of $f\circ\proj\circ\gamma$ follows from $\lambda$-concavity of $f$.
In the former case, the concavity has to be checked only at the point of crossing;
we may assume that it happens at $x=\gamma(0)$.
Since $\nabla_xf\in\partial\T_x$ for any $x\in\partial\spc{A}$, the $f$-gradient flows on $\spc{A}_0$ and $\spc{A}_1$ agree on the common boundary; so they induce a continuous flow on $\spc{W}$.

Assume $f\circ\proj\circ\gamma$ is not $\lambda$-concave at $0$.
Apply the constructed flow on $\spc{W}$ to shorten $\gamma$ keeping its ends as in the proof of \ref{ex:liberman},
and arrive at a contradiction.

\parbf{\ref{ex:Hilbert/G}.} Read \cite[Section 4]{terng-thorbergsson} and/or the solution for ``Quotient of the Hilbert space'' in \cite{petrunin2020}.

\parbf{\ref{ex:sumbetries(S^2)}}; \ref{SHORT.ex:sumbetries(S^2):1}.
Choose an isometric $\SSS^1$-action on $\SSS^2$ that fixes the poles of the sphere.
Consider the projection to the quotient space $\sigma_1\:\SSS^2\z\to \SSS^2/\SSS^1=[0,\pi]$.

\parit{\ref{SHORT.ex:sumbetries(S^2):2}.}
Take a semicircle $\gamma$ on $\SSS^2$ and define
$\sigma_2(x)\df\distfun_\gamma(x)_{\SSS^2}$.

\parit{\ref{SHORT.ex:sumbetries(S^2):n}.}
Subdivide $\SSS^2$ into $\SSS^1$-orbits of the action from~\ref{SHORT.ex:sumbetries(S^2):1}.
Cut $\SSS^2$ into two hemispheres by meridians, rotate one hemisphere by an angle $\alpha=\pi/n$ and glue it back.
Observe that there is a submetry $\sigma_n$ such that the inverse image $\sigma_n^{-1}\{y\}$ is a union of the arcs from the original $\SSS^1$-orbits.

Note that for $n=2$ we get the solution in \ref{SHORT.ex:sumbetries(S^2):2}.

\parbf{\ref{ex:sumbetries(E^2)}.}
Show that for any $x\in\EE^2$ there is a half-line $H\ni x$ such that 
the restriction $\sigma|_H$ is an isometry.
Suppose such a half-line $H$ starts at $p$ and passes thru $q$.
Show that $\langle x-p,q-p \rangle\le 0$ for any $x\in \sigma^{-1}\{0\}$.
Conclude that $\sigma^{-1}\{0\}$ is a convex closed set.
Use the definition of submetry to show that $\sigma^{-1}\{0\}$ has no interior points.
Make a conclusion.

\parbf{\ref{ex:S^3/S^1}};
\ref{SHORT.ex:S^3/S^1:pq}.
Our $\SSS^1$ is a commutative subgroup of $\SO(4)$.
Therefore it is a subgroup of a maximal torus in $\SO(4)$.
Show that the described torus action is induced by a maximal torus in $\SO(4)$.
Use that maximal tori in $\SO(4)$ are conjugate.

\parit{\ref{SHORT.ex:S^3/S^1:sphere}.}
Cut $\SSS^3$ along the Clifford torus $\tfrac1{\sqrt2}\cdot (\SSS^1\times \SSS^1)$.
Observe that the quotient of each half is a disk;
conclude that $\Sigma_{p,q}$ is a sphere.
The torus action on $\SSS^3$ induces the needed $\SSS^1$-action on $\Sigma_{p,q}$.

\parit{\ref{SHORT.ex:S^3/S^1:a}+\ref{SHORT.ex:S^3/S^1:b}+\ref{SHORT.ex:S^3/S^1:c}.} Calculations.

\parit{\ref{SHORT.ex:S^3/S^1:cc}.}
Consider the map $\Sigma_{p,q}\to\Sigma_{1,1}$ that sends poles to poles,
preserves the distance to the poles and respects the $\SSS^1$-actions.

\parbf{\ref{ex:extr-point(-1)}.}
Choose an $m$-dimensional $\Alex{-1}$ space $\spc{A}$.
Repeat the proof of \ref{thm:extr-point} for fixed $\alpha=\tfrac1{10}$ showing that any unit ball in $\spc{A}$ has at most $1000^m$ extremal points.
Then apply \ref{ex:tangent=Em:ball}+\ref{ex:pack-net} for $\eps=1$, $\kappa=-1$, and $r=\diam\spc{A}$.

\parbf{\ref{ex:number(m-1)}};
\ref{SHORT.ex:number(m-1):2}.
Suppose $\mathfrak{M}_{m-1}(\Gamma)\ge 3$;
that is, $\spc{A}=\EE^m/\Gamma$ has at least 3 boundary components.
Follow Case~3 in the proof of \ref{thm:hsiang-kleiner} to glue a train-like space from copies of $\spc{A}$ using two of these components.
Show that the obtained space splits and arrive at a contradiction.

(Alternatively, apply a similar construction to all components of the boundary.
Show that the obtained space has {}\emph{exponential volume growth};
that is, there is $a>1$ such that $\vol \oBall(p,r)>a^r$ for all large~$r$.
Arrive at a contradiction with the Bishop--Gromov inequality.)

\parit{\ref{SHORT.ex:number(m-1):1}.}
Apply the doubling theorem as in Case~2 in the proof of \ref{thm:hsiang-kleiner}.

\parbf{\ref{ex:S1actsS3}.}
Show that the quotient space $\Delta=\spc{A}/\mathbb{S}^1$ is an $\Alex1$ disk and $\gamma$ projects isometrically to its boundary $\partial\Delta$.
It remains to show that the perimeter of $\Delta$ cannot exceed $2\cdot\pi$.
The latter follows from Lytchak's problem \cite[3.3.5]{petrunin:survey};
it states that if $\Delta$ is an $m$-dimensional $\Alex1$ space, then $\vol_{m-1}\partial \Delta\z\le \vol_{m-1}\partial \mathbb{S}^{m-1}$.

\parbf{\ref{ex:geodesic-vertex}.}
Suppose a geodesic $\gamma$ passes thru a vertex $v$.
Denote by $\alpha$ and $\beta$ the angles that $\gamma$ cuts at $v$.
Since $v$ is essential, $\alpha+\beta<2\cdot\pi$.
Therefore $\alpha<\pi$ or $\beta<\pi$.
Arrive at a contradiction by showing that $\gamma$ is not length-minimizing.

\parbf{\ref{ex:gauss-bonnet}.}
Assume $\spc{P}$ has no boundary.
Denote by $k$, $l$, and $m$ the number of vertices, edges, and triangles, respectively, in a chosen triangulation of $\spc{P}$.
Note that
\[2\cdot l=3\cdot m
\qquad\text{and}\qquad
k-l+m=\chi(\spc{P}).\]
The first identity follows since each edge appears in two triangles and each triangle has 3 edges;
the second identity is Euler's formula.

Since each triangle contributes $\pi$ to the total sum of angles, we get that the total curvature is $2\cdot\pi\cdot k-\pi\cdot m$.
It remains to apply straightforward algebraic manipulations.

If $\spc{P}$ has nonempty boundary, then pass to its doubling, apply the formula and rewrite the result using inner turns.

\parbf{\ref{ex:poly-CBB}.}
We need to show that if a polyhedral surface is $\Alex0$, then the total angle $\theta$ at every vertex $p$ is at most $2\cdot\pi$.

Assume that $\theta>2\cdot\pi$.
Let $\phi=\min\{\,\pi,\tfrac13\cdot\theta\,\}$.
We can choose three points $x_1$, $x_2$, and $x_3$ close to $p$ such that
$\mangle \hinge p{x_i}{x_j}=\phi$ for $i\ne j$.
Since the points $x_i$ are close to $p$, we have $\mangle \hinge p{x_i}{x_j}=\angk p{x_i}{x_j}$.
The latter contradicts $\EE^2$-comparison.

\parbf{\ref{ex:construction}.}
\ref{SHORT.ex:approximation:nbhd}.
Apply \ref{cor:convex-nbhd}.

\parit{\ref{SHORT.ex:approximation:triangulation}.}
Observe that any chord divides a convex figure on $\spc{P}$ into two convex figures.
Use it to show that the union of two convex polygons can be triangulated by convex triangles.
Apply the last statement recursively to a finite cover of $\spc{P}$ by convex polygons.

\parit{\ref{SHORT.ex:approximation:poly}.}
Let $\Delta_1,\dots,\Delta_n$ be the triangles incident to $v$, listed in cyclic order, and let $\alpha_i$ be the angle of $\Delta_i$ at $v$.
Denote by $I_i\subset\Sigma_v$ the set of directions pointing into $\Delta_i$.
Note and use that (1) each $I_i$ is an interval of length $\alpha_i\le\tfrac12\cdot\length\Sigma_v$ and (2)
the intervals $I_i$ intersect only at their common endpoints.

\parbf{\ref{ex:approximation}}; \ref{SHORT.ex:approximation:excess}.
Show that the total angle around each vertex of the triangulation is at most $2\cdot\pi$ and argue as in \ref{ex:gauss-bonnet}.

\parit{\ref{SHORT.ex:approximation:length}.}
Denote by $\alpha$, $\beta$, and $\gamma$ the angles of $\Delta$, and let
$\tilde\alpha$, $\tilde\beta$, and $\tilde\gamma$ be the corresponding model angles.
Since $\tilde\alpha+\tilde\beta+\tilde\gamma=\pi$, we have
\[
\excess\Delta
=
\alpha+\beta+\gamma
-\tilde\alpha-\tilde\beta-\tilde\gamma.
\]
By comparison, $\tilde\alpha\le\alpha$, $\tilde\beta\le\beta$, and $\tilde\gamma\le\gamma$.
Hence
\[\alpha\le\tilde\alpha+\excess\Delta,\quad
\beta\le\tilde\beta+\excess\Delta,\quad
\gamma\le\tilde\gamma+\excess\Delta.
\]
It remains to apply the hinge comparison (\ref{angle}).

\parit{\ref{SHORT.ex:approximation:area}.}
The first inequality follows from the generalized Kirszbraun theorem proved by Urs Lang and Viktor Schroeder \cite{lang-schroeder:kirszbraun, alexander-kapovitch-kirszbraun, alexander-kapovitch-petrunin2024};
we will not discuss its proof since this inequality is not used in the proof.
In the proof of the second inequality, we use Alexandrov's idea \cite[X §~1]{alexandrov-1948}.

Suppose that a convex triangle $\Delta$ has vertices $x$, $y$, and $z$, with opposite side lengths $a$, $b$ and $c$ and angles $\alpha$, $\beta$, and $\gamma$, respectively.
Let $\tilde\Delta$ be its solid model triangle and let $\tilde \alpha$, $\tilde \beta$, and $\tilde \gamma$ be the corresponding angles.

Use that $\log_x$ is a noncontracting map to show that
\[\area\Delta-\area\tilde\Delta\le \tfrac12\cdot(\alpha-\tilde\alpha)\cdot (\diam \Delta)^2+o(a);\]
moreover, $o(a)/a\to0$ uniformly for all the triangles considered below.

Now let us subdivide $\Delta$ into triangles $\Delta_1,\dots,\Delta_n$ with common vertex $x$ and small opposite sides $a_1,\dots,a_n$; so, $a=a_1+\dots+a_n$.
The model triangles $\tilde\Delta_1,\dots,\tilde\Delta_n$ can be arranged on the plane with common corresponding sides.
This produces a convex region $\tilde\Lambda$ bounded by two straight sides $b$ and $c$ and a polygonal side of total length $a$.
Summing up the above inequalities for $\Delta_1,\dots,\Delta_n$, we get
\[\area\Delta-\area\tilde\Lambda\le \tfrac12\cdot(\alpha-\tilde\alpha)\cdot (\diam\Delta)^2+o(a_1)+\dots +o(a_n).\]
It follows that given $\eps>0$, we can choose the subdivision so that
\[\area\Delta-\area\tilde\Lambda\le \tfrac12\cdot(\alpha-\tilde\alpha)\cdot (\diam\Delta)^2+\eps.\]

By the hinge comparison (\ref{angle}), the angles adjacent to polygonal sides of $\tilde\Lambda$ do not exceed $\beta$ and $\gamma$.
It follows that
\[\area\tilde\Lambda-\area\tilde\Delta\le \tfrac12\cdot(\beta+\gamma-\tilde\beta-\tilde\gamma)\cdot (\diam \Delta)^2.\]
Therefore
\[\area\Delta-\area\tilde\Delta\le \tfrac12\cdot\excess\Delta\cdot (\diam \Delta)^2+\eps\]
for arbitrary $\eps>0$, hence the result.

\parbf{\ref{ex:surf-S2}.}
We can assume that the origin lies in the interior of the convex body.
Consider the central projection from its surface, say $\Sigma$, to the sphere $\SSS^2$ centered at the origin.
Show that this projection $\Sigma\to \SSS^2$ is a homeomorphism.

\begin{wrapfigure}{r}{25mm}
\vskip-0mm
\centering
\includegraphics{mppics/pic-1160}
\vskip-0mm
\end{wrapfigure}

\parbf{\ref{pr:tetrahedron}}; \ref{SHORT.pr:tetrahedron:=}.
Cut the surface of $T$ along three edges coming from one vertex $v_1$ and unfold the obtained surface onto the plane.
Show that this way we get a triangle; the three vertices correspond to $v_1$ and the midpoints of sides correspond to the remaining three vertices.
Make a conclusion.

\parit{\ref{SHORT.pr:tetrahedron:perp}}.
Let $v_1,v_2,v_3,v_4\in\RR^3$ be the vertices of $T$.
From \ref{SHORT.pr:tetrahedron:=}, we have that 
\[|v_1-v_2|=|v_3-v_4|,\quad |v_1-v_3|=|v_2-v_4|,\quad|v_1-v_4|=|v_2-v_3|.\]
Use these equalities to show that $\langle v_1-v_2,v_1+v_2-v_3-v_4\rangle=0$.
Make a conclusion.

\parbf{\ref{ex:surface-covergence}.}
We will use that the nearest-point projection $x\mapsto \bar x$ from the Euclidean space $\EE^m$ to a convex body is short; that is, distance-nonexpanding (see for example \cite[13.3]{petrunin-zamora}).
This is a consequence of the following observation: $\langle x-\bar x,\bar y-\bar x\rangle\le 0$ and $\langle y-\bar y,\bar x-\bar y\rangle\le 0$ for any $x,y\in \EE^m$.

Assume $K_\infty$ is nondegenerate.
Without loss of generality, we may assume that 
\[\cBall(0,r)\subset K_\infty\subset\cBall(0,1)\]
for some $r>0$.
Then there is a sequence $\eps_n\to 0$ such that 
\[K_\infty\subset(1+\eps_n)\cdot K_n
\quad\text{and}\quad
K_n\subset(1+\eps_n)\cdot K_\infty\]
for each large $n$.

Given $x\in K_n$, denote by $g_n(x)$ the nearest-point projection of $(1+\eps_n)\cdot x$ to $K_\infty$.
Similarly, given $x\in K_\infty$, denote by $h_n(x)$ the nearest-point projection of $(1+\eps_n)\cdot x$ to $K_n$.
Note that 
\begin{align*}
\dist{g_n(x)}{g_n(y)}{}&\le (1+\eps_n)\cdot\dist{x}{y}{}
\intertext{and}
\dist{h_n(x)}{h_n(y)}{}&\le (1+\eps_n)\cdot\dist{x}{y}{}.
\end{align*}

Denote by $\Sigma_\infty$ and $\Sigma_n$ the surface of $K_\infty$ and $K_n$, respectively.
The inclusions above imply that
\[
g_n(\Sigma_n)=\Sigma_\infty
\quad\text{and}\quad
h_n(\Sigma_\infty)=\Sigma_n.
\]
The preceding inequalities therefore imply
\begin{align*}
\dist{g_n(x)}{g_n(y)}{\Sigma_\infty}
&\le (1+\eps_n)\cdot\dist{x}{y}{\Sigma_n}
\intertext{for any $x,y\in\Sigma_n$, and}
\dist{h_n(x)}{h_n(y)}{\Sigma_n}
&\le (1+\eps_n)\cdot\dist{x}{y}{\Sigma_\infty}
\end{align*}
for any $x,y\in\Sigma_\infty$.
It remains to apply \ref{ex:GH-po}.

Alternatively, since the nearest-point projection cannot increase the length of a curve, we also get
\begin{align*}
\dist{x}{h_n\circ g_n(x)}{\Sigma_n}&\le 10\cdot \eps_n
\\
\dist{y}{g_n\circ h_n(y)}{\Sigma_\infty}&\le 10\cdot \eps_n
\end{align*}
for all large $n$.
Therefore, $g_n$ is a $\delta_n$-isometry $\Sigma_n\to\Sigma_\infty$ for a sequence $\delta_n\to 0$.

\parit{Comment.}
More generally, if a sequence of $m$-dimensional $\Alex\kappa$ spaces $\spc{A}_1,\spc{A}_2,\dots$ converges to $\spc{A}_\infty$ and $\dim \spc{A}_\infty=m<\infty$,
then $\partial \spc{A}_n$ equipped with the induced length metrics converge to $\partial \spc{A}_\infty$.
This statement is a special case of the theorem about extremal subsets proved by the third author \cite[1.2]{petrunin1997}.

\parbf{\ref{ex:alexandrov=<4}.}
\ref{SHORT.ex:alexandrov=<4:>=3}.
Observe that the curvature of any essential vertex is less than $2\cdot\pi$ and apply the Gauss--Bonnet formula \ref{ex:gauss-bonnet}.

\parit{\ref{SHORT.ex:alexandrov=<4:=3}.} Show that geodesics between essential vertices divide the surface into two flat triangles, which have to be isometric since their sides are equal.
Make a conclusion.

\parit{\ref{SHORT.ex:alexandrov=<4:4}.}
Show that geodesics between the essential vertices can be chosen so that they do not cross each other; that is, two such geodesics intersect only at a common endpoint, if they have one.
In this case, they divide the surface into Euclidean triangles.

Since the curvature is nonnegative, the sum of three angles of the triangles at each vertex is at most $2\cdot\pi$.
Show that the three angles at one (and therefore any) vertex satisfy the triangle inequality.
Conclude that these triangles form faces of a tetrahedron (possibly degenerate to a quadrangle).

\parbf{\ref{ex:convex}.} Assume the contrary.
Then there is a minimizing geodesic $\gamma\z{\not\subset}\Delta$ with endpoints $p$ and $q$ in $\Delta$.

Without loss of generality, we may assume that only one arc of $\gamma$ lies outside of $\Delta$.
Reflection of this arc with respect to $\Pi$ together with the remaining part of $\gamma$ forms another curve $\hat\gamma$ from $p$ to $q$;
it runs partly along the surface
and partly outside $K$,
but does not enter the interior of $K$.
Note that
\[\length\hat\gamma=\length\gamma.\]

Denote by $\bar\gamma$ the nearest-point projection of $\hat\gamma$ to $K$.
Since $K$ is convex, the projection does not increase the length.
Show that it decreases the length strictly on the reflected arc.
Therefore,
\[\length\bar\gamma<\length\gamma.\]
This means that $\gamma$ is not length-minimizing
--- a contradiction.

\parbf{\ref{ex:milka}.}
If the plane $py_1y_2$ supports $K$, then 
$\mangle\hinge p{y_1}{y_2}_{\EE^3}=\mangle\hinge p{x_1}{x_2}_{\spc{P}}$.
In this case, the statement follows from \ref{prop:conv-surf-CBB(0)}.

Now suppose that the line segment $[y_1y_2]_{\EE^3}$ intersects $K$.
Choose a geodesic $[y_1y_2]_W$;
note that it contains a point of $K$, say $z$.
Now consider the one-parameter family of points 
\[y_i(t)\df \gamma_i(t)+\gamma_i^+(t)\z\cdot (1-t).\]
This family is not necessarily continuous; note that $y_i(0)=y_i$ and $y_i(1)=x_i$.

Show that for any point $q\in K$, the function $t\mapsto \dist{q}{y_i(t)}{\EE^3}$ is nonincreasing.
Conclude that the function $t\mapsto \dist{q}{y_i(t)}{W}$ is nonincreasing for any $q\in \spc{P}$.
Therefore, 
\begin{align*}
\dist{y_1}{y_2}{W}
&=\dist{y_1(0)}{y_2(0)}{W}=
\\
&=\dist{y_1(0)}{z}{W}+\dist{y_2(0)}{z}{W}\ge
\\
&\ge\dist{y_1(1)}{z}{W}+\dist{y_2(1)}{z}{W}\ge
\\
&\ge\dist{x_1}{x_2}{\spc{P}}.
\end{align*}
The last inequality follows since the nearest-point projection $W\to \spc{P}$ is short.

\begin{wrapfigure}{o}{28mm}
\vskip-3mm
\centering
\includegraphics{mppics/pic-1111}
\vskip-0mm
\end{wrapfigure}

It remains to consider the case when the plane $py_1y_2$ does not support $K$,
and $[y_1y_2]_{\EE^3}$ does not intersect $K$.
In this case the plane $py_1y_2$ intersects $K$ along a convex figure $F$ that lies in the solid triangle 
$py_1y_2$ and contains its vertex $p$.

Choose points $y_1'\in [py_1]_{\EE^3}$ and $y_2'\in [py_2]_{\EE^3}$ such that $[y_1'y_2']$ touches~$F$.
Denote by $x_1'\z\in [px_1]_{\spc{P}}$ and $x_2'\in [px_2]_{\spc{P}}$ the corresponding points;
that is, $\dist{p}{y_1'}{\EE^3}=\dist{p}{x_1'}{\spc{P}}$ and $\dist{p}{y_2'}{\EE^3}=\dist{p}{x_2'}{\spc{P}}$.
From above, we have that $\dist{y_1'}{y_2'}{\EE^3}\ge\dist{x_1'}{x_2'}{\spc{P}}$;
in other words, 
\[\angk p{y_1'}{y_2'}\ge \angk p{x_1'}{x_2'};\]
here we think of $[p{y_1'}{y_2'}]$ as a triangle in $\EE^3$, and $[p{x_1'}{x_2'}]$ as a triangle in $\spc{P}$.
Note that 
\[\angk p{y_1'}{y_2'}=\angk p{y_1}{y_2}
\quad\text{and}\quad
\angk p{x_1}{x_2}\le \angk p{x_1'}{x_2'};
\]
the second inequality follows from \ref{ex:noncreasing}.
Hence the remaining case follows.